\documentclass[11 pt]{amsart}

\usepackage{amsmath,amssymb,amsthm,amsfonts,setspace,hyperref,dsfont,yhmath,euscript}

\usepackage{color}

\usepackage[nameinlink]{cleveref}
\usepackage{hyperref}\hypersetup{colorlinks=true, citecolor=green, linkcolor=blue, hypertexnames=false}

\usepackage[utf8]{inputenc}

\newcommand{\NN}{\mathbb{N}}
\newcommand{\RR}{\mathbb{R}}

\newcommand{\CC}{\mathbb{C}}

\newcommand{\TT}{\mathbb{T}}
\newcommand{\ZZ}{\mathbb{Z}}

\newcommand{\GG}{\mathbb{G}}

\newcommand{\norm}[1]{\lVert#1\rVert}
\newcommand{\abs}[1]{\lvert#1\rvert}

\newtheorem{theorem}{Theorem}[section]
\newtheorem{corollary}[theorem]{Corollary}

\newtheorem{lemma}[theorem]{Lemma}
\newtheorem{proposition}[theorem]{Proposition}
\newtheorem{example}{Example}
\newtheorem{observation}[theorem]{Observation}

\newcommand{\comment}[1]{}
\theoremstyle{definition}
\newtheorem{definition}[theorem]{Definition}
\newtheorem{remark}[theorem]{Remark}

\newtheorem{nsect}{}
\newtheorem{sect}{}

\numberwithin{equation}{section}

\DeclareMathOperator{\Ad}{Ad}

\makeatletter
\renewcommand*\env@matrix[1][*\c@MaxMatrixCols c]{%
  \hskip -\arraycolsep
  \let\@ifnextchar\new@ifnextchar
  \array{#1}}
\makeatother

\begin{document}
\title[Local rigidity]
{Local rigidity of higher rank partially hyperbolic algebraic actions}

\thanks{ $^1$ Based on research supported by NSF grant DMS-1845416}


\author[]{ Zhenqi Jenny Wang$^1$ }

\address{Department of Mathematics\\
        Michigan State University\\
        East Lansing, MI 48824,USA}
\email{wangzq@math.msu.edu}

\begin{abstract} We give a complete solution to the local classification  program of higher rank partially hyperbolic algebraic actions. We show $C^\infty$ local rigidity of abelian ergodic algebraic actions for symmetric space examples, twisted symmetric space examples and automorphisms on nilmanifolds.  The method is a combination of representation theory, harmonic analysis and a KAM iteration. A striking feature of the method is no specific information from representation theory is needed. It is the first time local rigidity for non-accessible partially hyperbolic actions has ever been obtained  other than torus examples.  Even for Anosov actions, our results are new: it is the first time twisted spaces with non-abelian nilradical have been treated in the literature.
\end{abstract}


\maketitle

\setcounter{tocdepth}{2}

\tableofcontents

\section{Introduction}
\subsection{Abelian algebraic actions}\label{sec:11} Let $\GG$ be a
connected Lie group, $\Upsilon$ a (cocompact) torsion free lattice in
$G$, and $A\subseteq \GG$ a closed abelian subgroup which
is isomorphic to $\ZZ^k\times \RR^l$. Let $L$ be a compact subgroup of the centralizer $Z(A)$ of
$A$. Then $A$ acts by left translation on the compact space
$\mathcal{X}=L\backslash\GG/\Upsilon$, which is called \emph{an algebraic $A$-action} and is denoted by $\alpha_A$. $\alpha_A$ is higher-rank if $k+\ell\geq2$. The linear part $d\alpha_A$ of $\alpha_A$ is induced by the adjoint representation of $A$ on $\mathfrak{G}=\text{Lie}(\GG)$. Define the Lyapunov exponents of $\alpha_A$ as the log's of the absolute values of the eigenvalues of $d\alpha_A$. We get linear functionals $\chi:\,A\to\RR$, which  are called \emph{Lyapunov functionals} of $d\alpha_A$.

\begin{itemize}
  \item $\alpha_A$  is \emph{hyperbolic} if $0$-Lyapunov functionals of $d\alpha_A$ only appear inside the orbit distribution.

  \item $\alpha_A$  is \emph{partially hyperbolic} if $d\alpha_A$ has non-zero Lyapunov functionals, and the $0$-Lyapunov distributions appear outside the orbit distribution. $\alpha_A$ is \emph{higher-rank} partially hyperbolic if $d\alpha_A$ has at least two linearly indepent Lyapunov functionals.

  \item $\alpha_A$ is \emph{accessible} if hyperbolic directions of various elements generate the whole tangent space;
  $\alpha_A$ is \emph{genuine}  higher rank if  the Lyapunov functionals restricted on each simple factors of $G$ has rank $\geq2$; $\alpha_A$ has \emph{strong} hyperbolicity
  if $0$-Lyapunov functionals of $d\alpha_A$ only appear inside the orbit distribution of the full Cartan action. We note that
  genuine higher rank and strong hyperbolic actions are accessible; strong hyperbolic actions are also semisimple.

\end{itemize}

There are three main classes of algebraic abelian
actions described by Katok and Spatzier in their seminal work \cite{Spatzier1}:

\begin{enumerate}
  \item for the symmetric space examples take $\GG$ a semisimple Lie group of the non-compact
type;

\smallskip
  \item for the twisted symmetric space examples take $\GG=G\ltimes_\rho V$, a semidirect
product of a semisimple Lie group of the non-compact
type with a simply connected nilpotent group $V$. The
multiplication of elements in $G_\rho$ is given by
\begin{align}\label{for:12}
(g_1,x_1)\cdot(g_2,x_2)=(g_1g_2,\rho(g_2^{-1})(x_1)x_2).
\end{align}
\item for suspensions of partially hyperbolic automorphisms of nilmanifolds take $\GG=\RR^m\ltimes \RR^N$ or $\GG=\RR^m\ltimes V$, where $V$ is a simply connected nilpotent Lie
group.
\end{enumerate}

\subsection{Rigidity of actions and related notions}\label{sec:59}
 We say that $\alpha_A$ is $C^{p,\ell,r}$ {\em locally
rigid} if for any $C^p$ perturbation $\ZZ^k\times \RR^l$-action $\tilde{\alpha}$ which is
sufficiently $C^{\ell}$ close to $\alpha_{A}$, there is a closed abelian subgroup $A'\subseteq Z(L)$ which
is isomorphic to $\ZZ^k\times \RR^\ell$ and $h\in \emph{\small Diff}^{\,\,r}(\mathcal{X})$ such that for any $x\in \mathcal{X}$ and $a\in A$ we have
\begin{align}\label{for:195}
 h(\tilde{\alpha}_a(x))=\alpha_{i_0(a)}(h(x)),
\end{align}
where $i_0:A\to A'$ is an isomorphism.

If $A'=A$ and $i_0=I_{id}$ we say that $\alpha_A$ is {\em strongly locally
rigid}.

\begin{definition}
An action $\alpha'$ of $\ZZ^k\times \RR^l$ on a manifold $\mathcal{X}'$ is a factor of $\alpha_A$ if there exists an
epimorphism $p:\mathcal{X}\rightarrow \mathcal{X}'$ such that $p\circ\alpha_A=\alpha'\circ p$.

An action $\alpha'$ is a rank one factor if it is an factor and if $\alpha'(\ZZ^k\times \RR^l)$
contains a rank one subgroup of finite index.
\end{definition}

\subsection{History and motivation} (Partially) hyperbolic actions have rich and interesting properties, finding
applications and perspectives in dynamical systems, number theory, PDE and geometry. Motivated by the Zimmer program, the study of smooth local rigidity  of higher rank (partially) hyperbolic actions has become one of the most active areas. The main goal of local classification is completely understanding the dynamics of smooth actions which are small perturbations of the given
action.

After the seminal work of Katok and Spatzier on Anosov actions \cite{Spatzier},  partially hyperbolic actions have since been extensively studied
over the past decades; some of the highlights are
\cite{Damjanovic3}, \cite{Damjanovic4}, \cite{Damjanovic2},  \cite{Zhenqi1}, \cite{zwang-2}, \cite{vinhage15}, \cite{vw}.
Among them, the study of genuine higher rank actions (see \eqref{sec:11}) has had the most progress. Currently, two different approaches
have been developed: geometric approach and KAM approach.

The first one was by Damjanovic and Katok based on a study of symmetric space exmaples \cite{Damjanovic1,Damjanovic2}.
They considered actions on $G/\Gamma$, where $G=SL(n,k)$, $n\geq 3$, $k=\RR$ or $\CC$. Damjanovic, Vinhage and the author \cite{Damjanovic3}, \cite{vw} generalized the results to a broad class of genuine higher rank actions.
All of these works rely on careful analysis of geometric obstructions to a H\"{o}lder conjugacy. Proving vanishing of the obstructions needs
algebraic information extracted from generating relations of the group \cite{Deodhar}, \cite{Matsumoto}.
The smoothness of the conjugacy then follows from a priori regularity method using the accessibility condition as in \cite{Spatzier2}, \cite{Spatzier}.
Since non-genuine actions allow additional geometric obstructions that can not be trivialized by algebraic data, applications of the geometric approach are not straightforward. As a consequence, no local
classification results have been obtained for broader class of partially hyperbolic systems.  Even for genuine higher rank actions,
the local rigidity problem is not complectly  resolved. The lack of precise description of
the needed algebraic information is a main difficulty for all actions on double coset spaces ($L\backslash G/\Gamma$), as well as
a large class of actions on twisted spaces.
Another issue is that the smoothing techniques \cite{Spatzier} are especially poorly
suited to non-semisimple actions. This implies that additional strong hyperbolicity assumption (see \eqref{sec:11}) is essential for $C^\infty$ rigidity.

The second one is also by Damjanovic and Katok based on a study of toral automorphisms \cite{Damjanovic4}.
They proved $C^\infty$ local rigidity for higher rank partially hyperbolic automorphisms on
torus. The approach is an adapted KAM type iteration
used by Moser in \cite{Moser}. We emphasize that the action considered in \cite{Damjanovic4} is geometrically more subtle than the genuine higher rank actions, because the accessibility condition fails, which was an essential ingredient for the geometric approach. This adapted KAM approach was further used for elliptic and parabolic systems \cite{DK-parabolic},  \cite{DT}, \cite{Petkovic}, \cite{W5}. The scheme has two essential ingredients \cite{Damjanovic5}: quantifying the error between the algebraic action and its perturbation and an application of the KAM iteration. The key step of the quantifying procedure is solving
the linearized conjugacy equation over the unperturbed action approximately by suitably smooth functions; then a $C^\infty$ conjugacy can be obtained by the KAM iteration.
Since the linearized conjugacy equation
is precisely the first cohomology group, this step relies on constructing tame inverses for both first
and second coboundary operators. Tameness means that there is a finite loss of regularity
between the solution and the given data.

Unfortunately, carrying out this scheme to partially hyperbolic actions other than torus examples seems to be a ``mission impossible", even to
automorphisms of step 1 nilmanifolds, and to $\ZZ^2$ actions on
$SL(2,\RR)\times SL(2,\RR)/\Gamma$, the most basic ones.
The following remarks may illustrate the difficulties. For partially hyperbolic actions, the quantifying procedure involves the study of  twisted cohomological equations. The twist  comes from eigenvalues of the linear part. Super-exponential mixing was an essential tool in the
study of twisted cohomology \cite{Damjanovic4}. However, it is not available to
smooth functions on nilmanifolds or homogeneous spaces $G/\Gamma$. We stress that even when a satisfactory understanding of the cohomology (untwisted case) is possible,
the obstructions to solving the coboundaries may lack explicit descriptions, in which case the construction of inverse coboundary operators may fail to be carried out.
This would result the KAM scheme may not apply. So while the obstructions to solving coboundaries (untwisted case) over algebraic partially hyperbolic actions have been completely classified \cite{Spatzier1}, \cite{W4},
it is unclear how to actually construct the inverse coboundary operators.

In \cite{Damjanovic4}, the quantifying procedure heavily relies upon Fourier analysis. However, the harmonic analysis for nilmanifolds is complicated
and not easily used.
For algebraic actions of sesmimple type, the quantifying procedure
usually needs the representation theory, which is usually hard to perform. So far, only representation theory of $SL(2,\RR)$ is applicable for the study of local rigidity,  which explains why current  results are restricted to actions generated by
elements inside various $SL(2,\RR)$ subgroups \cite{DK-parabolic}, \cite{W5}. It should be noted that one may need to carry out the quantifying procedure for different algebraic actions
other than the original one if coordinate changes appear in the iteration. This poses
new challenges to partially hyperbolic actions other than torus examples, for which coordinates change is usually inevitable in each iterative step. This implies that we are forced to carry out representation theory of higher rank semisimple
Lie groups, even if the original action has some element sitting inside  a $SL(2,\RR)$ subgroup. This makes the quantifying procedure substantially harder. In general the unitary dual of many higher rank simple Lie groups is not completely classified,
and even when the classification is known, it is too complicated to apply. It should be noted that the representation theory even for
$SL(2,\RR)\times SL(2,\RR)$ is hard to apply (whose unitary dual is  well-understood \cite{tan}), which is in many ways the next easiest for a higher rank simple Lie group (after $SL(2,\RR)$).  So as a consequence, despite the success of torus examples, in the past 20 years it has remained mysterious
whether the KAM scheme has broader applications in partially hyperbolic systems.



\subsection{Result of the paper}\label{sec:2} The aim of
this paper is to give a complete solution to the local classification program of higher rank partially hyperbolic algebraic actions.
In order to state our main theorems precisely,we need to set up some
notations. In this paper, $G$ denotes a higher-rank connected  semisimple Lie group with finite center and without compact factors and $\Gamma$ a cocompact irreducible lattice in
$G$. Let $V$ be a simply connected nilpotent group and $\Lambda$ be a cocomapct lattice of $V$.

Firstly, we show results for \emph{symmetric and twisted symmetric space examples}.

\begin{definition} Let $[V,V]$ be the commutator subgroup of $V$.  A representation $\rho: G\to GL(V )$ is called \emph{excellent} if the following are satisfied:

1. $\rho(G)$-fixed points in $V/[V,V]$ are $\{0\}$;

2. if $G_i$ is a simple factor of $G$ and is isomorphic to
$SO(n,1)$ or $SU(n,1)$, then $\rho(G_i)$-fixed points in $V/[V,V]$ are $\{0\}$.
\end{definition}
Suppose that $\Gamma$
acts on a compact nilmanifold $V/\Lambda$ such that its linearization extends to a homomorphism of $G$, which induces a representation $\rho: G\to GL(V )$.  Then $\Gamma\ltimes \Lambda$ is a lattice in $G\ltimes_\rho V$.

Set $\GG=G$ (resp. $G\ltimes_\rho V$) and $\Upsilon=\Gamma$ (resp. $\Gamma\ltimes \Lambda$). For any closed abelian subgroup $A$ of $\GG$ and $L$ a compact subgroup of the centralizer $Z(A)$ of
$A$, we use $\alpha_A$ to denote the the action of $A$ by left translations
on $\mathcal{X}=L\backslash\mathbb{G}/\Upsilon$. We say that $\GG$ is \emph{standard} if: $\GG=G$ or $\GG=G\ltimes_\rho V$, where $\rho$ is excellent.

\begin{theorem}\label{th:7} If $G$ is standard and $\alpha_A$ is a higher-rank partially hyperbolic action on $\mathcal{X}$, then there is $\ell_0\in\NN$ such that $\alpha_A$ is $C^{\infty,\ell_0,\infty}$ locally rigid.
\end{theorem}
\begin{remark} The excellence of $\rho$ is the weakest assumption put on twisted spaces so far. For the definition, the first condition is natural. It is used to guarantee ergodicity of any non-compact elements on the twisted space.
The second one is used to guarantee exponential mixing property for any partially hyperbolic elements on the twisted space.  We do not know if the property holds if we get rid of the restrictions on simple factors isomorphic to
$SO(n,1)$ or $SU(n,1)$ (this is true for symmetric space examples).
\end{remark}
The next result is for partially
hyperbolic diffeomorphisms on \emph{nilmanifolds}.

\begin{theorem}\label{th:5} Let $\alpha: A\times V/\Lambda\to V/\Lambda$, where $A=\ZZ^k$,  $k\geq2$ be a higher rank action by automorphisms of $V$. If $\alpha_A$ has no rank one factors, then there is $\ell_0\in\NN$ such that $\alpha_A$ is $C^{\infty,\ell_0,\infty}$ locally rigid.

\end{theorem}
We point out that unlike torus examples, where $\alpha_A$ is strongly locally
rigid (see Section \ref{sec:59}), for
nilmanifold examples, \eqref{for:388} of Section \ref{sec:19} shows that $\tilde{\alpha}_A$ may conjugate to an affine
  action other than $\alpha_A$.

Hence we complete the local classification of higher rank partially hyperbolic algebraic actions. It is the first time local rigidity for non-accessible actions has been obtained  other than torus examples.  Even for Anosov actions, our results are new: it is the first time twisted spaces with non-abelian nilradical have been treated in the literature.
Either the Fourier analysis method in \cite{Spatzier1} or the geometric approach in \cite{vw} does not seem to be
easily adaptable for non-abelian cases, which is why the results are all restricted to abelian nilradical.

\subsection{Proof strategy}  Our method is a combination of harmonic analysis, representation theory and a KAM iteration.
Unlike the previous work, no specific information from representation theory is needed for the new method. This enables the complete local classification.
Our method is similar to the KAM scheme in spirit insofar as the quantifying procedure is also needed. However,
it differs significantly from the framework in both main ingredients:
study of the (twisted)-cohomology and construction of a $C^\infty$ approximate
solution.

We use an approach that is quite different from existing work in the study of
both the untwisted and twisted cohomology.  The innovation in our approach is using
the system of imprimitivity to obtain super-exponential mixing, derived from the deformation of projection-valued measure under the natural dual action of partially hyperbolic elements. This idea is motivated by works \cite{tan} and \cite{Zhenqi2}, where the deformation of $K$-orbits under the dual action of split Cartan
was used to obtain upperbounds of matrix coefficients of simidirect products.

The super-exponential mixing allows us to get complete (twisted) cohomological information. However, it is interpreted in an abstract way.
So while the study is complete, little information is obtained in the
``actual" construction of inverse coboundary
operators. Explicit interpretation usually needs the full machinery of representation theory, which is usually hard to carry out.
This is the main difficulty  in obtaining a $C^\infty$ approximate
solution. So we abandon the framework of trying to construct inverse (twisted)-coboundary
operators. We instead use a quite different approach. A basic observation is: for any (twisted) cohomological equation, exponential mixing gives us a distributional solution. Then we construct a $C^\infty$ ``sufficiently tame" (see Remark \ref{re:3}) approximation of the distributional solution, which is exactly a $C^\infty$ approximate
solution. However,  this construction seems  quite unreasonable both from the general  functional analysis point of view and from the KAM scheme point of view. For the former, the (twisted) coboundary operator has infinite-dimensional
co-kernel and is without a spectral gap. For the latter,  the construction challenges essential ingredients of the framework.  Firstly, even if no error involves, the approximation  is not the
real solution. Thus it is unlikely to define inverse operators.  Secondly, the construction is not preserved by the group action, which means that
we can not make use of the higher rank condition. Both result that we face significant
technique issues in obtaining suitable estimates to ensure convergence of the KAM iterative
process. While our starting point is from  the abstract interpretation of invariant distributions,
the information we derive from it is quite distinct from that is
previously considered. The key innovation here is: we can precisely measure how far the distributional solution is away from a $C^\infty$ one
by using invariant distributions and the system of imprimitivity.
So we have enough
freedom in choosing the approximation in such a way that the estimates are sufficiently nice for the iteration.

Sections \ref{sec:13} and \ref{sec:14}, how the invariant distributions  are constructed and ``sufficiently tame"  estimates are obtained, Sections \ref{sec:20} and \ref{sec:7}, how the approximated solution is constructed,  are the heart and technical core of the present paper.



\subsection{Outline of this paper} We begin with the proof of automorphisms on nilmanifolds in Section \ref{sec:23} for the sake of transparency of ideas and exposition.
With proper modifications we extend the strategy to symmetric space and twisted space examples  in Section \ref{sec:24}.
In Section \ref{sec:13} and \ref{sec:14} we make a detailed study of the (twisted) cohomological  equation: we classify the obstructions to solving the equation, obtain
tame estimates for the solution, as well as sufficiently tame estimates for the obstructions. In Section \ref{sec:25} and \ref{sec:14} we study the almost twisted cocycle equations. We show that the norms of the invariant distributions are comparable to norms of the error. In Sections \ref{sec:20} and \ref{sec:7} we construct the approximated solution and obtain sufficiently tame estimates. Finally, in Sections \ref{sec:26} to \ref{sec:28} and Sections \ref{sec:29} to \ref{sec:30} we set up the iterative step, show convergence of the KAM iteration and thus Theorem \ref{th:5} and \ref{th:7} are proved respectively.

\section{Preliminaries on unitary representation theory}\label{sec:15}

\subsection{Sobolev space and elliptic regularity theorem}\label{sec:17} Let $\pi$ be a unitary representation of a Lie group $S$ with Lie algebra $\mathfrak{S}$ on a
Hilbert space $\mathcal{H}=\mathcal{H}(\pi)$. Fix an inner product $|\cdot|$ on $\mathfrak{S}=\text{Lie}(S)$.
\begin{definition}\label{de;1}
For $k\in\NN$, $\mathcal{H}^k(\pi)$ consists of all $v\in\mathcal{H}(\pi)$ such that the
$\mathcal{H}$-valued function $g\rightarrow \pi(g)v$ is of class $C^k$ ($\mathcal{H}^0=\mathcal{H}$). For $X\in\mathfrak{S}$, $d\pi(X)$ denotes the infinitesimal generator of the
one-parameter group of operators $t\rightarrow \pi(\exp tX)$, which acts on $\mathcal{H}$ as an essentially skew-adjoint operator. For any $v\in\mathcal{H}$, we also write $Xv:=d\pi(X)v$.
\end{definition}
We shall call $\mathcal{H}^k=\mathcal{H}^k(\pi)$ the space of $k$-times differentiable vectors for $\pi$ or the \emph{Sobolev space} of order $k$. The
following basic properties of these spaces can be found, e.g., in \cite{Nelson} and \cite{Goodman}:
\begin{enumerate}
  \item $\mathcal{H}^k=\bigcap_{m\leq k}D(d\pi(Y_{j_1})\cdots d\pi(Y_{j_m}))$, where $\{Y_j\}$ is a basis for $\mathfrak{S}$, and $D(T)$
denotes the domain of an operator on $\mathcal{H}$.

  \item $\mathcal{H}^k$ is a Hilbert space, relative to the inner product
  \begin{align*}
    \langle v_1,\,v_2\rangle_{S,k}:&=\sum_{1\leq m\leq k}\langle Y_{j_1}\cdots Y_{j_m}v_1,\,Y_{j_1}\cdots Y_{j_m}v_2\rangle+\langle v_1,\,v_2\rangle
  \end{align*}
  \item The spaces $\mathcal{H}^k$ coincide with the completion of the
subspace $\mathcal{H}^\infty\subset\mathcal{H}$ of \emph{infinitely differentiable} vectors with respect to the norm
\begin{align*}
    \norm{v}_{S,k}=\bigl\{\norm{v}^2+\sum_{1\leq m\leq k}\norm{Y_{j_1}\cdots Y_{j_m}v}^2\bigl\}^{\frac{1}{2}}.
  \end{align*}
induced by the inner product in $(2)$. The subspace $\mathcal{H}^\infty$
coincides with the intersection of the spaces $\mathcal{H}^k$ for all $k\geq 0$.

\item $\mathcal{H}^{-k}$, defined as the Hilbert space duals of
the spaces $\mathcal{H}^{k}$, are subspaces of the space $\mathcal{E}(\mathcal{H})$ of distributions, defined as the
dual space of $\mathcal{H}^\infty$.
  \end{enumerate}
We write $\norm{v}_{k}:=\norm{v}_{S,k}$ and $ \langle v_1,\,v_2\rangle_{k}:= \langle v_1,\,v_2\rangle_{S,k}$ if there is no confusion. Otherwise,
we use subscripts to emphasize that the regularity is measured with respect to $S$. If we want to consider the restricted representation on a subgroup $H$ of $S$
we use $\mathcal{H}^k_{H}$ to denote the Sobolev space of order $k$ with respect to $H$. $\norm{v}_{H,k}$ and $\langle v_1,\,v_2\rangle_{H,k}$ are defined accordingly.

We list the well-known elliptic regularity theorem which will be frequently
used in this paper (see \cite[Chapter I, Corollary 6.5 and 6.6]{Robinson}):
\begin{theorem}\label{th:4}
Fix a basis $\{Y_j\}$ for $\mathfrak{S}$ and set $L_{2m}=\sum Y_j^{2m}$, $m\in\NN$. Then
\begin{align*}
    \norm{v}_{2m}\leq C_m(\norm{L_{2m}v}+\norm{v}),\qquad \forall\, m\in\NN
\end{align*}
where $C_m$ is a constant only dependent on $m$ and $\{Y_j\}$.
\end{theorem}

We end this part with a technical result, which shows that if a distribution is sufficiently smooth, then it is in fact a vector of high Sobolev regularity.
\begin{lemma}\label{le:3} Fix a basis $\{Y_j\}$ for $\mathfrak{S}$ and set $L_{m}=\sum Y_j^{m}$, $m\in\NN$. Suppose $0<r<4m$.  If $D\in \mathcal{E}(\mathcal{H})$ satisfying
\begin{align*}
 \norm{(L_{4m}D)\xi}\leq E_{4m}\norm{\xi}_r,\quad\text{and}\quad \norm{D(\xi)}\leq E_0\norm{\xi}_r
\end{align*}
for any $\xi\in \mathcal{H}^r$, then $D\in \mathcal{O}^{4m-r}$ with estimates
\begin{align*}
 \norm{D}_{4m-r}\leq C_{m,r}(E_{4m}+E_0).
\end{align*}
\end{lemma}
\begin{proof} For any $\xi\in \mathcal{H}^\infty$ we have
\begin{align}\label{for:29}
 \langle L_{4m}\xi,\xi\rangle= \langle L_{2m}\xi,L_{2m}\xi\rangle=\norm{L_{2m}\xi}^2\overset{\text{(1)}}{\geq} \frac{1}{2C_m^2}\norm{\xi}^2_{2m}-\norm{\xi}^2
\end{align}
Here in $(1)$ we use Theorem \ref{th:4}. This implies that the operator $L_{4m}$ satisfies the $\pi$-coercive condition in \cite{Goodman1} (see definition 3.1). Then by
Theorem $3.1$ of \cite{Goodman1}, $D\in \mathcal{H}^{4m-r}$. Lemma $3.1$ of \cite{Goodman1} shows that there exist positive constants $c_{m}$ and $c_{r,m}$ such that \begin{align*}
  \norm{(L^{4m}+c_{r,m})D}_{-r}\geq c_m\norm{D}_{4m-r}.
\end{align*}
We note that the left side of the above inequality is bounded by $E_{4m}+c_{r,m} E_0$.   This implies the result.

\end{proof}

\subsection{Useful results}
We review several important results which will serve as ready
references later. Suppose $G$ denotes a semisimple Lie group of non-compact type with finite center and $\Gamma$ is an irreducible lattice of $G$. The following result is quoted from \cite{Mieczkowski1}, which is derived from \cite{Clozel}, \cite{Kleinbock1} and \cite{Shalom}.
\begin{theorem}\label{th:10} Suppose $G=G_1\times \cdots\times G_k$ where $G_i$, $1\leq i\leq k$ is a simple factor of $G$ and $H$ is a nonamenable closed subgroup of $G$. Then the restriction of $L^2(G/\Gamma)$ to $H$ has a spectral gap (outside a fixed neighborhood of the trivial
representation of $H$ in the Fell topology).
\end{theorem}
Kleinbock and Margulis (see \cite[appendix]{Kleinbock}) extended the matrix coefficient decay result about smooth vectors in \cite{Spatzier1}
to Lipschitz vectors. Fix a maximal compact subgroup $\mathcal{K}$ of $G$ and a Riemannian metric $d$ on $G$ which is bi-invariant with respect to $\mathcal{K}$.
\begin{theorem}\label{th:2} Let $(\pi,\,\mathcal{H})$ be a unitary representation of
$G$  such that the restriction of $\pi$ to any
simple factor of $G$ is isolated from the trivial representation. Then there exist constant $C,\,\sigma>0$  dependent only on $G$ and the spectral gap of each simple factor of $G$ such that if $v_i\in \mathcal{H}^1$, $i=1,2$, then for any $g\in G$
\begin{align*}
\abs{\langle\pi(g)v_1,v_2\rangle}&\leq
C\norm{v_1}_1\norm{v_2}_1e^{-\sigma d(e,g)}.
\end{align*}
\end{theorem}
Next, we formulate an important fact, which is derived from Theorem \ref{th:2}.  The proof is left for Appendix \ref{sec:8}. We assume notations in Section \ref{sec:2}.
\begin{theorem}\label{th:3}
Fix a spit Caran subalgebra $\mathcal{C}$ of $G$. Denote by $\vartheta$ the
sum of roots on $\mathcal{C}$.  There exist constants $0<\gamma<1$ and $E>0$, dependent only on $\GG$, $\Upsilon$ and $\rho$ such that for any
$f_i\in \mathcal{O}^1$, $i=1,2$ and  any $Y\in \mathcal{C}$, we have
\begin{align*}
\abs{\langle f_1(\exp(tY)),f_2\rangle}&\leq
E\norm{f_1}_{1}\norm{f_2}_{1}e^{-|t|\gamma\vartheta(Y)},\qquad \forall\,t\in\RR.
\end{align*}
\end{theorem}
Gorodnik and Spatizer (see \cite{Gorodnik}) obtain the following exponential mixing result for commuting groups of
ergodic automorphisms of nilmanifolds:
\begin{theorem}\label{th:1}
Suppose $\alpha_0$ is a $\ZZ^k$ action on a nilmanifold $V/\Lambda$ by nil-automorphisms such
the all nontrivial elements of $\ZZ^k$ act ergodically wrt Haar measure. Then there exists a constant
$\gamma$ only dependent on $\alpha_0$ such that for any $f,\,g\in C^1(V/\Lambda)$ orthogonal to constants and any $z\in\ZZ^k$ we have
\begin{align*}
\abs{\langle f(\alpha_0(z)x), g\rangle}&\leq
E\norm{f}_{C^1}\norm{g}_{C^1}e^{-\gamma\norm{z}}.
\end{align*}
\end{theorem}

\begin{remark} Theorem \ref{th:3} and \ref{th:1} give exponential mixing results for all partially hyperbolic algebraic actions. It will
allows us to define converging series in a suitable space of distributions,  which will be the basis for
the quantifying procedure.
\end{remark}

For the classical horocycle flow defined by the $\mathfrak{sl}(2,\RR)$-matrix $U=\begin{pmatrix}
  0 & 1 \\
  0 & 0
\end{pmatrix}$, Flaminio and Forni established the study to the solution of the cohomological
equation \cite{Forni}. Let $(\pi,\,\mathcal{H})$ be a unitary representation of
$S$. Suppose $\text{Lie}(S)=\mathfrak{sl}(2,\RR)$.

\begin{theorem}\label{th:8} Suppose $\pi$ has a spectral gap of $u_0$. For any $\omega\in \mathcal{H}^\infty$, if the equation $U\vartheta=\omega$ has a solution $\vartheta\in \mathcal{H}^{\infty}$, then
\begin{align*}
 \norm{\vartheta}_m\leq C_{m,u_0}\norm{\omega}_{m+2},\qquad \forall\,m\geq0.
\end{align*}
\end{theorem}

\section{Directional smoothing operators}\label{sec:5}

Here we summarize the results from \cite{W1}, \cite{Zhenqi2} and \cite{W5}. In this part we show a general construction of  smoothing operators by using group algebra. These smoothing operators play a crucial role in  constructing of $C^\infty$ approximations in Section \ref{sec:20}.
Even though all results in this section
can be found in \cite{W5} we include the proofs in Appendix \ref{sec:35} to \ref{sec:34}, for completeness.
If readers are not interested in these proofs, they can skip them without affecting their understanding of this section.

\subsection{Truncation by Fourier transform}\label{sec:37} Before we introduce the smoothing operators, it is worth reviewing the truncation by Fourier transform, which serves as the prototypical foundation for understanding these operators.

We denote by $W^{2,p}(\RR)$ the Sobolev space of $L^2$
functions with $L^2$ weak derivatives up to order $p$. Fix a bump function $f$. We define smoothing operators $\pi(f\circ \tau^{-1})$, $\tau>0$ on $W^{2,0}(\RR)=L^2(\RR)$ as follows:
\begin{align}\label{for:62}
 \pi(f\circ \tau^{-1})(g)(x)=\text{\tiny$\frac{1}{\sqrt{2\pi}}$}\int_{\RR} f(\text{\tiny$\frac{\chi}{\tau}$})\hat{g}(\chi)e^{\textrm{i}\chi\cdot x} d\chi.
\end{align}
where $\hat{g}(\chi)=\text{\tiny$\frac{1}{\sqrt{2\pi}}$}\int_{\RR} g(x)e^{-\textrm{i}\chi\cdot x} dx$.

 It is easy to check that the following properties hold:
\begin{enumerate}
  \item\label{for:45} $ \pi(f_1\circ \tau^{-1})\pi(f_2\circ \tau^{-1})=\pi\big((f_1f_2)\circ \tau^{-1} \big)$;
  \smallskip
  \item $\big\langle \pi(f\circ \tau^{-1})(g),\,g_1\big\rangle=\big\langle g,\,\pi(\bar{f}\circ \tau^{-1})(g_1)\big\rangle$, where $\bar{f}$ is the complex conjugate of $f$;

  \smallskip
  \item\label{for:2006} we have
\begin{align*}
  &\text{\tiny$\frac{\partial}{\partial x^n}$}\big(\pi(f\circ \tau^{-1})(g)(x)\big)=\text{\tiny$\frac{1}{\sqrt{2\pi}}$}\int_{\RR} f(\text{\tiny$\frac{\chi}{\tau}$})(\chi\textrm{i})^{n}\hat{g}(\chi)e^{\textrm{i}\chi\cdot x} d\chi\notag\\
  &=\tau^{n}\pi(f^*\circ \tau^{-1})(g)(x)
\end{align*}
where $f_n^*(x)=f(x)(x\textrm{i})^{n}$;

\smallskip

\item\label{for:40} for any $g\in W^{2,p}$, we have $\pi(f\circ \tau^{-1})(g)\in W^{2,\infty}$; and the following estimates hold:
\begin{align*}
 \norm{\pi(f\circ \tau^{-1})(g)}_{W^{2,q}}&\leq C_{p,q} \tau^{q-p}\norm{g}_{W^{2,p}}, \qquad \forall\,q\geq p; \\
 \norm{g-\pi(f\circ \tau^{-1})(g)}_{W^{2,k}}&\leq C_{k,p} \tau^{k-p}\norm{g}_{W^{2,p}}\qquad \forall\,k\leq p.
\end{align*}
\end{enumerate}

\subsection{Notations}\label{sec:36} Throughout this section, we fix a Lie group $S$ and set $\mathfrak{S}=\text{Lie}(S)$. Let $(\pi,\mathcal{H})$ be a unitary representation of $S$.

\begin{enumerate}
  \item We say that a vector $u\in \mathfrak{S}$ is \emph{nilpotent} if $\text{ad}_u$ is
nilpotent. Suppose $\mathfrak{u}\in \mathfrak{S}$ is nilpotent. We use $\{\mathfrak{u}\}$ to denote the one parameter subgroup of $\mathfrak{u}$.

\smallskip
  \item\label{de:4} We say that $f\in C^\infty(\RR)$ is of $(1,2)$-type if $f\geq0$,  $f(t)=1$ for $|t|\leq 1$ and $f(t)=0$ for $|t|\geq 2$.
For any $\tau>0$ denote $f(\frac{t}{\tau})$ by $f\circ \tau^{-1}$.

\smallskip
  \item We use $\mathcal{S}(\RR)$ to denote the Schwartz space of $\RR$. For $f\in C^\infty(\RR)$, $f^{(n)}$ denotes the $n$-th derivative of $f$. We set
\begin{align}\label{for:239}
 \tilde{\mathcal{S}}(\RR)=\{f\in \text{\small$C^{\infty}(\RR)$}: f^{(n)}\in \text{\small$L^\infty(\RR)$}, \forall\,n\geq 0\}
\end{align}
and we define the norm as
\begin{align*}
 \norm{f}_{\tilde{\mathcal{S}}(\RR),n}=\max_{0\leq i\leq n}\{ \norm{f^{(i)}}_{\text{\small$L^\infty(\RR)$}}\},\quad \forall\,f\in \tilde{\mathcal{S}}(\RR).
\end{align*}
\end{enumerate}
\subsection{Main results}

In Section \ref{sec:27},  for any $f\in L^\infty(\RR)$ we define a linear operator $\pi_{\mathfrak{u}}(f)$ on $\mathcal{H}$ satisfying the following properties:
\begin{enumerate}
  \item For any $f_1,\,f_2\in L^\infty(\RR)$
\begin{align*}
 \pi_{\mathfrak{u}}(f_1)\pi_{\mathfrak{u}}(f_2)=\pi_{\mathfrak{u}}(f_1f_2);
\end{align*}
and
\begin{align*}
 \bigl\langle \pi_{\mathfrak{u}}(f)\xi,\eta\bigl\rangle=\bigl\langle \xi,\pi_{\mathfrak{u}}(\bar{f})\eta\bigl\rangle,\qquad \xi,\,\eta\in \mathcal{H},
\end{align*}
where $\bar{f}$ is the complex conjugate of $f$.

From the the above relations, we see that: if $X\subseteq \RR$ is a Borel set and $I_X$ denotes the characteristic function of $X$, then $\pi_{\mathfrak{u}}(I_X)$ is idempotent and self-adjoint, i.e., an orthogonal projection onto a subspace of $\mathcal{H}$. Thus the assignment $X\to \pi_{\mathfrak{u}}(I_X)$ is a \emph{projection-value measure};

\smallskip
  \item (Lemma \ref{cor:1}) Suppose  $f\in \tilde{\mathcal{S}}(\RR)$ and $\tau>0$. Then:
\begin{enumerate}
  \item If $\xi\in \mathcal{H}$. Then $\pi_{\mathfrak{u}}(f\circ \tau^{-1})\xi\in \mathcal{H}_{\{\mathfrak{u}\}}^\infty$ with estimates
\begin{align*}
 \big\|\mathfrak{u}^l\big(\pi_{\mathfrak{u}}(f\circ \tau^{-1})\xi\big)\big\|\leq C_{f,l}\tau^l\norm{\xi},\qquad l\geq0.
\end{align*}
This result  shows that the $\pi_{\mathfrak{u}}(f\circ \tau^{-1})$ operator provides smoothness along the $\mathfrak{u}$-direction. This is the reason to call it a ``directional" smoothing operator.

\smallskip
  \item If $\tau\geq1$ and $\xi\in \mathcal{H}^s$, $s\geq 0$. Then
$\pi_{\mathfrak{u}}(f\circ \tau^{-1})\xi\in \mathcal{H}^s$ with estimates
\begin{align}\label{for:491}
 \norm{\pi_{\mathfrak{u}}(f\circ \tau^{-1})\xi}_l\leq C_{f,l}\norm{\xi}_l, \qquad \forall\, \,0\leq l\leq s.
\end{align}
This result shows that $\pi_{\mathfrak{u}}(f\circ \tau^{-1})$ is a global smoothing operator, i.e., $\pi_{\mathfrak{u}}(f\circ \tau^{-1})(\mathcal{H}^s)\subseteq \mathcal{H}^s$.

\end{enumerate}
  \item (Lemma \ref{cor:4}) Suppose $f$ is of $(1,2)$-type and $\xi\in \mathcal{H}_{\{\mathfrak{u}\}}^s$. If $\tau>0$, then
\begin{align}\label{for:492}
 \norm{\xi-\pi_{\mathfrak{u}}(f\circ \tau^{-1})\xi}=\norm{\pi_{\mathfrak{u}}\big((1-f)\circ \tau^{-1}\big)\xi}\leq C_{f}\tau^{-s}\norm{\xi}_{\{\mathfrak{u}\},s}.
\end{align}
This result provides the estimate for the error coming from the smoothing.

\eqref{for:491} and \eqref{for:492} show that the estimates of the directional smoothing operators
are similar
to those of the standard smoothing operators (see \eqref{sect:6} of Section \ref{sec:19}).

\smallskip

\item (Lemma \ref{le:4})  Suppose $h\in S$ and we denote by $\mathfrak{u}_1=\norm{\text{Ad}_{h}(\mathfrak{u})}^{-1}\text{Ad}_{h}(\mathfrak{u})$.
Then for any $f\in L^\infty(\RR)$ we have
\begin{align*}
 \pi(h)\pi_{\mathfrak{u}}(f)\pi(h^{-1})=\pi_{\mathfrak{u}_1}(f\circ \norm{\text{Ad}_{h}(\mathfrak{u})}).
\end{align*}
This property illustrates how the  spectrum changes under the conjugacy,  which will be invaluable to us.

\smallskip

\item (Proposition \ref{cor:6}) Suppose $Q_0$ and $Q'$ are subgroups of $S$ and $\mathfrak{S}=\text{Lie}(Q_0)+\text{Lie}(Q')$. Also suppose $Q'$ is nilpotent.  We choose a basis
$\{\mathfrak{u}_{1},\cdots, \mathfrak{u}_{n_0}\}$ of $\text{Lie}(Q')$ with a suitable order. Then for any distribution $\mathfrak{D}$ which is only smooth on $Q_0$,
\begin{align*}
 \pi_{\mathfrak{u}_{n_0}}(f\circ \tau^{-1})\pi_{\mathfrak{u}_{n_0-1}}(f\circ \tau^{-1})\cdots \pi_{\mathfrak{u}_{1}}(f\circ \tau^{-1})\mathfrak{D}\in \mathcal{H}^\infty.
\end{align*}
The proposition shows how to obtain a smooth vector from a distributional one.

\end{enumerate}

\subsection{Basic facts of group algebra}\label{sec:27} In this section, we will generalize the construction of smoothing operators by truncation using Fourier transform to a one parameter subgroup  $\{\mathfrak{u}\}$ in a general Lie group $S$, where
$\mathfrak{u}$ is nilpotent.

For any $\xi,\,\eta\in \mathcal{H}$, consider
the corresponding matrix coefficients of $\pi\mid_{\{\mathfrak{u}\}}$:
\begin{align*}
\phi_{\xi,\eta}(t)=\langle \pi(t)\xi,\,\eta\rangle,\qquad t\in \RR,
\end{align*}
where $\pi(t)=\pi(\exp(t\mathfrak{u}))$. There exists a regular Borel measure $\mu$ on $\widehat{\RR}$, called \emph{the associated measure of $\pi$ (with respect to $\widehat{\RR}$)}, such that $\xi=\int_{\widehat{\RR}}\xi_\chi d\mu(\chi)$, and
\begin{align}\label{for:84}
    \phi_{\xi,\eta}(t)&=\int_{\widehat{\RR}}\chi(t)\langle \xi_\chi,\,\eta_\chi\rangle d\mu(\chi).
\end{align}
Here $\chi(t)=e^{\textrm{i}\chi t}$ (we identify $\RR$ and $\widehat{\RR}$).

 The representation $\pi\mid_{\{\mathfrak{u}\}}$ extends to a $^\ast$-representation on $\mathcal{S}(\RR)$: for any $f(x)\in \mathcal{S}(\RR)$, $\pi_{\mathfrak{u}}(f)$ is the operator on $\mathcal{H}$ for which
\begin{align*}
\bigl\langle \pi_\mathfrak{u}(f)\xi,\eta\bigl\rangle=\text{\tiny$\frac{1}{\sqrt{2\pi}}$}\int_{\RR} \hat{f}(t)\phi_{\xi,\eta}(t) dt, \qquad \forall \xi,\,\eta\in \mathcal{H},
\end{align*}
where $\hat{f}(t)=\text{\tiny$\frac{1}{\sqrt{2\pi}}$}\int_{\RR} f(x)e^{-\textrm{i}tx} dx$.

It is easy to check that if $\mathcal{T}$ is linear operator on $\mathcal{H}$ and commutes with $\pi$ then
\begin{align}\label{for:73}
 \mathcal{T}\pi_\mathfrak{u}(f)=\pi_\mathfrak{u}(f)\mathcal{T}.
\end{align}
Computations show that
\begin{align}\label{for:93}
 \bigl\langle \pi_{\mathfrak{u}}(f)\xi,\eta\bigl\rangle=\bigl\langle \xi,\pi_{\mathfrak{u}}(\bar{f})\eta\bigl\rangle=\int_{\widehat{\RR}}f(\chi)\langle \xi_\chi,\,\eta_\chi \rangle d\mu(\chi),
\end{align}
where $\bar{f}$ is the complex conjugate of $f$. The proof is left for Appendix \ref{sec:31}.

Since
\begin{align}\label{for:215}
 \norm{\pi_{\mathfrak{u}}(f)}\leq \norm{f}_\infty
\end{align}
for any $f\in \mathcal{S}(\RR)$, we can extend $\pi_{\mathfrak{u}}$ from $\mathcal{S}(\RR)$ to $L^\infty(\RR)$ (see Appendix \ref{sec:35} for a detailed argument).

It is easy to check that the following property holds:
\begin{align}\label{for:2}
 \pi_{\mathfrak{u}}(f_1)\pi_{\mathfrak{u}}(f_2)=\pi_{\mathfrak{u}}(f_1f_2),\qquad\forall\,f_1,\,f_2\in L^\infty(\RR).
\end{align}\label{for:31}
If $f\in L^\infty(\RR)$ and $\xi\in \mathcal{H}$, then
\begin{align}\label{for:41}
  \mathfrak{u}^{k}(\pi_{\mathfrak{u}}(f)\xi)=\pi_{\mathfrak{u}}(f^*)\xi
\end{align}
where $f^*(x)=f(x)(x\textrm{i})^{k}$ for any $k\geq0$.

If $f\in \tilde{\mathcal{S}}(\RR)$ and $\xi\in \mathcal{H}^1$, for any $v\in \mathfrak{S}$ we have
\begin{align}\label{for:234}
  \pi_{\mathfrak{u}}(f)(v\xi)&=\text{\tiny$\sum_{\substack{0\leq \text{\tiny$j$}\leq\dim\mathfrak{S}-1}}$}d_{\text{\tiny$j$}}\text{ad}^{\text{\tiny$j$}}_{\mathfrak{u}}(v)\pi_{\mathfrak{u}}(f^{(j)})
  \xi,\quad\text{where }d_j\in\CC.
\end{align}
Similarly, the following relation holds:
\begin{align}
  v(\pi_{\mathfrak{u}}(f)\xi)&=\text{\tiny$\sum_{\substack{0\leq \text{\tiny$j$}\leq\dim\mathfrak{S}-1}}$}c_{\text{\tiny$j$}}\pi_{\mathfrak{u}}(f^{(j)})
  (\text{ad}^{\text{\tiny$j$}}_{\mathfrak{u}}(v)\xi),\quad\text{where }c_j\in\CC.  \label{for:98}
\end{align}
The proof is left for Appendix \ref{sec:31}. As a direct consequence, we have
\begin{align}\label{for:10}
  v(\pi_{\mathfrak{u}}(f)\xi)=\pi_{\mathfrak{u}}(f)(v\xi),\qquad\text{if }[v,\mathfrak{u}]=0.
\end{align}

\begin{lemma}\label{cor:1} Suppose  $f\in \tilde{\mathcal{S}}(\RR)$ and $\tau>0$. Then:
\begin{enumerate}
  \item\label{for:128} If $\xi\in \mathcal{H}$, then $\pi_{\mathfrak{u}}(f\circ \tau^{-1})\xi\in \mathcal{H}_{\{\mathfrak{u}\}}^\infty$ with estimates
\begin{align*}
 \big\|\mathfrak{u}^l\big(\pi_{\mathfrak{u}}(f\circ \tau^{-1})\xi\big)\big\|\leq C_{f,l}\tau^l\norm{\xi},\qquad l\geq0.
\end{align*}

\item\label{for:292} If $\xi\in \mathcal{H}^s$, $s\geq 0$, then
\begin{align*}
 \norm{\pi_{\mathfrak{u}}(f)\xi}_l\leq C_{l}\norm{f}_{\tilde{\mathcal{S}}(\RR),l\dim\mathfrak{S}}\norm{\xi}_l, \qquad \forall\, \,0\leq l\leq s.
\end{align*}

  \item\label{for:489} If $\xi\in \mathcal{H}^s$, $s\geq 0$, then $\pi_{\mathfrak{u}}(f\circ \tau^{-1})\xi\in \mathcal{H}^s$ with estimates:
  \begin{align}\label{re:7}
 \norm{\pi_{\mathfrak{u}}(f\circ \tau^{-1})\xi}_l\leq C_{f,l}\tau^{-l\dim\mathfrak{S}}\norm{\xi}_l, \qquad \forall\, \,0\leq l\leq s.
\end{align}
If $\tau\geq1$, the above estimates can be simplified as
   \begin{align}\label{for:7}
 \norm{\pi_{\mathfrak{u}}(f\circ \tau^{-1})\xi}_l\leq C_{f,l}\norm{\xi}_l, \qquad \forall\, \,0\leq l\leq s.
\end{align}

\end{enumerate}
\end{lemma}
The proof is left for Appendix \ref{sec:33}.

The following lemma provides the estimate for the error coming from the smoothing.
\begin{lemma}\label{cor:4} Suppose $f$ is of $(1,2)$-type and $\xi\in \mathcal{H}_{\{\mathfrak{u}\}}^s$. If $\tau>0$, then
\begin{align*}
 \norm{\xi-\pi_{\mathfrak{u}}(f\circ \tau^{-1})\xi}=\norm{\pi_{\mathfrak{u}}\big((1-f)\circ \tau^{-1}\big)\xi}\leq C_{f}\tau^{-s}\norm{\xi}_{\{\mathfrak{u}\},s}.
\end{align*}

\end{lemma}
The proof is left for Appendix \ref{sec:34}.

\subsection{Applications of group algebra} We list two important applications. The first is a generalization of the dual actions on projection-valued measure \cite{tan}, \cite{Zhenqi2}. It illustrates how the  spectrum changes under the conjugacy,  which will be invaluable to us:

 \begin{lemma}\label{le:4} Suppose $h\in S$ and we denote by $\mathfrak{u}_1=\norm{\text{Ad}_{h}(\mathfrak{u})}^{-1}\text{Ad}_{h}(\mathfrak{u})$.
Then for any $f\in L^\infty(\RR)$ we have
\begin{align}\label{for:130}
 \pi(h)\pi_{\mathfrak{u}}(f)\pi(h^{-1})=\pi_{\mathfrak{u}_1}(f\circ \norm{\text{Ad}_{h}(\mathfrak{u})}).
\end{align}
\end{lemma}
\begin{proof}
For any $f\in \mathcal{S}(\RR)$ we have
\begin{align*}
&\bigl\langle \pi(h)(\pi_{\mathfrak{u}}(f)\xi),\eta\bigl\rangle\overset{\text{(1)}}{=}\text{\tiny$\frac{1}{\sqrt{2\pi}}$}\int_{\RR}\hat{f}(t)\langle \pi(h)\pi(\exp(t\cdot \mathfrak{u}))\xi,\eta\rangle dt\\
&=\text{\tiny$\frac{1}{\sqrt{2\pi}}$}\int_{\RR}\hat{f}(t)\big\langle \pi\big(\exp(t\norm{\text{Ad}_{h}(\mathfrak{u})}\cdot \mathfrak{u}_1)\big)\pi(h)\xi,\eta\big\rangle dt\\
&=\text{\tiny$\frac{1}{\sqrt{2\pi}}$}\int_{\RR}\hat{f}(t\norm{\text{Ad}_{h}(\mathfrak{u})}^{-1})\langle \pi(\exp( t\cdot \mathfrak{u}_1))\pi(h)\xi,\eta\rangle \norm{\text{Ad}_{h}(\mathfrak{u})}^{-1}dt\\
&\overset{\text{(1,2)}}{=}\bigl\langle \pi_{\mathfrak{u}_1}(f\circ \norm{\text{Ad}_{h}(\mathfrak{u})})(\pi(h)\xi),\eta\bigl\rangle.
\end{align*}
In $(1)$ we use definition;  and in $(2)$ we note that
\begin{align*}
 \norm{\text{Ad}_{h}(\mathfrak{u})}^{-1}\hat{f}\circ \norm{\text{Ad}_{h}(\mathfrak{u})}^{-1}=\hat{g}
\end{align*}
where $g=f\circ \norm{\text{Ad}_{h}(\mathfrak{u})}$.

The above discussion shows that  \eqref{for:130} holds for any $f\in \mathcal{S}(\RR)$. By arguments at the end of Appendix \ref{sec:35}, we can extend \eqref{for:130} from $\mathcal{S}(\RR)$ to $L^\infty(\RR)$.
\end{proof}

The next  result shows how to obtain a smooth vector from a distributional one. From \eqref{for:128} of Lemma \ref{cor:1} we see that $\pi_{\mathfrak{u}}$ serves as a smoothing operator along direction $\mathfrak{u}$. Then for a distribution, it is natural to hope that if we apply $\pi_{\mathfrak{u}}$ successively along all its non-smooth directions, it should becomes a bona fide smooth vector. This trick was also used in \cite{W5} for a special case when $Q'$ is abelian. The proof is left for Appendix \ref{sec:18}. Let $\mathfrak{S}^1$ denote the set of the unit vectors in $\mathfrak{S}$.

\begin{proposition}\label{cor:6} Suppose $Q_0$ is a subgroup of $S$ and $\mathfrak{u}_i\in \mathfrak{S}$ is nilpotent, $1\leq i\leq n_0$. Let $Q_i=\{\mathfrak{u}_{i}\}$ for each $i$. For any $\mathfrak{D}\in \mathcal{E}(\mathcal{H})$ and $\tau\geq1$, define
\begin{align*}
 \mathcal{P}_{f,\tau}\stackrel{\text{def}}{=}\pi_{\mathfrak{u}_{1}}(f\circ \tau^{-1})\pi_{\mathfrak{u}_{2}}(f\circ \tau^{-1})\cdots \pi_{\mathfrak{u}_{n_0}}(f\circ \tau^{-1})
\end{align*}
and
\begin{align*}
 (\mathcal{P}_{f,\tau}\mathfrak{D})(\omega)\stackrel{\text{def}}{=}\mathfrak{D}\big(\pi_{\mathfrak{u}_{1}}(f\circ \tau^{-1})\pi_{\mathfrak{u}_{2}}(f\circ \tau^{-1})\cdots \pi_{\mathfrak{u}_{n_0}}(f\circ \tau^{-1})\omega \big),\quad \forall\,\omega\in \mathcal{H}^\infty.
\end{align*}
Denote by $Q'$ the subgroup generated by $Q_i$, $1\leq i\leq n_0$.  Fix $f$ which is of $(1,2)$ type. Suppose the following are satisfied:
\begin{enumerate}
  \item $\text{Lie}(Q')$ is linearly spanned by  $\mathfrak{u}_{i}$, $1\leq i\leq n_0$;

 \smallskip
  \item\label{for:72}  $\mathfrak{S}$ is linearly spanned by $\text{Lie}(Q_0)$ and $\text{Lie}(Q')$;

  \smallskip
  \item\label{for:71} for any $1\leq i,\,j\leq n_0$,  if $[\mathfrak{u}_{i},\mathfrak{u}_{j}]\neq0$, then $[\mathfrak{u}_{i},\mathfrak{u}_{j}]$ is inside the subspace linearly spanned by  $\mathfrak{u}_{k}$, $k>\max\{i,j\}$.

  \smallskip
  \item there are $r_0,\,s\in\NN$ with $s\geq 4+r_0$ such that
\begin{align}\label{for:75}
 \norm{v_1\cdots v_m\mathfrak{D}}_{-r_0}\leq E_m,\qquad \forall\,\, 0\leq m\leq s
\end{align}
for any $v_i\in \text{Lie}(Q_0)\bigcap \mathfrak{S}^1$, $1\leq i\leq m$.
\end{enumerate}
Then $\mathcal{P}_{f,\tau}\mathfrak{D}\in \mathcal{H}^{s-4-r_0}$ with estimates
\begin{align}\label{for:70}
 \norm{\mathcal{P}_{f,\tau}\mathfrak{D}}_m\leq C_{m,r_0}(\tau ^{m+4+r_0}E_0+E_{m+4+r_0}),
\end{align}
for any $0\leq m\leq s-4-r_0$.

\end{proposition}

\section{Proof of Theorem \ref{th:5}}\label{sec:23}

\subsection{Basic notations}\label{sec:19} We will use these notations throughout Section \ref{sec:23}. So the reader
should consult this section  if an unfamiliar symbol appears.

In what follows, $C$ will denote any constant that depends only
  on the given  action  $\alpha_A$ and the manifolds $\mathcal{X}=V/\Gamma$. $C_{x,y,z,\cdots}$ will denote any constant that in addition to the
above depends also on parameters $x, y, z,\cdots$.

\begin{sect}\label{sect:5}

We use $dx$ to denote the Harr measure on $\mathcal{X}=V/\Gamma$.
Set $L^2_0(\mathcal{X})$ to be the subspace of $L^2(\mathcal{X})$ orthogonal to constants.
We use $(\pi, \mathcal{O})$ to denote the regular representation of $L^2_0(\mathcal{X})$.
Let $\mathcal{V}=\text{Lie}(V)$ and $\mathcal{V}^1$ be the set of unit vectors in $\mathcal{V}$.
Denote by $C^r(\mathcal{X},\,\mathcal{V})$ the set of $C^r$ maps from $\mathcal{X}$ to $\mathcal{V}$.

\end{sect}

\begin{sect}\label{sect:1}
Suppose $\mathfrak{F}$ is a subspace of $\mathcal{V}$. Fix a basis $\{v_1,v_2$, $\cdots, v_{\dim\mathfrak{F}}\}$ of $\mathfrak{F}$.
$\mathfrak{F}(\mathcal{O}^s)$ denotes the set of maps from $\mathcal{X}$ to
$\mathfrak{F}$ with coordinate functions in $\mathcal{O}^s$ (see Section \ref{sec:17}).
For any  map $\mathcal{F}\in \mathfrak{F}(\mathcal{O}^s)$ with coordinate functions $f_i$, define:
\begin{align*}
 \norm{\mathcal{F}}_s:=\max_{i}\norm{f_i}_s\quad\text{and}\quad \int_{\mathcal{X}} \mathcal{F}=(\int_\mathcal{X}
f_1dx,\cdots,\int_\mathcal{X} f_{\dim \mathfrak{F}}dx).
\end{align*}
Then the regular representation $(\pi, \mathcal{O})$ naturally extends to the representation $(\pi, \mathfrak{F}(\mathcal{O}))$. $\norm{\mathcal{F}}_{C^s}$ is defined similarly.

For any distribution $\mathcal{D}=(\mathcal{D}_1,\cdots,\mathcal{D}_{\dim \mathfrak{F}})\in \mathfrak{F}(\mathcal{O}^{-s})$, $s>0$, and any
$\mathcal{F}=(f_1,\cdots,f_{\dim \mathfrak{F}})\in \mathfrak{F}(\mathcal{O}^s)$ define
\begin{align*}
 \langle \mathcal{D}, \mathcal{F}\rangle:=\big(\mathcal{D}_1(f_1),\mathcal{D}_2(f_2),\cdots, \mathcal{D}_{\dim \mathfrak{F}}(f_{\dim \mathfrak{F}}) \big)\in \mathfrak{F}.
\end{align*}
\end{sect}

\begin{sect}\label{sect:6}
There exists a collection of smoothing operators $\mathfrak{s}_b:C^\infty(\mathcal{X},\,\mathfrak{F})\to C^\infty(\mathcal{X},\,\mathfrak{F})$, $b>0$, such that for any $s, s_1,s_2\geq0$, the following holds:
\begin{align}
 \norm{\mathfrak{s}_b\mathcal{F}}_{C^{s+s_1}}&\leq C_{s,s_1}b^{s_1}\norm{\mathcal{F}}_{C^{s}},\quad \text{and}\label{for:45}\\
 \norm{(I-\mathfrak{s}_b)\mathcal{F}}_{C^{s-s_2}}&\leq C_{s,s_2}b^{-s_2}\norm{\mathcal{F}}_{C^{s}},\quad \text{if }s\geq s_2,\label{for:37}
\end{align}
see \cite{Hamilton}.

The next result follows directly from Sobolev embedding theorem on compact manifolds.  For any $\mathcal{F}\in \mathfrak{F}(\mathcal{O}^\infty)$ and $s\geq0$ the following hold:
\begin{align}\label{for:179}
  \norm{\mathcal{F}}_{s}\leq C_s\norm{\mathcal{F}}_{C^{s}},\quad \norm{\mathcal{F}}_{C^{s}}\leq C_{s}\norm{\mathcal{F}}_{s+\beta},
\end{align}
where $\beta>0$ is a constant only dependent on $\mathcal{X}$.
\end{sect}

\begin{sect}\label{sect:7}
 An  action $\alpha_A$ by nil-automorphisms has no rank one factors
if $\alpha_A(\ZZ^k)$ contains a subgroup isomorphic to $\ZZ^2$ such that all nontrivial
elements in this subgroup act ergodically wrt $dx$. Denote by
\begin{align*}
 ax:=\alpha_a(x)\quad\text{and}\quad da:=d\alpha_a,\quad \forall\,a\in \ZZ^k
\end{align*}
for any $x\in \mathcal{X}$, where $d\alpha$ is the induced action of $\alpha$ on $\mathcal{V}$.

We have the Lyapunov subspace decomposition:
\begin{align}\label{for:349}
 \mathcal{V}= \sum_{\lambda\in \Delta}\mathfrak{g}_{\lambda}
\end{align}
such that for any $\lambda\in \Delta$, $\lambda(dz)$ is the absolute value of the eigenvalue of $dz$ on $\mathfrak{g}_\lambda$, for any $z\in \ZZ^k$.

We consider two ergodic generators: $\textbf{a}_1$ and $\textbf{a}_2$ such that
\begin{enumerate}
  \item all nontrivial
elements in the $\ZZ^2$ subgroup generated by $\textbf{a}_1$ and $\textbf{a}_2$ act ergodically;

\smallskip
  \item for any $1\neq \lambda\in \Delta$, $\lambda(d\textbf{a}_i)\neq1$, $i=1,\,2$.
\end{enumerate}
It is clear that
\begin{align}\label{for:168}
 [\mathfrak{g}_{\lambda},\,\mathfrak{g}_{\mu}]\subseteq \mathfrak{g}_{\lambda\mu},\qquad \forall\,\lambda,\,\mu\in \Delta.
\end{align}
We have the following estimates: for any $\lambda\in \Delta$
\begin{align}\label{for:334}
 \norm{d\textbf{a}_i^n|_{\mathfrak{g}_{\lambda}}}\leq C\lambda(d\textbf{a}_i)^n(1+|n|)^{\dim\mathcal{V}},\quad \forall\,\,n\in\ZZ,\,\,i=1,2.
\end{align}
Set
\begin{gather*}
 \Delta_{\textbf{a}_1}^{+}=\{\mu\in\Delta:\mu(d\textbf{a}_1)>1\}, \qquad\Delta_{\textbf{a}_1}^{-}=\{\mu\in\Delta:\mu(d\textbf{a}_1)<1\},\\
 \Delta_{\textbf{a}_1}^{0}=\{\mu\in\Delta:\mu(d\textbf{a}_1)=1\}.
\end{gather*}
We also set
\begin{align*}
   \mathfrak{g}^+_{\textbf{a}_1}=\sum_{\chi\in \Delta_{\textbf{a}_1}^{+}}\mathfrak{g}_{\chi},\quad \mathfrak{g}^-_{\textbf{a}_1}=\sum_{\chi\in \Delta_{\textbf{a}_1}^{-}}\mathfrak{g}_{\chi},\quad \mathfrak{g}^0_{\textbf{a}_1}=\sum_{\chi\in \Delta_{\textbf{a}_1}^{0}}\mathfrak{g}_{\chi}.
\end{align*}
It is clear that $\mathfrak{g}^0_{\textbf{a}_1}=\mathfrak{g}_1$.

We fix a compact generating set $E\subset \ZZ^k\times \RR^l$ that contains $\textbf{a}_1$ and $\textbf{a}_2$.

\end{sect}


\begin{sect}\label{sect:3}

We say that $0\neq v\in \mathfrak{g}_{\mu}$ is a \emph{stable} (resp. \emph{unstable} or \emph{neutral}) direction of $\textbf{a}_1$ if $\mu\in \Delta_{\textbf{a}_1}^{-}$ (resp. $\mu\in \Delta_{\textbf{a}_1}^{+}$ or $\mu\in \Delta_{\textbf{a}_1}^{0}$).

If $\mu\in \Delta_{\textbf{a}_1}^{+}\cup\Delta_{\textbf{a}_1}^{0}$, let $L_\mu=\ZZ^-\cup\{0\}$ and $\delta(\mu)=-$. If $\mu\in \Delta_{\textbf{a}_1}^{-}$, let $L_\mu=\ZZ^+\cup\{0\}$ and $\delta(\mu)=+$.

Let  $L_\mu^c$ denote the complement of $L_\mu$ in $\ZZ$. For any $\mu,\,\lambda\in \Delta$, we say that $(\mu,\,\lambda)$ is a \emph{poor pair} if: either
$\mu\in \Delta_{\textbf{a}_1}^{+}\cup\Delta_{\textbf{a}_1}^{0}$ and $\lambda\in \Delta_{\textbf{a}_1}^{+}$;  or $\mu,\lambda\in \Delta_{\textbf{a}_1}^{-}$.

\end{sect}

\begin{sect}\label{sect:2} Set
\begin{align*}
 \kappa_1=\min_{\stackrel{\lambda\in\Delta,}
 {\big|\lambda(d\textbf{a}_1)|\neq1}}\{|\log\lambda(d\textbf{a}_1)\big|\}\quad\text{and}\quad \kappa_2=\max_{\mu\in \Delta}\{\big|\log\lambda(d\textbf{a}_1) \big|\}.
\end{align*}
Choose $\kappa\geq \beta+1$ (see \eqref{for:179} of \ref{sect:6}) such that the followings are satisfied:
\begin{align}
 &\kappa \kappa_1>2\kappa_2,\label{for:139}\\
 &\kappa\gamma_1\min_{\stackrel{\lambda,\mu\in\Delta,}
 {\lambda(d\textbf{a}_i)\neq1}}\{\frac{\big|\log\mu(d\textbf{a}_i)\big|}{\big|\log\lambda(d\textbf{a}_i)\big|}\}
 >4\max_{\mu\in \Delta}\{\big|\log\lambda(d\textbf{a}_i)\big|\},\quad i=1,2 \label{for:146}
\end{align}
($\gamma_1$ is defined in \eqref{sec:45} of Section \ref{sec:19}).

\end{sect}

\begin{sect}\label{sect:4}
We fix $f$ which is of $(1,2)$-type (see \eqref{de:4} of Section \ref{sec:36}). Let $f_1=1-f$.

\end{sect}

\begin{sect}\label{sect:8} Suppose $w_i\in V$, $i=1,\,2$. We denote by
\begin{align}\label{for:336}
  z_ix:=w_i\cdot \alpha_{\textbf{a}_i}(x),\qquad \forall\, x\in \mathcal{X},\,\,i=1,2,
\end{align}
where $\cdot $ denotes the left translation of $w_i$ on $\alpha_{\textbf{a}_i}(x)$. It is clear that  $z_i$, $i=1,2$ are affine actions.

 Then $dz_i$, the induced action of $z_i$ on $\mathcal{V}$ is $\text{Ad}_{w_i}\cdot d\textbf{a}_i$, $i=1,\,2$, where $\cdot $ denotes the composition
 of two linear maps.

  For any map $\mathcal{F}$ defined on $\mathcal{X}$
denote $\mathcal{F}(z_ix)$ by $\mathcal{F}\circ z_i$.

\smallskip

\noindent\textbf{\emph{Note}}. We point that probably $z_1$  and $z_2$ don't commute.

\smallskip

For any $dz_1$,  we have the Lyapunov space decomposition
\begin{align}
\mathcal{V}=\sum_{\chi\in\Delta_{z_1}}\mathfrak{g}_{\chi,z_1}
\end{align}
such that for any $\chi\in\Delta_{z_1}$, $\chi(dz_1)$ is the absolute value of the eigenvalue of $dz_1$ on $\mathfrak{g}_{\chi,z_1}$. Similar to notations in \eqref{sect:7}, we can also define $\Delta_{z_1}^+$, $\Delta_{z_1}^-$ and $\Delta_{z_1}^0$, as well as $\mathfrak{g}^+_{z_1}$, $\mathfrak{g}^-_{z_1}$ and $\mathfrak{g}^0_{z_1}$.

\end{sect}
\begin{sect}\label{sect:10}

\eqref{for:168} shows that $\mathfrak{g}_1$ is a subalgebra. Denote by $V^o$ the subgroup with its Lie algebra $\mathfrak{g}_1$. Throughout Section \ref{sec:23}, we center on the case of $z_i=w_i\cdot \textbf{a}_i$, where $w_i\in V^o$, $i=1,\,2$.

For any $\lambda\in \Delta$ and $v\in E$ (see the end of \eqref{sect:7}), if $z_v=w_v\cdot \alpha_v$, $w_v\in V^o$, then
each $\mathfrak{g}_{\lambda}$, $\lambda\in \Delta$ is invariant under $dz_v$, $v\in E$. In fact, both the Lyapunov
exponents and the Lyapunov subspace of $z_v$ are determined by $\alpha_v$, the automorphism part. Moreover, if we denote by $\tilde{A}$ the subgroup of affine actions on $V$ generated by $z_v$, $v\in E$, even though $\tilde{A}$ is probably not abelian, \eqref{for:349} of Section \ref{sec:19} is also the Lyapunov subspace decomposition for $\tilde{A}$:
\begin{align*}
 \mathcal{V}= \sum_{\lambda\in \Delta}\mathfrak{g}_{\lambda}
\end{align*}
such that for any $\lambda\in \Delta$ and any $\tilde{a}\in \tilde{A}$, $\lambda(da)$ is the absolute value of the eigenvalue of $d\tilde{a}$ on $\mathfrak{g}_\lambda$, where $a$ is the automorphism part of $\tilde{a}$.

Further,
the following estimates hold: for any $v\in E$
\begin{align}\label{for:138}
 \norm{dz_v^n|_{\mathfrak{g}_{\lambda}}}\leq C_{\norm{w_v}}\lambda(dv)^n(1+|n|)^{c_1},\quad \forall\,\,n\in\ZZ,
\end{align}
where  and $c_1=2\dim\mathcal{V}+(\dim\mathcal{V})^2$.  The proof is left for Appendix \ref{sec:43} and Remark \ref{re:12} of Appendix \ref{sec:58}.

\end{sect}

\begin{sect}\label{sec:45} Then there exists a constant
$\gamma_1$ only dependent on $\alpha_A$ such that for any $z_i=w_i\cdot \textbf{a}_1$, where $w_i\in V^o$, $i=1,\,2$,  $f,\,g\in C^1(V/\Lambda)$ orthogonal to constants and any $n,\,m\in\ZZ$ we have
\begin{align*}
\abs{\langle f\circ (z_1^nz_2^m), g\rangle}&\leq
C_{\norm{w_1},\norm{w_2}}\norm{f}_{C^1}\norm{g}_{C^1}e^{-\gamma_1(|n|+|m|)}.
\end{align*}
The proof is left for Appendix \ref{sec:58}.
\end{sect}

\begin{sect}\label{sect:9}

Fix a right invariant metric $\textbf{\emph{d}}$ on $V$.  $B_r(x)\stackrel{\text{def}}{=}\{y\in
V|\textbf{\emph{d}}(x,y)\leq r\}$. There are $0<\delta<1$ such that:
\begin{enumerate}
  \item\label{for:365} for any $u\in \mathcal{V}$ if $\norm{u}\leq \delta$, then $\textbf{\emph{d}}(\exp u,e)\leq C\norm{u}$;

  \smallskip
  \item\label{for:364} for any $w\in B_{\delta}(e)$, we have
\begin{align*}
 \norm{\log w}+\norm{\text{Ad}_{w}-I}\leq C\textbf{\emph{d}}(w,e)\leq C_1\norm{\log w}.
\end{align*}
\end{enumerate}
Suppose $w_v\in V^o$, $v\in E$. Let $z_v=w_v\cdot v$, $v\in E$.
It is easy to check that for any $v_1,v_2\in E$,  there is $\rho(v_1,v_2)\in V^o$ such that
\begin{align*}
 z_{v_1}z_{v_2}(x)=\rho\cdot z_{v_2}z_{v_1}(x),\qquad \forall\,x\in \mathcal{X}.
\end{align*}
We denote $\rho$ by $[z_{v_1},z_{v_2}]$. We fix $\eta(E)>0$ , such that for if each $w_v\in B_{2\eta}(e)$, $v\in E$ we have
\begin{align*}
 \big\|\log[z_{v_1},z_{v_2}]\big\|\leq \delta ^2
\end{align*}
for any $v_1,v_2\in E$.

\smallskip

\noindent\emph{Note}. We point that probably the set $\{z_v: v\in E\}$ is not commuting.

\end{sect}

\begin{sect}\label{for:388} We denote by $V^c$ the subgroup of $V$ such that
\begin{align*}
 a(v)=v,\qquad \forall\,a\in A\quad\text{and}\quad \forall\,v\in V^c.
\end{align*}
We say that $\alpha_A$ is \emph{non-weakly partially hyperbolic} if $V^c$ is trivial. If $V=\RR^N$, i.e., $\mathcal{X}=\TT^N$, $\alpha_A$ is always non-weakly partially hyperbolic. However, in the case of non-abelian nilmanifolds,  $\alpha_A$ is not necessarily non-weakly partially hyperbolic.

If $\alpha_A$ is not non-weakly partially hyperbolic, we construct an abelian affine $\ZZ^k$ action $\tilde{\alpha}_A$ as follows: for any $a\in A$ set
\begin{align*}
 \tilde{\alpha}_a(x):=i_0(a)\cdot \alpha_a(x),\qquad \forall x\in \mathcal{X},
\end{align*}
where $i_0:\,A\to V^c$ is a group homeomorphism. If $i_0$ is close to the trivial homeomorphism, then $\tilde{\alpha}_A$ is a small perturbation of $\alpha_A$.
This means that  $\alpha_A$ is not structure stable among automorphisms on $V$: we can find abelian affine
  actions other than automorphisms which are close to $\alpha_A$.
\end{sect}

\subsection{Constants in Section \ref{sec:23}}\label{sec:60} In this part we list all the constants that will be used in Section \ref{sec:23}.
\begin{enumerate}
  \item $\beta$ is defined in \eqref{sect:6} of Section \ref{sec:19};
  \item $\kappa_1$,  $\kappa_2$ and  $\kappa$ are defined in \eqref{sect:2} of Section \ref{sec:19};
  \item $c_1$ is defined in \eqref{sect:10} of Section \ref{sec:19};
  \item $\gamma_1$ is defined in \eqref{sec:45} of Section \ref{sec:19};
  \item $\delta$ and $\eta$ are defined in \eqref{sect:9} of Section \ref{sec:19};
  \item\label{for:351} $\gamma_2=\min\{\frac{1}{2}\gamma_1,\,\frac{1}{4}\kappa_1 \kappa\}$;
  \item $\delta_m$ is from Lemma \ref{le:1};
  \item $\gamma_3$ and $c_0$ are defined in \eqref{for:304} of Section \ref{sec:40};
  \item $\varrho$ is defined in Section \ref{for:504}.
\end{enumerate}

\subsection{The perturbation of $\alpha_A$ }\label{sec:22}
We fix a compact generating set $E\subset \ZZ^k$ that contains $\textbf{a}_1$ and $\textbf{a}_2$. We consider a $C^r$-small $C^\infty$ perturbation
$\tilde{\alpha}_A$ of $\alpha_A$ on $\mathcal{X}=V/\Gamma$. By this we mean for any $v\in E$, the $C^r$ distance between
$\tilde{\alpha}_v$ and $\alpha_v$ is small.

For any $\varepsilon>0$ we say that $E'$ is a \emph{$\varepsilon$-standard   perturbation of $E$}  if $E'=\{w_v\cdot \alpha_v:\,v\in E,\,w_v\in V^o\cap B_{\varepsilon}(e)\}$.
For any $v\in E$, we use $z_v=w_v\cdot \alpha_v$ to denote the corresponding element in $E'$.
We emphasize that elements in $E'$ do not need to commute.
 By using \eqref{sect:10} of Section \ref{sec:19}, we see that
 $\mathfrak{g}_{\mu}$, $\mu\in \Delta$ are invariant subspaces for $\text{Ad}_{z_v}$, $v\in E$.

Unlike torus examples, where the perturbation $\tilde{\alpha}_A$ is conjugate to the original action $\alpha_A$, for
nilmanifold examples, \eqref{for:388} of Section \ref{sec:19} shows that $\tilde{\alpha}_A$ may conjugate to an affine
  action other than $\alpha_A$. As a consequence,  we need to consider the coordinate change in each iterative step. This implies that
  we write
\begin{align*}
 \tilde{\alpha}_v(x)=\exp(\mathfrak{p}_v(x))\cdot z_v(x)\qquad \forall\,v\in E,
\end{align*}
for any $x\in \mathcal{X}$, where $\mathfrak{p}_z:\mathcal{X}\to \mathcal{V}$ are $C^\infty$ small maps and $z_v=w_v\cdot
\alpha_v(x)$, where $w_v\in V^o$ and $\textbf{\emph{d}}(w_v,e)$ are small (see \eqref{sect:9} of Section \ref{sec:19}).

  \smallskip

\emph{Note}. Even though $\tilde{\alpha}_A$ is a $\ZZ^k$ action, we do not require that the set $\{z_v: v\in E\}$ commute. This is one of the key difference with previous work.  During the iterative step, the new set of generators determined by the coordinate change may not commute. We point out that it is not necessarily possible to find a set of commuting generators in sufficiently small neighborhoods of a set of non-commuting generators.
\smallskip

Let
\begin{gather*}
  r_v=\int_{\mathcal{X}}\mathfrak{p}_v,\quad \mathfrak{p}'_v=\mathfrak{p}_{v}-r_v, \qquad\text{and}\\
   \mathbb{D}_{u,v}=\mathfrak{p}'_{u}\circ z_v-dz_v \mathfrak{p}'_{u}-\big(\mathfrak{p}'_{v}\circ z_u-dz_u\mathfrak{p}'_{v}\big)
\end{gather*}
for any $v,\,u\in E$.

Set $\textbf{\emph{d}}(w,e)=\max_{v\in E}\textbf{\emph{d}}(w_v,e)$ and set
\begin{align*}
 \norm{\mathfrak{p}}_{C^m}=\max_{z\in E}\{\norm{\mathfrak{p}_z}_{C^m}\}\text{ and }\norm{\mathbb{D}}_{C^m}=\max_{u,v\in E}\{\norm{\mathbb{D}_{u,v}}_{C^m}\}\quad \forall\,m\geq0.
\end{align*}
 \emph{Note}. The constants $r_v$, $v\in E$  determines the coordinate change in the iterative step.  $\mathbb{D}$ measures how far these $\mathfrak{p}_v'$ are away from twisted cocycles over the new action generated by $\{z_v, \, v\in E\}$.

\smallskip

Even though $\{z_v: v\in E\}$ is not necessarily an abelian set, the next result shows  that,  similar to the abelian cases (see \eqref{for:411} of Remark \ref{re:11}),  $\mathbb{D}$ is quadratically small with respect to $\mathfrak{p}$. The proof is left for Appendix \ref{sec:12}.
\begin{lemma}\label{le:1} For any $m\geq0$,  there is $\delta_m>0$ such that if $\norm{\mathfrak{p}}_{C^m}\leq\delta_m$ and $\textbf{\emph{d}}(w,e)<\delta_0$, then
\begin{align}
 \max_{u,v\in E}\norm{\log[z_{v}&,z_{u}]}\leq C\norm{\mathfrak{p}}_{C^0}, \qquad\text{ and }\label{for:347}\\
 \norm{\mathbb{D}}_{C^{m}}
 &\leq C\norm{\mathfrak{p}}_{C^{m+1}}\norm{\mathfrak{p}}_{C^{m}}.\label{for:337}
\end{align}
(see \eqref{sect:9} of Section \ref{sec:19}).  Moreover, for any $1\neq \mu\in \Delta$, we have
\begin{align}\label{for:348}
  \big\|(r_v-d z_ur_v-r_u+dz_vr_u)\big|_{\mathfrak{g}_{\mu}}\big\|\leq C\norm{\mathfrak{p}}_{C^1}\norm{\mathfrak{p}}_{C^0}
\end{align}
for any $u,v\in E$.
\end{lemma}
\begin{remark}\label{re:11} 1. Abelianess is used as an essential condition in previous work  to
obtain quadratical smallness of $\mathbb{D}$ in \eqref{for:337}. However, after the coordinate change,
the new action generated by $\{z_v, \, v\in E\}$ is likely to lose  abelian property. This challenge compels us to develop new techniques for handling non-abelian behavior, see Appendix \ref{sec:12}.

2. If  $\alpha_A$ is structure stable among automorphisms on $V$,  we don't need to worry about the coordinate change in the iterative step. Then:
\begin{enumerate}
  \item $z_v=v$ for each $v\in E$;

  \smallskip
  \item the map $\mathfrak{p}$ extends to a twisted cocycle over $\tilde{\alpha}_A$, for example see \cite{Damjanovic4}, \cite{vw}. We define
  \begin{align*}
  \mathbb{E}_{u,v}=\mathfrak{p}_{u}\circ z_v-dz_v \mathfrak{p}_{u}-\big(\mathfrak{p}_{v}\circ z_u-dz_u\mathfrak{p}_{v}\big).
  \end{align*}
  Usually, $\mathbb{E}$ is non-trivial as $\mathfrak{p}$ is not a twisted cocycle over $\alpha_A$;

  \smallskip

  \item \label{for:411} of Lemma \ref{le:1} becomes
  \begin{align}\label{for:338}
   \norm{\mathbb{E}}_{C^{m}}
 &\leq C\norm{\mathfrak{p}}_{C^{m+1}}\norm{\mathfrak{p}}_{C^{m}}
  \end{align}
  (see Appendix \ref{sec:12}). Obtaining estimates of  $\mathbb{E}$ follows closely to the proof of Lemma $4.7$ in \cite{Damjanovic4}.

\end{enumerate}

3. As we need to take into account  the coordinate change in each iterative step, we study the twisted equation \eqref{for:135} and the almost twisted cocycle equation \eqref{for:143} in the subsequent sections,  where $z_1$ and $z_2$ are small affine perturbations of $\textbf{a}_1$ and $\textbf{a}_2$,  respectively,  along $V^o$.

\end{remark}

\subsection{Change of spectrum for affine actions} We now adapt Lemma \ref{le:4} to the specific case concerning affine actions on nilmanifolds. Suppose $w_i\in V$, $i=1,\,2$. We denote by $z_i=w_i\cdot \textbf{a}_i $, $i=1,\,2$. We recall that $dz_i=\text{Ad}_{w_i}\cdot d\textbf{a}_i$, $i=1,\,2$ (see \eqref{sect:8} of Section \ref{sec:19}).

\begin{corollary}\label{le:15} Suppose $f\in L^\infty(\RR^m)$ and $\xi\in \mathfrak{F}(\mathcal{O})$, where $\mathfrak{F}$ is a subspace of $\mathcal{V}$. Then for any  $0\neq\mathfrak{u}\in \mathcal{V}$ and any $n,l\in\ZZ$, we have
\begin{align*}
 \pi_{\mathfrak{u}}(f)\big(\xi\circ (z_1^nz_2^l)\big)=\big(\pi_{\mathfrak{u}_{n,m}}(f\circ \norm{dz_1^ndz_2^l(\mathfrak{u})})\xi\big)\circ (z_1^nz_2^l),
\end{align*}
where
\begin{align*}
\mathfrak{u}_{n,l}=\norm{dz_1^ndz_2^l(\mathfrak{u})}^{-1}dz_1^ndz_2^l(\mathfrak{u})
\end{align*}

\end{corollary}
\begin{proof}  Let $S=\ZZ^2\ltimes V$ where $\ZZ^2$ is the subgroup generated by $\textbf{a}_1$ and $\textbf{a}_2$. The
multiplication of elements in $S$ is given by
\begin{align}\label{for:150}
(g_1,v_1)\cdot(g_2,v_2)=(g_1g_2,\, g_2^{-1}(v_1)v_2).
\end{align}
Let $(\bar{\pi},\,\mathfrak{F}(\mathcal{O}))$ denote the unitary representation of $S$ on $\mathfrak{F}(\mathcal{O})$:
\begin{align}\label{for:151}
 \bar{\pi}(g,v)\xi=\xi\circ (g,v)^{-1}.
\end{align}
\emph{Note}. We see that
\begin{align*}
 \big(\xi\circ (g,v)^{-1}\big)(x)=\xi\big(v^{-1}\cdot g^{-1}(x)\big),\qquad \forall\,x\in \mathcal{X}.
\end{align*}
Then
\begin{align*}
 z_i=(0,w_i)\cdot (\textbf{a}_i,e)=(\textbf{a}_i, \,\textbf{a}_i^{-1}w_i), \qquad i=1,\,2.
\end{align*}
Let $h=z_1^{n}z_2^{l}$.  From Lemma \ref{le:4}, letting $\mathfrak{u}_1=\norm{\text{Ad}_{h}(\mathfrak{u})}^{-1}\text{Ad}_{h}(\mathfrak{u})$, we have
\begin{align*}
 \bar{\pi}_{\mathfrak{u}}(f)\bar{\pi}(h^{-1})\xi&=\bar{\pi}(h^{-1})\bar{\pi}_{\mathfrak{u}_1}(f\circ \norm{\text{Ad}_{h}(\mathfrak{u})})\xi\\
 \stackrel{(1)}{\Rightarrow}\pi_{\mathfrak{u}}(f)\bar{\pi}(h^{-1})\xi&=\bar{\pi}(h^{-1})\pi_{\mathfrak{u}_1}(f\circ \norm{\text{Ad}_{h}(\mathfrak{u})})\xi\\
 \stackrel{(2)}{\Rightarrow}\pi_{\mathfrak{u}}(f)\bar{\pi}(h^{-1})\xi&=\bar{\pi}(h^{-1})\pi_{\mathfrak{u}_{n,l}}(f\circ \norm{dz_1^ndz_2^l(\mathfrak{u})})\xi\\
 \stackrel{(3)}{\Rightarrow}\pi_{\mathfrak{u}}(f)\big(\xi\circ (z_1^nz_2^l)\big)&=\big(\pi_{\mathfrak{u}_{n,l}}(f\circ \norm{dz_1^ndz_2^l(\mathfrak{u})})\xi\big)\circ (z_1^nz_2^l).
\end{align*}
Here in $(1)$ we use $\{\mathfrak{u}\}\in V$ and $\bar{\pi}|_{V}=\pi$; in $(2)$ we use $\text{Ad}_{h}(\mathfrak{u})=dz_1^ndz_2^l(\mathfrak{u})$, which is from \eqref{for:150};
in $(3)$ we use \eqref{for:151}.

\end{proof}

\subsection{Study of twisted equation}\label{sec:13} In this section, we assume $z_1=w_1\cdot \textbf{a}_1$ where $w_1\in V^o\cap B_{\eta}(e)$ (see Section \ref{sec:60}). We study the twisted equation:
\begin{align}\label{for:135}
\varphi\circ z_1-dz_1\varphi=\omega,
\end{align}
where $\varphi,\,\omega\in \mathcal{V}(\mathcal{O})$ (see \eqref{sect:1} of Section \ref{sec:19}). Results in this section serve as the first step in constructing the approximated solution in Section \ref{sec:20}.

Since each $\mathfrak{g}_{\lambda}$, $\lambda\in \Delta$ is invariant under $dz_1$ (see \eqref{sect:10} of Section \ref{sec:19}), \eqref{for:135} decomposes into equations of the following forms:
\begin{align}\label{for:136}
\varphi_\mu\circ z_1-(dz_1|_{\mathfrak{g}_{\mu}})\varphi_\mu=\omega_\mu,\qquad \mu\in \Delta.
\end{align}
where $\varphi_\mu=\varphi|_{\mathfrak{g}_{\mu}}$ and $\omega_\mu=\omega|_{\mathfrak{g}_{\mu}}$. We note that $\varphi_\mu,\,\omega_\mu\in \mathfrak{g}_{\mu}(\mathcal{O})$.

It is clear that the study of equation \eqref{for:135} boils down to that of the decomposed equations \eqref{for:136}.
In this part, we define obstructions to solving  equation \eqref{for:136}. After that, we show that vanishing of these obstructions implies solvability of \eqref{for:136} with tame Sobolev estimates.

\subsubsection{Notations}\label{sec:39} The positive constants $\eta, c_1,\,\kappa,\,\gamma_1,\,\kappa_1,\,\gamma_2$ (see Section \ref{sec:60}) will appear in is section.
Suppose $\xi\in \mathfrak{g}_\mu(\mathcal{O}^m)$, $m\geq \kappa$ and $z_1=w_1\cdot \textbf{a}_1$ where $w_1\in V^o\cap B_{\eta}(e)$.
\begin{enumerate}

  \item\label{for:284} $\mathcal{P}_{\mu}(\xi)\stackrel{\text{def}}{=} \delta(\mu)\sum_{n\in L_\mu} (dz_1|_{\mathfrak{g}_{\mu}})^{n-1} \xi\circ z_1^{-n}$ (see \eqref{sect:3} of Section \ref{sec:19}).

    \smallskip
  \item\label{for:285} Suppose $(\mu,\,\lambda)$ is a poor pair (see \eqref{sect:3} of Section \ref{sec:19}). Then for any subset $\mathbb{S}\subseteq\ZZ$, any $u\in \mathfrak{g}_{\lambda}\bigcap \mathcal{V}^1$ ($\mathcal{V}^1$ is defined in \eqref{sect:5} of Section \ref{sec:19}) and any $\tau>0$
\begin{align*}
 \mathfrak{D}_{\mu,\lambda,u,\tau,\mathbb{S}}(\xi)\stackrel{\text{def}}{=}\sum_{n\in \mathbb{S}} (dz_1|_{\mathfrak{g}_{\mu}})^{n-1} \pi_{u}(f_1\circ \tau^{-1})(\xi\circ z_1^{-n}).
\end{align*}
$f_1$ is defined in \eqref{sect:4} of Section \ref{sec:19}.

\smallskip

  \item\label{for:294} For any subset $\mathbb{S}\subseteq\ZZ$ and any $\tau>0$
  \begin{align*}
 \norm{\mathfrak{D}_{\mu,\tau,\mathbb{S}}(\xi)}\stackrel{\text{def}}{=}\max_{\stackrel{\lambda\in\Delta,}
 {\text{$(\mu,\lambda)$ is a poor pair}}}\sup_{u\in \mathfrak{g}_{\lambda}\bigcap \mathcal{V}^1}\sup_{\stackrel{\psi\in \mathfrak{g}_\mu(\mathcal{O}^\kappa),}{\norm{\psi}_\kappa=1}} \big\|\langle \mathfrak{D}_{\mu,\lambda,u,\tau,\mathbb{S}}(\xi),\,\psi\rangle\big\|
 \end{align*}
  (recall \eqref{sect:1} of Section \ref{sec:19} for the definition of $\langle \cdot \rangle$).

\end{enumerate}

\subsubsection{Main results}\label{sec:61} Suppose $\xi,\,\omega_\mu\in \mathfrak{g}_\mu(\mathcal{O}^\infty)$ and suppose $u\in \mathfrak{g}_{\lambda}\bigcap \mathcal{V}^1$, $\lambda\in \Delta$.  Also suppose $\tau\geq1$. Then:
\begin{enumerate}
  \item\label{for:502} (\eqref{for:350} of Proposition \ref{po:6}) for any $n\in L_\mu^c$ we have
\begin{align*}
 &\big\|\langle (dz_1|_{\mathfrak{g}_{\mu}})^{n-1} \pi_{u}(f_1\circ \tau^{-1})(\xi\circ z_1^{-n}),\,\psi\rangle\big\|\notag\\
 &\leq C_{m}\tau^{-m}e^{-\frac{1}{4}\kappa_1 |n|m}\norm{\xi}_m\norm{\psi}
\end{align*}
for any $m\geq \kappa$ and any $\psi\in \mathfrak{g}_{\mu}(\mathcal{O})$;

  \smallskip
  \item\label{for:501} (\eqref{for:354} of Proposition \ref{po:6}) for any $n\in\ZZ$ we have
\begin{align*}
 &\big\|\langle (dz_1|_{\mathfrak{g}_{\mu}})^{n-1} \pi_{u}(f_1\circ \tau^{-1})(\xi\circ z_1^{-n}),\,\psi\rangle\big\|\notag\\
 &\leq Ce^{-\gamma_2 |n|}\norm{\xi}_\kappa\norm{\psi}_\kappa
\end{align*}
for any $\psi\in \mathfrak{g}_{\mu}(\mathcal{O}^\kappa)$;

\smallskip
   \item (Lemma \ref{le:5}) $\mathcal{P}_{\mu}(\omega_\mu)$ is a distributional solution of equation \eqref{for:136};

   \smallskip

   \item (Lemma \ref{le:5}) if equation \eqref{for:136} has a solution  $\varphi_\mu\in \mathfrak{g}_\mu(\mathcal{O}^\kappa)$, then $\varphi_\mu$ is unique and
  \begin{align*}
   \norm{\mathfrak{D}_{\mu,\tau,\ZZ}(\omega_\mu)}=0,\qquad \forall\,\tau>0;
  \end{align*}

  \smallskip
  \item (\eqref{for:293} of Proposition \ref{le:17}) if there is $\tau_0>0$ such that
 $\norm{\mathfrak{D}_{\mu,\tau_0,\ZZ}(\omega_\mu)}=0$, then equation \eqref{for:136} has a solution
$\varphi_\mu\in \mathfrak{g}_\mu(\mathcal{O}^\infty)$ with estimates
\begin{align*}
 \norm{\varphi_\mu}_{l}\leq C_{l}\norm{\omega_\mu}_{l+\kappa+4}
\end{align*}
for any $l\geq0$.
\end{enumerate}

\subsubsection{A typical example}\label{sec:38} We use a typical example to explain some frequently used nations in this part. Suppose $\Delta=\{\lambda_1,\lambda_2,\lambda_3\}$ and
$\lambda_1(d\textbf{a}_1)=\frac{1}{2}$, $\lambda_2(d\textbf{a}_1)=1$, $\lambda_3(d\textbf{a}_1)=2$. We also suppose $\xi\in \mathfrak{g}_\mu(\mathcal{O}^\infty)$. Then:

\begin{enumerate}
  \item $\Delta_{z_1}^{+}=\{\lambda_3\}$, $\Delta_{z_1}^{-}=\{\lambda_1\}$ and $\Delta_{z_1}^{0}=\{\lambda_2\}$ (see \eqref{sect:10} of Section \ref{sec:19}).

  \smallskip
  \item Any $0\neq v\in \mathfrak{g}_{\lambda_1}$ is a stable direction of $z_1$, any $0\neq v\in \mathfrak{g}_{\lambda_3}$ is an unstable direction of $z_1$ and any
  $0\neq v\in \mathfrak{g}_{\lambda_2}$ is a neutral direction of $z_1$.

  \smallskip
  \item If $\mu\in \{\lambda_3,\lambda_2\}$, $L_\mu=\ZZ^-\cup\{0\}$, $\delta(\mu)=-$ and $L_\mu^c=\ZZ^+$.

  \smallskip
  \item If $\mu=\lambda_1$, $L_\mu=\ZZ^+\cup\{0\}$, $\delta(\mu)=+$ and $L_\mu^c=\ZZ^-$.

  \smallskip

  \item If $\mu=\lambda_3$ or $\lambda_2$, then $(\mu,\,\lambda)$ is a poor pair if $\lambda=\lambda_3$.
  \smallskip

   \item If $\mu=\lambda_1$, then $(\mu,\,\lambda)$ is a poor pair if $\lambda=\lambda_1$.

   \smallskip

   \item\label{for:283} If $\mu=\lambda_3$
   \begin{align*}
    \mathcal{P}_{\mu}(\xi)&= -\sum_{n\leq 0} (dz_1|_{\mathfrak{g}_{\mu}})^{n-1} \xi\circ z_1^{-n}.
   \end{align*}
  \eqref{for:138} of Section \ref{sec:19} shows that $\mathcal{P}_{\mu}(\xi)\in \mathfrak{g}_{\mu}(\mathcal{O})$. Further,
for any $v_i\in \mathfrak{g}^-_{\textbf{a}_1}\oplus\mathfrak{g}^0_{\textbf{a}_1}=\mathfrak{g}^{\delta(\mu)}_{\textbf{a}_1}\oplus\mathfrak{g}^0_{\textbf{a}_1}
=\mathfrak{g}^{\delta(\mu)}_{z_1}\oplus\mathfrak{g}^0_{z_1}$, we see that
\begin{align*}
 v_1v_2\cdots v_l\mathcal{P}_{\mu}(\xi)\in \mathfrak{g}_{\mu}(\mathcal{O}).
\end{align*}
If $\mu=\lambda_2$, $\mathcal{P}_{\mu}(\xi)$ is a distribution. Further, for any $v_i\in \mathfrak{g}^-_{\textbf{a}_1}\oplus\mathfrak{g}^0_{\textbf{a}_1}=\mathfrak{g}^{\delta(\mu)}_{\textbf{a}_1}\oplus\mathfrak{g}^0_{\textbf{a}_1}
 =\mathfrak{g}^{\delta(\mu)}_{z_1}\oplus\mathfrak{g}^0_{z_1}$,
\begin{align*}
 v_1v_2\cdots v_l\mathcal{P}_{\mu}(\xi)
\end{align*}
is also a distribution (see \eqref{for:148} of Proposition \ref{le:17}).

For either case, we say that $\mathcal{P}_{\mu}(\xi)$ is $\infty$ differentiable along $\mathfrak{g}^{\delta(\mu)}_{\textbf{a}_1}\oplus\mathfrak{g}^0_{\textbf{a}_1}$. In other words, $\mathcal{P}_{\mu}(\xi)$ is $\infty$ differentiable along $\mathfrak{g}_{\lambda}$ if $(\mu,\lambda)$ is a NOT poor pair.

    However, $\mathcal{P}_{\mu}(\xi)$ probably loses differentiability along $\mathfrak{g}_{\lambda_3}=\mathfrak{g}^{-\delta(\mu)}_{\textbf{a}_1}$. In this case, $(\mu,\lambda_3)$ is a poor pair.

   \smallskip

   \item If $\mu=\lambda_1$
   \begin{align*}
    \mathcal{P}_{\mu}(\xi)= \sum_{n\geq0} (dz_1|_{\mathfrak{g}_{\mu}})^{n-1} \xi\circ z_1^{-n}.
   \end{align*}
   \eqref{for:138} of Section \ref{sec:19} shows that $\mathcal{P}_{\mu}(\xi)\in \mathfrak{g}_{\mu}(\mathcal{O})$. Further,
for any $v_i\in \mathfrak{g}^+_{\textbf{a}_1}\oplus\mathfrak{g}^0_{\textbf{a}_1}=\mathfrak{g}^{\delta(\mu)}_{\textbf{a}_1}\oplus\mathfrak{g}^0_{\textbf{a}_1}
=\mathfrak{g}^{\delta(\mu)}_{z_1}\oplus\mathfrak{g}^0_{z_1}$,
\begin{align*}
 v_1v_2\cdots v_l\mathcal{P}_{\mu}(\xi)\in \mathfrak{g}_{\mu}(\mathcal{O}).
\end{align*}
Then we see that $\mathcal{P}_{\mu}(\xi)$ is $\infty$ differentiable along $\mathfrak{g}_{\lambda}$ if $(\mu,\lambda)$ is a NOT poor pair.

However, $\mathcal{P}_{\mu}(\xi)$ probably loses differentiability along $\mathfrak{g}_{\lambda_1}=\mathfrak{g}^{-\delta(\mu)}_{\textbf{a}_1}$. In this case, $(\mu,\lambda_1)$ is a poor pair.
\end{enumerate}

\subsubsection{Basic ideas}\label{sec:41} In this part, we will first explain the construction of invariant distributions. Subsequently, we will show that the vanishing of
the invariant distributions implies a smooth solution.

  For the torus case, i.e., $V=\RR^m$, the solvability condition for equation \eqref{for:136} is the vanishing of the invariant distribution
\begin{align}\label{for:198}
 \mathfrak{D}:=\sum_{n\in \ZZ} (dz_1|_{\mathfrak{g}_{\mu}})^{n-1} \omega_\mu\circ z_1^{-n}.
\end{align}
Equivalently, obstructions are sums along the orbits of $z_1$. Even though the description  generalizes to partially hyperbolic algebraic actions when $|\mu(d\textbf{a}_1)|=1$
(see \cite{Spatzier1} and \cite{W4}), for the case of $|\mu(d\textbf{a}_1)|\neq1$ it is a rare
occurrence that such a description of obstructions is possible. The convergence in \eqref{for:198}
is due to the fact that super-exponential decay of matrix coefficients
for smooth functions on the torus, which guarantees that $\mathfrak{D}$ is a distribution. By contrast, for semisimple and non-abelian nilpotent cases,
matrix coefficients decay exponentially at a particular speed, no matter how smooth the functions are. In general the absolute values of eigenvalues of $dz_1$ beat $\gamma$ (see Theorem \ref{th:1}),  and thus the right side of \eqref{for:198} diverges too fast to produce even a distribution.

We solve the problem by using the system
of imprimitivity. Before explaining the construciton, we present a key observation, which plays a crucial role in understanding the main idea.

\smallskip
$(*)$ \emph{Observation: test space in Sobolev space.} We recall notations in Section \ref{sec:37}. Suppose $g\in L^2(\RR)$.  It is known that we can use functions in $W^{2,\infty}$ to test the smoothness of $g$. More precisely,  $g\in W^{2,\infty}$ if there is $p\geq0$ such that
\begin{align*}
 |\langle g,\,\text{\tiny$\frac{\partial}{\partial x^n}$}\psi\rangle|\leq C_{n,g}\norm{\psi}_{W^{2,p}}
\end{align*}
for any $n\in\NN$ and any $\psi\in W^{2,\infty}$.

Suppose $f$ is of $(1,2)$ type (see \eqref{de:4} of Section \ref{sec:36}). We decompose any $\psi\in W^{2,\infty}$:
\begin{align}\label{for:279}
 \psi=\pi(f\circ \tau^{-1})\psi+\pi(f_1\circ \tau^{-1})\psi,
\end{align}
where $f_1=1-f$.  We call $\pi(f\circ \tau^{-1})\psi$ the $\tau$-truncation of $\psi$ and $\pi(f_1\circ \tau^{-1})\psi$ the $\tau$-tail of $\psi$.  \eqref{for:279} shows shows that
$W^{2,\infty}$ is a direct sum of two subspaces: the space of $\tau$-truncations and the space of $\tau$-tails.

The key observation is: for any $\tau>0$, functions in $\pi(f_1\circ \tau^{-1})(W^{2,\infty})$ (the space of $\tau$-tails) are sufficient for serving as test functions.  More precisely,  $g\in W^{2,\infty}$ if there is $p\geq0$ such that
\begin{align*}
 |\langle g,\,\text{\tiny$\frac{\partial}{\partial x^n}$}\pi(f_1\circ \tau^{-1})\psi\rangle|\leq C_{n,g}\norm{\psi}_{W^{2,p}}
\end{align*}
for any $n\in\NN$ and any $\psi\in W^{2,\infty}$.

To see this, from \eqref{for:2006} of Section \ref{sec:37} we have
\begin{align*}
\norm{\text{\tiny$\frac{\partial}{\partial x^n}$}\pi(f\circ \tau^{-1})\psi}_{L^2(\RR)}\leq \tau^n C_n\norm{\psi}_{L^2(\RR)},
\end{align*}
which gives
\begin{align}\label{for:278}
 |\langle g,\,\text{\tiny$\frac{\partial}{\partial x^n}$}\pi(f\circ \tau^{-1})\psi\rangle|&\leq \norm{g}_{L^2(\RR)}\cdot\norm{\text{\tiny$\frac{\partial}{\partial x^n}$}\pi(f\circ \tau^{-1})\psi}_{L^2(\RR)}\notag\\
 &\leq \norm{g}_{L^2(\RR)}\cdot \tau^n C_n\norm{\psi}_{L^2(\RR)}
\end{align}
for any $n\in\NN$. It shows that when taking higher-order derivatives  the $\tau$-truncations, the upperbounds on the right side of \eqref{for:278} depend only  on
the $L^2$ norms of $g$ and $\psi$. This means for any fixed $\tau>0$, the space of $\tau$-truncations fails to test the smoothness of $g$. Thus, it claims the observation:
\begin{enumerate}
  \item [] \hypertarget{o.2}{(*)} it is, in fact,
the space of $\tau$-tails that performs the testing job.
\end{enumerate}
Below, we use a typical example described in Section \ref{sec:38} to explain the main idea for the sake of clarity and transparency. Throughout Section \ref{sec:23}, it is recommended
to keep this example in mind.

It is harmless to assume that $z_1=\textbf{a}_1$ and each $\mathfrak{g}_{\lambda_i}$, $1\leq i\leq 3$ is one-dimensional.  Then equation \eqref{for:136} for $\mu=\lambda_3$ has the form:
\begin{align}\label{for:291}
\varphi_{\lambda_3}\circ z_1-2\varphi_{\lambda_3}=\omega_{\lambda_3}
\end{align}
where $\omega_{\lambda_3}\in \mathcal{O}^\infty$ (see \eqref{sect:5} of Section \ref{sec:19}). It is clear that
\begin{align}\label{for:286}
   \varphi_{\lambda_3} = -\sum_{n\geq0} 2^{-(n+1)} \omega_{\lambda_3}\circ z_1^{n}\in \mathcal{O}
   \end{align}
is a solution. Another key observation is:
\begin{enumerate}
  \item [] \hypertarget{o.3}{(**)} $\varphi_{\lambda_3}$ is $\infty$ differentiable along stable and neutral directions of $z_1$ (in the sense we described in \eqref{for:283} of Section \ref{sec:38}). However, it loses differentiability along
the unstable direction of $z_1$, which is determined by $\mathfrak{g}_{\lambda_3}=\mathfrak{g}^{+}_{\textbf{a}_1}$.
\end{enumerate}
To define the invariant distributions, it is natural for us to consider the following
formal sum
\begin{align*}
 \mathfrak{D}:=\sum_{n=-\infty}^\infty 2^{-(n+1)} \omega_{\lambda_3}\circ z_1^{n}
\end{align*}
and show the convergence in the sense of a distribution: there is $p\geq0$ such that
\begin{align*}
 |\langle \mathfrak{D},\,\psi\rangle|\leq C\norm{\psi}_p,\qquad \forall\,\psi \in \mathcal{O}^\infty.
\end{align*}
Unfortunately, as we explained at the beginning of this part, carrying out this construction naively using decay of matrix coefficients
seems to be a nightmare. To get around this problem, we start from a new angle which is motivated by the observations \hyperlink{o.2}{(*)} and \hyperlink{o.3}{(**)}: firstly, instead of showing the convergence of $\mathfrak{D}$ against the entire space of $\mathcal{O}^\infty$ testing functions, we focus solely on its convergence against the space of $\tau$-tails. Secondly, we only need to
consider $\tau$-tails coming from $\mathfrak{g}_{\lambda_3}=\mathfrak{g}^{+}_{\textbf{a}_1}$, as differentiability of $\varphi_{\lambda_3}$ along the remaining directions inside $\mathfrak{g}^{-}_{\textbf{a}_1}\oplus\mathfrak{g}^{0}_{\textbf{a}_1}$ is already obtained.
More precisely, we show that:
there is $p\geq0$ such that for any $\tau>0$
\begin{align}\label{for:280}
 |\langle \mathfrak{D},\,\pi_{u}(f_1\circ \tau^{-1})\psi\rangle|\leq C_\tau\norm{\psi}_p,\qquad \forall\,\psi \in \mathcal{O}^\infty,
\end{align}
where $0\neq u\in \mathfrak{g}_{\lambda_3}$.

\eqref{for:280} implies that we should define the invariant distributions as follows:
\begin{align}\label{for:331}
 \mathfrak{D}_{\lambda_3,\tau}:=\sum_{n=-\infty}^\infty 2^{-(n+1)} \pi_{u}(f_1\circ \tau^{-1})(\omega_{\lambda_3}\circ z_1^{n}),\qquad \forall\,\tau>0.
\end{align}
$\mathfrak{D}_{\lambda_3,\tau}$ is exactly the $\mathfrak{D}_{\lambda_3,\lambda_3,u,\tau,\ZZ}(\omega_{\lambda_3})$ defined in  \eqref{for:285} of Section \ref{sec:39}.

Before we show the convergence of $\mathfrak{D}_{\lambda_3,\tau}$, we study the deformation of $\tau$-tails under the natural dual action
of $z_1$, which is important for this purpose.

From Corollary \ref{le:15}, we have:
\begin{align}\label{for:281}
 \pi_{u}(f_1\circ \tau^{-1})(\omega_{\lambda_3}\circ z_1^n)=\big(\pi_{u}(f\circ \tau_n^{-1})\omega_{\lambda_3}\big)\circ (z_1^n),\qquad \forall\,n\in\ZZ.
\end{align}
where $\tau_n=2^{-n}\tau$. \eqref{for:281} shows that $z_1^n$ deforms $\tau$-tails to $2^{-n}\tau$-tails.

Now we are ready to show the convergence of $\mathfrak{D}_{\lambda_3,\tau}$. We can write
\begin{align*}
 \mathfrak{D}_{\lambda_3,\tau}=P_{1,\lambda_3,\tau}+P_{2,\lambda_3,\tau}
\end{align*}
where
\begin{align}\label{for:310}
    P_{1,\lambda_3,\tau}&=\sum_{n\geq 0} 2^{-(n+1)} \pi_{u}(f_1\circ \tau^{-1})(\omega_{\lambda_3}\circ z_1^n)\qquad\text{and}\notag\\
    P_{2,\lambda_3,\tau}&=\sum_{n<0} 2^{-(n+1)} \pi_{u}(f_1\circ \tau^{-1})(\omega_{\lambda_3}\circ z_1^n).
\end{align}
Firstly, we show the convergence of $P_{1,\lambda_3,\tau}$: for any $\psi \in \mathcal{O}^\infty$, we have
\begin{align}\label{for:289}
 |\langle P_{1,\lambda_3,\tau},\,\psi\rangle|&=\Big| \sum_{n\geq 0} 2^{-(n+1)} \big\langle\pi_{u}(f_1\circ \tau^{-1})(\omega_{\lambda_3}\circ z_1^n),\, \psi \big\rangle\Big|\notag\\
 &\overset{\text{(1)}}{=}\Big| \sum_{n\geq 0} 2^{-(n+1)} \big\langle \omega_{\lambda_3}\circ z_1^n,\, \pi_{u}(f_1\circ \tau^{-1})\psi \big\rangle\Big|\notag\\
 &\overset{\text{(2)}}{\leq} \sum_{n\geq 0} 2^{-(n+1)} \norm{\omega_{\lambda_3}\circ z_1^n}\cdot \norm{\pi_{u}(f_1\circ \tau^{-1})\psi }\notag\\
 & \overset{\text{(3)}}{\leq} \sum_{n\geq 0} 2^{-(n+1)} \norm{\omega_{\lambda_3}}\cdot \norm{\psi }\notag\\
 &\leq C\norm{\omega_{\lambda_3}}\cdot \norm{\psi }.
\end{align}
Here in $(1)$ we use \eqref{for:93} of Section \ref{sec:27}; in $(2)$ we use cauchy schwarz inequality; in $(3)$ we use \eqref{for:215} of  Section \ref{sec:27} and the fact that
\begin{align*}
\norm{\varphi}=\norm{\varphi\circ z_1^n},\qquad \forall\,\varphi\in \mathcal{O};
\end{align*}
\emph{Note}. When $n\geq 0$, the twist $2^{-(n+1)}$ is not an obstacle for convergence; rather; it is helpful for convergence.

\smallskip

Secondly, we show the convergence of $P_{2,\lambda_3,\tau}$. Before that, we show the super-exponential decay when $n<0$, which guarantees the convergence of $P_{2,\lambda_3,\tau}$.
For any $\psi \in \mathcal{O}^\infty$ and $n<0$, we have
\begin{align}\label{for:282}
 &\Big|\big\langle\pi_{u}(f_1\circ \tau^{-1})(\omega_{\lambda_3}\circ z_1^n),\, \psi \big\rangle\Big|\notag\\
 &\overset{\text{(1)}}{=}\Big| \big\langle \big(\pi_{u}(f_1\circ \tau_n^{-1})\omega_{\lambda_3}\big)\circ (z_1^n),\, \psi \big\rangle\Big|\notag\\
 &\overset{\text{(2)}}{\leq} \norm{\pi_{u}(f_1\circ \tau_n^{-1})\omega_{\lambda_3}}\cdot \norm{\psi}\notag\\
 & \overset{\text{(3)}}{\leq} C\tau_n^{-s}\norm{\omega_{\lambda_3}}_{s}\cdot \norm{\psi }\notag\\
 &\overset{\text{(4)}}{=} C2^{ns}\tau^{-s}\norm{\omega_{\lambda_3}}_{s}\cdot \norm{\psi }.
\end{align}
Here in $(1)$ we use \eqref{for:281}; in $(2)$ we use cauchy schwarz inequality and the fact that
\begin{align*}
\norm{\varphi}=\norm{\varphi\circ z_1^n},\qquad \forall\,\varphi\in \mathcal{O};
\end{align*}
in $(3)$ we use Lemma \ref{cor:4}; in $(4)$ we recall that $\tau_n=2^{-n}\tau$.

\smallskip

\emph{Note}. As $n<0$, $2^{ns}=(\frac{1}{2})^{|n|s}$. Noting that $\omega_{\lambda_3}\in \mathcal{O}^\infty$, we can choose $s$ to be as large as we want. This means that
the matrix coefficients can decay as quickly as desired. This is how we obtain super-exponential decay for $\tau$-tails.
\smallskip

By choosing $s>1$, the convergence of $P_{2,\lambda_3,\tau}$ follows directly from \eqref{for:282}:
\begin{align}\label{for:290}
 |\langle P_{2,\lambda_3,\tau},\,\psi\rangle|&=\Big| \sum_{n<0} 2^{-(n+1)} \big\langle\pi_{u}(f_1\circ \tau^{-1})(\omega_{\lambda_3}\circ z_1^n),\, \psi \big\rangle\Big|\notag\\
 &\leq\sum_{n<0} 2^{-(n+1)}\Big|  \big\langle\pi_{u}(f_1\circ \tau^{-1})(\omega_{\lambda_3}\circ z_1^n),\, \psi \big\rangle\Big|\notag\\
 &\leq\sum_{n<0} 2^{-(n+1)}\cdot C2^{ns}\tau^{-s}\norm{\omega_{\lambda_3}}_{s}\cdot \norm{\psi }\notag\\
 &\leq C_{n,s}\tau^{-s}\norm{\omega_{\lambda_3}}_{s}\cdot \norm{\psi }.
\end{align}
Hence, we showed the convergence of $\mathfrak{D}_{\lambda_3,\tau}$.
Thus we finish the construction of invariant distributions.

\smallskip

\emph{Note}. When $n<0$, the twist is $2^{|n|-1}$ and the matrix coefficients decay at $(\frac{1}{2})^{|n|s}$. To guarantee the convergence, we choose $s>1$ such that
$(\frac{1}{2})^{s}\cdot 2<1$.

\begin{remark}\label{re:13} Comments from the above discussion:

\begin{enumerate}
\item Equation \eqref{for:291} has two formal solutions:
\begin{align*}
   \varphi_{\lambda_3}^+& = -\sum_{n\geq0} 2^{-(n+1)} \omega_{\lambda_3}\circ z_1^{n}\quad\text{and}\quad\varphi^-_{\lambda_3}=\sum_{n<0} 2^{-(n+1)} \omega_{\lambda_3}\circ z_1^{n}.
   \end{align*}
We set the solution $\varphi_{\lambda_3}=\varphi_{\lambda_3}^+$ as the other one $\varphi_{\lambda_3}^-$ is not even a distribution. Of course,  if $|\lambda_3(d\textbf{a}_1)|=1$
then we can choose either one.

We point out that $\varphi_{\lambda_3}$ is exactly the $\mathcal{P}_{\lambda_3}(\omega_{\lambda_3})$ defined in \eqref{for:284} of Section \ref{sec:39}, where
$L_{\lambda_3}$ determines on which subset of $\ZZ$ we define $\mathcal{P}_{\lambda_3}(\omega_{\lambda_3})$.

\smallskip

\item \eqref{for:289} and \eqref{for:290} show that both $P_{1,\lambda_3,\tau}$ and $P_{2,\lambda_3,\tau}$ are not only distributions but are, in fact, inside $ \mathcal{O}$.
Consequently, $\mathfrak{D}_{\lambda_3,\tau}$ is inside $ \mathcal{O}$.

If $|\lambda_3(d\textbf{a}_1)|=1$, we need to use decay of matrix coefficients to show the convergence of $P_{1,\lambda_3,\tau}$. In general, $P_{1,\lambda_3,\tau}\notin \mathcal{O}$. However, $P_{2,\lambda_3,\tau}$ is still inside $ \mathcal{O}$ by using super-exponential decay.

  \smallskip

\item\label{for:295} Computation in \eqref{for:282} shows that super-exponential decay occurs only when  calculating $P_{2,\lambda_3,\tau}$, where the sum is taken over $L_{\lambda_3}^c$.

\smallskip

  \item\label{for:532} If $\dim \mathfrak{g}_{\lambda_3}>1$, to define the invariant distributions, we need to fix a basis set $\{u_1,\cdots, u_{\dim \mathfrak{g}_{\lambda_3}}\}$ of $\mathfrak{g}_{\lambda_3}$. Thus the the invariant distributions are defined as:
      \begin{align*}
 \mathfrak{D}_{\lambda_3,u_i, \tau}:=\sum_{n=-\infty}^\infty 2^{-(n+1)} \pi_{u_i}(f_1\circ \tau^{-1})(\omega_{\lambda_3}\circ z_1^{n}),\qquad \forall\,\tau>0,
\end{align*}
   where $1\leq i\leq \dim \mathfrak{g}_{\lambda_3}$.

  \smallskip
  \item\label{for:312} Using  general notations, $\varphi_{\lambda_3}$ is $\infty$ differentiable along $\mathfrak{g}_{\lambda}$, if $\lambda\in \Delta_{z_1}^{-}\cup \Delta_{z_1}^{0}$, i.e., $(\lambda_3,\lambda)$ is a NOT poor pair (see \eqref{sect:3} of Section \ref{sec:19}). This shows that $\varphi_{\lambda_3}$ is $\infty$ differentiable along
      $\mathfrak{g}^{\delta(\mu)}_{z_1}\oplus\mathfrak{g}^0_{z_1}$.

$\varphi_{\lambda_3}$ loses differentiability along $\mathfrak{g}_{\lambda}$, if $\lambda\in \Delta_{z_1}^{+}$, i.e., $(\lambda_3,\lambda)$ is a poor pair. This shows that $\varphi_{\lambda_3}$ is $\infty$ loses differentiability  along
      $\mathfrak{g}^{-\delta(\mu)}_{z_1}$. Thus the invariant distributions are defined as:
      \begin{align}\label{for:309}
 \mathfrak{D}_{\lambda_3,\lambda,u,\tau,\ZZ}(\omega_{\lambda_3}):=\sum_{n=-\infty}^\infty 2^{-(n+1)} \pi_{u}(f_1\circ \tau^{-1})(\omega_{\lambda_3}\circ z_1^{n})
\end{align}
  for any $\tau>0$ where $(\lambda_3,\lambda)$ is a poor pair and $u$ is inside $\mathfrak{g}_{\lambda}$. Of course, we need to choose sufficiently many $u$ such that their linear span is $\mathfrak{g}^{-\delta(\mu)}_{z_1}$.

  We point out that this is the reason we introduce the notation $\mathfrak{D}_{\mu,\lambda,u,\tau,\mathbb{S}}(\xi)$ in \eqref{for:285} of Section \ref{sec:39} for general cases.
  The convergence is proved by using \eqref{for:501} of Section \ref{sec:61}.

  \smallskip

 \item\label{for:311} We point out that
 \begin{align*}
  P_{2,\lambda_3,\tau}=\mathfrak{D}_{\lambda_3,\lambda_3,u,\tau,L_{\lambda_3}^c}(\omega_{\lambda_3}).
 \end{align*}
 Using \eqref{for:295}, we emphasize that super-exponential decay occurs only when  calculating $\mathfrak{D}_{\mu,\lambda,u,\tau,\mathbb{S}}(\xi)$ for $\mathbb{S}\subseteq L_\mu^c$ (see \eqref{for:502} of Section \ref{sec:61}).
\end{enumerate}

\end{remark}

Below, we will show that the vanishing of the invariant
distributions implies a smooth solution. Suppose $\mathfrak{D}_{\lambda_3,\tau}=0$ (see \eqref{for:331}) for some $\tau>0$.
The observation \hyperlink{o.3}{(**)} suggests that it suffices to show
the differentiability of $\varphi_{\lambda_3}$ along $\mathfrak{g}_{\lambda_3}$. More precisely, fixing $0\neq u\in \mathfrak{g}_{\lambda_3}$, we need to show that there is $p>0$ such that
for any $n>0$ and any $\psi \in \mathcal{O}^\infty$, we have
\begin{align*}
 |\langle\varphi_{\lambda_3},\, u^n\psi\rangle|\leq C_{n,\varphi_{\lambda_3}}\norm{\psi}_p.
\end{align*}
Similar to \eqref{for:279} we decompose $\psi$ as follows:
\begin{align*}
 \psi=\pi_u(f\circ \tau^{-1})\psi+\pi_u(f_1\circ \tau^{-1})\psi.
\end{align*}
Next, we will estimate $\Big|\big\langle\varphi_{\lambda_3},\, u^n\pi_{u}(f\circ \tau^{-1})\psi \big\rangle\Big|$ and $\Big|\big\langle\varphi_{\lambda_3},\, u^n\pi_{u}(f_1\circ \tau^{-1})\psi \big\rangle\Big|$, respectively.

For the former, we have
\begin{align}\label{for:288}
 &\Big|\big\langle\varphi_{\lambda_3},\, u^n\pi_{u}(f\circ \tau^{-1})\psi \big\rangle\Big|\notag\\
 &\overset{\text{(1)}}{\leq}\norm{\omega_{\lambda_3}}\cdot \norm{u^n\pi_{u}(f\circ \tau^{-1})\psi}\notag\\
  &\overset{\text{(2)}}{\leq}\norm{\omega_{\lambda_3}}\cdot C_{n}\tau^n\norm{\psi}.
\end{align}
Here in $(1)$ we use cauchy schwarz inequality;   in $(2)$ we use \eqref{for:128} of Lemma \ref{cor:1}.

\smallskip
\emph{Note}. Similar to \eqref{for:278}, \eqref{for:288} also uses the fact that the space of $\tau$-truncations fails to perform the testing job.

\smallskip

For the latter, we have
\begin{align*}
 \big\langle\varphi_{\lambda_3}&,\, u^n\pi_{u}(f_1\circ \tau^{-1})\psi \big\rangle\overset{\text{(1)}}{=}\big\langle \pi_{u}(f_1\circ \tau^{-1})\varphi_{\lambda_3},\, u^n\psi \big\rangle\notag\\
 &\overset{\text{(2)}}{=}- \sum_{m\geq0} 2^{-(m+1)}\big\langle \pi_{u}(f_1\circ \tau^{-1}) (\omega_{\lambda_3}\circ z_1^{m}),\, u^n\psi \big\rangle.
  \end{align*}
Here in $(1)$ we use \eqref{for:10} and \eqref{for:93} of Section \ref{sec:27}; in $(2)$ we use \eqref{for:286}.

Since $\mathfrak{D}_{\lambda_3,\tau}=0$, it follows that
\begin{align*}
 &- \sum_{m\geq0} 2^{-(m+1)}\big\langle \pi_{u}(f_1\circ \tau^{-1}) (\omega_{\lambda_3}\circ z_1^{m}),\, u^n\psi \big\rangle\\
 &=\sum_{m<0} 2^{-(m+1)}\big\langle \pi_{u}(f_1\circ \tau^{-1}) (\omega_{\lambda_3}\circ z_1^{m}),\, u^n\psi \big\rangle.
\end{align*}
Now we recall \eqref{for:310} and \eqref{for:311} of Remark \ref{re:13}. It reminds us to use super-exponential decay again to estimate the above sum.

It follows that
\begin{align*}
 &\Big|\big\langle\varphi_{\lambda_3},\, u^n\pi_{u}(f_1\circ \tau^{-1})\psi \big\rangle\Big|\notag\\
  &=\Big| \sum_{m<0} 2^{-(m+1)}\big\langle \pi_{u}(f_1\circ \tau^{-1}) (\omega_{\lambda_3}\circ z_1^{m}),\, u^n\psi \big\rangle\Big|\notag\\
   &\overset{\text{(1)}}{=}\Big| \sum_{m<0} 2^{-(m+1)}\big\langle \pi_{u}(f_1\circ \tau^{-1})u^n (\omega_{\lambda_3}\circ z_1^{m}),\, \psi \big\rangle\Big|\notag\\
   &\overset{\text{(2)}}{=}\Big| \sum_{m<0} 2^{-(m+1)}2^{mn}\big\langle \pi_{u}(f_1\circ \tau^{-1}) \big((u^n\omega_{\lambda_3})\circ z_1^{m}\big),\, \psi \big\rangle\Big|\notag\\
  &\overset{\text{(3)}}{=}\Big| \sum_{m<0} 2^{-(m+1)}2^{mn}\big\langle \big(\pi_{u}(f_1\circ \tau_m^{-1})(u^n\omega_{\lambda_3})\big)\circ (z_1^m),\, \psi \big\rangle\Big|
  \end{align*}
  \begin{align}\label{for:287}
  & \overset{\text{(4)}}{\leq} \sum_{m<0} 2^{-(m+1)}2^{mn}\Big\|\pi_{u}(f_1\circ \tau_m^{-1})(u^n\omega_{\lambda_3})\Big\|\cdot \norm{\psi}\notag\\
 & \overset{\text{(5)}}{\leq} \sum_{m<0} 2^{-(m+1)}2^{mn} \cdot C\tau_m^{-s}\norm{u^n\omega_{\lambda_3}}_{s}\cdot \norm{\psi }\notag\\
 &=\sum_{m<0} 2^{-(m+1)}2^{mn} \cdot C2^{ms}\tau^{-s}\norm{u^n\omega_{\lambda_3}}_{s}\cdot \norm{\psi }
\end{align}
for any $s>0$.

Here in $(1)$ we use \eqref{for:10} of Section \ref{sec:27}; in $(2)$ we use the fact that $\lambda_3(dz_1)=2$; in $(3)$ we use \eqref{for:281} and note that $\tau_m=2^{-m}\tau$; in $(4)$ we use cauchy schwarz inequality and the fact that
\begin{align*}
\norm{\varphi}=\norm{\varphi\circ z_1^n},\qquad \forall\,\varphi\in \mathcal{O};
\end{align*}
in $(5)$ we use Lemma \ref{cor:4}.

Choosing $s>1$ it follows from \eqref{for:287} that
\begin{align*}
 \Big|\big\langle\varphi_{\lambda_3}&,\, u^n\pi_{u}(f_1\circ \tau^{-1})\psi \big\rangle\Big|<C_{n,s}\tau^{-s}\norm{\omega_{\lambda_3}}_{n+s}\cdot \norm{\psi }.
\end{align*}
This together with  \eqref{for:288} show that
\begin{align*}
|\langle\varphi_{\lambda_3},\, u^n\psi\rangle|\leq C_{n,\tau,s}\norm{\omega_{\lambda_3}}_{n+s}\cdot \norm{\psi }
\end{align*}
for any $n>0$. This gives the differentiability of $\varphi_{\lambda_3}$ along $\mathfrak{g}_{\lambda_3}$.

\subsubsection{Invariant distributions} In next proposition the crucial upperbounds for matrix coefficients along individual orbits are obtained, which justifies the convergence of
specifically defined distributions.

The following positive constants $c_1,\,\kappa,\,\gamma_1,\,\kappa_1,\,\gamma_2$ (see Section \ref{sec:60}) will appear in the following proposition.
\begin{proposition}\label{po:6}
Suppose $\xi\in \mathfrak{g}_\mu(\mathcal{O}^m)$, $m\geq 0$. Also suppose $u\in \mathfrak{g}_{\lambda}\bigcap \mathcal{V}^1$, $\lambda\in \Delta$ and $\tau\geq1$.
  Then:
\begin{enumerate}
  \item\label{for:152} for any $n,\,l\in\ZZ$ we have
\begin{align*}
 &\Big\|\pi_{u}(f_1\circ \tau^{-1})\big(\xi\circ (z_1^nz_2^l)\big)\Big\|\\
 &\leq C\tau^{-m}\lambda(d\textbf{a}_1)^{nm}\lambda(d\textbf{a}_2)^{lm} (1+|n|)^{mc_1}(1+|l|)^{mc_1}\norm{\xi}_m
\end{align*}
for any $m>0$;

\smallskip

\item\label{for:158} for any $n,\,l\in\ZZ$, any $g\in \tilde{\mathcal{S}}(\RR)$ and any $\psi\in \mathfrak{g}_{\mu}(\mathcal{O}^\kappa)$ we have
\begin{align*}
\Big\|\big\langle  \xi\circ (z_1^{-n}z_2^{-l}),\,\pi_{u}(g\circ \tau^{-1})\psi\big\rangle\Big\|\leq C_ge^{-\gamma_1 (|n|+|l|)}\norm{\xi}_{\kappa}\norm{\psi}_{\kappa};
\end{align*}

\item\label{for:140} for any $n\in L_\mu$ and any $g\in \tilde{\mathcal{S}}(\RR)$ we have
\begin{align*}
 &\Big\|\big\langle (dz_1|_{\mathfrak{g}_{\mu}})^{n-1} \pi_{u}(g\circ \tau^{-1})(\xi\circ z_1^{-n}),\,\psi\big\rangle\Big\|\\
 &\leq C_g(1+|n-1|)^{\dim\mathcal{V}} e^{-\gamma_1 |n|}\norm{\xi}_{\kappa}\norm{\psi}_{\kappa}
\end{align*}
for any $\psi\in \mathfrak{g}_{\mu}(\mathcal{O}^\kappa)$;

  \smallskip

  \item\label{for:133} $\mathcal{P}_{\mu}(\xi)$ (see \eqref{for:284} of Section \ref{sec:39}) is a distribution  satisfying
\begin{align*}
 \|\langle \mathcal{P}_{\mu}(\xi),\,\psi\rangle\|\leq C\norm{\xi}_{\kappa}\norm{\psi}_{\kappa}
\end{align*}
for any $\psi\in \mathfrak{g}_{\mu}(\mathcal{O}^\kappa)$;

\smallskip

\item\label{for:350} for any $n\in L_\mu^c$ we have
\begin{align*}
 &\big\|\langle (dz_1|_{\mathfrak{g}_{\mu}})^{n-1} \pi_{u}(f_1\circ \tau^{-1})(\xi\circ z_1^{-n}),\,\psi\rangle\big\|\notag\\
 &\leq C_{m}\tau^{-m}e^{-\frac{1}{4}\kappa_1 |n|m}\norm{\xi}_m\norm{\psi}
\end{align*}
for any $\psi\in \mathfrak{g}_{\mu}(\mathcal{O})$ and $m\geq \kappa$;

\smallskip
\item\label{for:354} for any $n\in\ZZ$ we have
\begin{align*}
 &\big\|\langle (dz_1|_{\mathfrak{g}_{\mu}})^{n-1} \pi_{u}(f_1\circ \tau^{-1})(\xi\circ z_1^{-n}),\,\psi\rangle\big\|\notag\\
 &\leq Ce^{-\gamma_2 |n|}\norm{\xi}_\kappa\norm{\psi}_\kappa
\end{align*}
for any $\psi\in \mathfrak{g}_{\mu}(\mathcal{O}^\kappa)$;

\smallskip

 \item\label{for:153}  suppose $(\mu,\,\lambda)$ is a poor pair. Then for any subset $\mathbb{S}\subseteq\ZZ$, $\mathfrak{D}_{\mu,\lambda,u,\tau,\mathbb{S}}(\xi)$ (see \eqref{for:285} of Section \ref{sec:39}) is a distribution satisfying
\begin{align*}
 \|\langle \mathfrak{D}_{\mu,\lambda,u,\tau,\mathbb{S}}(\xi),\,\psi\rangle\|\leq C\norm{\xi}_{\kappa}\norm{\psi}_{\kappa}
\end{align*}
for any $\psi\in \mathfrak{g}_{\mu}(\mathcal{O}^\kappa)$. Moreover, we have
\begin{align*}
 \norm{\mathfrak{D}_{\mu,\tau,\mathbb{S}}(\xi)}\leq C\norm{\xi}_{\kappa}
\end{align*}
(see \eqref{for:294} of Section \ref{sec:39});

\smallskip
\item\label{for:156} for any $\psi\in \mathfrak{g}_\mu(\mathcal{O})$ we have
\begin{align*}
 &\big\|\langle \mathfrak{D}_{\mu,w,\tau,L_\mu^c}(\xi),\,\psi\rangle\big\|\leq C_m\tau^{-m}\norm{\xi}_m\norm{\psi},\qquad \forall\,m\geq\kappa.
\end{align*}

\end{enumerate}

\end{proposition}
\begin{remark}\label{re:8} Comparing \eqref{re:7} and \eqref{for:7} of Lemma \ref{cor:1}, we conclude that
 if we remove the assumption that $\tau\geq1$, the results in Proposition \ref{po:6} still hold, although the constants may now depend on $\tau$.

\end{remark}
\begin{observation}\label{ob:2} In \eqref{for:152} letting $l=0$ we have
\begin{align*}
 \big\|\pi_{u}(f_1\circ \tau^{-1})\big(\xi\circ z_1^n\big)\big\|&\leq C\tau^{-m}\lambda(d\textbf{a}_1)^{nm} (1+|n|)^{mc_1}\norm{\xi}_m
 \end{align*}
for any $m>0$.

Hence, we have super exponential decay when $n<0$ (resp. $n>0$)   if $\lambda(d\textbf{a}_1)>1$ (resp. $\lambda(d\textbf{a}_1)<1$). This shows that if $(\mu,\,\lambda)$ is a poor pair, then super exponential decay occurs when $n\in L_\mu^c$.

\end{observation}
\begin{proof} We recall that the Lyapunov
exponents and the Lyapunov subspace of $z_1$ match those of $\textbf{a}_1$ (see \eqref{sect:10} of Section \ref{sec:19}). We only consider the case of $\mu\in \Delta_{\textbf{a}_1}^{+}\cup\Delta_{\textbf{a}_1}^{0}$. The proof for the case of $\mu\in \Delta_{\textbf{a}_1}^{-}$ follows in a similar way.

\smallskip
\eqref{for:152}: We have
\begin{align*}
&\Big\|\pi_{u}(f_1\circ \tau^{-1})\big(\xi\circ (z_1^nz_2^l)\big)\Big\|\\
&\overset{\text{(1)}}{=}\Big\|\pi_{\mathfrak{u}_{n,l}}\big(f_1\circ (\tau^{-1}\norm{dz_1^ndz_2^l(u)})\big)\xi\Big\|\\
&\overset{\text{(2)}}{\leq} C(\tau\norm{dz_1^ndz_2^l(u)}^{-1})^{-m}\norm{\xi}_m\\
&\overset{\text{(3)}}{\leq} C\tau^{-m}|\lambda(d\textbf{a}_1)|^{nm}|\lambda(d\textbf{a}_2)|^{lm}(1+|n|)^{mc_1}(1+|l|)^{mc_1}\norm{\xi}_m.
\end{align*}
Here in $(1)$ we use Corollary \ref{le:15}; in $(2)$ we use Lemma \ref{cor:4}; in $(3)$ we use \eqref{for:138} of Section \ref{sec:19}.  Then we get the result.

\smallskip
\eqref{for:158}: For any $\psi\in \mathfrak{g}_{\mu}(\mathcal{O}^\kappa)$ we have
\begin{align*}
 &\Big\|\big\langle  \xi\circ (z_1^{-n}z_2^{-l}),\,\pi_{u}(g\circ \tau^{-1})\psi\big\rangle\Big\|\\
 &\overset{\text{(1)}}{\leq} Ce^{-\gamma_1 (|n|+|l|)}\norm{\xi}_{C^1}\norm{\pi_{u}(g\circ \tau^{-1})\psi}_{C^1}\\
 &\overset{\text{(2)}}{\leq} C_1 e^{-\gamma_1 (|n|+|l|)}\norm{\xi}_{\kappa}\norm{\pi_{u}(g\circ \tau^{-1})\psi}_{\kappa}\\
 &\overset{\text{(3)}}{\leq} C_ge^{-\gamma_1 (|n|+|l|)}\norm{\xi}_{\kappa}\norm{\psi}_{\kappa}.
\end{align*}
Here in $(1)$ we use \eqref{sec:45} of Section \ref{sec:19}; in $(2)$ we use Sobolev embedding inequality  \eqref{for:179} of Section \ref{sec:19} and the fact $\kappa\geq \beta+1$ (see \eqref{sect:2} of of Section \ref{sec:19}); in $(3)$ we use \eqref{for:7} of Lemma \ref{cor:1}.

\smallskip

\eqref{for:140}: We note that $L_\mu=\ZZ^-\cup\{0\}$. Then for any $n\in L_\mu$, i.e., $n\leq0$ we have
\begin{align*}
  &\big\|\langle (dz_1|_{\mathfrak{g}_{\mu}})^{n-1} \pi_{u}(g\circ \tau^{-1})(\xi\circ z_1^{-n}),\,\psi\rangle\big\|\notag\\
  &\leq \norm{(dz_1|_{\mathfrak{g}_{\mu}})^{n-1}}\big\|\langle  \pi_{u}(g\circ \tau^{-1})(\xi\circ z_1^{-n}),\,\psi\rangle\big\|\notag\\
  &\overset{\text{(1)}}{=} \norm{(dz_1|_{\mathfrak{g}_{\mu}})^{n-1}}\big\|\langle  \xi\circ z_1^{-n},\,\pi_{u}(g\circ \tau^{-1})\psi\rangle\big\|\notag\\
  &\overset{\text{(2)}}{\leq} C_g|\mu(d\textbf{a}_1)|^{n-1}(1+|n-1|)^{c_1} e^{-\gamma_1 |n|}\norm{\xi}_{\kappa}\norm{\psi}_{\kappa}\notag\\
  &\overset{\text{(3)}}{\leq} C_g(1+|n-1|)^{c_1} e^{-\gamma_1 |n|}\norm{\xi}_{\kappa}\norm{\psi}_{\kappa}.
\end{align*}
Here in $(1)$ we use \eqref{for:93} of Section \ref{sec:27}; $(2)$ we use \eqref{for:138} of Section \ref{sec:19} and \eqref{for:158}; in $(3)$ we use $|\mu(d\textbf{a}_1)|\geq1$.

\smallskip

\eqref{for:133}: In \eqref{for:140}  by letting $g=1$ and noting that $\pi_{u}(1\circ \tau^{-1})=I$, we have
\begin{align*}
 \big\|\langle (dz_1|&_{\mathfrak{g}_{\mu}})^{n-1} \xi\circ z_1^{-n},\,\psi\rangle\big\|\leq C(1+|n-1|)^{\dim\mathcal{V}} e^{-\gamma_1 |n|}\norm{\xi}_{\kappa}\norm{\psi}_{\kappa}
\end{align*}
for any $n\in L_\mu$ and any $\psi\in \mathfrak{g}_{\mu}(\mathcal{O}^\kappa)$. The above inequality implies the result.

\smallskip
\eqref{for:350}: For any $n\in L_\mu^c$, , i.e., $n>0$ we have
\begin{align}\label{for:155}
 &\big\|\langle (dz_1|_{\mathfrak{g}_{\mu}})^{n-1} \pi_{u}(f_1\circ \tau^{-1})(\xi\circ z_1^{-n}),\,\psi\rangle\big\|\notag\\
 &\overset{\text{(1)}}{\leq} \big\|(dz_1|_{\mathfrak{g}_{\mu}})^{n-1} \pi_{u}(f_1\circ \tau^{-1})(\xi\circ z_1^{-n})\big\|\norm{\psi}\notag\\
&\overset{\text{(2)}}{\leq} C|\mu(d\textbf{a}_1)|^{n-1}(1+|n-1|)^{c_1}\big\|\pi_{u}(f_1\circ \tau^{-1})(\xi\circ z_1^{-n})\big\|\norm{\psi}\notag\\
 &\overset{\text{(3)}}{\leq} C_1\tau^{-m}|\mu(d\textbf{a}_1)|^{n-1}|\lambda(d\textbf{a}_1)|^{-nm}(1+|n-1|)^{(m+1)c_1}\norm{\xi}_m\norm{\psi}\notag\\
 &\overset{\text{(4)}}{\leq} C_1\tau^{-m}e^{\kappa_2(n-1)}e^{-\kappa_1 nm}(1+|n-1|)^{(m+1)c_1}\norm{\xi}_m\norm{\psi}\notag\\
  &\overset{\text{(5)}}{\leq} C_{m}\tau^{-m}e^{-\frac{1}{4}\kappa_1 nm}\norm{\xi}_m\norm{\psi}
\end{align}
for any $m\geq\kappa$.

Here in $(1)$ we use cauchy schwarz inequality; in $(2)$ use \eqref{for:138} of Section \ref{sec:19}; in $(3)$ we use \eqref{for:152} and choose $m\geq\kappa$; in $(4)$ we recall notations in \eqref{sect:2} of Section \ref{sec:19} and the fact that $|\lambda(d\textbf{a}_1)|>e^{\kappa_1}>1$; in $(5)$ we note that
\begin{align*}
 &e^{\kappa_2(n-1)}e^{-\kappa_1nm}(1+|n-1|)^{(m+1)c_1}\\
 &=e^{\kappa_2(n-1)}e^{-\frac{1}{2}\kappa_1 nm} e^{-\frac{1}{2}\kappa_1 nm}(1+|n-1|)^{(m+1)c_1}\\
 &\leq (e^{\kappa_2(n-1)}e^{-\frac{1}{2}\kappa_1 n \kappa})e^{-\frac{1}{2}\kappa_1 nm}(1+|n-1|)^{(m+1)c_1}\\
 &\overset{\text{(a)}}{\leq} Ce^{-\frac{1}{2}\kappa_1 nm}(1+|n-1|)^{(m+1)c_1}\\
 &\leq C_me^{-\frac{1}{4}\kappa_1 nm}.
\end{align*}
Here in $(a)$ we use \eqref{for:139} of Section \ref{sec:19}. Then we get the result.

\smallskip

\eqref{for:354}: If $n\in L_\mu$, by  \eqref{for:140}, letting $g=f_1$, we have
\begin{align*}
 &\Big\|\big\langle (dz_1|_{\mathfrak{g}_{\mu}})^{n-1} \pi_{u}(f_1\circ \tau^{-1})(\xi\circ z_1^{-n}),\,\psi\big\rangle\Big\|\\
 &\leq C_g(1+|n-1|)^{\dim\mathcal{V}} e^{-\gamma_1 |n|}\norm{\xi}_{\kappa}\norm{\psi}_{\kappa}
\end{align*}
for any $\psi\in \mathfrak{g}_{\mu}(\mathcal{O}^\kappa)$.

If $n\in L_\mu^c$, by \eqref{for:350}, letting $\tau=1$ and $m=\kappa$, we have
\begin{align*}
 &\big\|\langle (dz_1|_{\mathfrak{g}_{\mu}})^{n-1} \pi_{u}(f_1\circ \tau^{-1})(\xi\circ z_1^{-n}),\,\psi\rangle\big\|\notag\\
 &\leq Ce^{-\frac{1}{4}\kappa_1 |n|\kappa}\norm{\xi}_\kappa\norm{\psi}.
\end{align*}
Recall $\gamma_2=\min\{\frac{1}{2}\gamma_1,\,\frac{1}{4}\kappa_1 \kappa\}$. The above discussion implies the result.

\smallskip

\eqref{for:153}: It follows directly from \eqref{for:354}.

\smallskip

\eqref{for:156} is a direct
 consequence of \eqref{for:350}.
\end{proof}

The implications of specifically defined distributions in Proposition \ref{po:6} are clear from the next result: $\mathcal{P}_{\mu}(\omega_\mu)$ is a distributional solution of equation \eqref{for:136} and $\mathfrak{D}_{\mu,\lambda,u,\tau,\ZZ}(\omega_\mu)$ are invariant distributions for solving \eqref{for:136}.
\begin{lemma}\label{le:5} Suppose $\omega_\mu\in \mathfrak{g}_\mu(\mathcal{O}^\kappa)$. Then:

\begin{enumerate}
  \item\label{for:159} $\mathcal{P}_{\mu}(\omega_\mu)$ is a distributional solution of equation \eqref{for:136};

  \smallskip
  \item\label{for:141} if equation \eqref{for:136} has a solution $\varphi_\mu\in \mathfrak{g}_\mu(\mathcal{O}^\kappa)$,  then $\varphi_\mu$ is unique;

  \smallskip
  \item\label{for:142} if equation \eqref{for:136} has a solution $\varphi_\mu\in \mathfrak{g}_\mu(\mathcal{O}^\kappa)$,  then
  \begin{align*}
   \norm{\mathfrak{D}_{\mu,\tau,\ZZ}(\omega_\mu)}=0,\qquad \forall\,\tau>0.
  \end{align*}
\end{enumerate}
\end{lemma}
\begin{proof}
\eqref{for:159}: A direct consequence \eqref{for:133} of Proposition \ref{po:6}.

\smallskip
\eqref{for:141}: It suffices to show that if: $\xi\in \mathfrak{g}_\mu(\mathcal{O}^\kappa)$ satisfies the equation
\begin{align}\label{for:167}
 \xi\circ z_1-(dz_1|_{\mathfrak{g}_{\mu}})\xi=0,
\end{align}
then $\xi=0$. From \eqref{for:167} we have
\begin{align*}
 \xi\circ z_1^n=(dz_1|_{\mathfrak{g}_{\mu}})^n\xi,\qquad \forall\,n\in\ZZ.
\end{align*}
Suppose $\mu\in \Delta_{\textbf{a}_1}^{+}\cup\Delta_{\textbf{a}_1}^{0}$. For any $n\leq0$ the above equation implies that
\begin{align*}
 \|\langle \xi,\,\xi\rangle\|=\big\|\langle (dz_1|_{\mathfrak{g}_{\mu}})^{-n}(\xi\circ z_1^n),\,\xi\rangle\big\|
 \overset{\text{(1)}}{\leq}
  C(1+|n-1|)^{\dim\mathcal{V}} e^{-\gamma_1 |n|}\norm{\xi}_{\kappa}^2
\end{align*}
Here in $(1)$ we use \eqref{for:140} of Proposition \ref{po:6}  by letting $g=1$ and noting that $\pi_{u}(1\circ \tau^{-1})=I$.

Letting $n\to-\infty$, we have $\langle \xi,\,\xi\rangle=0$. Then we get $\xi=0$. The proof for the case of $\mu\in \Delta_{z_1}^{-}$ follows in a similar way.

\smallskip
\eqref{for:142}: It suffices to show that: for any $\lambda\in \Delta$ such that $(\mu,\,\lambda)$ is a poor pair and any $u\in \mathfrak{g}_{\lambda}\bigcap \mathcal{V}^1$,
$\mathfrak{D}_{\mu,\lambda,u,\tau,\ZZ}(\omega_\mu)=0$.

\eqref{for:153} of Proposition \ref{po:6} and Remark \ref{re:8} show that both $\mathfrak{D}_{\mu,\lambda,u,\tau,\ZZ}(\omega_\mu)$ and
$\mathfrak{D}_{\mu,\lambda,u,\tau,\ZZ}(\varphi_\mu)$ are well defined. Applying the operator $\mathfrak{D}_{\mu,\lambda,u,\tau,\ZZ}$ on both sides of equation \eqref{for:136}
we have
\begin{align*}
 &\mathfrak{D}_{\mu,\lambda,u,\tau,\ZZ}(\varphi_\mu\circ z_1)-\mathfrak{D}_{\mu,\lambda,u,\tau,\ZZ}\big((dz_1|_{\mathfrak{g}_{\mu}})\varphi_\mu\big)=\mathfrak{D}_{\mu,\lambda,u,\tau,\ZZ}(\omega_\mu).
 \end{align*}
Recalling \eqref{for:285} of Section \ref{sec:39}, the left side is:
 \begin{align*}
  &\sum_{n\in \ZZ} (dz_1|_{\mathfrak{g}_{\mu}})^{n-1} \pi_{u}(f_1\circ \tau^{-1})(\varphi_\mu\circ z_1^{-n+1})\\
  &-\sum_{n\in \ZZ} (dz_1|_{\mathfrak{g}_{\mu}})^{n-1} \pi_{u}(f_1\circ \tau^{-1})\big(((dz_1|_{\mathfrak{g}_{\mu}})\varphi_\mu)\circ z_1^{-n}\big)\\
  &\overset{\text{(1)}}{=}\sum_{n\in \ZZ} (dz_1|_{\mathfrak{g}_{\mu}})^{n-1} \pi_{u}(f_1\circ \tau^{-1})(\varphi_\mu\circ z_1^{-n+1})\\
  &-\sum_{n\in \ZZ} (dz_1|_{\mathfrak{g}_{\mu}})^{n} \pi_{u}(f_1\circ \tau^{-1})(\varphi_\mu\circ z_1^{-n})\\
  &=0.
\end{align*}
Here in $(1)$ we use the fact that for any $n\in\ZZ$
\begin{align}\label{for:296}
 \pi_{u}(f_1\circ \tau^{-1})&\big(((dz_1|_{\mathfrak{g}_{\mu}})\varphi_\mu)\circ z_1^{-n}\big)=\pi_{u}(f_1\circ \tau^{-1})\big((dz_1|_{\mathfrak{g}_{\mu}})(\varphi_\mu\circ z_1^{-n})\big)\notag\\
 &\overset{\text{(1)}}{=}(dz_1|_{\mathfrak{g}_{\mu}})\pi_{u}(f_1\circ \tau^{-1})(\varphi_\mu\circ z_1^{-n}).
\end{align}
Here in $(1)$ we use \eqref{for:73} of Section \ref{sec:27}. We thus finish the proof.

\end{proof}
The following result  shows that the vanishing of these invariant distributions  implies the solvability of \eqref{for:136}. Tame Sobolev estimates of the solution are also formulated.
\begin{proposition}\label{le:17}  Suppose $\omega_\mu\in \mathfrak{g}_\mu(\mathcal{O}^m)$, $m\geq \kappa$.
\begin{enumerate}
  \item\label{for:148} For any $u_i\in (\mathfrak{g}^{\delta(\mu)}_{\textbf{a}_1}\oplus\mathfrak{g}^0_{\textbf{a}_1})\bigcap \mathcal{V}^1$ (see \eqref{sect:3} of Section \ref{sec:19}), $1\leq i\leq l$ where $0\leq l\leq m-\kappa$,  we have
  \begin{align*}
   \big\|\langle u_1u_2\cdots u_l\mathcal{P}_{\mu}(\omega_\mu),\,\psi\rangle\big\|\leq C_l \norm{\omega_\mu}_{l+\kappa}\norm{\psi}_{\kappa}
  \end{align*}
  for any $\psi\in \mathfrak{g}_{\mu}(\mathcal{O}^\kappa)$;

  \smallskip
  \item\label{for:293} Suppose $m\geq 2\kappa+4$. If there is $\tau_0>0$ such that
 $\norm{\mathfrak{D}_{\mu,\tau_0,\ZZ}(\omega_\mu)}=0$, then
 \begin{align*}
 \norm{\mathfrak{D}_{\mu,\tau,\ZZ}(\omega_\mu)}=0,\qquad \forall\,\tau>0.
\end{align*}
 Moreover, equation \eqref{for:136} has a solution
$\varphi_\mu\in \mathfrak{g}_{\mu}(\mathcal{O}^{m-\kappa-4})$ with estimates
\begin{align*}
 \norm{\varphi_\mu}_{l}\leq C_{l}\norm{\omega_\mu}_{l+\kappa+4}
\end{align*}
for any $0\leq l\leq m-\kappa-4$;

\end{enumerate}
\end{proposition}
\begin{proof} We only consider the case of $\mu\in \Delta_{\textbf{a}_1}^{+}\cup\Delta_{\textbf{a}_1}^{0}$. The proof for the case of $\mu\in \Delta_{\textbf{a}_1}^{-}$ follows in a similar way.

Let $\varphi_\mu=\mathcal{P}_{\mu}(\omega_\mu)$. \eqref{for:159} of Lemma \ref{le:5} shows that $\varphi_\mu$ is a distributional solution.

\smallskip
\eqref{for:148}: We note that
\begin{align}\label{for:317}
 \norm{dz_1^{-n}|_{\mathfrak{g}^{-}_{\textbf{a}_1}\oplus\mathfrak{g}^0_{\textbf{a}_1}}}\leq C(1+|n|)^{c_1},\quad \forall\,\,n\in L_\mu=\ZZ^-\cup\{0\}
\end{align}
by using  \eqref{for:138} of Section \ref{sec:19}.

For any $u_i\in (\mathfrak{g}^{-}_{\textbf{a}_1}\oplus\mathfrak{g}^0_{\textbf{a}_1})\bigcap \mathcal{V}^1$ and any $\psi\in \mathfrak{g}_{\mu}(\mathcal{O}^\kappa)$ we have
\begin{align}\label{for:257}
 &\big\|\langle u_1u_2\cdots u_l\mathcal{P}_{\mu}(\omega_\mu),\,\psi\rangle\big\|\notag\\
 &\leq \sum_{n\in L_\mu}\Big\|\big\langle (dz_1|_{\mathfrak{g}_{\mu}})^{n-1} u_1u_2\cdots u_l(\omega_\mu\circ z_1^{-n}),\,\psi\big\rangle\Big\|\notag\\
 &\overset{\text{(1)}}{=}\sum_{n\in L_\mu}\norm{dz_1^{-n} u_1}\cdots \norm{dz_1^{-n} u_l}\Big\|\big\langle (dz_1|_{\mathfrak{g}_{\mu}})^{n-1} \textbf{u}_n(\omega_\mu)\circ z_1^{-n},\,\psi\big\rangle\Big\|\notag\\
 &\overset{\text{(2)}}{\leq}C_l\sum_{n\in L_\mu}(1+|n|)^{l\dim\mathcal{V}}\Big\|\big\langle (dz_1|_{\mathfrak{g}_{\mu}})^{n-1} \textbf{u}_n(\omega_\mu)\circ z_1^{-n},\,\psi\big\rangle\Big\|\notag\\
 &\overset{\text{(3)}}{\leq}\sum_{n\in L_\mu}C_{l,1}|\lambda_1(1+|n|)^{l\dim\mathcal{V}}\notag\\
 &\quad\cdot(1+|n-1|)^{\dim\mathcal{V}}e^{-\gamma_1 |n|}\norm{\textbf{u}_{n}\omega_\mu}_{\kappa}\norm{\psi}_{\kappa}\notag\\
&\leq C_{l,2} \norm{\omega_\mu}_{l+\kappa}\norm{\psi}_{\kappa}
\end{align}
for any $0\leq l\leq m-\kappa$. Here in $(1)$
\begin{align*}
\textbf{u}_{n}=(\norm{dz_1^{-n}(u_1)}\cdots \norm{dz_1^{-n} u_l})^{-1}dz_1^{-n}(u_1)\cdots dz_1^{-n}(u_l);
\end{align*}
in $(2)$ we use \eqref{for:317}; in $(3)$ we use \eqref{for:140} of Proposition \ref{po:6}.

As a result of \eqref{for:257}, we get \eqref{for:148}.

\smallskip

\eqref{for:293}: \eqref{for:148} shows the smoothness of $\varphi_\mu$ along stable or neutral directions of $\textbf{a}_1$.
Next, we firstly show differentiability of
$\varphi_\mu$ along unstable directions of $\textbf{a}_1$.

\smallskip

\noindent\emph{Step 1: Smoothness of $\mathcal{P}_{\mu}(\omega_\mu)$ along unstable directions of $\textbf{a}_1$. }

Suppose $|\lambda(d\textbf{a}_1)|>1$. Then $(\mu,\lambda)$ is a poor pair. For any $u\in \mathfrak{g}_{\lambda}\bigcap \mathcal{V}^1$, any $\psi\in \mathfrak{g}_{\mu}(\mathcal{O}^\kappa)$, we can write
\begin{align*}
 \psi=\pi_{u}(f\circ \tau_0^{-1})\psi+\pi_{u}(f_1\circ \tau_0^{-1})\psi.
\end{align*}
Hence for any $l$, we have
\begin{align*}
\langle \mathcal{P}_{\mu}(\omega_\mu),\,u^l\psi\rangle&=\mathcal{E}_1+\mathcal{E}_2,
\end{align*}
where
\begin{align*}
 \mathcal{E}_1&=\big\langle \mathcal{P}_{\mu}(\omega_\mu),\,u^l\big(\pi_{u}(f\circ \tau_0^{-1})\psi\big) \big\rangle\qquad\text{and}\\
 \mathcal{E}_2&=\big\langle \mathcal{P}_{\mu}(\omega_\mu),\,u^l\big(\pi_{u}(f_1\circ \tau_0^{-1})\psi \big) \big\rangle.
\end{align*}
Next, we estimate $\norm{\mathcal{E}_1}$ and $\norm{\mathcal{E}_2}$ respectively.

 For the former,  we have
\begin{align}\label{for:157}
 \norm{\mathcal{E}_1}&\overset{\text{(1)}}{=} \tau_0^l\big\|\langle \mathcal{P}_{\mu}(\omega_\mu),\,\pi_{u}(f^*_l\circ \tau_0^{-1})\psi\rangle\big\|\notag\\
 &\overset{\text{(2)}}{\leq}C_{\tau_0}\tau_0^l\norm{\omega_\mu}_\kappa\big\|\pi_{u}(g_l\circ \tau_0^{-1})\psi\big\|_\kappa\notag\\
 &\overset{\text{(3)}}{\leq} C_{l,\tau_0}\tau_0^l\norm{\omega_\mu}_\kappa\|\psi\|_\kappa.
\end{align}
Here in $(1)$ we use \eqref{for:41} of Section \ref{sec:27}, where $f^*_{l}(t)=f(t)(t\textrm{i})^{l}$; in $(2)$ we use \eqref{for:133} of Proposition \ref{po:6} and Remark \ref{re:8}; in $(3)$ we use \eqref{re:7} of Lemma \ref{cor:1}.

For the latter, we note that
\begin{align*}
\mathcal{E}_2&\overset{\text{(1)}}{=}\big\langle \mathcal{P}_{\mu}(\omega_\mu),\,\pi_{u}(f_1\circ \tau_0^{-1})(u^l\psi)\big\rangle\\
&\overset{\text{(2)}}{=}\big\langle \pi_{u}(f_1\circ \tau_0^{-1})\mathcal{P}_{\mu}(\omega_\mu),\,u^l\psi\big\rangle\\
&\overset{\text{(3)}}{=}\delta(\mu)\big\langle \mathfrak{D}_{\mu,\lambda,u,\tau_0,L_\mu}(\omega_\mu),\,u^l\psi\big\rangle.
\end{align*}
Here in $(1)$ we use \eqref{for:10} of Section \ref{sec:27}; in $(2)$ we use \eqref{for:93} of Section \ref{sec:27}; in $(3)$ we recall notations in Section \ref{sec:39}
and recall the fact that $(\mu,\lambda)$ is a poor pair. Then we use \eqref{for:296}.

By the assumption $\mathfrak{D}_{\mu,\lambda,u,\tau_0,\ZZ}(\omega_\mu)=0$, it follows that
\begin{align*}
 \mathcal{E}_2=-\delta(\mu)
 \langle \mathfrak{D}_{\mu,\lambda,u,\tau_0,L_\mu^c}(\omega_\mu),\,u^l\psi\rangle.
\end{align*}
We recall Observation \ref{ob:2}. It reminds us to use
super-exponential decay again to estimate $ \langle \mathfrak{D}_{\mu,\lambda,u,\tau_0,L_\mu^c}(\omega_\mu),\,u^l\psi\rangle$. Consequently,
\begin{align*}
\|\mathcal{E}_2\|&=\big\|
 \langle \mathfrak{D}_{\mu,\lambda,u,\tau_0,L_\mu^c}(\omega_\mu),\,u^l\psi\rangle\big\|\\
 &\overset{\text{(1)}}{=}\Big\|
 \sum_{n\in L_\mu^c}\big\langle (dz_1|_{\mathfrak{g}_{\mu}})^{n-1} \pi_{u}(f_1\circ \tau_0^{-1})u^l(\omega_\mu\circ z_1^{-n}),\,\psi\big\rangle\Big\|\\
 &\overset{\text{(2)}}{=}\Big\|
 \sum_{n\in L_\mu^c}\norm{dz_1^{-n} u}^l\big\langle (dz_1|_{\mathfrak{g}_{\mu}})^{n-1} \pi_{u}(f_1\circ \tau_0^{-1})
 \big((u_n^l\omega_\mu)\circ z_1^{-n}\big),\,\psi\big\rangle\Big\|.
\end{align*}
Here in $(1)$ we use \eqref{for:10} of Section \ref{sec:27}; in $(2)$
$u_{n}=\norm{dz_1^{-n}(u)}^{-1}dz_1^{-n}(u)$.

It follows that for any $0\leq l\leq m-\kappa$ we have
\begin{align*}
 \norm{\mathcal{E}_2}&\overset{\text{(1)}}{\leq}
 \sum_{n\geq1}\norm{dz_1^{-n} u}^l\Big\|(dz_1|_{\mathfrak{g}_{\mu}})^{n-1} \pi_{u}(f_1\circ \tau_0^{-1})\big((u_n^l\omega_\mu)\circ z_1^{-n}\big)\Big\|\norm{\psi}\notag\\
 &\leq
 \sum_{n\geq1}\norm{dz_1^{-n}|_{\mathfrak{g}_{\lambda}}}^l\norm{(dz_1|_{\mathfrak{g}_{\mu}})^{n-1}}\Big\| \pi_{u}(f_1\circ \tau_0^{-1})\big((u_n^l\omega_\mu)\circ z_1^{-n}\big)\Big\|\norm{\psi}\notag\\
 &\overset{\text{(2)}}{\leq}C\sum_{n\geq1}|\lambda(d\textbf{a}_1)|^{-ln}(1+|-n|)^{lc_1}|\mu(d\textbf{a}_1)|^{n-1}(1+|n-1|)^{c_1}\notag\\
 &\quad\quad\cdot \Big\| \pi_{u}(f_1\circ \tau_0^{-1})\big((u_n^l\omega_\mu)\circ z_1^{-n}\big)\Big\|\norm{\psi}\notag\\
 &\overset{\text{(3)}}{\leq}C_{\tau_0}\sum_{n\geq1}|\lambda(d\textbf{a}_1)|^{-ln}(1+|-n|)^{lc_1}|\mu(d\textbf{a}_1)|^{n-1}(1+|n-1|)^{c_1}\notag\\
 &\quad\quad\cdot \tau_0^{-\kappa}|\lambda(d\textbf{a}_1)|^{-n\kappa}(1+|-n|)^{\kappa c_1}\norm{u_n^l\omega_\mu}_\kappa\norm{\psi}\notag\\
 &\overset{\text{(4)}}{\leq}C_{l,\tau_0}\tau_0^{-\kappa}\norm{\omega_\mu}_{l+\kappa}\norm{\psi}.
\end{align*}
Here in $(1)$ we use $L_\mu^c=\NN$ and cauchy schwarz inequality; in $(2)$ we use \eqref{for:138} of Section \ref{sec:19};
in $(3)$ we use \eqref{for:152} of Proposition \ref{po:6} and Remark \ref{re:8}; in $(4)$ recalling $|\lambda(d\textbf{a}_1)|>1$ we have
\begin{align*}
|\lambda(d\textbf{a}_1)|^{-ln}(1+|-n|)^{(l+\kappa)c_1}(1+|n-1|)^{c_1}\leq C_l,\qquad \forall\,n\geq1;
\end{align*}
and by \eqref{for:139} of Section \ref{sec:19} we have
\begin{align*}
 |\mu(d\textbf{a}_1)|^{n-1}|\lambda(d\textbf{a}_1)|^{-n\kappa}\leq C|\lambda(d\textbf{a}_1)|^{-\frac{1}{2}n\kappa}.
\end{align*}
It follows from the above discussion that
\begin{align}\label{for:160}
&\big\|\langle \mathcal{P}_{\mu}(\omega_\mu),\,u^l\psi\rangle\big\|\leq \norm{\mathcal{E}_1}+\norm{\mathcal{E}_2}\notag\\
&\leq C_{l,\tau_0}\tau_0^l\norm{\omega_\mu}_\kappa\|\psi\|_\kappa+C_{l,\tau_0}\tau_0^{-\kappa}\norm{\omega_\mu}_{l+\kappa}\norm{\psi}
\end{align}
for any $0\leq l\leq m-\kappa$.

\smallskip

\noindent\emph{Step 2: We show that $\mathcal{P}_{\mu}(\omega_\mu)\in \mathfrak{g}_{\mu}(\mathcal{O}^{m-\kappa-4})$}.

Combining \eqref{for:148} and \eqref{for:160}, then using  Lemma \ref{le:3},  we conclude that $\mathcal{P}_{\mu}(\omega_\mu)\in \mathfrak{g}_{\mu}(\mathcal{O}^{m-\kappa-4})$ with estimates
\begin{align}\label{for:161}
 \norm{\mathcal{P}_{\mu}(\omega_\mu)}_{l}\leq C_{l,\tau_0}\tau_0^{l+4}\norm{\omega_\mu}_{\kappa}+C_{l,\tau_0}\tau_0^{-\kappa}\norm{\omega_\mu}_{l+\kappa+4}
\end{align}
for any $0\leq l\leq m-\kappa-4$.

\smallskip

\noindent\emph{Step 3: We show that}
\begin{align}\label{for:298}
 \norm{\mathfrak{D}_{\mu,\tau,\ZZ}(\omega_\mu)}=0,\qquad \forall\,\tau>0.
\end{align}
Since $m\geq 2\kappa+4$, it follows from \emph{Step 2} that $\varphi_\mu=\mathcal{P}_{\mu}(\omega_\mu)\in \mathfrak{g}_{\mu}(\mathcal{O}^{\kappa})$. Then by \eqref{for:142} of Lemma \ref{le:5}, we get the result.

\smallskip

\noindent\emph{Conclusion}: \eqref{for:298} gives  $\norm{\mathfrak{D}_{\mu,1,\ZZ}(\omega_\mu)}=0$. This  means we can replace $\tau_0$ by $1$
in \emph{Step 1-2}. Then \eqref{for:161} becomes
\begin{align*}
 \norm{\mathcal{P}_{\mu}(\omega_\mu)}_{l}\leq C_{l,1}\norm{\omega_\mu}_{\kappa}+C_{l,1}\norm{\omega_\mu}_{l+\kappa+4}\leq C_{l,2}\norm{\omega_\mu}_{l+\kappa+4}
\end{align*}
for any $0\leq l\leq m-\kappa-4$. Hence we finish the proof.

\end{proof}

\begin{corollary}\label{cor:3} Suppose $\omega\in \mathcal{V}(\mathcal{O}^m)$, $m\geq 2\kappa+4$. If equation \eqref{for:135} has a solution $\varphi\in \mathcal{V}(\mathcal{O}^\kappa)$,  then $\varphi\in \mathcal{V}(\mathcal{O}^{m-\kappa-4})$ with estimates
\begin{align*}
 \norm{\varphi}_{l}\leq C_{l}\norm{\omega}_{l+\kappa+4}
\end{align*}
for any $0\leq l\leq m-\kappa-4$.
\end{corollary}

\begin{proof}
By \eqref{for:142} of  Lemma \ref{le:5} $\norm{\mathfrak{D}_{\mu,\tau,\ZZ}(\omega_\mu)}=0$ for any $\tau>0$ and any $\mu\in\Delta$. By \eqref{for:141} of  Lemma \ref{le:5} $\varphi_\mu$ is the unique solution of \eqref{for:135}. Using \eqref{for:293} of Proposition \ref{le:17}, we get the estimates of $\varphi_\mu$. This implies the estimates of $\varphi$.
\end{proof}

\subsection{Study of almost cocycle equation}\label{sec:25} In this section, we assume $z_i=w_i\cdot \textbf{a}_i$ where $w_i\in V^o\cap B_{\eta}(e)$ (see Section \ref{sec:60}), $i=1,\,2$. We emphasize that we don't require that $z_1$ and $z_2$ commute.

  We study the twisted equation:
\begin{align}\label{for:143}
\varphi\circ z_1-dz_1\varphi=\theta\circ z_2-dz_2\theta+\omega.
\end{align}
where $\varphi,\,\,\theta,\,\omega\in \mathcal{V}(\mathcal{O})$.

Since each $\mathfrak{g}_{\lambda}$, $\lambda\in \Delta$ is invariant under $dz_i$ (see \eqref{sect:10} of Section \ref{sec:19}), $i=1,\,2$,  \eqref{for:143} decomposes into equations of the following forms:
\begin{align}\label{for:164}
\varphi_\mu\circ z_1-(dz_1|_{\mathfrak{g}_{\mu}})\varphi_\mu=\theta_\mu\circ z_2-(dz_2|_{\mathfrak{g}_{\mu}})\theta_\mu+\omega_\mu, \quad\mu\in\Delta
\end{align}
where $\varphi_\mu=\varphi|_{\mathfrak{g}_\mu}$,  $\theta_\mu=\theta|_{\mathfrak{g}_\mu}$ and $\omega_\mu=\omega|_{\mathfrak{g}_\mu}$.

In this part, we show that the norms of $\mathfrak{D}_{\mu,\tau,\ZZ}(\theta_\mu)$ are bounded by the Sobolev norm of $\omega$ of order $\kappa$ (see Proposition \ref{po:2}). These estimates are then used in Section \ref{sec:20}
as a base for constructing the approximate solution.

\subsubsection{Notations}\label{sec:40} The positive constants $\eta, c_1,\,\kappa,\,\gamma_1,\,\gamma_2$ (see Section \ref{sec:60}) will appear in is section. We assume $z_i=w_i\cdot \textbf{a}_i$ where $w_i\in V^o\cap B_{\eta}(e)$.
\begin{enumerate}
  \item\label{for:301} For any $\mu\in\Delta$ set
\begin{align*}
 \mathbb{S}_{\mu}=\{(m,n)\in\ZZ^2: \mu(d\textbf{a}_1)^{n}\mu(d\textbf{a}_2)^{m}
 \leq e^{\frac{1}{2}\gamma_1(|n|+|m|)}\}.
\end{align*}

  \item For any poor pair $(\mu,\lambda)$, we define $L_{\mu,\lambda}$,  a subset of $\ZZ$ as follows:
  \begin{enumerate}
    \item Suppose $\mu(d\textbf{a}_1)=1$. If $\mu(d\textbf{a}_2)\geq 1$ (resp. $\mu(d\textbf{a}_2)< 1$), then  set
\begin{align*}
 L_{\mu,\lambda}=\ZZ^-\bigcup \{0\}\quad (\text{resp. }L_{\mu,\lambda}=\ZZ^+\bigcup \{0\}).
\end{align*}

    \item\label{for:302} Suppose  $\mu(d\textbf{a}_1)\neq1$. So there is $c(\mu,\lambda)>0$ such that $\lambda(d\textbf{a}_1)=\mu(d\textbf{a}_1)^c$. If $\lambda(d\textbf{a}_2)\geq \mu(d\textbf{a}_2)^c$ (resp. $\lambda(d\textbf{a}_2)< \mu(d\textbf{a}_2)^c$),
\begin{align*}
 L_{\mu,\lambda}=\ZZ^+\bigcup \{0\}\quad (\text{resp. }L_{\mu,\lambda}=\ZZ^-\bigcup \{0\}).
\end{align*}

\end{enumerate}

\item\label{for:304} Set $\gamma_3=\min_{(\mu,\lambda)\text{ a poor pair}}\{\frac{1}{4}\gamma_1,\,\frac{1}{5}c(\mu,\lambda)\gamma_1 \kappa\}$ and $c_0=\min\{\frac{4}{7},\,\frac{\gamma_2}{2},\,\frac{\gamma_3}{2}\}$. It is clear that both $\gamma_3$ and $c_0$ are positive.

\end{enumerate}

\subsubsection{Main result}  The main purpose of this section is to prove the following result:
\begin{proposition}\label{po:2} Suppose $\varphi,\,\,\theta,\,\omega\in \mathcal{V}(\mathcal{O}^\kappa)$. Then
\begin{align*}
  \norm{\mathfrak{D}_{\mu,\tau,\ZZ}(\theta_\mu)}\leq C\norm{\log[z_1,z_2]}^{c_0}\max\{\norm{\theta_\mu}_{\kappa+1},\,\norm{\varphi_\mu}_{\kappa}\}+C\norm{\omega_\mu}_\kappa
\end{align*}
for any $\mu\in\Delta$.
\end{proposition}

\begin{remark}\label{re:10} In \eqref{for:143}, if $\omega=0$ and $z_1$ and $z_2$ commute, then it is a a twisted cocycle over the $\ZZ^2$ action generated by $z_1$ and $z_2$.
$\omega$  measures how far the pair $(\varphi,\,\theta)$ is away from a twisted cocycle over $z_1$ and $z_2$.

$\norm{\mathfrak{D}_{\mu,\tau,\ZZ}(\theta_\mu)}$ measures how far $\theta_\mu$ is away from a twisted coboundary. We note that  Sobolev norms of $\omega_\mu$ are usually  ``quadratically small" with respect to those of $\theta_\mu$. This means the estimate in Proposition \ref{po:2} is usually sharper than that  of Proposition \ref{po:6} (see \eqref{for:153}).

In Section \ref{sec:46}, we provide the proof of Proposition \ref{po:2} for the case when $z_1$ and $z_2$ commute.
The complete proof is left for  Appendix \ref{sec:44}.
\end{remark}

\subsubsection{A preparatory step} The following result serves as the first step in proving the Proposition \ref{po:2}:
\begin{lemma}\label{le:2}  For any $\tau\geq1$, any $\lambda\in\Delta$ with $(\mu,\lambda)$ a poor pair, any $u\in \mathfrak{g}_{\lambda}\bigcap \mathcal{V}^1$, any $n\in\ZZ$ and any $m\in L_{\mu,\lambda}$, we have
  \begin{align*}
 \Big\|(dz_2|&_{\mathfrak{g}_{\mu}})^{m-1}(dz_1|_{\mathfrak{g}_{\mu}})^{n-1}\big\langle \xi\circ (z_1^{-n}z_2^{-m}), \,\pi_{u}(f_1\circ \tau^{-1})\psi\big\rangle\Big\|\notag\\
 &\leq Ce^{-\gamma_3(|n|+|m|)}\norm{\xi}_{\kappa}\norm{\psi}_\kappa;
\end{align*}
for any $\xi,\, \psi\in \mathfrak{g}_\mu(\mathcal{O}^\kappa)$.

\end{lemma}

\emph{Note}. We emphasize that Lemma \ref{le:2} holds,  even though $z_1$ and $z_2$ do \textbf{not} commute.

\begin{proof} Set
\begin{align*}
 \mathcal{E}_{n,m}=(dz_2|_{\mathfrak{g}_{\mu}})^{m-1}(dz_1|_{\mathfrak{g}_{\mu}})^{n-1}\big\langle \xi\circ (z_1^{-n}z_2^{-m}), \,\pi_{u}(f_1\circ \tau^{-1})\psi\big\rangle.
\end{align*}
For any $n,\,m\in\ZZ$ we have
\begin{align}\label{for:162}
  \norm{\mathcal{E}_{n,m}} &\leq \norm{(dz_2|_{\mathfrak{g}_{\mu}})^{m-1}}\norm{(dz_1|_{\mathfrak{g}_{\mu}})^{n-1}}
  \Big\|\big\langle \xi\circ (z_1^{-n}z_2^{-m}), \,\pi_{u}(f_1\circ \tau^{-1})\psi\big\rangle\Big\|\notag\\
  &\overset{\text{(a)}}{\leq} C|\mu(d\textbf{a}_1)|^{n-1}|\mu(d\textbf{a}_2)|^{m-1}(1+|n-1|)^{\dim\mathcal{V}}(1+|m-1|)^{\dim\mathcal{V}}\notag\\
  &\cdot\Big\|\big\langle \xi\circ (z_1^{-n}z_2^{-m}), \,\pi_{u}(f_1\circ \tau^{-1})\psi\big\rangle\Big\|.
  \end{align}
Here in $(a)$ we use \eqref{for:138} of Section \ref{sec:19}.

\smallskip

\noindent\emph{Step 1: We show that}: if $(m,n)\in \mathbb{S}_\mu$ (see \eqref{for:301} of Section \ref{sec:40}), then
\begin{align*}
 \norm{\mathcal{E}_{n,m}}\leq Ce^{-\frac{1}{4}\gamma_1(|n|+|m|)}\norm{\xi}_\kappa\norm{\psi}_\kappa.
\end{align*}
For any $(m,n)\in \mathbb{S}_\mu$, by using \eqref{for:162} and \eqref{for:158} of Proposition \ref{po:6} we have
\begin{align}\label{for:144}
 \norm{\mathcal{E}_{n,m}}&\leq C|\mu(d\textbf{a}_1)|^{n-1}|\mu(d\textbf{a}_2)|^{m-1}(1+|n-1|)^{c_1}(1+|m-1|)^{c_1}\notag\\
  &\cdot e^{-\gamma_1 (|n|+|m|)}\norm{\xi}_{\kappa}\norm{\psi}_{\kappa}\notag\\
&\overset{\text{(a)}}{\leq} C_1e^{-\frac{1}{2}\gamma_1(|n|+|m|)}(1+|n|)^{c_1}(1+|m|)^{c_1}\norm{\xi}_\kappa\norm{\psi}_\kappa\notag\\
&\leq C_2e^{-\frac{1}{4}\gamma_1(|n|+|m|)}\norm{\xi}_\kappa\norm{\psi}_\kappa.
\end{align}
Here in $(a)$ we use the assumption for $\mathbb{S}_\mu$.

Hence, we get the result.

\smallskip
\noindent\emph{Step 2: We show that}: if $\mu(d\textbf{a}_1)=1$, then for any $n\in\ZZ$ and $m\in L_{\mu,\lambda}$ we have
\begin{align*}
 \norm{\mathcal{E}_{n,m}}\leq Ce^{-\frac{1}{4}\gamma_1(|n|+|m|)}\norm{\xi}_\kappa\norm{\psi}_\kappa.
\end{align*}
We note that if $\mu(d\textbf{a}_1)=1$, then for any $n\in\ZZ$ and $m\in L_{\mu,\lambda}$, $(m,n)\in \mathbb{S}_\mu$. Then the result follows from \emph{Step 1}.

\smallskip

\emph{Note}. From  now on, we only consider the case of $(m,n)\notin \mathbb{S}_\mu$. Otherwise, the result follows directly from \emph{Step 1}.

\smallskip
\noindent\emph{Step 3: We show that}: if $\mu(d\textbf{a}_1)\neq1$, then for any $n\in\ZZ$ and $m\in L_{\mu,\lambda}$ we have
\begin{align}\label{for:163}
 &\Big\|\pi_{u}(f_1\circ \tau^{-1})\big(\xi\circ (z_1^{-n}z_2^{-m})\big)\Big\|\notag\\
&\leq Ce^{-\frac{1}{2}c\gamma_1\kappa (|n|+|m|)}(1+|n|)^{\kappa c_1}(1+|m|)^{\kappa c_1}\norm{\xi}_\kappa
\end{align}
($c$ is defined in \eqref{for:302} of Section \ref{sec:40}).

From the definition of $L_{\mu,\lambda}$ (see \eqref{for:302} of Section \ref{sec:40}) we note that in this case $(\lambda(d\textbf{a}_2)\mu(d\textbf{a}_2)^{-c})^{m}>1$, which further gives
\begin{align}\label{for:145}
\lambda(d\textbf{a}_1)^{n}&\lambda(d\textbf{a}_2)^{m}=\mu(d\textbf{a}_1)^{cn}\mu(d\textbf{a}_2)^{cm}(\lambda(d\textbf{a}_2)\mu(d\textbf{a}_2)^{-c})^{m}\notag\\
&\geq \mu(d\textbf{a}_1)^{cn}\mu(d\textbf{a}_2)^{cm}\overset{\text{(a)}}{\geq} e^{\frac{1}{2}c\gamma_1 (|n|+|m|)}.
\end{align}
Here in $(a)$ we use $(m,n)\notin \mathbb{S}_\mu$.

 It follows from  \eqref{for:152} of Proposition \ref{po:6} that
\begin{align}
 &\Big\|\pi_{u}(f_1\circ \tau^{-1})\big(\xi\circ (z_1^{-n}z_2^{-m})\big)\Big\|\notag\\
&\leq C\tau^{-\kappa}\lambda(d\textbf{a}_1)^{-n\kappa}\lambda(d\textbf{a}_2)^{-m\kappa}(1+|n|)^{\kappa c_1}(1+|m|)^{\kappa c_1}\norm{\xi}_\kappa\notag\\
&\overset{\text{(a)}}{\leq}Ce^{-\frac{1}{2}c\gamma_1\kappa (|n|+|m|)}(1+|n|)^{\kappa c_1}(1+|m|)^{\kappa c_1}\norm{\xi}_\kappa.
\end{align}
Here in $(a)$ we use \eqref{for:145} and the fact $\tau\geq1$. Then we get the result.

\smallskip
\noindent\emph{Step 4: We show that}:
\begin{align}\label{for:303}
\mu(d\textbf{a}_1)^{n-1}\mu(d\textbf{a}_2)^{m-1}\leq Ce^{\frac{1}{4}c\gamma_1 \kappa(|n|+|m|)},\qquad \forall\,m,n\in\ZZ.
\end{align}
We recall that
\begin{align*}
 \lambda(d\textbf{a}_1)=\mu(d\textbf{a}_1)^c\Rightarrow c=\frac{\log \lambda(d\textbf{a}_1)}{\log \mu(d\textbf{a}_1)}>0.
\end{align*}
From \eqref{for:146} of Section \ref{sec:19} we have
\begin{align*}
 e^{\frac{1}{4}\kappa \gamma_1 c}>\max\{\mu(d\textbf{a}_1), \mu(d\textbf{a}_1)^{-1}\}.
\end{align*}
This implies the result.

\smallskip
\noindent\emph{Step 5: We show that}: if $|\mu(d\textbf{a}_1)|\neq1$, then for any $n\in\ZZ$ and $m\in L_{\mu,\lambda}$ we have
\begin{align*}
 \norm{\mathcal{E}_{n,m}}\leq Ce^{-\frac{1}{5}c\gamma_1 \kappa(|n|+|m|)}\norm{\xi}_\kappa\norm{\psi}.
\end{align*}
We have
\begin{align}\label{for:147}
  \norm{\mathcal{E}_{n,m}}&\overset{\text{(a)}}{\leq}C\mu(d\textbf{a}_1)^{n-1}\mu(d\textbf{a}_2)^{m-1}(1+|n-1|)^{c_1}(1+|m-1|)^{c_1}\notag\\
  &\cdot\Big\|\big\langle \pi_{u}(f_1\circ \tau^{-1})\big(\xi\circ (z_2^{-m}z_1^{-n})\big), \,\psi\big\rangle\Big\|\notag\\
  &\overset{\text{(b)}}{\leq} C_1e^{\frac{1}{4}c\gamma_1 \kappa(|n|+|m|)}(1+|n|)^{c_1}(1+|m|)^{c_1}\notag\\
  &\cdot \big\|\pi_{u}(f_1)\big(\xi\circ (z_2^{-m}z_1^{-n})\big)\big\|\norm{\psi}\notag\\
  &\overset{\text{(d)}}{\leq} C_1e^{\frac{1}{4}c\gamma_1 \kappa(|n|+|m|)}(1+|n|)^{c_1}(1+|m|)^{c_1}\notag\\
  &\cdot Ce^{-\frac{1}{2}c\gamma_1\kappa (|n|+|m|)}(1+|n|)^{\kappa c_1}(1+|m|)^{\kappa c_1}\norm{\xi}_\kappa\norm{\psi}\notag\\
 &\leq  C_2e^{-\frac{1}{5}c\gamma_1 \kappa(|n|+|m|)}\norm{\xi}_\kappa\norm{\psi}.
\end{align}
Here in $(a)$ we use \eqref{for:162} and \eqref{for:93} of Section \ref{sec:27}; in $(b)$ we use
cauchy schwarz inequality and \eqref{for:303}; in  $(d)$ we use \eqref{for:163}.

We recall the definition of $\gamma_3$ (see \eqref{for:304} of Section \ref{sec:40}). Then  \eqref{for:144} and \eqref{for:147} imply the result.
Hence we finish the proof.

\end{proof}

\subsubsection{Proof of Proposition \ref{po:2} if $z_1$ and $z_2$ commute}\label{sec:46}

To show that $\mathfrak{D}_{\mu,\tau,\ZZ}(\theta_\mu)$ can be bounded
by Sobolev norms of $\omega_\mu$, we need to express $\mathfrak{D}_{\mu,\tau,\ZZ}(\theta_\mu)$
in terms of $\omega_\mu$.

Suppose $(\mu,\,\lambda)$ is a poor pair and $u\in \mathfrak{g}_{\lambda}\bigcap \mathcal{V}^1$.
From equation \ref{for:164}, we have
\begin{align}\label{for:166}
0&\overset{\text{(1)}}{=}\mathfrak{D}_{\mu,\lambda,u,\tau,\ZZ}\big(\varphi_\mu\circ z_1-(dz_1|_{\mathfrak{g}_{\mu}})\varphi_\mu\big)\notag\\
&=\mathfrak{D}_{\mu,\lambda,u,\tau,\ZZ}(\theta_\mu\circ z_2)-\mathfrak{D}_{\mu,\lambda,u,\tau,\ZZ}\big((dz_2|_{\mathfrak{g}_{\mu}})\theta_\mu\big)+
\mathfrak{D}_{\mu,\lambda,u,\tau,\ZZ}(\omega_\mu).
\end{align}
Here in $(1)$ we use \eqref{for:142} of Lemma \ref{le:5}.

\eqref{for:166} implies that
\begin{align*}
 &\mathfrak{D}_{\mu,\lambda,u,\tau,\ZZ}(\theta_\mu\circ z_2^{-l})-(dz_2|_{\mathfrak{g}_{\mu}})\mathfrak{D}_{\mu,\lambda,u,\tau,\ZZ}\big(\theta_\mu\circ z_2^{-1-l}\big)\\
 &=-
\mathfrak{D}_{\mu,\lambda,u,\tau,\ZZ}(\omega_\mu\circ z_2^{-1-l})
\end{align*}
for any $l\geq0$.  Iterating this equation with respect to $z_2$ backwards and using $z_1z_2=z_2z_1$ we obtain
\begin{align}\label{for:297}
 &\mathfrak{D}_{\mu,\lambda,u,\tau,\ZZ}(\theta_\mu)-(dz_2|_{\mathfrak{g}_{\mu}})^{l+1}\mathfrak{D}_{\mu,\lambda,u,\tau,\ZZ}\big(\theta_\mu\circ z_2^{-l-1}\big)\notag\\
 &=-
\sum_{i=0}^l(dz_2|_{\mathfrak{g}_{\mu}})^{i}\mathfrak{D}_{\mu,\lambda,u,\tau,\ZZ}(\omega_\mu\circ z_2^{-1-i})
\end{align}
for any $l\geq0$.

\eqref{for:166} also implies that
\begin{align*}
 &\mathfrak{D}_{\mu,\lambda,u,\tau,\ZZ}(\theta_\mu\circ z_2^l)-(dz_2|_{\mathfrak{g}_{\mu}})^{-1}\mathfrak{D}_{\mu,\lambda,u,\tau,\ZZ}(\theta_\mu\circ z_2^{l+1})\\
 &=
(dz_2|_{\mathfrak{g}_{\mu}})^{-1}\mathfrak{D}_{\mu,\lambda,u,\tau,\ZZ}(\omega_\mu\circ z_2^l).
\end{align*}
for any $l\geq0$. Iterating this equation we obtain
\begin{align}\label{for:299}
 &\mathfrak{D}_{\mu,\lambda,u,\tau,\ZZ}(\theta_\mu)-(dz_2|_{\mathfrak{g}_{\mu}})^{-l-1}\mathfrak{D}_{\mu,\lambda,u,\tau,\ZZ}\big(\theta_\mu\circ z_2^{l+1}\big)\notag\\
 &=
\sum_{i=0}^{-l}(dz_2|_{\mathfrak{g}_{\mu}})^{-1+i}\mathfrak{D}_{\mu,\lambda,u,\tau,\ZZ}(\omega_\mu\circ z_2^{-i})
\end{align}
for any $l\geq0$.

\eqref{for:297} and \eqref{for:299} suggest that we can express $\mathfrak{D}_{\mu,\tau,\ZZ}(\theta_\mu)$
in terms of $\omega_\mu$ if at least one of the limits
\begin{align*}
 \lim_{l\to+\infty}(dz_2|_{\mathfrak{g}_{\mu}})^{l+1}\mathfrak{D}_{\mu,\lambda,u,\tau,\ZZ}\big(\theta_\mu\circ z_2^{-l-1}\big)
\end{align*}
or
\begin{align*}
 \lim_{l\to-\infty}(dz_2|_{\mathfrak{g}_{\mu}})^{l+1}\mathfrak{D}_{\mu,\lambda,u,\tau,\ZZ}\big(\theta_\mu\circ z_2^{-l-1}\big)
\end{align*}
is $0$ as a distribution.

Suppose $L_{\mu,\lambda}=\ZZ^+\bigcup \{0\}$. Then we show that
\begin{align}\label{for:355}
 \lim_{l\to\infty}(dz_2|_{\mathfrak{g}_{\mu}})^{l+1}\mathfrak{D}_{\mu,\lambda,u,\tau,\ZZ}\big(\theta_\mu\circ z_2^{-l-1}\big)=0
\end{align}
as a distribution.

Firstly,  using $z_1z_2=z_2z_1$ and $dz_1dz_2=dz_2dz_1$ we note that
\begin{align}\label{for:306}
&(dz_2|_{\mathfrak{g}_{\mu}})^{l-1}\mathfrak{D}_{\mu,\lambda,u,\tau,\ZZ}(\xi\circ z_2^{-l})\notag\\
&=\sum_{n}(dz_1|_{\mathfrak{g}_{\mu}})^{n-1}(dz_2|_{\mathfrak{g}_{\mu}})^{l-1}\pi_{u}(f_1\circ \tau^{-1})\big(\xi\circ (z_2^{-l}z_1^{-n})\big)
\end{align}
for any $l\in\ZZ$.

Then by using \eqref{for:306} and Lemma \ref{le:2}  we have
\begin{align}\label{for:356}
 &\big\|\big\langle (dz_2|_{\mathfrak{g}_{\mu}})^{l-1}\mathfrak{D}_{\mu,\lambda,u,\tau,\ZZ}(\xi\circ z_2^{-l}), \,\psi\big\rangle\big\|\notag\\
  &\leq\sum_{n} Ce^{-\gamma_3(|n|+|l|)}\norm{\theta_\mu}_{\kappa}\norm{\psi}_\kappa\notag\\
  &\leq C_1e^{-\gamma_3|l|}\norm{\theta_\mu}_{\kappa}\norm{\psi}_\kappa.
\end{align}
Letting $l\to\infty$, we get  \eqref{for:355}.

It follows from  \eqref{for:297} that
\begin{align*}
 &\mathfrak{D}_{\mu,\lambda,u,\tau,\ZZ}(\theta_\mu)=-
\sum_{i=0}^\infty(dz_2|_{\mathfrak{g}_{\mu}})^{i}\mathfrak{D}_{\mu,\lambda,u,\tau,\ZZ}(\omega_\mu\circ z_2^{-1-i}),
\end{align*}
which gives
\begin{align*}
 \big\|\big\langle \mathfrak{D}&_{\mu,\lambda,u,\tau,\ZZ}(\theta_\mu), \,\psi\big\rangle\big\|
 =\Big\|\sum_{i\geq0}\big\langle (dz_2|_{\mathfrak{g}_{\mu}})^{i}\mathfrak{D}_{\mu,\lambda,u,\tau,\ZZ}(\omega_\mu\circ z_2^{-1-i}), \,\psi\big\rangle\Big\|\notag\\
 & \overset{\text{(a)}}{\leq }\sum_{i\geq0} Ce^{-c_0|i+1|}\norm{\omega_\mu}_{\kappa}\norm{\psi}_\kappa\leq C_1\norm{\omega_\mu}_{\kappa}\norm{\psi}_\kappa
\end{align*}
for any $\psi\in \mathfrak{g}_\mu(\mathcal{O}^\kappa)$. Here in $(a)$ we use \eqref{for:356}.

Suppose $L_{\mu,\lambda}=\ZZ^-\bigcup \{0\}$. From \eqref{for:299} applying the
same reasoning, we obtain
\begin{align*}
 \mathfrak{D}_{\mu,\lambda,u,\tau,\ZZ}(\theta_\mu)=
\sum_{i\leq0}(dz_2|_{\mathfrak{g}_{\mu}})^{-1+i}\mathfrak{D}_{\mu,\lambda,u,\tau,\ZZ}(\omega_\mu\circ z_2^{-i});
\end{align*}
and thus we also have
\begin{align*}
 \big\|\big\langle \mathfrak{D}_{\mu,\lambda,u,\tau,\ZZ}(\theta_\mu), \,\psi\big\rangle\big\|\leq C\norm{\omega_\mu}_{\kappa}\norm{\psi}_\kappa
\end{align*}
for any $\psi\in \mathfrak{g}_\mu(\mathcal{O}^\kappa)$. Hence we finish the proof.

\subsection{Construction of approximate $C^\infty$ solution}\label{sec:20}
We recall nations in Section \ref{sec:22}. In this part we give an approximate $C^\infty$ solution $\Theta$ of the almost twisted cocycle equations (see \eqref{for:173}) involving every element in the compact generating set $E$. The entire section is dedicated to proving Theorem \ref{th:6}.

 \subsubsection{Basic notations}\label{for:504} The following notations will be used in this part:
 \begin{enumerate}

\item The positive constants $\eta$, $\kappa$, $c_0$ and $\beta$ (see Section \ref{sec:60}) will appear in this section. Let
\begin{align}\label{for:171}
 \varrho=\max\{4+\kappa+\beta,\,9+2\kappa+2\beta,\,\text{\small$\frac{4}{3}$}(2\beta+\kappa+4)+4+\kappa \}.
\end{align}

   \item\label{for:506} Suppose $E'=\{z_v:\,v\in E\}$ is an $\eta$-standard   perturbation of $E$. We recall that we can write $z_v=w_v\cdot
\alpha_v(x)$,  where $w_v\in V^o\cap B_{\eta}(e)$, $v\in E$. We emphasize that we do not require the set $\{z_v: v\in E\}$ to be abelian.

\smallskip

\item
For any $\mathfrak{p}_v\in \mathcal{V}(\mathcal{O}^\infty)$, $v\in E$, we set
\begin{enumerate}
  \item$\norm{\log[z,z]}=\max_{u,v\in E}\norm{\log[z_{v},z_{u}]}$ (see \eqref{sect:9} of Section \ref{sec:19});

  \smallskip
  \item $\norm{\mathfrak{p}}_{C^m}=\max_{v\in E}\{\norm{\mathfrak{p}_{v}}_{C^m}\}$;
  \item for any $u,v\in E$
  \begin{align}\label{for:173}
 \mathbb{D}_{v,u}=\mathfrak{p}_u\circ z_v-dz_v\mathfrak{p}_u-\big( \mathfrak{p}_u\circ z_v-dz_v\mathfrak{p}_u\big);
\end{align}

\item $\norm{\mathbb{D}}_{C^m}=\max_{u,\,v\in E}\{\norm{\mathbb{D}_{v,u}}_{C^m}\}$, for any $m\geq0$. $\norm{\mathfrak{p}}_{m}$ and $\norm{\mathbb{D}}_{m}$ are defined accordingly.
\end{enumerate}

 \end{enumerate}

\subsubsection{Notations in the proof}\label{sec:42}
\begin{enumerate}
  \item For each $\lambda\in \Delta$, we fix a basis
$\{u_{(\lambda,1)},\cdots, u_{(\lambda,x_\chi)}\}$ of $\mathfrak{g}_{\lambda}$, where $x_\lambda=\dim \mathfrak{g}_{\lambda}$.

\smallskip

\item  Set
\begin{align*}
 \Lambda^+&=\{u_{\phi,i}:\phi\in \Delta_{\textbf{a}_1}^{+},\,\,1\leq i\leq x_\phi\}\qquad\text{and}\\
 \Lambda^-&=\{u_{\phi,i}:\phi\in \Delta_{\textbf{a}_1}^{-},\,\,1\leq i\leq x_\phi\}.
\end{align*}
Then $\Lambda^+$ is a basis of $\mathfrak{g}_{\textbf{a}_1}^+$ and $\Lambda^-$ is a basis of $\mathfrak{g}_{\textbf{a}_1}^-$.

\smallskip

\item\label{for:313} We can write
$\Lambda^+=\{\nu_1,\cdots, \nu_{p_+}\}$ (resp. $\Lambda^-=\{\nu_1,\cdots, \nu_{p_-}\}$), where $\nu_n=u_{(\phi_{i_n},j_n)}$, $1\leq n\leq p_+$ (resp. $1\leq n\leq p_-$) in such a way that if $n\leq m$ then
\begin{align*}
 \phi_{i_n}(d\textbf{a}_1)\leq \phi_{i_m}(d\textbf{a}_1)\quad\big(\text{resp.}\quad\phi_{i_n}(d\textbf{a}_1)\geq \phi_{i_m}(d\textbf{a}_1)\big).
\end{align*}
It is clear that
\begin{align}\label{for:507}
 [\nu_n,\nu_m] &\text{ is inside the subspace
linearly spanned by }\notag\\
&\nu_k, \,\, k>\max\{n,\,m\}.
\end{align}

\item \label{for:307} We define:
\begin{align*}
 \mathcal{E}_{+,\tau}&\stackrel{\text{def}}{=}\pi_{\nu_1}(f\circ \tau^{-1})\pi_{\nu_2}(f\circ \tau^{-1})\cdots \pi_{\nu_{p_{+}}}(f\circ \tau^{-1})\\
 (resp. \quad \mathcal{E}_{-,\tau}&\stackrel{\text{def}}{=}\pi_{\nu_1}(f\circ \tau^{-1})\pi_{\nu_2}(f\circ \tau^{-1})\cdots \pi_{\nu_{p_{-}}}(f\circ \tau^{-1})).
\end{align*}
(see Section \ref{sec:27} and \eqref{sect:4} of Section \ref{sec:19}).

We also define:
\begin{align*}
 \mathcal{E}_{+,\tau}'&\stackrel{\text{def}}{=}\pi_{\nu_{p_{+}}}(f\circ \tau^{-1})\cdots\pi_{\nu_2}(f\circ \tau^{-1})\pi_{\nu_1}(f\circ \tau^{-1}), \\
 \mathcal{E}_{+,\tau,k}&\stackrel{\text{def}}{=}\pi_{\nu_k}(f\circ \tau^{-1})\pi_{\nu_2}(f\circ \tau^{-1})\cdots \pi_{\nu_{p_{+}}}(f\circ \tau^{-1}),
\end{align*}
for any $1\leq k\leq p_{+}$.

$\mathcal{E}_{-,\tau}'$ and $\mathcal{E}_{-,\tau,k}$ are defined accordingly.

\smallskip

\item\label{for:314} For any $\mu\in\Delta$,
\begin{align*}
 \mathcal{E}_{\mu,\tau}\stackrel{\text{def}}{=} \mathcal{E}_{-\delta(\mu),\tau},\quad \mathcal{E}'_{\mu,\tau}\stackrel{\text{def}}{=} \mathcal{E}'_{-\delta(\mu),\tau}
 \,\,\text{ and }\,\,\mathcal{E}_{\mu,\tau,k}\stackrel{\text{def}}{=} \mathcal{E}_{-\delta(\mu),\tau,k}
\end{align*}
(see \eqref{sect:3}  of Section \ref{sec:19}). Set
$p_\mu=p_{-\delta(\mu)}$.

\end{enumerate}

 \subsubsection{Main result}

Theorem \ref{th:6} is the central part of the proof of Theorem \ref{th:5}. The entire section is dedicated to proving Theorem \ref{th:6}. We recall notations in Section \ref{for:504}.

\begin{theorem}\label{th:6} Suppose $\mathfrak{p}_v\in \mathcal{V}(\mathcal{O}^\infty)$, $v\in E$. Also suppose $E'=\{z_v:\,v\in E\}$ is an $\eta$-standard   perturbation of $E$.
Then for any $\tau>1$ there exist $\Theta\in \mathcal{V}(\mathcal{O}^\infty)$ and $\mathcal{R}_v\in \mathcal{V}(\mathcal{O}^\infty)$ for each $v\in E$ such that
\begin{align*}
 \mathfrak{p}_v=\Theta\circ z_v-dz_v\Theta+\mathcal{R}_v
\end{align*}
with estimates: for any $v\in E$
\begin{align}\label{for:206}
 \max\{\norm{\Theta}_{C^m},\,\norm{\mathcal{R}_v}_{C^m}\}
 \leq C_{m}(\tau^{m+\varrho}\norm{\mathfrak{p}}_{C^\varrho}+\norm{\mathfrak{p}}_{C^{m+\varrho}})
\end{align}
for any $m\geq0$;  and
\begin{align}\label{for:188}
 \norm{\mathcal{R}_v}_{C^0}&\leq C\tau^{\varrho}\norm{\log[z,z]}^{\frac{c_0}{8}}\norm{\mathfrak{p}}_{C^\varrho}+C\tau^{\varrho}
(\norm{\mathbb{D}}_{C^\varrho})^{\frac{1}{8}}(\norm{\mathfrak{p}}_{{C^\varrho}})^{\frac{7}{8}}\notag\\
 &+C\tau^{\varrho}
\big(C_m\tau^{-m}\norm{\mathfrak{p}}_{C^m}\big)^{\frac{1}{8}}(\norm{\mathfrak{p}}_{{C^\varrho}})^{\frac{7}{8}}+C\norm{\mathbb{D}}_{{C^\varrho}}
\end{align}
for any $m\geq \varrho$.
\end{theorem}

\subsubsection{Main ideas}

In previous work, once the invariant distributions are classified,  a
``good" projection $\mathfrak{P}$ on $\mathcal{V}(\mathcal{O}^\infty)$ is constructed such that  for any $\mathfrak{p}_v$, $v\in E$ there exists a splitting:
$\mathfrak{p}_v=\mathfrak{P}(\mathfrak{p}_v)+\mathfrak{R}(\mathfrak{p}_v)$,
where $\mathfrak{P}(\mathfrak{p}_v)$ vanishes on the set of invariant distributions. Both $\mathfrak{P}(\mathfrak{p}_v)$ and the error part $\mathfrak{R}(\mathfrak{p}_v)$ have good estimates. By  ``good", we mean that one can reasonably hope to employ the KAM iterative method.

Since $\mathfrak{P}(\mathfrak{p}_v)$ is a twisted coboundary, its solution is an approximate solution. In all these cases construction of the tame projection
needs explicit description of the invariant distributions using representation theory, which heavily relies on special  features of $\GG$ (see Section \ref{sec:11}).
However,
it seems not easy to obtain similar explicit descriptions for partially hyperbolic algebraic actions other than the torus examples.

We solve the problem in a fairly general way  using the system of imprimitivity. Instead of obtaining an approximate solution by constructing a tame splitting,
we start from the reverse way: we construct an approximate solution directly from the distributional solution, then the approximate solution naturally gives a good splitting.

 Below, we continue to use the typical example previously used in Section \eqref{sec:41} to explain the main idea. We recall equation \eqref{for:291}
 \begin{align}\label{for:505}
 \varphi_{\lambda_3}\circ z_1-2\varphi_{\lambda_3}=\omega_{\lambda_3},
 \end{align}
as well as the solution
 $\varphi_{\lambda_3}$ defined in \eqref{for:286}. In particular,
$\varphi_{\lambda_3}$ satisfies the following properties:
\begin{enumerate}
  \item $\varphi_{\lambda_3}$ is smooth along $\mathfrak{g}^-_{\textbf{a}_1}\oplus\mathfrak{g}^0_{\textbf{a}_1}$, i.e., the stable and neutral directions of $z_1$;

  \smallskip
  \item $\varphi_{\lambda_3}$ probably loses differentiability along $\mathfrak{g}^+_{\textbf{a}_1}$, i.e., the unstable direction of $z_1$.
\end{enumerate}
Suppose $\mathfrak{g}^+_{\textbf{a}_1}=\{u\}$.   Set
\begin{align*}
 \Theta=\pi_{u}(f\circ \tau^{-1})\varphi_{\lambda_3}.
\end{align*}
Proposition \ref{cor:6} shows that $\Theta$ is smooth along $\mathfrak{g}^+_{\textbf{a}_1}$, while maintaining the smoothness of $\varphi_{\lambda_3}$ along $\mathfrak{g}^-_{\textbf{a}_1}\oplus\mathfrak{g}^0_{\textbf{a}_1}$. Thus $\Theta$ is in fact a $C^\infty$ approximation of $\varphi_{\lambda_3}$.
Let
\begin{align}\label{for:308}
\mathcal{R}_{\lambda_3}=\omega_{\lambda_3}-(\Theta\circ z_1-2\Theta).
\end{align}
As a result, $\Theta$ is a $C^\infty$ approximate solution of equation \ref{for:291} with the error part $\mathcal{R}_{\lambda_3}$. Comparing \eqref{for:505} and \eqref{for:308}, we see that
\begin{align*}
\mathcal{R}_{\lambda_3}=(\varphi_{\lambda_3}-\Theta)\circ z_1-2(\varphi_{\lambda_3}-\Theta).
\end{align*}
To show that the splitting in \eqref{for:308} is good, we need to show that both $\Theta$ and $\mathcal{R}_{\lambda_3}$ have good estimates.
The goodness of $\Theta$ comes from \eqref{for:70} of Proposition \ref{cor:6}. The key step is to show that $\varphi_{\lambda_3}-\Theta$ is good, which consequently leads to the goodness of  $\mathcal{R}_{\lambda_3}$.  We have
\begin{align*}
 \varphi_{\lambda_3}-\Theta&=\varphi_{\lambda_3}-\pi_{u}(f\circ \tau^{-1})\varphi_{\lambda_3}\overset{\text{(1)}}{=}\pi_{u}(f_1\circ \tau^{-1})\varphi_{\lambda_3}\\
 &\overset{\text{(2)}}{=}-\sum_{n\geq0} 2^{-(n+1)}\pi_{u}(f_1\circ \tau^{-1})\big(\omega_{\lambda_3}\circ z_1^{n}\big)\\
 &=-\sum_{n=-\infty}^\infty 2^{-(n+1)}\pi_{u}(f_1\circ \tau^{-1})\big(\omega_{\lambda_3}\circ z_1^{n}\big)\notag\\
 &+\sum_{n<0} 2^{-(n+1)}\pi_{u}(f_1\circ \tau^{-1})\big(\omega_{\lambda_3}\circ z_1^{n}\big)\\
 &\overset{\text{(3)}}{=}-\mathfrak{D}_{\lambda_3,\lambda_3,u,\tau,\ZZ}(\omega_{\lambda_3})+P_{2,\lambda_3,\tau}.
 \end{align*}
Here in $(1)$ we recall $f_1=1-f$; in $(2)$ we use \eqref{for:286} of Section \ref{sec:41}; $(3)$ we recall \eqref{for:331}, as well as the comments below \eqref{for:331}  and \eqref{for:310} of Section \ref{sec:41}.

We note that $\mathfrak{D}_{\lambda_3,\lambda_3,u,\tau,\ZZ}(\omega_{\lambda_3})$ is an invariant distribution, whose size is ``quadratically" small with respect to that of $\omega_{\lambda_3}$ (see proposition \ref{po:2} and Remark \ref{re:10}). \eqref{for:290} of Section \ref{sec:41} shows that
\begin{align*}
  \norm{P_{2,\lambda_3,\tau}}\leq C_{n,s}\tau^{-s}\norm{\omega_{\lambda_3}}_{s},\qquad \forall\,s>1.
\end{align*}
 For well chosen $\tau$ and sufficiently  large $s$, we see that $\tau^{-s}\norm{\omega_{\lambda_3}}_{s}$ is also quadratically small with respect to
$\norm{\omega_{\lambda_3}}_{C^0}$. Hence, we showed the goodness of $\varphi_{\lambda_3}-\Theta$, which consequently gives the goodness of $\mathcal{R}_{\lambda_3}$.
As a result, we have successfully shown the goodness the splitting.

\begin{remark} Comments from the above discussion.
\begin{enumerate}
\item If we let $\Theta_1=\mathfrak{s}_b\varphi_{\lambda_3}$ where $\mathfrak{s}_b\varphi_{\lambda_3}$ is the standard smoothing operator defined in \eqref{sect:6} of
Section \ref{sec:19}, then $\Theta_1$ is also a $C^\infty$ function. Moreover, $\Theta_1$ is also a good approximation of $\varphi_{\lambda_3}$. However, showing the goodness of
\begin{align*}
\mathcal{R}_{\lambda_3,1}=\omega_{\lambda_3}-(\Theta_1\circ z_1-2\Theta_1)
\end{align*}
is challenging.

\smallskip
  \item We point out that the directional smoothing operators and the description of invariant distributions play a key role in showing the splitting is good;

  \smallskip
  \item we recall that $P_{2,\lambda_3,\tau}=\mathfrak{D}_{\lambda_3,\lambda_3,u,\tau,L_{\lambda_3}^c}(\omega_{\lambda_3})$ (see \eqref{for:311} of Remark \ref{re:13}). In general cases, to estimate $ \norm{\mathfrak{D}_{\lambda_3,\lambda,u,\tau,L_{\lambda_3}^c}(\omega_{\lambda_3})}$ where $(\lambda_3, \lambda)$ is a poor pair, we use \eqref{for:156} of Proposition \ref{po:6};

  \smallskip
  \item in general cases (we assume notations in \eqref{for:532} of Remark \ref{re:13}),  $\Theta$ is naturally constructed as follows:
  \begin{align*}
   \Theta=\pi_{u_1}(f\circ \tau^{-1})\cdots \pi_{u_n}(f\circ \tau^{-1})\omega_{\lambda_3}
  \end{align*}
  where $\{u_1,\cdots, u_n\}$ is a basis set of $\mathfrak{g}^+_{\textbf{a}_1}$. Of course,  we need to arrange $u_i$ in a good way to guarantee the smoothness of $\Theta$. This is the reason we introduce notations $\Lambda^+$ and $\Lambda^-$, set orders on both of them (see \eqref{for:313} of Section \ref{sec:42}),  and then define the operators $\mathcal{E}_{\mu,\tau}$ and $\mathcal{E}_{\mu,\tau}'$ (see \eqref{for:314} of Section \ref{sec:42}).

\end{enumerate}

\end{remark}
\subsubsection{Proof of Theorem \ref{th:6}}\label{sec:57}

\emph{Step 1: Construction of $\Theta$ and $R_v$, $v\in E$. }

Let
\begin{align}\label{for:203}
 \mathfrak{p}_{v,\mu}=\mathfrak{p}_{v}|_{\mathfrak{g}_{\mu}},\quad
 v\in E\text{ and }\mu\in \Delta.
\end{align}
We denote by $z_i=z_{\textbf{a}_i}$, $i=1,\,2$,  to simplify notations.

Set
\begin{align}\label{for:204}
 \Theta_{\mu}\stackrel{\text{def}}{=} \mathcal{E}'_{\mu,\tau}\big(\mathcal{P}_{\mu}(\mathfrak{p}_{\textbf{a}_1,\mu})\big), \quad \mu\in \Delta
\end{align}
(recall \eqref{for:314} of Section \ref{sec:42} and \eqref{for:284} of Section \ref{sec:39}) and set $\Theta$ to be the  map with coordinate maps $\Theta_{\mu}$.
Let
\begin{align}\label{for:184}
 \mathcal{R}_v&\stackrel{\text{def}}{=}\mathfrak{p}_v-(\Theta\circ z_v-dz_v\Theta),\quad\forall\,v\in E\quad\text{and}\notag\\
 \mathcal{R}_{\textbf{a}_1,\mu}&\stackrel{\text{def}}{=}\mathcal{R}_{\textbf{a}_1}|_{\mathfrak{g}_{\mu}},\quad\quad \mu\in\Delta.
\end{align}
It follows from \eqref{for:184} that
\begin{align}\label{for:318}
\mathcal{R}_{\textbf{a}_1,\mu}&\overset{\text{(1)}}{=}\mathfrak{p}_{z_1,\mu}-\big(\Theta_\mu\circ z_1-(dz_1|_{\mathfrak{g}_{\mu}})\Theta_\mu\big),\qquad \mu\in\Delta.
\end{align}
Here in $(1)$ we use the fact that each $\mathfrak{g}_{\mu}$ is invariant under $dz_1$ (see \eqref{sect:10} of Section \ref{sec:19}).

\smallskip

\noindent\emph{Step 2: Estimates for $\Theta$ and $\mathcal{R}_v$, $v\in E$. }

We recall that $z_{i}=w_{\textbf{a}_i}\cdot
\textbf{a}_i$,  where $w_{\textbf{a}_i}\in V^o\cap B_{\eta}(e)$, $i=1,\,2$ (see \eqref{for:506} of Section \ref{for:504}).
\eqref{for:159} of Lemma \ref{le:5} shows that $\mathcal{P}_{\mu}(\mathfrak{p}_{z_1,\mu})$ is a distributional solution of the equation
\begin{align}\label{for:172}
\mathcal{P}_{\mu}(\mathfrak{p}_{\textbf{a}_1,\mu})\circ z_1-(dz_1|_{\mathfrak{g}_{\mu}})\mathcal{P}_{\mu}(\mathfrak{p}_{\textbf{a}_1,\mu})=\mathfrak{p}_{\textbf{a}_1,\mu},\quad \mu\in\Delta.
\end{align}
Moreover, \eqref{for:148} of Proposition \ref{le:17} shows that
\begin{align}\label{for:508}
   \big\|u_1u_2\cdots u_l\mathcal{P}_{\mu}(\mathfrak{p}_{\textbf{a}_1,\mu})\|_{-\kappa}\leq C_l \norm{\mathfrak{p}_{\textbf{a}_1,\mu}}_{l+\kappa}
  \end{align}
for any $l\geq0$, any $u_i\in (\mathfrak{g}^{\delta(\mu)}_{\textbf{a}_1}\oplus\mathfrak{g}^0_{\textbf{a}_1})\bigcap \mathcal{V}^1$, $1\leq i\leq l$.

We note that for any $\psi\in \mathfrak{g}_\mu(\mathcal{O}^\infty)$, keeping using \eqref{for:93} of Section \ref{sec:27} we see that
 \begin{align}\label{for:316}
  \langle \Theta_{\mu}, \, \psi\rangle=\langle\mathcal{P}_{\mu}(\mathfrak{p}_{\textbf{a}_1,\mu}), \, \mathcal{E}_{\mu,\tau}\psi\rangle.
 \end{align}
 We will use Proposition \ref{cor:6} to estimate $\Theta_{\mu}$. To applying Proposition \ref{cor:6},
firstly, we  define the following sets:
 \begin{enumerate}
   \item  Let $Q_i=\{\nu_i\}$, $1\leq i\leq p_{-\delta(\mu)}=p_\mu$ (recall \eqref{for:313} and \eqref{for:314} of Section \ref{sec:42}).

   \smallskip
   \item Let $Q'(\mu)$ be the subgroup with Lie algebra $\mathfrak{g}_{\textbf{a}_1}^{-\delta(\mu)}$.

   \smallskip
   \item Let
$Q_0(\mu)$ be the subgroup with Lie algebra $\mathfrak{g}_{\textbf{a}_1}^{\delta(\mu)}\oplus\mathfrak{g}_{\textbf{a}_1}^{0}$.
 \end{enumerate}
 Secondly, we need to show the following conditions are satisfied:
 \begin{enumerate}
  \item by definition, we see that $\text{Lie}(Q')$ is linearly spanned by  $\nu_i$, $1\leq i\leq p_{\mu}$;

 \smallskip
  \item it is clear that $\mathcal{V}=\mathfrak{g}_{\textbf{a}_1}^{-\delta(\mu)}\oplus\mathfrak{g}_{\textbf{a}_1}^{\delta(\mu)}\oplus\mathfrak{g}_{\textbf{a}_1}^{0}$;

  \smallskip
  \item by \eqref{for:507} of Section \ref{sec:42} we have
\begin{align*}
 [\nu_n,\nu_m] &\text{ is inside the subspace
linearly spanned by }\notag\\
&\nu_k, \,\, k>\max\{n,\,m\};
\end{align*}
 \item  estimates of derivatives of $\mathcal{P}_{\mu}(\mathfrak{p}_{\textbf{a}_1,\mu})$ along $\mathfrak{g}_{\textbf{a}_1}^{\delta(\mu)}\cup\mathfrak{g}_{\textbf{a}_1}^{0}$ are available (see \eqref{for:508}).
\end{enumerate}
 \eqref{for:316} shows that
 \begin{align*}
 \Theta_{\mu}=\mathcal{P}_{f,\tau}\big(\mathcal{P}_{\mu}(\mathfrak{p}_{\textbf{a}_1,\mu})\big)
 \end{align*}
 in the sense described  in Proposition \ref{cor:6}.

The above discussion justifies the application of Proposition \ref{cor:6} to $\Theta_{\mu}$.  As a consequence, we have
\begin{align*}
  \norm{\Theta_{\mu}}_m&\leq  C_{m}(\tau^{m+4+\kappa}\norm{\mathfrak{p}_{\textbf{a}_1,\mu}}_\kappa+\norm{\mathfrak{p}_{\textbf{a}_1,\mu}}_{m+4+\kappa})
 \end{align*}
for any $m\geq0$. As an immediate consequence we have
\begin{align*}
  \norm{\Theta}_m\leq  C_{m}(\tau^{m+4+\kappa}\norm{\mathfrak{p}}_\kappa+\norm{\mathfrak{p}}_{m+4+\kappa}),\qquad \forall\,m\geq0.
\end{align*}
It follows that: for any $v\in E$ and any $m\geq0$
\begin{align}\label{for:181}
  \norm{\mathcal{R}_v}_m&\leq \norm{\mathfrak{p}_v}_m+\norm{\Theta\circ z_v-dz_v\Theta}_m\notag\\
  &\leq  C_{m}(\tau^{m+4+\kappa}\norm{\mathfrak{p}}_\kappa+\norm{\mathfrak{p}}_{m+4+\kappa}).
\end{align}
By Sobolev embedding inequality \eqref{for:179} of Section \ref{sec:19} we have
\begin{align}
 \max\{\norm{&\Theta}_{C^m},\,\norm{\mathcal{R}_v}_{C^m}\}\leq C_{m}\max\{\norm{\Theta}_{m+\beta},\,\norm{\mathcal{R}_v}_{m+\beta}\}\notag\\
 &\leq C_{m}(\tau^{m+4+\kappa+\beta}\norm{\mathfrak{p}}_\kappa+\norm{\mathfrak{p}}_{m+4+\kappa+\beta})\notag\\
  &\leq C_{m}(\tau^{m+4+\kappa+\beta}\norm{\mathfrak{p}}_{C^\kappa}+\norm{\mathfrak{p}}_{C^{m+4+\kappa+\beta}})\label{for:374}\\
 &\overset{\text{(1)}}{\leq} C_{m}(\tau^{m+\varrho}\norm{\mathfrak{p}}_{C^\varrho}+\norm{\mathfrak{p}}_{C^{m+\varrho}})\label{for:373}.
\end{align}
for any $v\in E$ and any $m\geq0$. Here in $(1)$ we use the definition of $\varrho$ in \eqref{for:171}.

Thus we get \eqref{for:206}.

\smallskip
\noindent\emph{Step 3:  In this part we show that:}
\begin{align}\label{for:178}
 \langle \mathcal{R}_{\textbf{a}_1,\mu}, \, \psi\rangle=\big\langle (I-\mathcal{E}_{\mu,\tau}')\mathcal{P}_{\mu}(\mathfrak{p}_{\textbf{a}_1,\mu}), \, \psi\circ z_1^{-1}-(dz_1|_{\mathfrak{g}_{\mu}})^*\psi\big\rangle
\end{align}
 for any $\mu\in\Delta$ and  any $\psi\in \mathfrak{g}_\mu(\mathcal{O})$. Here $(dz_1|_{\mathfrak{g}_{\mu}})^*$ is the transpose of $dz_1|_{\mathfrak{g}_{\mu}}$.

\smallskip
\emph{Note}. \emph{Steps} $3-7$ aim to obtain a sharper estimate of  $\norm{\mathcal{R}_{z_1,\mu}}$ than the one in \eqref{for:181}. This is the reason we rewrite $\mathcal{R}_{z_1,\mu}$
in this step.

\smallskip

It follows from \eqref{for:318} that
\begin{align*}
 \mathcal{R}_{\textbf{a}_1,\mu}&=\mathfrak{p}_{\textbf{a}_1,\mu}-\big(\Theta_\mu\circ z_1-(dz_1|_{\mathfrak{g}_{\mu}})\Theta_\mu\big)\notag\\
 &\overset{\text{(1)}}{=}\big(\mathcal{P}_{\mu}(\mathfrak{p}_{\textbf{a}_1,\mu})-\Theta_\mu\big)\circ z_1-(dz_1|_{\mathfrak{g}_{\mu}})
 \big(\mathcal{P}_{\mu}(\mathfrak{p}_{\textbf{a}_1,\mu})-\Omega_\mu\big)\notag\\
 &\overset{\text{(2)}}{=}\big((I-\mathcal{E}'_{\mu,\tau})\mathcal{P}_{\mu}(\mathfrak{p}_{\textbf{a}_1,\mu})\big)\circ z_1-(dz_1|_{\mathfrak{g}_{\mu}})
 \big((I-\mathcal{E}'_{\mu,\tau})\mathcal{P}_{\mu}(\mathfrak{p}_{\textbf{a}_1,\mu})\big)
\end{align*}
for any $\mu\in\Delta$.

Here in $(1)$ we use \eqref{for:172}; in $(2)$ we recall \eqref{for:204}.

Then for any $\psi\in \mathfrak{g}_\mu(\mathcal{O})$ we have
\begin{align*}
 \langle \mathcal{R}_{z_1,\mu}, \, \psi\rangle&=\big\langle\big((I-\mathcal{E}'_{\mu,\tau})\mathcal{P}_{\mu}(\mathfrak{p}_{\textbf{a}_1,\mu})\big)\circ z_1, \, \psi\big\rangle\\
 &-\big\langle(dz_1|_{\mathfrak{g}_{\mu}})
 \big((I-\mathcal{E}'_{\mu,\tau})\mathcal{P}_{\mu}(\mathfrak{p}_{\textbf{a}_1,\mu})\big), \, \psi\big\rangle\\
 &=\big\langle(I-\mathcal{E}'_{\mu,\tau})\mathcal{P}_{\mu}(\mathfrak{p}_{\textbf{a}_1,\mu}), \, \psi\circ z_1^{-1}\big\rangle\\
 &-\big\langle(I-\mathcal{E}'_{\mu,\tau})\mathcal{P}_{\mu}(\mathfrak{p}_{\textbf{a}_1,\mu}), \, (dz_1|_{\mathfrak{g}_{\mu}})^*\psi\big\rangle\\
 &=\big\langle (I-\mathcal{E}_{\mu,\tau}')\mathcal{P}_{\mu}(\mathfrak{p}_{\textbf{a}_1,\mu}), \, \psi\circ z_1^{-1}-(dz_1|_{\mathfrak{g}_{\mu}})^*\psi\big\rangle.
\end{align*}
We thus finish the proof.

\smallskip
\noindent\emph{Step 4:  In this part we show that:}
\begin{align}\label{for:205}
   I-\mathcal{E}_{\mu,\tau}=\sum_{k=1}^{p_\mu}\pi_{\nu_k}(f_1\circ \tau^{-1})\mathcal{E}_{\mu,\tau,k+1},\qquad \mu\in\Delta.
  \end{align}
 From \eqref{for:178}, we see that to estimate $\mathcal{R}_{\textbf{a}_1,\mu}$, we should estimate $(I-\mathcal{E}'_{\mu,\tau})\mathcal{P}_{\mu}(\mathfrak{p}_{\textbf{a}_1,\mu})$ first.
We also note that
\begin{align}\label{for:319}
 \big\langle (I-\mathcal{E}'_{\mu,\tau})&\mathcal{P}_{\mu}(\mathfrak{p}_{\textbf{a}_1,\mu}), \, \psi\big\rangle=\big\langle \mathcal{P}_{\mu}(\mathfrak{p}_{\textbf{a}_1,\mu}), \, \psi\big\rangle-\big\langle \mathcal{E}'_{\mu,\tau}\mathcal{P}_{\mu}(\mathfrak{p}_{\textbf{a}_1,\mu}), \, \psi\big\rangle\notag\\
 &\overset{\text{(1)}}{=}\big\langle \mathcal{P}_{\mu}(\mathfrak{p}_{\textbf{a}_1,\mu}), \, \psi\big\rangle-\big\langle \mathcal{P}_{\mu}(\mathfrak{p}_{\textbf{a}_1,\mu}), \, \mathcal{E}_{\mu,\tau}\psi\big\rangle\notag\\
 &=\big\langle \mathcal{P}_{\mu}(\mathfrak{p}_{\textbf{a}_1,\mu}), \, (I-\mathcal{E}_{\mu,\tau})\psi\big\rangle.
\end{align}
Here in $(1)$ we keep using \eqref{for:93} of Section \ref{sec:27}.

Then it is natural for us to rewrite $I-\mathcal{E}_{\mu,\tau}$ in this step.

\smallskip
Recalling the notation in \eqref{for:314} of Section \ref{sec:42} we see that
  \begin{align}
   I-\mathcal{E}_{\mu,\tau}&=I-\mathcal{E}_{\mu,\tau,1}\overset{\text{(1)}}{=}\sum_{k=1}^{p_\mu}(\mathcal{E}_{\mu,\tau,k+1}-\mathcal{E}_{\mu,\tau,k})\notag\\
   &\overset{\text{(2)}}{=}\sum_{k=1}^{p_\mu}\big(\pi_{\nu_k}(1\circ \tau^{-1})\mathcal{E}_{\mu,\tau,k+1}-\pi_{\nu_k}(f\circ \tau^{-1})\mathcal{E}_{\mu,\tau,k+1}\big)\notag\\
   &=\sum_{k=1}^{p_\mu}\pi_{\nu_k}(f_1\circ \tau^{-1})\mathcal{E}_{\mu,\tau,k+1}.
  \end{align}
  Here in $(1)$ we set $\mathcal{E}_{\mu,\tau,p_\mu+1}=I$; in $(2)$ we note that $\pi_{\nu_{k}}(1\circ \tau^{-1})=I$ for each $k$.

\smallskip
\noindent\emph{Step 5:  In this part we show that:}
\begin{align}\label{for:180}
 &\Big\|\big\langle (I-\mathcal{E}'_{\mu,\tau})\mathcal{P}_{\mu}(\mathfrak{p}_{\textbf{a}_1,\mu}), \, \psi\big\rangle\Big\|\notag\\
 &\leq C(\norm{\log[z,z]}^{c_0}\norm{\mathfrak{p}}_{\kappa+1}+\norm{\mathbb{D}}_\kappa
 +C_m\tau^{-m}\norm{\mathfrak{p}}_m)\|\psi\|_\kappa
\end{align}
for any $\mu\in\Delta$, any $m\geq\kappa$ and any $\psi\in \mathfrak{g}_\mu(\mathcal{O}^\kappa)$.

\medskip
For any $\psi\in \mathfrak{g}_\mu(\mathcal{O}^\kappa)$, we have
\begin{align*}
 &\big\langle (I-\mathcal{E}'_{\mu,\tau})\mathcal{P}_{\mu}(\mathfrak{p}_{\textbf{a}_1,\mu}), \, \psi\big\rangle\notag\\
 &\overset{\text{(1)}}{=}\big\langle \mathcal{P}_{\mu}(\mathfrak{p}_{\textbf{a}_1,\mu}), \, (I-\mathcal{E}_{\mu,\tau})\psi\big\rangle\notag\\
 &\overset{\text{(2)}}{=}\sum_{k=1}^{p_\mu}\big\langle
 \mathcal{P}_{\mu}(\mathfrak{p}_{\textbf{a}_1,\mu}), \, \pi_{\nu_k}(f_1\circ \tau^{-1})\mathcal{E}_{\mu,\tau,k+1}\psi\big\rangle\notag\\
 &\overset{\text{(3)}}{=}\sum_{k=1}^{p_\mu}\big\langle
 \pi_{\nu_k}(f_1\circ \tau^{-1})\mathcal{P}_{\mu}(\mathfrak{p}_{\textbf{a}_1,\mu}), \, \mathcal{E}_{\mu,\tau,k+1}\psi\big\rangle\notag\\
 &\overset{\text{(4)}}{=}\sum_{k=1}^{p_\mu}\big\langle
 \mathfrak{D}_{\mu,\lambda_{i_k},\nu_k,\tau,L_\mu}(\mathfrak{p}_{\textbf{a}_1,\mu}), \, \mathcal{E}_{\mu,\tau,k+1}\psi\big\rangle\notag\\
 &\overset{\text{(5)}}{=}P_1-P_2.
\end{align*}
Here in $(1)$ we use \eqref{for:319}; in $(2)$ we use \eqref{for:205}; in $(3)$ we use \eqref{for:93} of Section \ref{sec:27};
in $(4)$ we recall notations in Sections \eqref{sec:39} and \eqref{sec:42}. We note that $(\mu,\lambda_{i_k})$ is a poor pair for each $k$ and
note that $\nu_k=u_{(\lambda_{i_k},j_k)}$; in $(5)$ we set
\begin{align*}
 P_1&=\sum_{k=1}^{p_\mu}\big\langle
 \mathfrak{D}_{\mu,\lambda_{i_k},\nu_k,\tau,\ZZ}(\mathfrak{p}_{\textbf{a}_1,\mu}), \, \mathcal{E}_{\mu,\tau,k+1}\psi\big\rangle,\qquad\text{and}\\
 P_2&=\sum_{k=1}^{p_\mu}\big\langle
 \mathfrak{D}_{\mu,\lambda_{i_k},\nu_k,\tau,L_\mu^c}(\mathfrak{p}_{\textbf{a}_1,\mu}), \, \mathcal{E}_{\mu,\tau,k+1}\psi\big\rangle.
\end{align*}
It follows that
\begin{align}\label{for:175}
 \Big\|\big\langle (I-\mathcal{E}'_{\mu,\tau})\mathcal{P}_{\mu}(\mathfrak{p}_{\textbf{a}_1,\mu}), \, \psi\big\rangle\Big\|\leq \norm{P_1}+\norm{P_2}.
\end{align}
Next, we estimate $P_1$ and $P_2$, respectively. For $P_1$, we recall that for each $k$, $\mathfrak{D}_{\mu,\lambda_{i_k},\nu_k,\tau,\ZZ}(\mathfrak{p}_{\textbf{a}_1,\mu})$ is an invariant distribution . Then it is natural for us to use Proposition \ref{po:2}. Thus we turn to the almost twisted cocycle equation equation \eqref{for:173} for $z_1$ and $z_2$:
\begin{align}\label{for:174}
 \mathfrak{p}_{\textbf{a}_1,\mu}\circ z_2-dz_2\mathfrak{p}_{\textbf{a}_1,\mu}&=\mathfrak{p}_{\textbf{a}_2,\mu}\circ z_1-dz_1\mathfrak{p}_{\textbf{a}_2,\mu}+\mathbb{D}_{\textbf{a}_2,\textbf{a}_1,\mu},
\end{align}
where $\mathbb{D}_{\textbf{a}_2,\textbf{a}_1,\mu}=\mathbb{D}_{\textbf{a}_2,\textbf{a}_1}|_{\mathfrak{g}_{\mu}}$.

Applying Proposition \ref{po:2} to \eqref{for:174}, we have
\begin{align}\label{for:372}
 \norm{\mathfrak{D}_{\mu,\tau,\ZZ}(\mathfrak{p}_{\textbf{a}_1,\mu})}&\leq C\norm{\log[z_1,z_2]}^{c_0}\max\{\norm{\mathfrak{p}_{\textbf{a}_1,\mu}}_{\kappa+1},\,\norm{\mathfrak{p}_{\textbf{a}_2,\mu}}_{\kappa}\}\notag\\
 &+C\norm{\mathbb{D}_{\textbf{a}_2,\textbf{a}_1,\mu}}_\kappa\notag\\
 &\leq C\norm{\log[z,z]}^{c_0}\norm{\mathfrak{p}}_{\kappa+1}+C\norm{\mathbb{D}}_\kappa.
\end{align}
It follows that
\begin{align}\label{for:176}
 \norm{P_1}&\overset{\text{(1)}}{\leq} \sum_{k=1}^{p_\mu}\norm{\mathfrak{D}_{\mu,\tau,\ZZ}(\mathfrak{p}_{\textbf{a}_1,\mu})}
 \|\mathcal{E}_{\mu,\tau,k+1}\psi\|_\kappa\notag\\
 &\overset{\text{(2)}}{\leq}\sum_{k=1}^{p_\mu}\big(C\norm{\log[z,z]}^{c_0}\norm{\mathfrak{p}}_{\kappa+1}+C\norm{\mathbb{D}}_\kappa\big)
 \|\mathcal{E}_{\mu,\tau,k+1}\psi\|_\kappa\notag\\
 &\overset{\text{(3)}}{\leq} C_1\big(\norm{\log[z,z]}^{c_0}\norm{\mathfrak{p}}_{\kappa+1}+\norm{\mathbb{D}}_\kappa\big)
 \|\psi\|_\kappa.
\end{align}
Here in $(1)$ we recall the definition of $\norm{\mathfrak{D}_{\mu,\tau,\ZZ}(\mathfrak{p}_{z_1,\mu})}$ in Section \eqref{sec:39}; in $(2)$ we use \eqref{for:372}; in $(3)$ we use \eqref{for:7} of Lemma \ref{cor:1}.

Next, we estimate $P_2$.
\begin{align}\label{for:177}
 \norm{P_2}&\overset{\text{(1)}}{\leq} \sum_{k=1}^{p_\mu}\big\| \mathfrak{D}_{\mu,\lambda_{i_k},\nu_k,\tau,L_\mu^c}(\mathfrak{p}_{\textbf{a}_1,\mu})\big\|
 \|\mathcal{E}_{\mu,\tau,k+1}\psi\|\notag\\
 &\overset{\text{(2)}}{\leq} \sum_{k=1}^{p_\mu}C_m\tau^{-m}\norm{\mathfrak{p}_{\textbf{a}_1,\mu}}_m
 \|\mathcal{E}_{\mu,\tau,k+1}\psi\|\notag\\
 &\overset{\text{(3)}}{\leq} \sum_{k=1}^{p_\mu}C_{m,1}\tau^{-m}\norm{\mathfrak{p}_{\textbf{a}_1,\mu}}_m
 \|\psi\|
\end{align}
for any $m\geq\kappa$.

Here in $(1)$ we use cauchy schwarz inequality; in $(2)$ we use \eqref{for:156} of Proposition \ref{po:6}; in $(3)$ we use  \eqref{for:7} of Lemma \ref{cor:1}.

Hence, \eqref{for:180} follows from \eqref{for:175}, \eqref{for:176} and \eqref{for:177}.

\smallskip
\noindent\emph{Step 6:  In this part, we prove that}
\begin{align}\label{for:182}
 \norm{\mathcal{R}_{\textbf{a}_1}}\leq C\tau^{\frac{\varrho}{2}}\big(\norm{\log[z,z]}^{c_0}\norm{\mathfrak{p}}_{\varrho}+\norm{\mathbb{D}}_\varrho
 +C_m\tau^{-m}\norm{\mathfrak{p}}_m\big)^{\frac{1}{2}}(\norm{\mathfrak{p}}_{\varrho})^{\frac{1}{2}}
\end{align}
for any $m\geq\varrho$.

\medskip
For any $\mu\in\Delta$,  by using \eqref{for:178} we have
\begin{align*}
 \|\langle \mathcal{R}_{\textbf{a}_1,\mu},& \, \mathcal{R}_{\textbf{a}_1,\mu}\rangle\|=\Big\|\big\langle (I-\mathcal{E}'_{\mu,\tau})\mathcal{P}_{\mu}(\mathfrak{p}_{\textbf{a}_1,\mu}), \, \mathcal{R}_{\textbf{a}_1,\mu}\circ z_1^{-1}-(dz_1|_{\mathfrak{g}_{\mu}})^*\mathcal{R}_{\textbf{a}_1,\mu}\big\rangle\Big\|\\
 &\overset{\text{(1)}}{\leq}C(\norm{\log[z,z]}^{c_0}\norm{\mathfrak{p}}_{\kappa+1}+\norm{\mathbb{D}}_\kappa+C_m\tau^{-m}\norm{\mathfrak{p}}_m)\\
 &\quad\cdot\big\|\mathcal{R}_{\textbf{a}_1,\mu}\circ z_1^{-1}-(dz_1|_{\mathfrak{g}_{\mu}})^*\mathcal{R}_{\textbf{a}_1,\mu}\big\|_\kappa\\
 &\leq C_1(\norm{\log[z,z]}^{c_0}\norm{\mathfrak{p}}_{\kappa+1}+\norm{\mathbb{D}}_\kappa
 +C_m\tau^{-m}\norm{\mathfrak{p}}_m)\\
 &\quad\cdot\|\mathcal{R}_{\textbf{a}_1,\mu}\|_\kappa\\
 &\overset{\text{(2)}}{\leq}C_2(\norm{\log[z,z]}^{c_0}\norm{\mathfrak{p}}_{\kappa+1}+\norm{\mathbb{D}}_\kappa
 +C_m\tau^{-m}\norm{\mathfrak{p}}_m)\\
 &\quad\cdot(\tau^{2\kappa+4}\norm{\mathfrak{p}}_\kappa+\norm{\mathfrak{p}}_{2\kappa+4})\\
 &\overset{\text{(3)}}{\leq} 2C_2\tau^{\varrho}(\norm{\log[z,z]}^{c_0}\norm{\mathfrak{p}}_{\varrho}+\norm{\mathbb{D}}_\varrho
 +C_m\tau^{-m}\norm{\mathfrak{p}}_m)\norm{\mathfrak{p}}_{\varrho}
\end{align*}
for any $m\geq \varrho$.

Here in $(1)$ we use \eqref{for:180};  in $(2)$ we use \eqref{for:181}; in $(3)$ we use definition of $\varrho$ in  \eqref{for:171}.

The above estimate implies \eqref{for:182}.

\smallskip
\noindent\emph{Step 7:  Estimates for $\mathcal{R}_{v}$, $v\in E$}.

In this step it is standard to use the ``higher rank trick" (see for example \cite{Damjanovic4} and \cite{W5}). We use the almost twisted cocycle equation over $v$ and $z_1$ from \eqref{for:173}:
\begin{align*}
 \mathfrak{p}_v\circ z_1-dz_1\mathfrak{p}_v=\mathfrak{p}_{\textbf{a}_1}\circ z_v-dz_v\mathfrak{p}_{\textbf{a}_1}+\mathbb{D}_{\textbf{a}_1,v},\quad \forall \,v\in E.
\end{align*}
We substitute  the expressions for $\mathfrak{p}_{\textbf{a}_1}$ and $\mathfrak{p}_{v}$
from \eqref{for:184} respectively into the above equation. Then we get
\begin{align*}
&\big(\mathcal{R}_v+(\Theta\circ z_v-dz_v\Theta)\big)\circ z_1-dz_1\big(\mathcal{R}_v+(\Theta\circ z_v-dz_v\Theta)\big)\\
&=\big(\mathcal{R}_{\textbf{a}_1}+(\Theta\circ z_1-dz_1\Theta)\big)\circ z_v-dz_v\big(\mathcal{R}_{\textbf{a}_1}+(\Theta\circ z_1-dz_1\Theta)\big)+\mathbb{D}_{\textbf{a}_1,v}
\end{align*}
for each $v\in E$.

The above equation can be simplified to
\begin{align}\label{for:376}
\mathcal{R}_v\circ z_1-dz_1\mathcal{R}_v=\mathcal{R}_{\textbf{a}_1}\circ z_v-dz_v\mathcal{R}_{\textbf{a}_1}+\mathbb{D}_{z_1,v}+\Omega_{\textbf{a}_1,v},\quad \forall\,v\in E.
\end{align}
where
\begin{align*}
 \Omega_{\textbf{a}_1,v}=\Theta\circ (z_{1}z_v)-\Theta\circ (z_vz_{1})+\big(dz_vdz_{1}-dz_{1}dz_{v}\big)\Theta.
\end{align*}
For any $m\geq0$ we have
\begin{align}\label{for:375}
 \norm{\Omega_{\textbf{a}_1,v}}_{C^m}&\overset{\text{\tiny$(1)$}}{\leq} C_m\norm{\Theta}_{C^{m+1}}\textbf{\emph{d}}([z_v,z_1],e)+C_m\norm{d[z_v,z_1]-I}\norm{\Theta}_{C^{m}}\notag\\
 &\overset{\text{\tiny$(2)$}}{\leq} C_{m,1}\norm{\log[z_v,z_1]}\norm{\Theta}_{C^{m+1}}\notag\\
 &\overset{\text{\tiny$(3)$}}{\leq} C_{m,2}\norm{\log[z_v,z_1]}(\tau^{m+5+\kappa+\beta}\norm{\mathfrak{p}}_{C^\kappa}+\norm{\mathfrak{p}}_{C^{m+5+\kappa+\beta}}).
\end{align}
Here in $(1)$ we recall notations in \eqref{sect:9} of Section \ref{sec:19}; in $(2)$ we note that
\begin{align*}
 \textbf{\emph{d}}([z_v,z_1],e)\overset{\text{\tiny$(a)$}}{\leq} C\norm{\log[z_v,z_1]}\quad\text{and}\quad \norm{d[z_v,z_1]-I}\overset{\text{\tiny$(b)$}}{\leq} C\norm{\log[z_v,z_1]}.
\end{align*}
Here in $(a)$ we use \eqref{for:365} of \eqref{sect:9} of Section \ref{sec:19}; in $(b)$ we note that $d[z_v,z_1]=\text{Ad}_{[z_v,z_1]}$ and then we use  \eqref{for:364} of \eqref{sect:9} of Section \ref{sec:19}; in $(3)$ we use \eqref{for:374}.

Then we have: for any $m\geq \varrho$
\begin{align}
  \norm{\mathcal{R}_v}_{C^0}&\overset{\text{\tiny$(1)$}}{\leq} C\norm{\mathcal{R}_{\textbf{a}_1}\circ v-dv\mathcal{R}_{\textbf{a}_1}+\mathbb{D}_{z_1,v}+\Omega_{\textbf{a}_1,v}}_{C^{\beta+\kappa+4}}\notag\\
  &\overset{\text{\tiny$(2,3)$}}{\leq} C\norm{\mathcal{R}_{\textbf{a}_1}}_{2\beta+\kappa+4}+C\norm{\mathbb{D}}_{C^\varrho}+\norm{\Omega_{\textbf{a}_1,v}}_{C^{\beta+\kappa+4}}\notag\\
  &\overset{\text{\tiny$(4)$}}{\leq} C_1
\norm{\mathcal{R}_{\textbf{a}_1}}^{\frac{1}{4}}
(\norm{\mathcal{R}_{\textbf{a}_1}}_{\frac{4}{3}(2\beta+\kappa+4)})^{\frac{3}{4}}+C\norm{\mathbb{D}}_{C^\varrho}\notag\\
&+C\norm{\log[z_v,z_1]}(\tau^{9+2\kappa+2\beta}\norm{\mathfrak{p}}_{C^\kappa}+\norm{\mathfrak{p}}_{C^{9+2\kappa+2\beta}})\label{for:189}\\
 &\overset{\text{\tiny$(5)$}}{\leq} C_2\tau^{\frac{3\varrho}{4}}
\norm{\mathcal{R}_{\textbf{a}_1}}^{\frac{1}{4}}
(\norm{\mathfrak{p}}_\varrho)^{\frac{3}{4}}+C\norm{\mathbb{D}}_{C^\varrho}\notag\\
&+C_1\tau^{\varrho}\norm{\log[z_v,z_1]}\norm{\mathfrak{p}}_{C^\varrho}\notag\\
&\overset{\text{\tiny$(6)$}}{\leq} C_2\tau^{\frac{7\varrho}{8}}
\big(\norm{\log[z,z]}^{c_0}\norm{\mathfrak{p}}_{\varrho}+\norm{\mathbb{D}}_\varrho
 +C_m\tau^{-m}\norm{\mathfrak{p}}_m\big)^{\frac{1}{8}}(\norm{\mathfrak{p}}_{\varrho})^{\frac{7}{8}}\notag\\
&+C\norm{\mathbb{D}}_{\varrho}+C_1\tau^{\varrho}\norm{\log[z_v,z_1]}\norm{\mathfrak{p}}_{C^\varrho}\notag\\
&\overset{\text{\tiny$(7)$}}{\leq}C_2\tau^{\varrho}
\norm{\log[z,z]}^{\frac{c_0}{8}}\norm{\mathfrak{p}}_{\varrho}+ C_2\tau^{\varrho}
(\norm{\mathbb{D}}_\varrho)^{\frac{1}{8}}(\norm{\mathfrak{p}}_{\varrho})^{\frac{7}{8}}\notag\\
 &+C_2\tau^{\varrho}
\big(C_m\tau^{-m}\norm{\mathfrak{p}}_m\big)^{\frac{1}{8}}(\norm{\mathfrak{p}}_{\varrho})^{\frac{7}{8}}
+C\norm{\mathbb{D}}_{\varrho}\notag\\
&+C_1\tau^{\varrho}\norm{\log[z_v,z_1]}\norm{\mathfrak{p}}_{C^\varrho}\notag\\
&\overset{\text{\tiny$(3,8)$}}{\leq} C_2\tau^{\varrho}
\norm{\log[z,z]}^{\frac{c_0}{8}}\norm{\mathfrak{p}}_{\varrho}+C_2\tau^{\varrho}
(\norm{\mathbb{D}}_{C^\varrho})^{\frac{1}{8}}(\norm{\mathfrak{p}}_{{C^\varrho}})^{\frac{7}{8}}\notag\\
 &+C_2\tau^{\varrho}
\big(C_m\tau^{-m}\norm{\mathfrak{p}}_{C^m}\big)^{\frac{1}{8}}(\norm{\mathfrak{p}}_{{C^\varrho}})^{\frac{7}{8}}
+C\norm{\mathbb{D}}_{{C^\varrho}}.\notag
\end{align}
Here in $(1)$ we apply Corollary \ref{cor:3} to \eqref{for:376};
in $(2)$ we use definition of $\varrho$ in  \eqref{for:171}; in $(3)$ we use Sobolev embedding inequality \eqref{for:179} of Section \ref{sec:19}; in $(4)$ we use interpolation inequalities and \eqref{for:375}; in $(5)$ we use \eqref{for:181}:
\begin{align*}
\norm{\mathcal{R}_{z_1}}_{\frac{4}{3}(2\beta+\kappa+4)}&\leq C( \tau^{\frac{4}{3}(2\beta+\kappa+4)+4+\kappa}\norm{\mathfrak{p}}_\kappa+\norm{\mathfrak{p}}_{\frac{4}{3}(2\beta+\kappa+4)+4+\kappa})\notag\\
&\overset{\text{\tiny$(1)$}}{\leq} 2C\tau^{\varrho}\norm{\mathfrak{p}}_\varrho.
\end{align*}
Here in $(1)$ we use the definition of $\varrho$ in  \eqref{for:171}. Similarly, we have
\begin{align*}
 &\norm{\log[z_v,z_1]}(\tau^{9+2\kappa+2\beta}\norm{\mathfrak{p}}_{C^\kappa}+\norm{\mathfrak{p}}_{C^{9+2\kappa+2\beta}})\\
 &\leq 2C\tau^{\varrho}\norm{\log[z_v,z_1]}\norm{\mathfrak{p}}_\varrho;
\end{align*}
in $(6)$ we use \eqref{for:182}; in $(7)$ we use $\tau>1$ and the inequality:
\begin{align}\label{for:217}
 (x+y+a)^c\leq x^c+y^c+a^c,\qquad \forall\, x,\,y,\,a>0,\,\,\,0<c<1;
\end{align}
in $(8)$ we note that $c_0<1$ (see \eqref{for:304} of Section \ref{sec:40}).

Thus we get \eqref{for:188}. Hence, we finish the proof.

\begin{remark} In \eqref{for:189} we can choose any $r>2\beta+\kappa+4$ to estimate $\norm{\mathcal{R}_{z_1}}_{2\beta+\kappa+4}$ using interpolation inequalities:
 \begin{align*}
  \norm{\mathcal{R}_{\textbf{a}_1}}_{2\beta+\kappa+4}\leq C_r(\norm{\mathcal{R}_{\textbf{a}_1}})^{\frac{r-(2\beta+\kappa+4)}{r}}
(\norm{\mathcal{R}_{z_1}}_r)^{\frac{2\beta+\kappa+4}{r}}.
 \end{align*}
The reason we choose $r=\text{\small$\frac{4}{3}$}(2\beta+\kappa+4)$ is to obtain uniform estimates for both the nilmanifolds and (twisted) symmetric spaces, see Theorem \ref{th:9}.
\end{remark}
\subsubsection{Revised version of Theorem \ref{th:6}}
From Theorem \ref{th:6}, we see that there is a fixed loss of regularity in solving the twisted almost cocycle equations.
To overcome the fixed loss of
regularity at each step of the iteration, it is standard (see  \cite{Zehnder})
to use the smoothing operators $\mathfrak{s}_{b}$ (see \eqref{sect:6} of Section \ref{sec:19}). As a result, we introduce two more parameter $b$ and $\ell$.
\begin{corollary}\label{cor:8} Suppose $\mathfrak{p}_v\in \mathcal{V}(\mathcal{O}^\infty)$, $v\in E$. Also suppose $E'=\{z_v:\,v\in E\}$ is an $\eta$-standard   perturbation of $E$. Set
\begin{align*}
 \mathfrak{p}_u\circ z_v-dz_v\mathfrak{p}_u=\mathfrak{p}_v\circ z_u-dz_u\mathfrak{p}_v+\mathbb{D}_{u,v}, \quad \forall\,u,v\in E.
\end{align*}
Let  $\norm{\mathbb{D}}_{{C^m}}$, $\norm{\mathfrak{p}}_{{C^m}}$ and $\norm{\log[z,z]}$ be as defined in Theorem \ref{th:6} (see Section \ref{for:504}). Then for any $b,\,\tau>1$ there exist $\Theta\in \mathcal{V}(\mathcal{O}^\infty)$ and $\mathcal{R}_v\in \mathcal{V}(\mathcal{O}^\infty)$, $v\in E$ such that
\begin{align*}
 \mathfrak{p}_v=\Theta\circ z_v-dz_v\Theta+\mathcal{R}_v
\end{align*}
with estimates: for any $v\in E$
\begin{align}\label{for:250}
 \max\{\norm{\Theta}_{C^m},\,\norm{\mathcal{R}_v}_{C^m}\}
 \leq C_{m}b^\varrho(\tau^{m}\norm{\mathfrak{p}}_{C^\varrho}+\norm{\mathfrak{p}}_{C^{m}})
\end{align}
for any $m\geq\varrho$;  and
\begin{align}\label{for:252}
 \norm{\mathcal{R}_v}&_{C^0}\leq C\tau^{\varrho}\norm{\log[z,z]}^{\frac{c_0}{8}}\norm{\mathfrak{p}}_{C^\varrho}+C\tau^{\varrho}
(\norm{\mathbb{D}}_{C^\varrho})^{\frac{1}{8}}(\norm{\mathfrak{p}}_{{C^\varrho}})^{\frac{7}{8}}\notag\\
& +C\tau^{\varrho}
\big(C_m\tau^{-m}\norm{\mathfrak{p}}_{C^m}\big)^{\frac{1}{8}}(\norm{\mathfrak{p}}_{{C^\varrho}})^{\frac{7}{8}}\notag\\
&+C\norm{\mathbb{D}}_{{C^\varrho}}+C_{\ell} b^{-\ell+\varrho}(\tau^{\ell}\norm{\mathfrak{p}}_{C^\varrho}+\norm{\mathfrak{p}}_{C^{\ell}})
\end{align}
for any $m\geq \varrho$ and $\ell\geq\varrho$.
\end{corollary}
\begin{proof}
By Theorem \ref{th:6}, there exist $\Theta'\in \mathcal{V}(\mathcal{O}^\infty)$ and $\mathcal{R}'_v\in \mathcal{V}(\mathcal{O}^\infty)$, $v\in E$ such that
\begin{align}\label{for:249}
 \mathfrak{p}_v=\Theta'\circ z_v-dz_v\Theta+\mathcal{R}'_v, \quad\forall\,v\in E.
\end{align}
and the estimates \eqref{for:206} and \eqref{for:188} are satisfied for $\Theta'$ and $\mathcal{R}'_v$, $v\in E$.  Let
\begin{align}\label{for:320}
  \Theta=\mathfrak{s}_{b}\Theta'\quad \text{and}\quad
  \mathcal{R}_v=\mathfrak{p}_v-(\Theta\circ z_v-dz_v\Theta), \quad v\in E
\end{align}
(see \eqref{sect:6} of Section \ref{sec:19}). Then we have
\begin{align*}
 \norm{\Theta}_{C^m}&\overset{\text{\tiny$(1)$}}{\leq} C_{m}b^\varrho\norm{\Theta'}_{C^{m-\varrho}}\notag\\
 &\overset{\text{\tiny$(2)$}}{\leq} C_{m,1}b^\varrho(\tau^{m}\norm{\mathfrak{p}}_{C^\varrho}+\norm{\mathfrak{p}}_{C^{m}})
 \end{align*}
for any $m\geq\varrho$.  Here in $(1)$ we use \eqref{for:45} of Section \ref{sec:19}; in $(2)$ we use \eqref{for:206}.

It follows that: for any $v\in E$ and any $m\geq\varrho$
\begin{align*}
  \norm{\mathcal{R}_v}_{C^m}&\leq \norm{\mathfrak{p}_v}_{C^m}+\norm{\Theta\circ z_v-dz_v\Theta}_m\notag\leq
  \norm{\mathfrak{p}_v}_{C^m}+C\norm{\Theta}_{C^m}\\
  &\leq  C_{m}b^\varrho(\tau^{m}\norm{\mathfrak{p}}_{C^\varrho}+\norm{\mathfrak{p}}_{C^{m}}).
\end{align*}
Hence we get \eqref{for:250}.

Next, we estimate $\norm{\mathcal{R}_v}_{C^0}$, $v\in E$. From \eqref{for:249}, we have
\begin{align}\label{for:321}
 \mathfrak{p}_v=(\mathfrak{s}_{b}\Theta')\circ z_v-dz_v\mathfrak{s}_{b}\Theta'+T_v, \quad\forall\,v\in E
\end{align}
where
\begin{align*}
 T_v=\mathcal{R}'_v+\big((I-\mathfrak{s}_{b})\Theta'\big)\circ z_v-dz_v\big((I-\mathfrak{s}_{b})\Theta'\big),\quad \forall\,v\in E.
\end{align*}
Comparing \eqref{for:320} and \eqref{for:321} we see that
\begin{align*}
 \mathcal{R}_v=T_v,
\end{align*}
which gives
\begin{align}\label{for:251}
 \norm{\mathcal{R}_v}_{C^0}&=\norm{T_v}_{C^0}\leq C\norm{\mathcal{R}'_v}_{C^0}+C\norm{(I-\mathfrak{s}_{b})\Theta'}_{C^0}\notag\\
 &\overset{\text{\tiny$(1)$}}{\leq} C\norm{\mathcal{R}'_v}_{C^0}+C_\ell b^{-(\ell-\varrho)}\norm{\Theta'}_{C^{\ell-\varrho}}\notag\\
 &\overset{\text{\tiny$(2)$}}{\leq} C\norm{\mathcal{R}'_v}_{C^0}+C_{\ell,1} b^{-\ell+\varrho}(\tau^{\ell}\norm{\mathfrak{p}}_{C^\varrho}+\norm{\mathfrak{p}}_{C^{\ell}})
\end{align}
for any $\ell\geq\varrho$.

Here in $(1)$ we use \eqref{for:37} of Section \ref{sec:19}; in $(2)$ we use \eqref{for:206}.

 Hence \eqref{for:252} follows from \eqref{for:188} and \eqref{for:251}.
\end{proof}

\subsection{ Iterative step and error estimate}\label{sec:26} In this part for any given  sufficiently small $C^\infty$ perturbation $\tilde{\alpha}_A$
of the action $\alpha_A$, which is also close to $\alpha_{E'(z)}$ on $E$ where $E'(z)$ is a standard perturbation of $E$,
we construct a $C^\infty$ diffeomorphism $h$ and a new standard perturbation $E'(z^{(1)})$ of $E$;  and then estimate the closeness between the new action $\tilde{\alpha}'_A=h^{-1}\tilde{\alpha}h$ and $\alpha_{E'(z^{(1)})}$ on $E$.

We recall notations in
Section \ref{sec:22}. We point out that the positive constants $\varrho$, $c_0$ and $\delta_\varrho$ (see Section \ref{sec:60}) will be used in this part.

\begin{proposition}\label{po:7} There exists $0<\bar{c}<1$ such that the following holds: if
$\tilde{\alpha}_A$ is a $C^\infty$ small perturbation of $\alpha_A$ which can be written as
\begin{align*}
 \tilde{\alpha}_v(x)=\exp(\mathfrak{p}_v(x))\cdot
z_v(x),\qquad \forall\,v\in E,\,x\in \mathcal{X}
\end{align*}
where $\mathfrak{p}_v\in C^\infty(\mathcal{X},\,\mathcal{V})$ (\eqref{sect:5} of Section \ref{sec:19}) and $z_v=w_v\cdot \alpha_v$, $w_v\in V^o$ for each  $v\in E$  satisfying
\begin{align*}
 \norm{\mathfrak{p}}_{C^{\varrho+1}}:=\max_{v\in E}\norm{\mathfrak{p}_v}_{C^{\varrho+1}}\leq \bar{c}\quad\text{and}\quad\textbf{d}(w,e)=\max_{v\in E}\textbf{d}(w_v,e)\leq \bar{c},
\end{align*}
then for any $b,\,\tau>1$,  there is $h\in C^\infty(\mathcal{X})$
such that the following hold:
\begin{enumerate}

  \item\label{for:211}  for any $m\geq\varrho$
  \begin{align*}
   d(h, I)_{C^m}\leq C_{m}b^\varrho(\tau^{m}\norm{\mathfrak{p}}_{C^\varrho}+\norm{\mathfrak{p}}_{C^{m}})
  \end{align*}
  where $\norm{\mathfrak{p}}_{C^m}=\max_{v\in E}\{\norm{\mathfrak{p}_{v}}_{C^m}\}$.

If $b,\,\tau$ are well chosen such that
 \begin{align*}
   b^\varrho\tau^{\varrho}\norm{\mathfrak{p}}_{C^{\varrho}}< \bar{c},
 \end{align*}
then $h$ is invertible;

\smallskip
\item\label{for:209} for $\tilde{\alpha}_A'=h^{-1}\tilde{\alpha}_Ah$, there exist $\mathfrak{p}_v^{(1)}\in C^\infty(\mathcal{X},\,\mathcal{V})$ and $w_v^{(1)}\in V^o$ for each $v\in E$ such that
\begin{align*}
 \tilde{\alpha}'_v(x)=\exp(\mathfrak{p}^{(1)}_v(x))\cdot z_v^{(1)}(x),\qquad \forall\,v\in E,\,x\in \mathcal{X}
\end{align*}
where $z_v^{(1)}=w_v^{(1)}\cdot z_v$ with the estimates:
\begin{align}\label{for:518}
  \norm{\mathfrak{p}^{(1)}}_{C^m}\leq C_{m}(\tau^{m}b^\varrho\norm{\mathfrak{p}}_{C^\varrho}+b^\varrho\norm{\mathfrak{p}}_{C^{m}}+1)
 \end{align}
for any $m\geq\varrho$ and
\begin{align*}
 \max_{z\in E} \textbf{d}(w_v^{(1)},\,e)<C\norm{\mathfrak{p}}_{C^0},
\end{align*}
where $\norm{\mathfrak{p}^{(1)}}_{C^m}:=\max_{z\in E}\norm{\mathfrak{p}^{(1)}_z}_{C^m}$;

\smallskip
\item \label{for:210}   the following estimate for $\norm{\mathfrak{p}^{(1)}}_{C^0}$ holds:
 \begin{align}\label{for:197}
\norm{\mathfrak{p}^{(1)}}&_{C^0}\leq  Cb^{2\varrho}\tau^{2\varrho}
(\norm{\mathfrak{p}}_{{C^{\varrho+1}}})^{1+y_0}\notag\\
 &+C\tau^{\varrho}
\big(C_m\tau^{-m}\norm{\mathfrak{p}}_{C^m}\big)^{\frac{1}{8}}(\norm{\mathfrak{p}}_{{C^\varrho}})^{\frac{7}{8}}\notag\\
&+C_{\ell} b^{-\ell+\varrho}(\tau^{\ell}\norm{\mathfrak{p}}_{C^\varrho}+\norm{\mathfrak{p}}_{C^{\ell}})
\end{align}
for any $m,\,\ell\geq\varrho$, where $ y_0=\frac{c_0}{8}$.

 \end{enumerate}

\end{proposition}

\begin{remark}\label{re:3}
In \eqref{for:197}, $C^m$ and $C^\ell$ norms of $\mathfrak{p}$ are used simultaneously to bound $\norm{\mathfrak{p}^{(1)}}_{C^0}$, which are more complex than classical KAM estimates. In \eqref{for:518} the term $\tau^{m}b^\varrho \norm{\mathfrak{p}}_{C^{\varrho}}$ is not found in previous KAM works.
  If the term $\tau^{m}b^\varrho \norm{\mathfrak{p}}_{C^{\varrho}}$ could be ignored, then \eqref{for:518} would be
\begin{align}\label{for:519}
  \norm{\mathfrak{p}^{(1)}}_{C^r}\leq C_r(b^\varrho\norm{\mathfrak{p}}_{C^r}+1),\qquad \forall\,r\geq\varrho.
 \end{align}
If we let $m=\ell$ and $b=\tau^2$ in \eqref{for:197}, then $\norm{\mathfrak{p}^{(1)}}_{C^0}$ can be simplified as (we note that $\varrho>2$):
 \begin{align}\label{for:520}
\norm{\mathfrak{p}^{(1)}}&_{C^0}\leq  C\tau^{6\varrho}
(\norm{\mathfrak{p}}_{{C^{\varrho+1}}})^{1+y_0}+C_\ell\tau^{\varrho}
\big(\tau^{-\ell}\norm{\mathfrak{p}}_{C^\ell}\big)^{\frac{1}{8}}(\norm{\mathfrak{p}}_{{C^\varrho}})^{\frac{7}{8}}\notag\\
&+C_{\ell} (\tau^{-\ell+2\varrho}\norm{\mathfrak{p}}_{C^\varrho}+\tau^{-2\ell+2\varrho}\norm{\mathfrak{p}}_{C^{\ell}}).
\end{align}
The proof of convergence would be standard for the KAM iteration if we had \eqref{for:519} and \eqref{for:520} (see \cite{Fayad}).
However, the presence of the additional term $\tau^{m}b^\varrho \norm{\mathfrak{p}}_{C^{\varrho}}$ directly arises from our method and cannot be ignored, as it results from the use of directional smoothing operators.

Since the term $\tau^{m\mu_1}\norm{\mathfrak{p}}_{C^{0}}$ with $\mu_1=1$ appears for the estimate of  $\norm{\mathfrak{p}^{(1)}}_{C^r}$ (see \eqref{for:518}) and
the term $\tau^{-m\mu_2}\norm{\mathfrak{p}}_{C^m}$ with $\mu_2=1$ appears for the estimate of $\norm{\mathfrak{p}^{(1)}}_{C^0}$ (see \eqref{for:197}), the KAM iteration may diverge in $C^r$ topology for large $r$. The general KAM scheme needs $\frac{\mu_2}{\mu_1}>2$ to guartentee the convergence in $C^r$ topology for each $r$ (in fact $\frac{\mu_2}{\mu_1}>1$ is sufficient)  (see \cite{Fayad}).

To overcome the difficulty,   two parameters $m$  and $\ell$ are used simultaneously to estimate $\norm{\mathcal{R}}_{C^0}$. This is the main difference with the
traditional KAM scheme where only the case of $m=\ell$ is considered.   The two-parameters trick was first introduced in \cite{W5}. The application in \cite{W5} relies on the fact that the estimate for $\norm{\mathfrak{p}^{(1)}}_{C^0}$ is independent of the extra parameter $m$. However, our case is quite different (see \eqref{for:197}). To overcome this difficulty, we modify the trick. The role of the two-parameters trick is explained precisely in  Remark \ref{re:6}.

\end{remark}

\subsubsection{Proof of Proposition \ref{po:7}} \emph{Basic notations}:

Let
\begin{gather}
  r_v=\int_{\mathcal{X}}\mathfrak{p}_v,\quad \mathfrak{p}'_v=\mathfrak{p}_{v}-r_v, \qquad\text{and}\notag\\
    \mathbb{D}_{u,v}(x)=\mathfrak{p}'_{u}(z_vx)-dv \mathfrak{p}'_{u}(x)-\big(\mathfrak{p}'_{v}(z_ux)-dz_u\mathfrak{p}'_{v}(
x)\big)\label{for:170}
\end{gather}
for any $u,v\in E$. Set $\norm{\mathbb{D}}_{C^m}=\max_{u,\,v\in E}\{\norm{\mathbb{D}_{v,u}}_{C^m}\}$.

We can assume $\bar{c}<\delta_\varrho$ (see Lemma \ref{le:1}). By Lemma \ref{le:1} we have
\begin{align}
 \max_{u,v\in E}\norm{\log[z_{v}&,z_{u}]}\leq C\norm{\mathfrak{p}}_{C^0}, \qquad\text{ and }\label{for:377}\\
\norm{\mathbb{D}}_{C^\varrho}
 &\leq C\norm{\mathfrak{p}}^2_{C^{\varrho+1}}.\label{for:190}
\end{align}
Let
\begin{align*}
 \mathfrak{c}_{\mu}=\big((I-dz_{\textbf{a}_1})|_{\mathfrak{g}_{\mu}}\big)^{-1}(r_{\textbf{a}_1}|_{\mathfrak{g}_{\mu}});\quad 1\neq\mu\in \Delta, \,\,\text{and}\,\,
 \mathfrak{c}_{1}=0;
\end{align*}
and let $\mathfrak{c}_{\textbf{a}_1}\in \mathcal{V}$ be the vector with coordinates $\mathfrak{c}_{\mu}$, $\mu\in \Delta$.

\smallskip

\emph{Construction of $\Theta$ and $\mathcal{R}_v$, $v\in E$}:

We apply Corollary \ref{cor:8} to almost twisted cocycle equation \eqref{for:170}. For any $b,\,\tau>1$ there exist $\Theta\in C^\infty(\mathcal{X},\,\mathcal{V})$ and $\mathcal{R}_v\in C^\infty(\mathcal{X},\,\mathcal{V})$ for each $v\in E$ such that
\begin{align}\label{for:255}
 \mathfrak{p}'_v=\Theta\circ z_v-dz_v\Theta+\mathcal{R}_v
\end{align}
with estimates: for any $v\in E$
\begin{align}\label{for:193}
 \max\{\norm{\Theta}_{C^m},\,\norm{\mathcal{R}_v}_{C^m}\}
 \leq C_{m}b^\varrho(\tau^{m}\norm{\mathfrak{p}}_{C^\varrho}+\norm{\mathfrak{p}}_{C^{m}})
\end{align}
for any $m\geq\varrho$;  and
\begin{align}\label{for:194}
\norm{\mathcal{R}_v}&_{C^0}\leq C\tau^{\varrho}\norm{\log[z,z]}^{\frac{c_0}{8}}\norm{\mathfrak{p}}_{C^\varrho}+C\tau^{\varrho}
(\norm{\mathbb{D}}_{C^\varrho})^{\frac{1}{8}}(\norm{\mathfrak{p}}_{{C^\varrho}})^{\frac{7}{8}}\notag\\
 &+C\tau^{\varrho}
\big(C_m\tau^{-m}\norm{\mathfrak{p}}_{C^m}\big)^{\frac{1}{8}}(\norm{\mathfrak{p}}_{{C^\varrho}})^{\frac{7}{8}}\notag\\
&+C\norm{\mathbb{D}}_{{C^\varrho}}+C_{\ell} b^{-\ell+\varrho}(\tau^{\ell}\norm{\mathfrak{p}}_{C^\varrho}+\norm{\mathfrak{p}}_{C^{\ell}})\notag\\
&\overset{\text{\tiny$(1)$}}{\leq} C\tau^{\varrho}\norm{\log[z,z]}^{\frac{c_0}{8}}\norm{\mathfrak{p}}_{C^\varrho}+C\tau^{\varrho}
(\norm{\mathfrak{p}}_{{C^{\varrho+1}}})^{\frac{1}{4}}(\norm{\mathfrak{p}}_{{C^\varrho}})^{\frac{7}{8}}\notag\\
 &+C\tau^{\varrho}
\big(C_m\tau^{-m}\norm{\mathfrak{p}}_{C^m}\big)^{\frac{1}{8}}(\norm{\mathfrak{p}}_{{C^\varrho}})^{\frac{7}{8}}\notag\\
&+C\norm{\mathfrak{p}}^2_{C^{\varrho+1}}+C_{\ell} b^{-\ell+\varrho}(\tau^{\ell}\norm{\mathfrak{p}}_{C^\varrho}+\norm{\mathfrak{p}}_{C^{\ell}})\notag\\
&\overset{\text{\tiny$(2)$}}{\leq} 3C\tau^{\varrho}
(\norm{\mathfrak{p}}_{{C^{\varrho+1}}})^{1+\frac{c_0}{8}}
 +C\tau^{\varrho}
\big(C_m\tau^{-m}\norm{\mathfrak{p}}_{C^m}\big)^{\frac{1}{8}}(\norm{\mathfrak{p}}_{{C^\varrho}})^{\frac{7}{8}}\notag\\
&+C_{\ell} b^{-\ell+\varrho}(\tau^{\ell}\norm{\mathfrak{p}}_{C^\varrho}+\norm{\mathfrak{p}}_{C^{\ell}})
\end{align}
for any $m,\,\ell\geq \varrho$.

Here in $(1)$ we use \eqref{for:377} and \eqref{for:190}; in $(2)$ we use $\norm{\mathfrak{p}}_{C^{\varrho+1}}<\bar{c}<1$, $c_0<1$ (see \eqref{for:304} of Section \ref{sec:40}) and $\tau>1$.

\smallskip

\emph{Construction of $h$ and $z^{(1)}_v$, $v\in E$}:

Let $\mathfrak{h}=\Theta+c_{\textbf{a}_1}$
and let
\begin{gather*}
 h(x)=\exp(\mathfrak{h}(x))x,\quad \forall\,x\in \mathcal{X},\\
 w_v^{(1)}=\exp(r_v|_{\mathfrak{g}_{1}}),\quad z_v^{(1)}=w_v^{(1)}\cdot z_v,\,\,v\in E.
\end{gather*}
\emph{Note}.  $w_v^{(1)}$, $v\in E$ is the coordinate change, which is determined by $r_v$, $v\in E$.

\smallskip
\emph{Estimates of $c_{\textbf{a}_1}$ and $w^{(1)}_v$, $v\in E$}:

From \eqref{for:138} of Section \ref{sec:19} we see that $dz_{\textbf{a}_1}$ has the same Lyapunov exponents as $d\textbf{a}_1$. Hence, we have
\begin{align}\label{for:378}
 \norm{c_{\textbf{a}_1}}\leq C\norm{\mathfrak{p}}_{C^0}.
\end{align}
From \eqref{for:365} of \eqref{sect:9} of Section \ref{sec:19} we have
\begin{align}\label{for:379}
 \max_{z\in E} \textbf{\emph{d}}(w_v^{(1)},\,e)\leq C\max_{z\in E}\norm{r_v|_{\mathfrak{g}_{1}}}\leq C_1\norm{\mathfrak{p}}_{C^0}.
\end{align}

\eqref{for:211}: We have
\begin{align}\label{for:212}
 d(h, I)_{C^m}&\leq C_m\norm{\mathfrak{h}}_{C^m} \leq C_m(\norm{\Theta}_{C^m}+\norm{c_{\textbf{a}_1}})\notag\\
 &\overset{\text{\tiny$(1)$}}{\leq} C_{m,1}b^\varrho(\tau^{m}\norm{\mathfrak{p}}_{C^\varrho}+\norm{\mathfrak{p}}_{C^{m}})
 \end{align}
for any $m\geq\varrho$.  Here in $(1)$ we use \eqref{for:193}, \eqref{for:378} and the assumption $b>1$.

In \eqref{for:212} letting $m=\varrho$, we have
\begin{align}\label{for:322}
 d(h, I)_{C^\varrho}&\leq Cb^\varrho(\tau^{\varrho}\norm{\mathfrak{p}}_{C^\varrho}+\norm{\mathfrak{p}}_{C^{\varrho}})\leq 2Cb^\varrho\tau^{\varrho}\norm{\mathfrak{p}}_{C^\varrho}<2C\bar{c}.
 \end{align}
We can assume that $\bar{c}$ is sufficiently small such that $h$ is invertible; moreover,
\begin{align}\label{for:208}
  d(h^{-1}, I)_{C^m}\leq C_md(h, I)_{C^m},\qquad \forall\,m\geq 0,
\end{align}
see \cite[Lemma AII.26]{Lazutkin}.

\medskip
\eqref{for:209}: \eqref{for:208} justifies the definition of $\tilde{\alpha}_A'$. For any $v\in E$ we have
\begin{align*}
  d(\tilde{\alpha}'_v,& \,z_v^{(1)})_{C^0}=d(h^{-1}\circ \tilde{\alpha}_v\circ h, \,z_v^{(1)})_{C^0} \\
  &\leq d(h^{-1}\circ \tilde{\alpha}_v\circ h, \tilde{\alpha}_v\circ h)_{C^0}+d(\tilde{\alpha}_v\circ h, z_v\circ h)_{C^0}\\
  &+d(z_v\circ h, z_v)_{C^0}+d(z_v, z_v^{(1)})_{C^0}\\
  &\leq d(h^{-1}, I)_{C^0}+d(\tilde{\alpha}_v, z_v)_{C^0}+\norm{z_v}_{C^1}d(h, I)_{C^0}\\
  &+C\textbf{\emph{d}}(w_v^{(1)},\,e)\\
  &\overset{\text{\tiny$(1)$}}{\leq} Cd(h, I)_{C^0}+C\norm{\mathfrak{p}}_{C^0}+Cd(h, I)_{C^0}+C_1\norm{\mathfrak{p}}_{C^0}\\
  &\overset{\text{\tiny$(2)$}}{\leq} C_2\bar{c}.
\end{align*}
Here in $(1)$ we use \eqref{for:208} and \eqref{for:379};  $(2)$ we use \eqref{for:322}.

Hence, if $\bar{c}$ is sufficiently small there exists $\mathfrak{p}_v^{(1)}\in C^\infty(\mathcal{X},\,\mathcal{V})$ for each $v\in E$ such that
\begin{align*}
 \tilde{\alpha}'_v(x)=\exp(\mathfrak{p}^{(1)}_v(x))\cdot
z_v^{(1)}(x),\qquad \forall\,v\in E,\,x\in \mathcal{X}.
\end{align*}
Then we have
\begin{align}
 \norm{\mathfrak{p}^{(1)}}_{C^m}&\leq C_m\max_{v\in E}d(\tilde{\alpha}'_v, \,z_v^{(1)})_{C^m}\leq C_m(\max_{v\in E}\norm{h^{-1}\tilde{\alpha}_vh}_{C^m}+1)\notag\\
 &\overset{\text{\tiny$(1)$}}{\leq} C_{m,1}(\norm{h^{-1}}_{C^m}+\norm{h}_{C^m}+\max_{v\in E}\norm{\tilde{\alpha}_v}_{C^m}+1)\notag\\
 &\overset{\text{\tiny$(2)$}}{\leq} C_{m,2}(\norm{h}_{C^m}+\norm{h}_{C^m}+1)\notag\\
 &\overset{\text{\tiny$(3)$}}{\leq}C_{m,3}(b^\varrho\tau^{m}\norm{\mathfrak{p}}_{C^\varrho}+b^\varrho\norm{\mathfrak{p}}_{C^{m}}+1)
\end{align}
for any $m\geq\varrho$. Here $(1)$ is a standard ``linear" estimate, see for example, \cite[Propositions 5.5]{Llave}; in $(2)$
we use \eqref{for:208}; in $(3)$ we use \eqref{for:212}.

\medskip

\eqref{for:210}: For each $v\in E$,  we have
\begin{gather*}
 h\big(\tilde{\alpha}'_v(x)\big)=\tilde{\alpha}_v(h(x))\\
 \Rightarrow \exp\big(\mathfrak{h}(\tilde{\alpha}'_v(x))\big)\cdot\tilde{\alpha}'_v(x)=\exp\big(\mathfrak{p}_v(h(x))\big)\cdot
z_v\big(\exp(\mathfrak{h}(x))\cdot x\big)\\
\Rightarrow \exp\big(\mathfrak{h}(\tilde{\alpha}'_v(x))\big)\cdot\exp(\mathfrak{p}^{(1)}_v(x))\cdot w_v^{(1)}\\
=\exp\big(\mathfrak{p}_v(h(x))\big)\cdot
\exp(dz_v\mathfrak{h}(x)).
\end{gather*}
By  using Baker-Campbell-Hausdorff formula, for each $v\in E$,   we have
\begin{align*}
  \mathfrak{p}^{(1)}_v&=-\mathfrak{h}\circ \tilde{\alpha}'(v)+\mathfrak{p}_v\circ h+dz_v\mathfrak{h}-r_v|_{\mathfrak{g}_{1}}+\Psi_v
\end{align*}
where
\begin{align}\label{for:196}
 \max_{v\in E}\big\|\Psi_v\big\|_{C^0}&\leq C\max_{v\in E}\{\norm{\mathfrak{h}}_{C^0}, \norm{\mathfrak{p}}_{C^0},\,\norm{\log w_v^{(1)}}\}^2\notag\\
 &\overset{\text{\tiny$(a)$}}{\leq}Cb^{2\varrho}\tau^{2\varrho}(\norm{\mathfrak{p}}_{C^\varrho})^2,
\end{align}
where in $(a)$ we use \eqref{for:212} and recall $\log w_v^{(1)}=r_v|_{\mathfrak{g}_{1}}$.
It follows that
\begin{align}\label{for:510}
 \mathfrak{p}^{(1)}_v&=\big(-\mathfrak{h}\circ \tilde{\alpha}'(v)+\mathfrak{h}\circ z^{(1)}_v\big)+\big(\mathfrak{h}\circ z_v-\mathfrak{h}\circ z^{(1)}_v\big)\notag\\
 &+\big(-\mathfrak{h}\circ z_v-r_v|_{\mathfrak{g}_{1}}+\mathfrak{p}_v+dz_v\mathfrak{h}\big)\notag\\
 &+\big(\mathfrak{p}_v\circ h-\mathfrak{p}_v\big)+\Psi_v.
\end{align}
We note that
\begin{align}\label{for:511}
 &-\mathfrak{h}\circ z_v-r_v|_{\mathfrak{g}_{1}}+\mathfrak{p}_v+dz_v\mathfrak{h}\notag\\
 &=-(\Theta+c_{\textbf{a}_1})\circ z_v-r_v|_{\mathfrak{g}_{1}}+\mathfrak{p}_v+dz_v(\Theta+c_{\textbf{a}_1})\notag\\
 &\overset{(1)}{=}\mathcal{R}_v+\big(r_v+(dz_v-I)c_{\textbf{a}_1}-r_v|_{\mathfrak{g}_{1}}\big)
\end{align}
Here in $(1)$ we recall that
\begin{align*}
 \mathcal{R}_v=\mathfrak{p}_v-r_v-(\Theta\circ z_v-dz_v\Theta)
\end{align*}
(see \eqref{for:255}).

Recall $c_{\textbf{a}_1}|_{\mathfrak{g}_{1}}=0$ and the fact that $\mathfrak{g}_{\mu}$, $\mu\in \Delta$ are invariant subspaces for $dz_v$, $v\in E$ (see \eqref{sect:10} of Section \ref{sec:19}).
Consequently, we see that: if $1\neq \mu\in \Delta$
\begin{align}\label{for:513}
 \big(r_v-r_v|_{\mathfrak{g}_{1}}+(dz_v-I)c_{\textbf{a}_1}\big)|_{\mathfrak{g}_{\mu}}
 =\big(r_v+(dz_v-I)c_{\textbf{a}_1}\big)|_{\mathfrak{g}_{\mu}}
\end{align}
and
\begin{align}\label{for:514}
 \big(r_v-r_v|_{\mathfrak{g}_{1}}+(dz_v-I)c_{\textbf{a}_1}\big)|_{\mathfrak{g}_{1}}=0.
\end{align}
It follows from \eqref{for:510} and \eqref{for:511} that
\begin{align}\label{for:512}
 \norm{\mathfrak{p}^{(1)}}_{C^0}\leq \mathcal{E}_1+\mathcal{E}_2+\mathcal{E}_3+\mathcal{E}_4,
\end{align}
where
\begin{gather*}
  \mathcal{E}_1=\max_{v\in E}\|-\mathfrak{h}\circ \tilde{\alpha}'(v)+\mathfrak{h}\circ z_v^{(1)}\|_{C^0},\\
  \mathcal{E}_2=\max_{v\in E}\|-\mathfrak{h}\circ z_v^{(1)}+\mathfrak{h}\circ z_v \|_{C^0}
  +\max_{v\in E}\|\mathfrak{p}_v\circ h-\mathfrak{p}_v \|_{C^0},\\
  \mathcal{E}_3=\max_{v\in E}\|\mathcal{R}_v\|_{C^0}
  +\max_{v\in E}\|\Psi_v \|_{C^0}\\
  \mathcal{E}_4=\max_{v\in E}\|r_v-r_v|_{\mathfrak{g}_{1}}+(dz_v-I)c_{\textbf{a}_1}\|.
\end{gather*}
We note that
\begin{align}\label{for:259}
\mathcal{E}_1&\leq C\norm{\mathfrak{h}}_{C^1} \norm{\mathfrak{p}^{(1)}}_{C^0}\overset{\text{\tiny$(1)$}}{\leq} C_2b^\varrho\tau^{\varrho}\norm{\mathfrak{p}}_{C^{\varrho}}\norm{\mathfrak{p}^{(1)}}_{C^0}\overset{\text{\tiny$(2)$}}{<}\text{\small$\frac{1}{2}$}\norm{\mathfrak{p}^{(1)}}_{C^0}.
\end{align}
Here in $(1)$ we use \eqref{for:212}; in $(2)$ by assumption $b^\varrho\tau^{\varrho}\norm{\mathfrak{p}}_{C^{\varrho}}<\bar{c}$. Thus
$C_1b^\varrho\tau^{\varrho}\norm{\mathfrak{p}}_{C^{\varrho}}<\frac{1}{2}$ providing $\bar{c}$ is sufficiently small.

It follows from \eqref{for:512} and \eqref{for:259} that
\begin{align}\label{for:384}
 \norm{\mathfrak{p}^{(1)}}_{C^0}\leq 2\mathcal{E}_2+2\mathcal{E}_3+2\mathcal{E}_4.
\end{align}
Firstly, we estimate $\mathcal{E}_2$:
\begin{align}\label{for:383}
  \mathcal{E}_2&\leq \norm{\mathfrak{h}}_{C^1}\max_{v\in E}\emph{\textbf{d}}(z_v^{(1)},z_v)+\norm{\mathfrak{p}}_{C^1}d(h,I)_{C^0}\notag\\
  &\overset{(1)}{\leq} Cb^\varrho\tau^{\varrho}\norm{\mathfrak{p}}_{C^{\varrho}}
 \norm{\mathfrak{p}}_{C^{0}}+Cb^\varrho\tau^{\varrho}\norm{\mathfrak{p}}_{C^{\varrho}}
 \norm{\mathfrak{p}}_{C^{1}}
\end{align}
Here in $(1)$ we use \eqref{for:379} and \eqref{for:212}.

Secondly, we estimate $\mathcal{E}_3$:
\begin{align}\label{for:516}
  \mathcal{E}_3&\overset{(1)}{\leq} C\tau^{\varrho}(\norm{\mathfrak{p}}_{C^{\varrho+1}})^{1+\frac{c_0}{8}}
+C\tau^{\varrho}(C_{m}\tau^{-m+\varrho}\norm{\mathfrak{p}}_{C^m})^{\frac{1}{8}}(\norm{\mathfrak{p}}_{C^\varrho})^{\frac{7}{8}}\notag\\
&+C_{\ell} b^{-\ell+\varrho}(\tau^{\ell}\norm{\mathfrak{p}}_{C^\varrho}+\norm{\mathfrak{p}}_{C^{\ell}})+Cb^{2\varrho}\tau^{2\varrho}(\norm{\mathfrak{p}}_{C^\varrho})^2.
\end{align}
for any $m,\,\ell\geq\varrho$.

Here in $(1)$ we use \eqref{for:194} and \eqref{for:196}.

Finally, we estimate $\mathcal{E}_4$: we show that
\begin{align}\label{for:515}
 \mathcal{E}_4\leq C(\norm{\mathfrak{p}}_{C^1})^2.
\end{align}
From \eqref{for:513} and  \eqref{for:514}, it suffices to prove: for any $1\neq\mu\in \Delta$
\begin{align}\label{for:382}
 \big\|r_v|_{\mathfrak{g}_{\mu}}-(I-dz_v)c_{\mu}\big\|\leq C(\norm{\mathfrak{p}}_{C^1})^2,\quad \forall\,v\in E.
\end{align}
We note that
\begin{align}\label{for:381}
  &(I-dz_{\textbf{a}_1})\big(r_v|_{\mathfrak{g}_{\mu}}-(I-dz_v)c_{\mu}\big)\notag\\
  &=(I-dz_{\textbf{a}_1})(r_v|_{\mathfrak{g}_{\mu}})-(I-dz_{v})(I-dz_{\textbf{a}_1})c_{\mu}\notag\\
  &+(dz_{v}dz_{\textbf{a}_1}-dz_{\textbf{a}_1}dz_{v})c_{\mu}\notag\\
  &\overset{(1)}{=}(I-dz_{\textbf{a}_1})(r_v|_{\mathfrak{g}_{\mu}})-(I-dz_{v})(r_{\textbf{a}_1}|_{\mathfrak{g}_{\mu}})\notag\\
  &+(dz_{v}dz_{\textbf{a}_1}-dz_{\textbf{a}_1}dz_{v})c_{\mu}\notag\\
  &\overset{(2)}{=}(r_v-dz_{\textbf{a}_1}r_v-r_{\textbf{a}_1}+dz_vr_{\textbf{a}_1})\big|_{\mathfrak{g}_{\mu}} \notag\\
  &+(dz_{v}dz_{\textbf{a}_1}-dz_{\textbf{a}_1}dz_{v})c_{\mu}.
\end{align}
Here in $(1)$ we recall definition of $c_{\mu}$; in $(2)$ we use the fact that each $\mathfrak{g}_{\mu}$ is invariant under $dz_v$, $v\in E$ (see \eqref{sect:10} of Section \ref{sec:19}).

It follows that: for any $1\neq\mu\in \Delta$
\begin{align*}
  &\big\|r_v|_{\mathfrak{g}_{\mu,A}}-(I-\text{Ad}_{z_v})c_{\mu}\big\|\notag\\
  &=\Big\|(I-dz_{\textbf{a}_1})^{-1}(I-dz_{\textbf{a}_1})\big(r_v|_{\mathfrak{g}_{\mu}}-(I-dz_v)c_{\mu}\big)\Big\|\notag\\
  &\leq C\Big\|(I-dz_{\textbf{a}_1})\big(r_v|_{\mathfrak{g}_{\mu}}-(I-dz_v)c_{\mu}\big)\Big\|\notag\\
  &\overset{\text{\tiny$(1)$}}{\leq} \big\|(r_v-dz_{\textbf{a}_1}r_v-r_{\textbf{a}_1}+dz_vr_{\textbf{a}_1})\big|_{\mathfrak{g}_{\mu}} \big\|\notag\\
  &+\big\|(dz_{v}dz_{\textbf{a}_1}-dz_{\textbf{a}_1}dz_{v})c_{\mu}\big\|\notag\\
  &\overset{\text{\tiny$(2)$}}{\leq} C\norm{\mathfrak{p}}_{C^1}\norm{\mathfrak{p}}_{C^0}+C\norm{\log[z_{\textbf{a}_1},z_v]}\norm{c_{\textbf{a}_1}}\notag\\
 &\overset{\text{\tiny$(3)$}}{\leq}  C\norm{\mathfrak{p}}_{C^1}\norm{\mathfrak{p}}_{C^0}+C\norm{\mathfrak{p}}_{C^0}^{2}.
\end{align*}
Here in $(1)$ we use \eqref{for:381}; in $(2)$ we use \eqref{for:348} of Lemma \ref{le:1};
in $(3)$ we use \eqref{for:347} of Lemma \ref{le:1} and \eqref{for:378}.  Hence we get \eqref{for:382}.

\medskip
Finally, it follows from \eqref{for:384}, \eqref{for:383}, \eqref{for:516}  and \eqref{for:515} that
\begin{align}\label{for:517}
\norm{\mathfrak{p}^{(1)}}_{C^0}&\leq Cb^\varrho\tau^{\varrho}\norm{\mathfrak{p}}_{C^{\varrho}}
 \norm{\mathfrak{p}}_{C^{0}}+Cb^\varrho\tau^{\varrho}\norm{\mathfrak{p}}_{C^{\varrho}}
 \norm{\mathfrak{p}}_{C^{1}}\notag\\
 &+C\tau^{\varrho}(\norm{\mathfrak{p}}_{C^{\varrho+1}})^{1+\frac{c_0}{8}}
+C\tau^{\varrho}(C_{m}\tau^{-m+\varrho}\norm{\mathfrak{p}}_{C^m})^{\frac{1}{8}}(\norm{\mathfrak{p}}_{C^\varrho})^{\frac{7}{8}}\notag\\
&+C_{\ell} b^{-\ell+\varrho}(\tau^{\ell}\norm{\mathfrak{p}}_{C^\varrho}+\norm{\mathfrak{p}}_{C^{\ell}})+Cb^{2\varrho}\tau^{2\varrho}(\norm{\mathfrak{p}}_{C^\varrho})^2\notag\\
&+C(\norm{\mathfrak{p}}_{C^1})^2
\end{align}
for any $m,\,\ell\geq\varrho$.

Noting that $\norm{\mathfrak{p}}_{C^{\varrho+1}}<1$, $b,\,\tau>1$ and $0<c_0<1$ (\eqref{for:304} of Section \ref{sec:40}),  we can simplify \eqref{for:517} as
\begin{align*}
\norm{\mathfrak{p}^{(1)}}_{C^0}
&\leq 5Cb^{2\varrho}\tau^{2\varrho}(\norm{\mathfrak{p}}_{C^{\varrho+1}})^{1+\frac{c_0}{8}}+C\tau^{\varrho}(C_{m}\tau^{-m+\varrho}\norm{\mathfrak{p}}_{C^m})^{\frac{1}{8}}
(\norm{\mathfrak{p}}_{C^\varrho})^{\frac{7}{8}}\notag\\
&+C_{\ell} b^{-\ell+\varrho}(\tau^{\ell}\norm{\mathfrak{p}}_{C^\varrho}+\norm{\mathfrak{p}}_{C^{\ell}}).
\end{align*}
Hence we finish the proof.

\subsection{The iterative process }\label{sec:21}  In this part we obtain a smooth
conjugacy between the initial perturbation $\tilde{\alpha}_A$ and $\alpha_A$ by using the KAM iteration.

\subsubsection{Constants in this part} We recall the following positive constants defined in previous sections:
\begin{enumerate}
  \item $\varrho$ is defined in Section \ref{for:504};

  \smallskip
  \item $\bar{c}$ is from Proposition \ref{po:7};

  \smallskip
  \item $y_0=\frac{c_0}{8}$, where $c_0$ is defined in \eqref{for:304} of Section \ref{sec:40}.
\end{enumerate}
To set up the iterative process we define positive constants $\gamma,\,\gamma_1,\,x_0,\,\delta$ and $\ell$ as follows:
we  first pick up
\begin{align*}
 1<\gamma<1+y_0,\quad \gamma_1>8,\quad 0<x_0<1\text{ satisfying }1-x_0<(\gamma-1)^2.
\end{align*}
Since $\gamma_1+8\gamma-7<\gamma\gamma_1+1$ we pick up $\nu$ such that
\begin{align*}
 \gamma_1+8\gamma-7<\nu<\gamma\gamma_1+1.
\end{align*}
 Fix
\begin{align}\label{for:326}
 0<\delta<\min\{\frac{1}{2},\,1-\text{\small$\frac{\gamma}{2}$},\,2-\gamma,\,\gamma\gamma_1+1-\nu,\,\frac{\nu-(\gamma_1+8\gamma-7)}{8},\,1-x_0\}.
\end{align}
Let $\ell$ be sufficiently large such that the following are satisfied:
\begin{gather}
1-\text{\small$\frac{(1+\gamma_1)(\varrho+1)}{\ell}$}>\text{\small$\frac{1}{2}$}+\delta,\quad
1-\text{\tiny$\text{\small$\frac{3\nu\varrho}{\ell}$}$}-\text{\small$\frac{(1+\gamma_1)\varrho}{\ell}$}> \text{\small$\frac{\gamma}{2}$}+\delta;\label{for:220}\\
\text{\small$\frac{(2\nu+1+\gamma_1)\varrho}{\ell}$}+\nu-1<\gamma\gamma_1-\delta;\label{for:221}\\
-\text{\small$\frac{6\nu}{\ell}$}+(1-\text{\small$\frac{(1+\gamma_1)
(\varrho+1)}{\ell}$})(1+y_0)>\gamma+\delta; \label{for:222}\\
\min\{\text{\small$\frac{2\nu(\ell-\varrho)}{\ell}$}-\nu+1-\text{\small$\frac{(1+\gamma_1)\varrho}{\ell}$},\, \text{\small$\frac{2\nu(\ell-\varrho)}{\ell}$}-\gamma_1\}>\gamma+\delta; \label{for:224}
\end{gather}
\begin{gather}
-\text{\small$\frac{\nu}{\ell}$}+(\nu-\gamma_1)\frac{1}{8}
+(1-\text{\small$\frac{(1+\gamma_1)\varrho}{\ell}$})\frac{7}{8}>\gamma+\delta;
\label{for:223}\\
\text{\small$\frac{2\nu\varrho}{\ell}$}+8(\gamma-x_0)+8\gamma<8\gamma^2; \label{for:213}\\
-\text{\small$\frac{\nu\varrho}{\ell}$}
+(8\gamma-8x_0+1-\text{\small$\frac{\varrho(1+\gamma_1)}{\ell}$})\frac{1}{8}+(1-\text{\small$\frac{\varrho(1+\gamma_1)}{\ell}$})
\frac{7}{8}
>\gamma+\delta\label{for:254}
\end{gather}
Next, we show that it is possible to choose an $\ell$ that satisfies all these constraints. By letting $\ell\to\infty$ inequalities \eqref{for:220} to \eqref{for:254} become
\begin{gather*}
 \eqref{for:220}\to\big(1>\text{\small$\frac{1}{2}$}+\delta,\quad 1>\text{\small$\frac{\gamma}{2}$}+\delta\big),\quad \eqref{for:221}\to\big(\nu-1<\gamma\gamma_1-\delta \big)\\
 \eqref{for:222}\to\big(1+y_0>\gamma+\delta \big) \quad  \eqref{for:224}\to\big(\min\{\nu+1,2\nu-\gamma_1\}>\gamma+\delta \big),\\
 \eqref{for:223}\to\big( \nu-\gamma_1+7>8\gamma+8\delta \big), \quad  \eqref{for:213}\to\big( 8(\gamma-x_0)+8\gamma<8\gamma^2 \big),\\
 \eqref{for:254}\to\big( \gamma-x_0+1>\gamma+\delta \big).
\end{gather*}
All the above  inequalities hold either automatically or as a direct consequence of assumptions. Thus inequalities \eqref{for:220} to \eqref{for:254} hold if we choose $\ell$ big enough.

We fix an increasing sequence $\beta_n\to \infty$ with $\beta_1>2\ell$.

\subsubsection{Inductive assumption}
We construct $\mathfrak{p}_v^{(n)}$, $v\in E$ and $h_n$ inductively as follows. We denote
\begin{align*}
 \tilde{\alpha}_v(x)=\exp(\mathfrak{p}_v(x))\cdot
\alpha_v(x),\qquad \forall\,v\in E,\,x\in \mathcal{X}.
\end{align*}
Set
\begin{align*}
 \tilde{\alpha}^{(0)}=\tilde{\alpha},\quad \mathfrak{p}_v^{(0)}=\mathfrak{p}_v,\quad z_v^{(0)}=v,\quad w_v^{(0)}=e, \quad h_{0}=I \quad\text{and}\quad\epsilon_n=\epsilon^{\text{\tiny$\gamma^n$}}
\end{align*}
where $0<\epsilon^{\text{\tiny$\frac{1}{4}$}}<\bar{c}$ is sufficiently small so that the following hold:
\begin{align*}
\norm{\mathfrak{p}^{(0)}}_{C^0}\leq \epsilon_0=\epsilon,\quad \norm{\mathfrak{p}^{(0)}}_{C^{\ell}}\leq \epsilon_0^{-\gamma_1}, \quad d(h_0, I)_{C^1}<\epsilon_0^{\text{\tiny$\frac{1}{2}$}},
\end{align*}
where $\norm{\mathfrak{p}^{(0)}}_{C^r}=\max_{v\in E}\norm{\mathfrak{p}^{(0)}_v}_{C^r}$ for any $r\geq0$.

Suppose  inductively that
\begin{align*}
 \tilde{\alpha}^{(n)}_v(x)=\exp(\mathfrak{p}^{(n)}_v(x))\cdot
z_v^{(n)}(x),\qquad \forall\,v\in E,\,x\in \mathcal{X}
\end{align*}
where
\begin{gather*}
 \tilde{\alpha}^{(n)}=h_{n}^{-1}\circ \tilde{\alpha}^{(n-1)}\circ h_{n},\quad z_v^{(n)}=w_v^{(n)}\cdot z_v^{(n-1)}
\end{gather*}
with $w_v^{(n)}\in V^o$, $v\in E$ and
\begin{align}\label{for:323}
 &\norm{\mathfrak{p}^{(n)}}_{C^0}\leq \epsilon_n, \quad \norm{\mathfrak{p}^{(n)}}_{C^{\ell}}\leq \epsilon^{-\gamma_1}_n,\quad \textbf{\emph{d}}(w^{(n)},e)\leq \epsilon_n\notag\\
&\norm{\mathfrak{p}^{(n)}}_{C^{\beta_q}}<(r_q)^n\epsilon^{\text{\tiny$-8\gamma$}}_n(\norm{\mathfrak{p}^{(q-1)}}_{C^{\beta_q}}+1),\quad d(h_{n},I)_{C^1}<\epsilon_n^{\frac{1}{2}}
 \end{align}
for any $1\leq q\leq n$; where $r_q$ is a constant dependent only on $q$ and
\begin{align*}
 \norm{\mathfrak{p}^{(n)}}_{C^m}=\max_{v\in E}\norm{\mathfrak{p}^{(n)}_v}_{C^m},\quad \textbf{\emph{d}}(w^{(n)},e)=\max_{v\in E}\textbf{\emph{d}}(w^{(n)}_v,e).
\end{align*}

\begin{remark}\label{re:6} Firstly, we explain briefly how the constants $\gamma,\,\gamma_1,\,\tau,\,b,\,m$ are chosen to ensure the induction works for $C^0$ and $C^{\ell}$ norms.
Let $m=\ell$ and $b=\tau^2$. From \eqref{for:520} we see that the new error $\norm{\mathfrak{p}^{(n+1)}}_{C^{0}}$ is at best to be $(\norm{\mathfrak{p}^{(n)}}_{{C^{\varrho+1}}})^{1+y_0}$ small. Thus we let $1<\gamma<1+y_0$. Moreover, from \eqref{for:520} we see that the main increasing parts for  $\norm{\mathfrak{p}^{(n+1)}}_{C^{0}}$ is $(\tau^{-\ell}\epsilon^{-\gamma_1}_n)^{\frac{1}{8}}\epsilon_n^{\frac{7}{8}}$. By  \eqref{for:209} of Proposition \ref{po:7}, the main increasing parts for $\norm{\mathfrak{p}^{(n+1)}}_{C^{\ell}}$ is
$\tau^{\ell}\epsilon_n$. To make the iteration work we should have
\begin{gather*}
 \tau^{\ell}\epsilon_n<\epsilon_{n+1}^{-\gamma_1}=\epsilon_{n}^{-\gamma\gamma_1}\quad\text{and}\quad (\tau^{-\ell}\epsilon^{-\gamma_1}_n)^{\frac{1}{8}}\epsilon_n^{\frac{7}{8}}<\epsilon_{n+1}=\epsilon_{n}^\gamma\\
 \Rightarrow \epsilon^{-(\gamma_1+8\gamma-7)}_n <\tau^{\ell}<\epsilon_{n}^{-\gamma\gamma_1-1}\\
 \Rightarrow \epsilon^{-(\gamma_1+8\gamma-7)}_n <\epsilon_{n}^{-\gamma\gamma_1-1}\\
 \Rightarrow \gamma_1+8\gamma-7<\gamma\gamma_1+1\\
 \Rightarrow \gamma>1\text{ and } \gamma_1>8.
\end{gather*}
As a result, we choose $\gamma_1+8\gamma-7<\nu<\gamma\gamma_1+1$ and let $b=\tau^2=\epsilon^{-\text{\small$\frac{2\nu}{\ell}$}}_n$.

To show the convergence in $C^{\beta_q}$ norm for any $q$. This means that we need to bound the term $\tau^{\beta_q}\norm{\mathfrak{p}^{(n)}}_{C^\varrho}$
 by $\norm{\mathfrak{p}^{(n)}}_{C^{\beta_q}}$. We introduce   parameters $\zeta_q$  (see \eqref{for:214}), which is used to measure
the ``competition'' between the increasing speed of $\norm{\mathfrak{p}^{(n)}}_{C^{\beta_q}}$ and $(b^{\frac{1}{2}})^{\beta_q}\norm{\mathfrak{p}^{(n)}}_{C^\varrho}$.

Suppose $\zeta=\zeta_q$. If $\zeta>b^{\frac{1}{2}}$, then $\norm{\mathfrak{p}^{(n)}}_{C^{\beta_q}}$ wins the game.  We let $\tau=b^{\frac{1}{2}}$ and $m=\ell$. Thus, $\tau^{\beta_q}\norm{\mathfrak{p}^{(n)}}_{C^\varrho}$ is bounded by  $\norm{\mathfrak{p}^{(n)}}_{C^{\beta_q}}$. If $\zeta\leq b^{\frac{1}{2}}$,
then $(b^{\frac{1}{2}})^{\beta_q}\norm{\mathfrak{p}^{(n)}}_{C^\varrho}$ wins. Thus, we have to let $\tau=\zeta$ so that $\tau^{\beta_q}\norm{\mathfrak{p}^{(n)}}_{C^\varrho}$ is also bounded by  $\norm{\mathfrak{p}^{(n)}}_{C^{\beta_q}}$. However, as $\tau\leq b^{\frac{1}{2}}$, we face the problem that
$\tau^{-\ell}\norm{\mathfrak{p}^{{(n)}}}_{C^{\ell}}$ needs not to be small when estimating $\norm{\mathfrak{p}^{{(n+1)}}}_{C^{0}}$. The key observation is in this case
$\zeta_q^{-\beta_q}\norm{\mathfrak{p}^{{(n)}}}_{C^{\beta_q}}$ is sufficiently small. Then we let $m=\beta_q$  instead of $\ell$ to bound $\norm{\mathfrak{p}^{{(n+1)}}}_{C^{0}}$.   This is where the two-parameters trick plays the crucial role.

\end{remark}

\subsection{Convergence}\label{sec:28} In this part, we use induction to prove that all the bounds in \eqref{for:323} remain valid for any $n\in\NN$.
\begin{proposition}\label{po:8}
Suppose $n\geq0$ and all the bounds in \eqref{for:323} hold for $n$.
Then there is $h_{n+1}\in C^\infty(\mathcal{X})$, $\mathfrak{p}_v^{(n+1)}\in C^\infty(\mathcal{X},\,\mathcal{V})$ and $w_v^{(n+1)}\in V^o$ for each $v\in E$
such that
for
\begin{align*}
 \tilde{\alpha}^{(n+1)}_v(x)=h_{n+1}^{-1}\circ \tilde{\alpha}^{(n)}\circ h_{n+1}=\exp(\mathfrak{p}^{(n+1)}_v(x))\cdot
z_v^{(n+1)}(x),\quad \forall\, v\in E
\end{align*}
where $z_v^{(n+1)}=w_v^{(n+1)}\cdot z_v^{(n)}$, the following estimates hold:
\begin{enumerate}

  \item\label{for:324} $d(h_{n+1},I)_{C^1}\leq \epsilon_{n+1}^{\text{\tiny$\frac{1}{2}$}}$;

  \smallskip

  \item\label{for:386} $\textbf{d}(w^{(n+1)},e)=\max_{v\in E}\textbf{d}(w^{(n+1)}_v,e)\leq \epsilon_{n+1}^{\text{\tiny$\frac{1}{2}$}}$

  \smallskip

  \item\label{for:325} $\norm{\mathfrak{p}^{(n+1)}}_{C^{\ell}}\leq \epsilon^{-\gamma_1}_{n+1}$;

  \smallskip

  \item\label{for:328} $\norm{\mathfrak{p}^{(n+1)}}_{C^0}\leq \epsilon_{n+1}$;

  \smallskip

  \item\label{for:327} for any $q\leq n+1$ we have
  \begin{align*}
   &\max\{\norm{\mathfrak{p}^{(n+1)}}_{C^{\beta_q}},\,\norm{h_{n+1}-I}_{C^{\beta_q}}\}\\
   &<r_q^{n+1}\epsilon^{\text{\tiny$-8\gamma$}}_{n+1}(\norm{\mathfrak{p}^{(q-1)}}_{C^{\beta_q}}+1)
  \end{align*}
  where $r_q$ is a constant dependent only on $q$;

  \smallskip

  \item \label{for:329} for any $q\leq n+1$ we have
  \begin{align*}
 d(h_{n+1},I)_{C^{\frac{\beta_q}{8}}}\leq C_q(\norm{\mathfrak{p}^{(q-1)}}_{C^{\beta_q}}+1)^{\frac{1}{8}}r_q^{\frac{(n+1)}{8}}
 \epsilon_{n+1}^{\frac{3.5-\gamma}{8}}.
\end{align*}

\end{enumerate}
\end{proposition}
\subsubsection{Proof of Proposition \ref{po:8}}
By interpolation inequalities we have: for any $0\leq m\leq \ell$
\begin{align}\label{for:183}
 \norm{\mathfrak{p}^{(n)}}_{C^m}&\leq C_{\ell} (\norm{\mathfrak{p}^{(n)}}_{C^0})^{\text{\small$\frac{\ell-m}{\ell}$}}(\norm{\mathfrak{p}^{(n)}}_{C^{\ell}})^{\text{\small$\frac{m}{\ell}$}}\notag\\
 &\leq C_{\ell}\epsilon_n^{1-\text{\small$\frac{(1+\gamma_1)m}{\ell}$}}.
\end{align}
By \eqref{for:183} we have
\begin{align}\label{for:76}
 \norm{\mathfrak{p}^{(n)}}_{C^{\varrho+1}}\leq C_{\ell}\epsilon_n^{1-\text{\small$\frac{(1+\gamma_1)(\varrho+1)}{\ell}$}}\overset{(1)}{<} C_\ell\epsilon_n^{\text{\tiny$\frac{1}{2}$}+\delta}<\epsilon_n^{\text{\tiny$\frac{1}{2}$}}<\bar{c},
\end{align}
Here in $(1)$ we use \eqref{for:220}.

We note that
\begin{align*}
z_v^{(n)}=w_v^{(n)}\cdot w_v^{(n-1)}\cdots w_v^{(1)}\cdot w_v^{(0)}\cdot \alpha_v.
\end{align*}
Then we have
\begin{align}\label{for:385}
 &\textbf{\emph{d}}(w_v^{(n)}\cdot w_v^{(n-1)}\cdots w_v^{(1)}\cdot w_v^{(0)},e)\notag\\
 &\overset{(1)}{\leq}  \sum_{j=0}^n \textbf{\emph{d}}(w_v^{(j)},e)\leq \sum_{j=0}^n\epsilon_{j}^{\text{\tiny$\frac{1}{2}$}}\leq C\epsilon^{\text{\tiny$\frac{1}{2}$}}<\epsilon^{\text{\tiny$\frac{1}{4}$}}<\bar{c}.
\end{align}
Here in $(1)$ we use right invariance of $\textbf{\emph{d}}$ (see \eqref{sect:9} of Section \ref{sec:19}).

\eqref{for:76} and \eqref{for:385} allow us to apply Proposition \ref{po:7} to obtain  the new iterates $\mathfrak{p}^{(n+1)}$, $h_{n+1}$ and $w_v^{(n+1)}$:
there are $h_{n+1}\in C^\infty(\mathcal{X})$, $\mathfrak{p}_v^{(n+1)}\in C^\infty(\mathcal{X},\,\mathcal{V})$ and $w_v^{(n+1)}\in V^o$ for each $v\in E$
such that
\begin{align*}
 \tilde{\alpha}^{(n+1)}_v(x)=h_{n+1}^{-1}\circ \tilde{\alpha}^{(n)}\circ h_{n+1}=\exp(\mathfrak{p}^{(n+1)}_v(x))\cdot
z_v^{(n+1)}(x)
\end{align*}
for any $v\in E$ and $x\in \mathcal{X}$, where $z_v^{(n+1)}=w_v^{(n+1)}\cdot z_v^{(n)}$.

Let $b=\epsilon^{-\text{\small$\frac{2\nu}{\ell}$}}_n$. Set
\begin{align}\label{for:214}
 \zeta_q=\big(C_{\beta_q}\norm{\mathfrak{p}^{(n)}}_{C^{\beta_q}}\norm{\mathfrak{p}^{(n)}}_{C^\varrho}^{-1}\epsilon^{\text{\tiny$-8(\gamma-x_0)$}}_n\big)^{\frac{1}{\beta_q}},\quad 1\leq q\leq n
\end{align}
 and $\zeta=\min_{1\leq q\leq n}\{\zeta_q\}$. It is clear that both $b$ and $\zeta$ are dependent on $n$.

The following two lemmas establish Proposition \ref{po:8} for the cases of $\zeta>b^{\frac{1}{2}}$ and  $\zeta\leq b^{\frac{1}{2}}$, respectively.

\begin{lemma} Proposition \ref{po:8} holds when $\zeta>b^{\frac{1}{2}}$.
\end{lemma}
\begin{proof}Set $\tau=b^{\frac{1}{2}}=\epsilon^{-\text{\small$\frac{\nu}{\ell}$}}_n$ and $m=\ell$.

\eqref{for:324}: By \eqref{for:211} of Proposition \ref{po:7}, we have
\begin{align}\label{for:219}
 d(h_{n+1},I)_{C^{1}}&\leq C b^{\text{\tiny$\varrho$}}\tau^{\text{\tiny$\varrho$}}\norm{\mathfrak{p}^{(n)}}_{C^{\varrho}}\overset{(1)}{\leq} C_\ell\epsilon^{-\text{\small$\frac{2\nu\varrho}{\ell}$}}_n\epsilon^{-\text{\small$\frac{\nu\varrho}{\ell}$}}_n
 \epsilon_n^{1-\text{\small$\frac{(1+\gamma_1)\varrho}{\ell}$}}\notag\\
 &\overset{(2)}{\leq} C_\ell(\epsilon_{n})^{\text{\tiny$\frac{\gamma}{2}+\delta$}}
 <(\epsilon_n^{\text{\tiny$\gamma$}})^{\text{\tiny$\frac{1}{2}$}}=\epsilon_{n+1}^{\text{\tiny$\frac{1}{2}$}}.
\end{align}
Here in $(1)$ we use \eqref{for:183} to estimate $\norm{\mathfrak{p}^{(n)}}_{C^{\varrho}}$; in $(2)$ we use \eqref{for:220}.

\smallskip

\eqref{for:386}: By \eqref{for:209} of Proposition \ref{po:7}, we have
\begin{align}\label{for:387}
 \textbf{\emph{d}}(w^{(n+1)},e)\leq C\norm{\mathfrak{p}}_{C^0}\leq C\epsilon_n\overset{(1)}{\leq} C(\epsilon_n^{\text{\tiny$\gamma$}+\delta})^{\text{\tiny$\frac{1}{2}$}}<(\epsilon_n^{\text{\tiny$\gamma$}})^{\text{\tiny$\frac{1}{2}$}}=\epsilon_{n+1}^{\text{\tiny$\frac{1}{2}$}}
\end{align}
Here in $(1)$ we use the assumption for $\gamma$ and $\delta$ (see \eqref{for:326}).

\smallskip

\eqref{for:325}: By \eqref{for:209} of Proposition \ref{po:7}, we have
\begin{align}\label{for:115}
 \norm{\mathfrak{p}^{(n+1)}}_{C^{\ell}}&\leq C_{\ell}(\tau^{\ell+2\varrho}\norm{\mathfrak{p}^{(n)}}_{C^\varrho}+\tau^{2\varrho}\norm{\mathfrak{p}^{(n)}}_{C^{\ell}}+1)\notag\\
 &\overset{(1)}{\leq} C_{\ell}(\epsilon^{-\text{\small$\frac{\nu(\ell+2\varrho)}{\ell}$}}_n
 \epsilon_n^{1-\text{\small$\frac{(1+\gamma_1)\varrho}{\ell}$}}+
\epsilon^{-\text{\small$\frac{2\nu\varrho}{\ell}$}}_n\epsilon^{-\gamma_1}_n+1)\notag\\
&\overset{(2)}{<} C_{\ell}(2
\epsilon^{-\text{\small$\frac{(2\nu+1+\gamma_1)\varrho}{\ell}$}}_n\epsilon^{1-\nu}_n+1)\notag\\
&\overset{(3)}{<}C_\ell(2\epsilon_n^{\text{\tiny$-\gamma\gamma_1$}+\delta}+1)\notag\\
&\overset{(4)}{<}3C_\ell\epsilon_n^{\text{\tiny$-\gamma\gamma_1$}+\delta}<\epsilon_n^{\text{\tiny$-\gamma\gamma_1$}}=(\epsilon_{n+1})^{-\gamma_1}.
\end{align}
Here in $(1)$ we use \eqref{for:183} to estimate $\norm{\mathfrak{p}^{(n)}}_{C^{\varrho}}$; in $(2)$ we note that $\epsilon^{-\gamma_1}_n<\epsilon^{1-\nu}_n$, which is implied by the assumption $\gamma_1+8\gamma-7<\nu$; in $(3)$ we use \eqref{for:221}; in $(4)$ from \eqref{for:326} we see that $\epsilon_n^{\text{\tiny$-\gamma\gamma_1$}+\delta}>1$.

\smallskip

\eqref{for:328}: By \eqref{for:197} of Proposition \ref{po:7}, we have
\begin{align}
\norm{\mathfrak{p}^{(n+1)}}&_{C^0}\leq  Cb^{2\varrho}\tau^{2\varrho}
(\norm{\mathfrak{p}^{(n)}}_{{C^{\varrho+1}}})^{1+y_0}\notag\\
&+C_{\ell} b^{-\ell+\varrho}(\tau^{\ell}\norm{\mathfrak{p}^{(n)}}_{C^\varrho}+\norm{\mathfrak{p}^{(n)}}_{C^{\ell}})
 \notag\\
&+C_\ell\tau^{\varrho}
\big(\tau^{-\ell}\norm{\mathfrak{p}^{(n)}}_{C^\ell}\big)^{\frac{1}{8}}(\norm{\mathfrak{p}^{(n)}}_{{C^\varrho}})^{\frac{7}{8}}\notag\\
 &\overset{(1)}{\leq}  C\epsilon^{-\text{\small$\frac{4\nu}{\ell}$}}_n\epsilon^{-\text{\small$\frac{2\nu}{\ell}$}}_n(\epsilon_n^{1-\text{\small$\frac{(1+\gamma_1)
(\varrho+1)}{\ell}$}})^{1+y_0}\notag\\
&+C_\ell \epsilon^{\text{\small$\frac{2\nu(\ell-\varrho)}{\ell}$}}_n\epsilon^{-\nu}_n\epsilon_n^{1-\text{\small$\frac{(1+\gamma_1)\varrho}{\ell}$}}+C_\ell \epsilon^{\text{\small$\frac{2\nu(\ell-\varrho)}{\ell}$}}_n\epsilon^{-\gamma_1}_n\notag\\
&+C_\ell\epsilon^{-\text{\small$\frac{\nu}{\ell}$}}_n(\epsilon^{\nu}_n\epsilon^{-\gamma_1}_n)^{\frac{1}{8}}
 (\epsilon_n^{1-\text{\small$\frac{(1+\gamma_1)\varrho}{\ell}$}})^{\frac{7}{8}}\notag\\
  &\overset{(2)}{\leq} 4C_\ell\epsilon_n^{\gamma+\delta}<\epsilon_n^{\text{\tiny$\gamma$}}=\epsilon_{n+1}. \label{for:256}
\end{align}
Here in $(1)$ we use \eqref{for:183} to estimate $\norm{\mathfrak{p}^{(n)}}_{C^{\varrho}}$ and $\norm{\mathfrak{p}^{(n)}}_{C^{\varrho+1}}$; in $(2)$ we use \eqref{for:222}, \eqref{for:224} and \eqref{for:223}.

\smallskip

\eqref{for:327}: By \eqref{for:211} and \eqref{for:209} of Proposition \ref{po:7} we have
\begin{align}\label{for:216}
 &\max\{\norm{\mathfrak{p}^{(n+1)}}_{C^{\beta_m}},\,d(h_{n+1},I)_{C^{\beta_m}}\}\notag\\
 &\leq C_m(\tau^{\beta_m}b^\varrho\norm{\mathfrak{p}^{(n)}}_{C^{\varrho}}+b^\varrho\norm{\mathfrak{p}^{(n)}}_{C^{\beta_m}}+1).
 \end{align}
 For $q=n+1$, let $r_{n+1}=2C_{n+1}\tau^{\beta_{n+1}}b^\varrho$, then
 \begin{align}\label{for:225}
 &\max\{\norm{\mathfrak{p}^{(n+1)}}_{C^{\beta_{n+1}}},\,d(h_{n+1},I)_{C^{\beta_{n+1}}}\}\leq r_{n+1}(\norm{\mathfrak{p}^{(n)}}_{C^{\beta_{n+1}}}+1).
 \end{align}
Next, we consider the case $1\leq q\leq n$. By using \eqref{for:214}, we have
\begin{align}\label{for:521}
\zeta_q^{\beta_q}\norm{\mathfrak{p}^{(n)}}_{C^\varrho}=C_{\beta_q}\norm{\mathfrak{p}^{(n)}}_{C^{\beta_q}}\epsilon^{\text{\tiny$-8(\gamma-x_0)$}}_n.
\end{align}
We recall that $\zeta_q\geq \zeta>\tau$. Then
\begin{align}\label{for:522}
\tau^{\beta_q}<\zeta_q^{\beta_q}.
\end{align}
Using \eqref{for:216} we have
 \begin{align}
 &\max\{\norm{\mathfrak{p}^{(n+1)}}_{C^{\beta_q}},\,d(h_{n+1},I)_{C^{\beta_q}}\}\notag\\
 &\overset{(1)}{\leq} C_{q}
(\zeta_q^{\beta_q}\epsilon^{-\text{\small$\frac{2\nu\varrho}{\ell}$}}_n\norm{\mathfrak{p}^{(n)}}_{C^{\varrho}}+1)+C_{q}
\epsilon^{-\text{\small$\frac{2\nu\varrho}{\ell}$}}_n\norm{\mathfrak{p}^{(n)}}_{C^{\beta_q}}\label{for:226}\\
&\overset{(2)}{\leq} C_{q}
(\epsilon^{-\text{\small$\frac{2\nu\varrho}{\ell}$}}_nC_{\beta_q}\norm{\mathfrak{p}^{(n)}}_{C^{\beta_q}}\epsilon^{\text{\tiny$-8(\gamma-x_0)$}}_n+1)+C_{q}
\epsilon^{-\text{\small$\frac{2\nu\varrho}{\ell}$}}_n\norm{\mathfrak{p}^{(n)}}_{C^{\beta_q}}\notag\\
&\leq 2C_{q}\epsilon^{-\text{\small$\frac{2\nu\varrho}{\ell}$}}_nC_{\beta_q}\epsilon^{\text{\tiny$-8(\gamma-x_0)$}}_n (\norm{\mathfrak{p}^{(n)}}_{C^{\beta_q}}+1)\notag\\
&\overset{(3)}{\leq} 4C_{q}C_{\beta_q}\epsilon^{-\text{\small$\frac{2\nu\varrho}{\ell}$}}_n\epsilon^{\text{\tiny$-8(\gamma-x_0)$}}_n (r_q)^n\epsilon^{\text{\tiny$-8\gamma$}}_n(\norm{\mathfrak{p}^{(q-1)}}_{C^{\beta_q}}+1)\notag\\
&\overset{(4)}{\leq} 4C_{q}C_{\beta_q}(r_q)^n\epsilon^{\text{\tiny$-8\gamma$}}_{n+1}(\norm{\mathfrak{p}^{(q-1)}}_{C^{\beta_q}}+1)\notag\\
&=(r_q)^{n+1}\epsilon^{\text{\tiny$-8\gamma$}}_{n+1}(\norm{\mathfrak{p}^{(q-1)}}_{C^{\beta_q}}+1).\label{for:192}
\end{align}
Here in $(1)$ we use \eqref{for:522} to replace $\tau^{\beta_q}$; in $(2)$ we use \eqref{for:521} to replace $\zeta_q^{\beta_q}\norm{\mathfrak{p}^{(n)}}_{C^\varrho}$; in $(3)$ we use the induction assumption (see \eqref{for:323}); in $(4)$ we use \eqref{for:213}.

\eqref{for:225} and \eqref{for:192} imply the result.

\smallskip

\eqref{for:329}: By interpolation inequalities we have
\begin{align}\label{for:218}
 d(h_{n+1},I)_{C^{\frac{\beta_q}{8}}}&\leq C_q\big(d(h_{n+1},I)_{C^{0}}\big)^{\frac{7}{8}}\big(d(h_{n+1},I)_{C^{\beta_q}}\big)^{\frac{1}{8}}\notag\\
 &\overset{(1)}{\leq}C_q(\epsilon_{n+1}^{\text{\tiny$\frac{1}{2}$}})^{\frac{7}{8}}
 \big((r_q)^{n+1}\epsilon^{\text{\tiny$-8\gamma$}}_{n+1}(\norm{\mathfrak{p}^{(q-1)}}_{C^{\beta_q}}+1)\big)^{\frac{1}{8}}\notag\\
 &\leq C_q(\norm{\mathfrak{p}^{(q-1)}}_{C^{\beta_q}}+1)^{\frac{1}{8}}r_q^{\frac{(n+1)}{8}}\epsilon_{n+1}^{\frac{3.5-\gamma}{8}}.
\end{align}
Here in $(1)$ we use \eqref{for:324} and \eqref{for:327}.

\end{proof}
\begin{lemma} Proposition \ref{po:8} holds when $\zeta\leq b^{\frac{1}{2}}$.
\end{lemma}
\begin{proof} Choose $1\leq q\leq n$ such that $\zeta_q=\zeta$.  Set $\tau=\zeta$ and $m=\beta_q$.

\eqref{for:324}, \eqref{for:386} and \eqref{for:325}: Since $\tau=\zeta\leq b^{\frac{1}{2}}= \epsilon^{-\text{\small$\frac{\nu}{\ell}$}}_n$, the estimates for $d(h_{n+1},I)_{C^1}$ (see \eqref{for:219}) and $\norm{\mathfrak{p}^{(n+1)}}_{C^{\ell}}$ (see \eqref{for:115}) still hold.

The estimate for $\textbf{\emph{d}}(w^{(n+1)},e)$ in \eqref{for:387} remains valid, as it is independent of the choice of $b$, $m$ or $\tau$.

\smallskip

\eqref{for:328}: From \eqref{for:197} of Proposition \ref{po:7}, replacing $m$ by $\beta_q$ we see that we need to estimate
$C_{\beta_q}\zeta_q^{-\beta_q}\norm{\mathfrak{p}^{(n)}}_{C^{\beta_q}}$ at first:
\begin{align}\label{for:330}
 C_{\beta_q}\zeta_q^{-\beta_q}&\norm{\mathfrak{p}^{(n)}}_{C^{\beta_q}}\overset{\text{\tiny$(1)$}}{=}\norm{\mathfrak{p}^{(n)}}_{C^\varrho}\epsilon^{\text{\tiny$8(\gamma-x_0)$}}_n
 \overset{\text{\tiny$(2)$}}{\leq} C_{\ell}\epsilon_n^{1-\text{\small$\frac{\varrho(1+\gamma_1)}{\ell}$}}\epsilon^{\text{\tiny$8(\gamma-x_0)$}}_n\notag\\
 &=C_{\ell}\epsilon^{\text{\tiny$8\gamma-8x_0+1-\text{\small$\frac{\varrho(1+\gamma_1)}{\ell}$}$}}_n.
\end{align}
Here in $(1)$ we use  \eqref{for:214} and in $(2)$ we use \eqref{for:183}.

It follows from \eqref{for:197} of Proposition \ref{po:7} that
\begin{align}\label{for:23}
\norm{\mathfrak{p}^{(n+1)}}&_{C^0}\leq  Cb^{2\varrho}\tau^{2\varrho}
(\norm{\mathfrak{p}^{(n)}}_{{C^{\varrho+1}}})^{1+y_0}\notag\\
&+C_{\ell} b^{-\ell+\varrho}(\tau^{\ell}\norm{\mathfrak{p}^{(n)}}_{C^\varrho}+\norm{\mathfrak{p}^{(n)}}_{C^{\ell}})
 \notag\\
&+C\tau^{\varrho}
\big(C_{\beta_q}\zeta_q^{-\beta_q}\norm{\mathfrak{p}^{(n)}}_{C^{\beta_q}}\big)^{\frac{1}{8}}(\norm{\mathfrak{p}^{(n)}}_{{C^\varrho}})^{\frac{7}{8}}\notag\\
&\overset{(1)}{\leq}  C\epsilon^{-\text{\small$\frac{4\nu}{\ell}$}}_n\epsilon^{-\text{\small$\frac{2\nu}{\ell}$}}_n(\epsilon_n^{1-\text{\small$\frac{(1+\gamma_1)
(\varrho+1)}{\ell}$}})^{1+y_0}\notag\\
&+C_\ell \epsilon^{\text{\small$\frac{2\nu(\ell-\varrho)}{\ell}$}}_n\epsilon^{-\nu}_n\epsilon_n^{1-\text{\small$\frac{(1+\gamma_1)\varrho}{\ell}$}}+C_\ell \epsilon^{\text{\small$\frac{2\nu(\ell-\varrho)}{\ell}$}}_n\epsilon^{-\gamma_1}_n\notag\\
&+C_\ell\epsilon^{-\text{\small$\frac{\nu}{\ell}$}}_n(\epsilon^{\text{\tiny$8\gamma-8x_0+1-\text{\small$\frac{\varrho(1+\gamma_1)}{\ell}$}$}}_n)^{\frac{1}{8}}
 (\epsilon_n^{1-\text{\small$\frac{(1+\gamma_1)\varrho}{\ell}$}})^{\frac{7}{8}}\notag\\
  &\overset{(2)}{\leq} 4C_\ell\epsilon_n^{\gamma+\delta}<\epsilon_n^{\text{\tiny$\gamma$}}=\epsilon_{n+1}
\end{align}
Here in $(1)$ we use \eqref{for:183} to estimate $\norm{\mathfrak{p}^{(n)}}_{C^{\varrho}}$ and $\norm{\mathfrak{p}^{(n)}}_{C^{\varrho+1}}$ and we use \eqref{for:330}; in $(2)$ we use \eqref{for:222}, \eqref{for:224} and \eqref{for:254}.

\smallskip
\eqref{for:327}: By \eqref{for:216}, for any $q\leq n+1$ we have
\begin{align*}
 &\max\{\norm{\mathfrak{p}^{(n+1)}}_{C^{\beta_q}},\,d(h_{n+1},I)_{C^{\beta_q}}\}\notag\\
 &\leq C_q(\tau^{\beta_q}b^\varrho\norm{\mathfrak{p}^{(n)}}_{C^{\varrho}}+b^\varrho\norm{\mathfrak{p}^{(n)}}_{C^{\beta_q}}+1)\\
 &\overset{(1)}{\leq}C_q(\zeta_q^{\beta_q}\epsilon^{-\text{\small$\frac{2\nu\varrho}{\ell}$}}_n\norm{\mathfrak{p}^{(n)}}_{C^{\varrho}}
 +\epsilon^{-\text{\small$\frac{2\nu\varrho}{\ell}$}}_n\norm{\mathfrak{p}^{(n)}}_{C^{\beta_q}}+1)\\
 &\overset{(2)}{=} C_{q}
(\zeta_q^{\beta_q}\epsilon^{-\text{\small$\frac{2\nu\varrho}{\ell}$}}_n\norm{\mathfrak{p}^{(n)}}_{C^{\varrho}}+1)+C_{q}
\epsilon^{-\text{\small$\frac{2\nu\varrho}{\ell}$}}_n\norm{\mathfrak{p}^{(n)}}_{C^{\beta_q}}.
 \end{align*}
Here in $(1)$ we use $\zeta_q\geq\zeta=\tau$; in $(2)$ we note that it is exactly \eqref{for:226}.
By the same arguments as in \eqref{for:192}, we still get the result.

 \smallskip
\eqref{for:329}: Once \eqref{for:324} and \eqref{for:327} are proved, the result follows in the same manner as obtaining \eqref{for:218}.


\end{proof}

\subsection{Complement of the proof  of Theorem \ref{th:5}} By  using Proposition \ref{po:8}, we can obtain infinite sequences $h_n$ and $\mathfrak{p}_v^{(n)}$, $v\in E$ inductively. Set
\begin{align*}
 H_n=h_0\circ \cdots \circ h_n\quad\text{and}\quad w_{n,v}=w_v^{(n)}\cdot w_v^{(n-1)}\cdots w_v^{(1)}\cdot w_v^{(0)},\,\,v\in E.
\end{align*}
\eqref{for:386} of Proposition \ref{po:8} shows $w_{n,v}$ converges to an element of $w_v$ in $V^o$ for each $v\in E$.  \eqref{for:324} of Proposition \ref{po:8} shows that $H_n$ converges in $C^1$ topology to a $C^1$ diffeomorphism $h$; moreover, \eqref{for:329} of Proposition \ref{po:8} shows that
 the convergence of the sequence $H_n$ in $C^{\frac{\beta_q}{8}}$ topology for any $q\in\NN$. Hence we see that
 $h$ is of class $C^{\infty}$. The convergence step shows that:
\begin{align*}
  h^{-1}\circ\tilde{\alpha}_v\circ h=w_v\cdot \alpha_{v},\qquad \forall\,v\in E.
\end{align*}
Since $\tilde{\alpha}$ is abelian, $\{w_v\cdot \alpha_{v}:\,v\in E\}$ is an abelian set. As $E$ is a generating set of the group $A$, the map $v\to w_v$, $v\in E$ extends to a map $i_0$ from $A$ to $V^o$ such that
\begin{align*}
  h^{-1}\circ\tilde{\alpha}_a(hx)=i_0(a)\cdot \alpha_a(x),\quad \text{for all }a\in A,\,x\in \mathcal{X}.
\end{align*}
Of course, $A\to i_0(a)\cdot \alpha_a$ is a group homeomorphism from $A$ to an abelian set of affine actions.

Then we complete the proof of Theorem \ref{th:5}.

\section{Proof of Theorem \ref{th:7}}\label{sec:24}

\subsection{Difference between nilmanifold examples and (twisted) symmetric space examples}\label{sec:32}

The first key difference between nilmanifold examples and (twisted) symmetric space examples  lies in the fact that,  for the former, the coordinate change along natural
directions preserves the Lyapunov subspace decomposition (see \eqref{sect:10} of Section \ref{sec:19}). However, for the latter, this property no longer holds, which poses new difficulties. The following example will help illustrate the distinction more clearly.
\begin{example}\label{ex:1} We consider the $\ZZ^2$ action $\alpha$ on $SL(2,\RR)\times SL(2,\RR)\times SL(2,\RR)/\Gamma$ with generators $\textbf{a}_1$ and $\textbf{a}_2$:
$\textbf{a}_1=\begin{pmatrix}
\sqrt{2} & 0 \\
0 & \frac{1}{\sqrt{2}}
\end{pmatrix}\times	\begin{pmatrix}
2 & 0 \\
0 & \frac{1}{2}
\end{pmatrix}\times	I$ and $\textbf{a}_2=\begin{pmatrix}
\sqrt{2} & 0 \\
0 & \frac{1}{\sqrt{2}}
\end{pmatrix}\times \begin{pmatrix}
\frac{1}{2} & 0 \\
0 & 2
\end{pmatrix}\times	I$. We note that $d\textbf{a}_1=\text{Ad}_{\textbf{a}_1}$. The eigenvalues of $\text{Ad}_{\textbf{a}_1}$ are $1$, $2$, $\frac{1}{2}$, $4$ and $\frac{1}{4}$.
Let
\begin{align*}
 z_1=\begin{pmatrix}
\sqrt{2}(1+\epsilon) & 0 \\
0 & \frac{1}{\sqrt{2}(1+\epsilon)}
\end{pmatrix}\times \begin{pmatrix}
2(1+\epsilon) & 0 \\
0 & \frac{1}{2(1+\epsilon)}
\end{pmatrix}\times	\begin{pmatrix}
1+\epsilon & 0 \\
0 & (1+\epsilon)^{-1}
\end{pmatrix}
\end{align*}
where $\epsilon>0$ and $z_2=\textbf{a}_2$. Then $z_1$ and $z_2$ generate a $\ZZ^2$ action $\alpha_\epsilon$ which is
close to $\alpha$ if $\epsilon$ is small. The eigenvalues of $dz_1=\text{Ad}_{z_1}$ are $1$, $2(1+\epsilon)^2$, $4(1+\epsilon)^2$, $\frac{1}{4(1+\epsilon)^2}$, $\frac{1}{2(1+\epsilon)^2}$, $(1+\epsilon)^2$ and $(1+\epsilon)^{-2}$.
It is clear that $d\textbf{a}_1$ and $dz_1$ have different set of Lyapunov exponents.

\end{example}
To study the following twisted equation
\begin{align}\label{for:389}
\varphi\circ z_1-2(1+\epsilon)^2\varphi=\omega,\qquad \omega\in \mathcal{O}^r.
\end{align}
(we can assume that $r$ is sufficiently large if needed).  Let
\begin{align*}
 u_{3,1,2}=\begin{pmatrix}
0 & 0 \\
0 & 0
\end{pmatrix}\times \begin{pmatrix}
0 & 0 \\
0 & 0
\end{pmatrix}\times \begin{pmatrix}
0 & 1 \\
0 & 0
\end{pmatrix}. \quad\text{Then } \text{Ad}_{z_1}u_{3,1,2}=(1+\epsilon)^2u_{3,1,2}.
\end{align*}
 $u$ is a unstable direction of $z_1$.

 We recall \eqref{for:285} of Section \ref{sec:39}. Naturally, we expect that
\begin{align}\label{for:391}
 \mathfrak{D}_{\mu,u_{3,1,2},\tau,\ZZ}(\omega)&\stackrel{\text{def}}{=}\sum_{n} (2(1+\epsilon)^2)^{n-1}\pi_{u_{3,1,2}}(f_1\circ \tau^{-1})\omega\circ z_1^{-n}\notag\\
 &=\sum_{n} (2(1+\epsilon)^2)^{-(n+1)}\pi_{u_{3,1,2}}(f_1\circ \tau^{-1})\omega\circ z_1^{n}
\end{align}
are invariant distributions for solving equation \eqref{for:389}.

We recall the construction of invariant distributions for the twisted equation \eqref{for:291} of Section \ref{sec:41}, which rely on the following ingredients:
\begin{enumerate}

  \item\label{for:390} when $n\geq0$, the exponential decay of the twist $2^{-(n+1)}$ and the uniform upperbound  of the operator $\pi_{u}(f_1\circ \tau^{-1})$, see \eqref{for:289} of Section \ref{sec:41};

  \smallskip
  \item\label{for:392} when $n<0$, the super-exponential decay of matrix coefficients of $z_1$, see \eqref{for:282} of Section \ref{sec:41}. We note that the decay rate $2^{ns}$,  is determined by the deformation of $dz_1^n$ along $u$, which is $2^{n}$,  and the $s$-Sobolev order of $\omega_{\lambda_3}$.
  Consequently, to beat $2^{-(n+1)}$,  the exponential growth of the
twist, one needs $s>1$,  see \eqref{for:290} of Section \ref{sec:41}. More precisely, the order of the invariant distribution is
\begin{align}\label{for:523}
 c(z_1,+)=\frac{\text{the largest eigenvalue of $dz_1|_{\mathfrak{g}^+_{z_1}}$} }{\text{the smallest eigenvalue of $dz_1|_{\mathfrak{g}^+_{z_1}}$}}.
\end{align}
\end{enumerate}
To show the convergence of \eqref{for:391}, we firstly note that the ingredients in \eqref{for:390} still hold.
When $n<0$, since the deformation of $\text{Ad}_{z_1^n}$ along $u_{3,1,2}$ is $(1+\epsilon)^{2n}$, the super-exponential decay rate is
$(1+\epsilon)^{2ns}$. To beat $(2(1+\epsilon)^2)^{-(n+1)}$, the growth of the twist,  one needs $s$ to satisfy
\begin{align*}
  (1+\epsilon)^{2ns}<(2(1+\epsilon)^2)^{n+1},\qquad n<0.
\end{align*}
We see that as $\epsilon\to 0$, $s\to \infty$. This implies that as the deformation becomes flatter, the convergence rate of \eqref{for:391} slows down.
As a consequence, as $\epsilon\to 0$,   the order of the invariant distribution $\mathfrak{D}_{\mu,u_{3,1,2},\tau,\ZZ}(\omega)$, which is $c(z_1,+)$ defined in \eqref{for:523},  goes to $\infty$.

The root cause of the difficulty  lies in $dz_1$ having eigenvalues very close to $1$, which also brings other serious problems. For the $\ZZ^2$ action of $z_1$ and $z_2$
we can also define $\kappa(z_1,z_2)$ (see \eqref{for:139} of Section \ref{sec:19}). As $\kappa(z_1,z_2)$ is not uniformly bounded in any small neighborhood of $\textbf{a}_1$ and $\textbf{a}_2$,  for the study of the twisted equation (see Corollary \ref{cor:3})
\begin{enumerate}
  \item\label{for:531} for any $k\in \NN$,  it is not necessarily that the existence of a solution of Sobolev order $k$ implies the existence of a solution of higher Sobolev order;

  \smallskip
  \item for any $k\in \NN$, it is not necessarily that the $r$-th Sobolev norm of the solution can be uniformly bounded by the $(r+k)$-th Sobolev norm of the (twisted)-coboundary. Uniform means the estimates are independence of the choice of $z_1$ and $z_2$.
\end{enumerate}
Either of the two problems mentioned above results in significant difficulty in obtaining the local rigidity. So as a consequence, even though we follow
the main proof line for the nilmanifold examples, new ingredients are often
required for the (twisted) symmetric space examples.

The second key difference lies in the study
of the cohomological equation
\begin{align}\label{for:420}
 v\varphi=\omega,
\end{align}
where $v$ is inside an eigenspace of $dz_1$, i.e., the problem of obtaining Sobolev estimates of a function $\varphi$ on $\mathcal{X}$ whose Lie derivative $v\varphi$ in
the direction of $v$ is a given function $\omega$ on $\mathcal{X}$.

If $\mathcal{X}$ is a nilmanifold, even though $dz_1$ is Anosov, it is known that we may not be able to compare Sobolev estimates between $\varphi$ and $\omega$ unless
additional assumptions are placed on $v$ (for example, $v$ is Diophantine).

For (twisted) symmetric space examples, Theorem \ref{th:8} shows that  for any $C^\infty$ function $\omega$ of average zero, if  there
exists a $C^\infty$ solution $\varphi$ of equation \eqref{for:420} with zero average, then for any $r\geq0$,  the $r$-th Sobolev norm of
$\varphi$ on a subgroup isomorphic to $SL(2,\RR)$ can be bounded by the $C^{r+2}$ norm of $\omega$. This implies that we can recover the
Sobolev norms of $\varphi$  after differentiating. This fact allows us to balance the twist by directional derivatives, which is the main
novelty for this case.

\subsection{Basic notations}\label{sec:3}

We will use notations from this section throughout subsequent
sections. So the reader should consult this section if an unfamiliar symbol
appears.

In what follows, $C$ will denote any constant that depends only
  on the given  group $\GG$, the manifolds $\mathcal{X}$ and the action $\alpha_A$. $C_{x,y,z,\cdots}$ will denote any constant that in addition to the
above depends also on parameters $x, y, z,\cdots$.

\begin{nsect}\label{for:67}
We recall notations in Section \ref{sec:2}.   $\mathbb{G}=G$ or $G_\rho$ where $G_\rho=G\ltimes_\rho V$.   Let $\mathcal{M}=\mathbb{G}/\Upsilon$ and let $\mathcal{L}=\text{Lie}(L)$.  We use $\mathfrak{G}$ (resp. $\mathfrak{g}$) to denote the Lie algebra of $\GG$ (resp. $G$). Denote by $U(\mathfrak{G})$ the universal enveloping algebra of $\mathfrak{G}$, with its usual filtration $\{U_{m}\}_{m\geq0}$.

We use $dx$ to denote the Harr measure on $\GG/\Gamma$. Set $L^2_0(\GG/\Gamma)$ to the subspace of $L^2(\GG/\Gamma)$ orthogonal to constants.
We use $(\pi, \mathcal{O})$ to denote the regular representation of $L^2_0(\GG/\Gamma)$.

Fix a right invariant metric $\textbf{\emph{d}}$ on $\GG$.  $B_\eta(x)\stackrel{\text{def}}{=}\{y\in
\GG|\textbf{\emph{d}}(x,y)\leq\eta\}$.

\end{nsect}

\begin{nsect}\label{for:52} For any subset $S$ of $\GG$, we use $Z(S)$ to denote its centralizer in $\GG$.
For any subalgebra $\mathfrak{F}_1, \mathfrak{F}_1,\cdots $ of $\mathfrak{G}$, we use $\exp\{\mathfrak{F}_1,\mathfrak{F}_2,\cdots\}$ to denote the subgroup with its Lie algebra generated by $\mathfrak{F}_1, \mathfrak{F}_1,\cdots $.

For any subalgebra $\mathfrak{F}$ of $\mathfrak{G}$, if it has a decomposition $\mathfrak{F}=\mathfrak{F}_1\oplus\mathfrak{F}_2$, where $\mathfrak{F}_i$ is a linear subspace of $\mathfrak{F}$, $i=1,2$, then there exist
bounded, open, connected neighborhoods $U_{\mathfrak{F}_1}$
and $U_{\mathfrak{F}_2}$ of $0$ in $\mathfrak{F}_1$
and $\mathfrak{F}_2$ respectively, such that the mapping:
\begin{gather*}
(U_{\mathfrak{F}_1},U_{\mathfrak{F}_2})\rightarrow
\exp(U_{\mathfrak{F}_1})\exp(U_{\mathfrak{F}_2})
\end{gather*}
is a
diffeomorphism. We call it \emph{the uniqueness decomposition
property}.

\end{nsect}

\begin{nsect}\label{for:17} Suppose $\mathfrak{F}$ is a subspace of $\mathfrak{G}$. $\mathfrak{F}(\mathcal{O}^s)$ denotes the set of maps defined on $\mathcal{M}$
taking values on $\mathfrak{F}$ with coordinate functions in $\mathcal{O}^s$ (see Section \ref{sec:17}).
 For any  map $\mathcal{F}\in \mathfrak{F}(\mathcal{O}^s)$, $\norm{\mathcal{F}}_s$ and $\int_{\mathcal{M}} \mathcal{F}$ are similarly defined as in
\eqref{sect:1} of Section \ref{sec:19}. Then the regular representation $(\pi, \mathcal{O})$ naturally extends to the representation $(\pi, \mathfrak{F}(\mathcal{O}))$.

For any distribution $\mathcal{D}\in \mathfrak{F}(\mathcal{O}^{-s})$, $s>0$, and any
$\mathcal{F}\in \mathfrak{F}(\mathcal{O}^s)$, $\langle \mathcal{D}, \mathcal{F}\rangle\in \mathfrak{F} $ is similarly defined as in
\eqref{sect:1} of Section \ref{sec:19}. Denote by $C^r(\mathcal{M},\,\mathfrak{F})$ the set of $C^r$ maps from $\mathcal{M}$ to $\mathfrak{F}$. $\norm{\mathcal{F}}_{C^s}$ is defined similarly.

For any map $\mathcal{F}$ defined on $\mathcal{M}$ and any $z\in \GG$, denote $\mathcal{F}(z\cdot x)$ by $\mathcal{F}\circ z$.
\end{nsect}

\begin{nsect}\label{for:191} Let $\mathcal{S}$ stands for $\mathcal{M}$ or $\mathcal{X}$.
There exists a collection of smoothing operators $\mathfrak{s}_b:C^\infty(\mathcal{S},\,\mathfrak{F})\to C^\infty(\mathcal{S},\,\mathfrak{F})$, $b>0$, such that for any $s, s_1,s_2\geq0$, the following holds:
\begin{align}
 \norm{\mathfrak{s}_b\mathcal{F}}_{C^{s+s_1}}&\leq C_{s,s_1}b^{s_1}\norm{\mathcal{F}}_{C^{s}},\quad \text{and}\\
 \norm{(I-\mathfrak{s}_b)\mathcal{F}}_{C^{s-s_2}}&\leq C_{s,s_2}b^{-s_2}\norm{\mathcal{F}}_{C^{s}},\quad \text{if }s\geq s_2,\label{for:471}
\end{align}
see \cite{Hamilton}.

The next result follows directly from Sobolev embedding theorem on compact manifolds.  For any $\mathcal{F}\in \mathfrak{F}(\mathcal{O}^\infty)$ and $s\geq0$ the following hold:
\begin{align}\label{for:453}
  \norm{\mathcal{F}}_{s}\leq C_s\norm{\mathcal{F}}_{C^{s}},\quad \norm{\mathcal{F}}_{C^{s}}\leq C_{s}\norm{\mathcal{F}}_{s+\beta},
\end{align}
where $\beta>0$ is a constant only dependent on $\mathcal{M}$.
\end{nsect}

\begin{nsect}\label{sect:11} For any $c\in G$ we have a corresponding Jordan-Chevalley normal form decomposition of $3$ commuting elements $c=s_ck_cn_c$ where $s_c$ is semisimple, $k_c$ is compact and $n_c$ is unipotent. Set $p(c)=s_c$.

For any $c\in G_\rho$ we have the expression $(g_c, v_c)$, where $g_c\in G$ and $v_c\in V$. We also use $s_c$ to denote $s_{g_c}$ and $p(c)$ to denote $p(g_c)$ if there is no confusion.

For any $c\in\GG$ let $y_c=k_cn_c$ if $\GG=G$ or $y_c=(k_{g_c}n_{g_c},v_c)$ if $\GG=G\ltimes_\rho V$.
Using \eqref{for:12}, we have a decomposition:
\begin{align*}
  c=s_cy_c,\qquad \forall\,c\in\GG.
\end{align*}
We say that $c$ is \emph{perfect} if $s_c$ and $y_c$ commute.

For any $c\in\GG$, there is $\delta(c)>0$ such that the following estimates hold:
\begin{align}\label{for:18}
  \norm{\text{Ad}_{(y_{c'})^j}}\leq C_c(\abs{j}+1)^{\dim\mathfrak{G}},\qquad \forall j\in\ZZ
\end{align}
for any $c'\in
  B_{\delta(c)}(c)$ (see Appendix \ref{sec:47}).

\smallskip

\noindent \emph{Note}. If $\GG=G$, then any $c$ is perfect. If $\GG=G\ltimes_\rho V$, generally, $c$ is not perfect.
However, Proposition \ref{po:3} shows that there exists a conjugate $c'$ of $c$ such that $c'$ is perfect and $s_{c'}=s_{c}=p(c)$.

\end{nsect}

\begin{nsect}\label{sec:48} For any $c\in\GG$, under the adjoint representation $\Ad_{c}$ we have the Lyapunov space decomposition
\begin{align}\label{for:436}
\mathfrak{G}=\sum_{\chi\in\Delta_c}\mathfrak{g}_{\chi,c}
\end{align}
such that for any $\chi\in\Delta_c$, $\chi(p(c))$ is the absolute value of the eigenvalue of $\Ad_{c}$ on $\mathfrak{g}_{\chi,c}$ (see Appendix \ref{ob:5}).

\eqref{for:18} of \eqref{sect:11} shows that if $c$ is perfect, then  the Lyapunov exponents, as well as the Lyapunov subspaces of $\Ad_{c}$   match  those of $\Ad_{p(c)}$.

Set
\begin{align}\label{for:273}
 \mathfrak{g}^+_c=\sum_{\stackrel{\chi\in \Delta_c,}
 {\chi(p(c))>1}}\mathfrak{g}_{\chi,c},\quad \mathfrak{g}^-_c=\sum_{\stackrel{\chi\in \Delta_c,}
 {\chi(p(c))<1}}\mathfrak{g}_{\chi,c},\quad \mathfrak{g}^0_c=\sum_{\stackrel{\chi\in \Delta_c,}
 {\chi(p(c))=1}}\mathfrak{g}_{\chi,c}.
\end{align}
\end{nsect}

\begin{nsect}\label{for:62}
Since $A$ is abelian, there is a $\RR$-split Cartan subgroup $A_0$ of $G$ such that $p(a)\in A_0$ for any $a\in A$ and
$p: A\to A_0$ is a group homeomorphism (see Proposition \ref{po:3}). Under the adjoint representation $\Ad|_{p(A)}$ we have the eigenspace decomposition
\begin{align}\label{for:6}
\mathfrak{G}=\sum_{\mu\in\Delta_A}\mathfrak{g}_{\mu,A}
\end{align}
such that for any $a\in A$ and $\mu\in\Delta_A$, $\mu(p(a))$ is the eigenvalue (resp. absolute value of the eigenvalue) of $\Ad_{p(a)}$ (resp. $\Ad_{a}$) on $\mathfrak{g}_\mu$. \textbf{\emph{We also write $\mu(a)$ for simplicity}}. Elements of
$A\setminus p^{-1}\left(\bigcup_{\phi\in\Delta_A\backslash  1}\ker(\phi)\right)$ are
\emph{regular} elements (for $A$).

It is clear that the following relations hold:
\begin{align}\label{for:475}
 [\mathfrak{g}_{\mu,A},\,\mathfrak{g}_{\lambda,A}]\subseteq \mathfrak{g}_{\mu\lambda,A},\qquad \forall\,\mu,\,\lambda\in \Delta_A.
\end{align}
We choose two regular elements $\textbf{a}_1$ and $\textbf{a}_2$ of $A$, which are as described in Proposition \ref{po:9}. Set
\begin{align*}
 \Delta_{\textbf{a}_1}^{+}=\{\mu\in\Delta_{A}:\mu(\textbf{a}_1)>1\}\quad\text{and}\quad\Delta_{\textbf{a}_1}^{-}=\{\mu\in\Delta_{A}:\mu(\textbf{a}_1)<1\}.
\end{align*}
We fix a compact generating set $E\subset A$ such that each element of $E$ is regular and $E$ contains $\textbf{a}_1$ and $\textbf{a}_2$.

\smallskip

\noindent\emph{Note}. If $A=A_0$ and $\GG=G$, then $\Delta_{A_0}$ is the set of roots of $A_0$ and $\mathfrak{g}_{\mu,A_0}$ is the root space of $\mu\in \Delta_{A_0}$. Since $p(A)\subset A_0$, we see that for any $\mu\in \Delta_{A}$,  $\mathfrak{g}_{\mu,A}$ is a union of $\mathfrak{g}_{\phi,A_0}$, $\phi\in \Delta_{A_0}$. More precisely, there exists a subset $s_\mu\subset \Delta_{A_0}$  such that
\begin{align*}
  \mathfrak{g}_{\mu,A}=\bigcup_{\lambda\in s_\mu}\mathfrak{g}_{\lambda,A_0},\qquad \forall\, \mu\in \Delta_{A}.
\end{align*}
Since $\textbf{a}_1$ and $\textbf{a}_2$ are both regular, we point out that
\begin{align}\label{for:11}
 \mathfrak{g}_{1,A}=\mathfrak{g}_{\textbf{a}_i}^{0}=Z(p(A)),\qquad i=1,\,2.
\end{align}
\eqref{for:475} shows that $\mathfrak{g}_{1,A}$ is a subalgebra. Then $\exp(\mathfrak{g}_{1,A})$ is a subgroup with its Lie algebra $\mathfrak{g}_{1,A}$.  For any $w\in\exp(\mathfrak{g}_{1,A})$, \eqref{for:475} shows that each $\mathfrak{g}_{\mu,A}$, $\mu\in \Delta_A$
in invariant under $\text{Ad}_{w}$. Consequently, we have
\begin{align}\label{for:490}
 \text{Ad}_{wa}(\mathfrak{g}_{\mu,A})\subseteq \mathfrak{g}_{\mu,A},\qquad \forall\,a\in A,\,\,\forall\,\mu\in \Delta_A.
\end{align}
Passing to a conjugation will not affect the local rigidity results in the paper. Then by Proposition \ref{po:3},  it is harmless to assume that $A$ satisfies the following property: \begin{enumerate}
  \item []\hypertarget{o.1}{(P)} For any $a,\,c\in A$,  $p(a)$ commutes with $y_c$ (see \eqref{sect:11}).
\end{enumerate}
\noindent\emph{Note}. Property  \hyperlink{o.1}{(P)} shows that any element in $A$ is perfect. This means  that \eqref{for:6} is also
the Lyapunov space decomposition of $\Ad_{A}$.

\smallskip

\end{nsect}

\begin{nsect}\label{for:54} We note that $\mathcal{L}\subseteq \mathfrak{g}_{1,A}$; and  $\mathfrak{g}_{\mu,A}$, $\mu\in\Delta_A$ are invariant under $\text{Ad}_L$.
Then there is a  positive definite $\text{Ad}_L$-invariant inner product $\langle\cdot\rangle$ on $\mathfrak{G}$ such that $\mathfrak{g}_{\mu,A}\subseteq \mathcal{L}^\bot$, for all $1\neq\mu\in\Delta_A$, where  $\mathcal{L}^\bot$ is the
orthogonal complement of $\mathcal{L}$.

A map $F:\mathcal{M}\to \mathcal{L}^\bot$ is called \emph{adjoint $L$-invariant} if
\begin{align}\label{for:117}
 \text{Ad}_{l}F\circ l^{-1}=F,\qquad \forall\,\,l\in L.
\end{align}
Set $J_\mu=\mathfrak{g}_{\mu,A}$, if $1\neq\mu\in\Delta_A$; and $J_1=\mathcal{L}^\bot\bigcap\mathfrak{g}_{1,A}$.
 Denote by $\mathfrak{i}_{\mathcal{L}^\bot}$ the projection from $\mathfrak{G}$ to $\mathcal{L}^\bot$.   For any $q\in \exp(\mathfrak{g}_{1,A})\cap Z(L)$, define
\begin{align*}
 \text{Ad}_{q}\nmid_{J_1}:=\mathfrak{i}_{\mathcal{L}^\bot}\circ\text{Ad}_{q}|_{J_1}
\end{align*}
We note that
$\mathcal{L}^\bot$ is not necessarily invariant under $\text{Ad}_{Z(L)}$. Then $J_1$ is not necessarily an invariant subspace for $\text{Ad}_{q}$, but it is invariant for
$\text{Ad}_{q}\nmid$. Moreover, for any $q_1,\,q_2\in \exp(\mathfrak{g}_{1,A})\cap Z(L)$, we have
\begin{align}\label{for:53}
(\text{Ad}_{q_1}\nmid_{J_1})(\text{Ad}_{q_2}\nmid_{J_1})=\text{Ad}_{q_1q_2}\nmid_{J_1}.
\end{align}
We will denote $\text{Ad}_{q}\nmid_{J_1}$ by $\text{Ad}_{q}|_{J_1}$ without further specification.

\end{nsect}

\begin{nsect}\label{for:35} Let $\mathfrak{G}^1$ be the set of unit vectors in $\mathfrak{G}$ with respect to the inner product defined in \eqref{for:54}.  There are $0<\delta<1$ such that the following are satisfied:
\begin{enumerate}
  \item for any $v\in \mathfrak{G}$,  if $\norm{v}\leq \delta$, then $\textbf{\emph{d}}(\exp v,e)\leq C\norm{v}$;

  \smallskip
  \item for any $g\in B_{\delta}(e)$, we have
\begin{align*}
 \norm{\log g}+\norm{\text{Ad}_{g}-I}\leq C\textbf{\emph{d}}(g,e).
\end{align*}
\end{enumerate}

\end{nsect}

\begin{nsect}\label{sec:49}
For any $c\in\GG$, there is $\delta(c)>0$ such that
\begin{align}\label{for:260}
  \norm{\text{Ad}_{(c')^j}|_{\mathfrak{g}_{\chi,c'}}}\leq C_c\chi(p(c'))^j(\abs{j}+1)^{\dim\mathfrak{G}}
\end{align}
for any $c'\in
  B_{\delta(c)}(c)$,  any $j\in\ZZ$, any $\chi\in \Delta_{c'}$ (see \eqref{sec:48}). The proof of \eqref{for:260} is left for Appendix \ref{cor:9}.
\end{nsect}

\begin{nsect}\label{for:1}
For any $\mu\in \Delta_A$ let $L_\mu=\ZZ^-\cup\{0\}$ (resp. $L_\mu=\NN$) and $\delta(\mu)=-$ (resp. $\delta(\mu)=+$) if $\mu\in \Delta_{\textbf{a}_1}^{+}\bigcup\{1\}$ (resp. $\mu\in \Delta_{\textbf{a}_1}^{-}$).

We say that $(\lambda_1,\cdots,\lambda_j)$ is a  \emph{poor row} if either $\lambda_1\in \Delta_{\textbf{a}_1}^{+}\bigcup\{1\}$ and $\lambda_l\in \Delta_{\textbf{a}_1}^{+}$, $2\leq l\leq j$; or
$\lambda_l\in \Delta_{\textbf{a}_1}^{-}$, $1\leq l\leq j$.

For any $\lambda\in \Delta_A$ we define a subset $\mathcal{U}_{\lambda,k}$  of $U_k$ (see \eqref{for:67}):
\begin{align*}
 \mathcal{U}_{\lambda,k}=\{u_1\cdots u_k:u_i\in \mathfrak{g}_{\lambda_i,A}\cap\mathfrak{G}^1,\,\,\,(\lambda,\lambda_1,\cdots, \lambda_k)\text{ is a poor row}\}.
\end{align*}
\end{nsect}

\subsection{Constants in Section \ref{sec:24}} \label{for:432} In this part, we list all the constants that will be used:
\begin{enumerate}
  \item Regarding $\textbf{a}_1$ and $\textbf{a}_2$, we define the following positive numbers:
   \begin{enumerate}
     \item $\emph{\textbf{l}}(\textbf{a}_1,\textbf{a}_2)$, $\emph{\textbf{l}}_1(\textbf{a}_1,\textbf{a}_2)=\emph{\textbf{l}}(\textbf{a}_1,\textbf{a}_2)^2$ (see Proposition \ref{po:9});

         \smallskip
     \item\label{for:488} $\eta(\textbf{a}_1,\textbf{a}_2)$ (see Proposition \ref{po:9}) and   $\eta_m(\textbf{a}_1)$ (see \eqref{for:426} of Corollary \ref{for:129});

     \smallskip
     \item We choose $\sigma(\textbf{a}_1,\textbf{a}_2)\geq4$ such that
\begin{align}
 \sigma(\textbf{a}_1,\textbf{a}_2)\emph{\textbf{l}}(\textbf{a}_1,\textbf{a}_2)&\geq\max_{\mu\in \Delta_A}\{1,\,|\log\mu(p(\textbf{a}_1))|,\,|\log\mu(p(\textbf{a}_2))|\}.\label{for:27}
\end{align}
We also write $\emph{\textbf{l}}$, $\emph{\textbf{l}}_1$, $\eta$, $\eta_m$ and $\sigma$ for simplicity if there is no confusion.
   \end{enumerate}

\smallskip
  \item $\beta$ is defined in \eqref{for:453} of \eqref{for:191}.

  \smallskip
  \item $0<\delta<1$ is defined in \eqref{for:35}.

  \smallskip

  \item $\delta_m$, $m\geq1$ is from Lemma \ref{le:8}.

  \smallskip

  \item\label{for:470} $0<c_0\leq\frac{4}{7}$ is defined in Proposition \ref{po:1}.

  \smallskip

  \item\label{for:236} $\varrho$ is defined in Section \ref{sec:56}.

\smallskip

\item $\bar{c}$ is from Proposition \ref{po:5}.

\smallskip

\item We choose $\delta_0$ such that $\delta_0<\eta_{\sigma+6+\beta}(\textbf{a}_1)$ (see \eqref{for:488}). Note that $\eta_{\sigma+6+\beta}(\textbf{a}_1)<\eta(\textbf{a}_1, \textbf{a}_2)$ (see \eqref{for:488}).
\end{enumerate}

\subsection{Useful results}\label{sec:50} In this part, we list some basic facts about the abelian group $A\subseteq \GG$. These facts will be frequently used later.

\subsubsection{Algebraic property of $A$} The following result shows that the Lyapunov exponents of $\alpha_A$ are determined  by the semisimple part of $A$.
\begin{proposition}\label{po:3}
There exists a split Cartan subgroup $A_0$ and  $u\in V$ such that:
\begin{enumerate}
  \item\label{for:393} $p(A)\subseteq A_0$ and $p: A\to A_0$ is a group homeomorphism (see \eqref{sect:11} of Section \ref{sec:3});

  \smallskip
  \item\label{for:394} $uau^{-1}\in Z(p(A))$, for any $a\in A$ (see \eqref{for:52} of Section \ref{sec:3}).
\end{enumerate}
\end{proposition}
\emph{Note}. Proposition \ref{po:3} is valid for any abelian group $A$.

\smallskip
\noindent The proof of Proposition \ref{po:3}   is left for Appendix \ref{sec:9}.

\medskip

\subsubsection{Asymptotic behaviours in small neighborhoods of $\textbf{a}_1$ and $\textbf{a}_2$}
In this part, we assume that $A$ satisfies property  \hyperlink{o.1}{(P)}.  We obtain uniform upperbounds of asymptotic behaviours in small neighborhoods of $\textbf{a}_1$ and $\textbf{a}_2$.

\begin{proposition}\label{po:9} We can choose two regular elements $\textbf{a}_1$ and $\textbf{a}_2$ of $A$ and find positive constants $\eta(\textbf{a}_1,\textbf{a}_2)$ and $\textbf{l}(\textbf{a}_1,\textbf{a}_2)$ such that the following are satisfied:
\begin{enumerate}
  \item\label{for:402} $p(\textbf{a}_1)$ and $p(\textbf{a}_2)$ are linearly independent;

  \smallskip

  \item\label{for:409} $\textbf{l}<\frac{1}{5}\min_{1\neq\mu\in \Delta_A}\{\big|\log \mu(p(\textbf{a}_1))\big|,\,\big|\log \mu(p(\textbf{a}_2))\big|\}$;

  \smallskip

  \item\label{for:50} for any $w_v\in B_{2\eta(\textbf{a}_1,\textbf{a}_2)}(e)$, $v\in E$ (see the end of \eqref{for:62} of Section \ref{sec:3}), we have
  \begin{align*}
    \norm{\log[w_{v_1}v_1, w_{v_2}v_2]}\leq \delta ^2,\qquad \forall\,v_1,v_2\in E
  \end{align*}
  ($\delta$ is defined in \eqref{for:35} of Section \ref{sec:3});

  \smallskip

  \item for
 any $w_i\in \exp(\mathfrak{g}_{1,A})\bigcap B_{\eta}(e)$, $i=1,2$, letting $z_i=w_i\textbf{a}_i$, $i=1,2$ then:
 \begin{enumerate}
   \item \label{for:403} for any $\xi_i\in\mathcal{O}^1$, $i=1,2$, we have
   \begin{align*}
\big|\langle \xi_1\circ (z_1^nz_2^m),\xi_2\rangle\big|&\leq
E\norm{\xi_1}_{1}\norm{\xi_2}_{1}e^{-4
\textbf{l}(|n|+|m|)}
\end{align*}
 for any $(m,n)\in\ZZ^2$;

\smallskip
   \item\label{for:404} $\mathfrak{g}_{\mu,A}$ are invariant subspaces for $\text{Ad}_{z_i}$, $i=1,2$ for any $\mu\in\Delta_A$. Moreover,  we have
   \begin{align}\label{for:447}
  \norm{\text{Ad}_{z_i^n}|_{\mathfrak{g}_{\mu,A}}}\leq C\mu(\textbf{a}_i)^ne^{\textbf{l}_1|n|}\leq C\mu(\textbf{a}_i)^ne^{\textbf{l}|n|},\quad \forall\,\,n\in\ZZ,\,\,i=1,2;
\end{align}
where $\textbf{l}_1=\textbf{l}^2$.

\begin{remark}Except in Section \ref{sec:16}, we always use the estimates
  \begin{align}\label{for:108}
  \norm{\text{Ad}_{z_i^n}|_{\mathfrak{g}_{\mu,A}}}\leq C\mu(\textbf{a}_i)^ne^{\emph{\textbf{l}}|n|},\quad \forall\,\,n\in\ZZ,\,\,i=1,2;
\end{align}
\end{remark}

\item\label{for:472} $I-\text{Ad}_{z_{1}}$ is inventible on each $\mathfrak{g}_{\mu,A}$, $1\neq\mu\in \Delta_A$ with estimates
\begin{align*}
  \Big\|\big((I-\text{Ad}_{z_{1}})|_{\mathfrak{g}_{\mu,A}}\big)^{-1}\Big\|\leq C,\qquad \forall\,1\neq\mu\in \Delta_A.
\end{align*}



 \item\label{for:407}  if $w_i\in \exp(\mathfrak{g}_{1,A})\bigcap B_{\eta}(e)\cap Z(L)$, $i=1,2$, $\mathfrak{g}_{\mu,A}$, as well as $J_{\mu}$ (see \eqref{for:54} of Section \ref{sec:3}) are still invariant subspaces for $\text{Ad}_{z_i}$,  for any $\mu\in\Delta_A$.
Hence
\eqref{for:108} still hold if we replace $\mathfrak{g}_{\mu,A}$ by $J_\mu$.

 \end{enumerate}
  \smallskip


\end{enumerate}
\end{proposition}
The proof is left for Appendix \ref{sec:10}.

\begin{corollary}\label{for:129}
We can assume that $\eta(\textbf{a}_1,\textbf{a}_2)$ (see Proposition \ref{po:9}) is sufficiently small such that for any $w_i\in \exp(\mathfrak{g}_{1,A})\bigcap B_{\eta}(e)$, $i=1,\,2$,  letting $z_i=w_i\textbf{a}_i$, then there exists $q(z_1)\in \GG$ with $\textbf{d}(q,e)+\textbf{d}(q^{-1},e)<1$ such that for $\mathfrak{z}_i=qz_iq^{-1}$, $i=1,\,2$  the following properties hold:
\begin{enumerate}
 \item \label{for:78} $\mathfrak{z}_1$ is perfect (see \eqref{sect:11} o Section \ref{sec:3}) and
  $p(\mathfrak{z}_1)\in A_0$ ($A_0$ is Proposition \ref{po:3});

  \smallskip

  \item \label{for:57} we have
  \begin{align*}
   \norm{\text{Ad}_{(y_{\mathfrak{z}_i})^j}}\leq C(\abs{j}+1)^{\dim\mathfrak{G}},\quad \forall\,\,n\in\ZZ,\,\,i=1,2
  \end{align*}
  (see \eqref{sect:11} of Section \ref{sec:3});

  \smallskip
  \item \label{for:8} for any $\chi\in \Delta_{\mathfrak{z}_i}$
  \begin{align*}
   \norm{\text{Ad}_{\mathfrak{z}_i^n}|_{\mathfrak{g}_{\chi,\mathfrak{z}_i}}}\leq C\chi(p(\mathfrak{z}_i))^j(\abs{j}+1)^{\dim\mathfrak{G}},\quad \forall\,\,n\in\ZZ,\,\,i=1,2
  \end{align*}
  (see \eqref{for:273} of Section \ref{sec:3});

\smallskip

\item\label{for:435} \eqref{for:403} and \eqref{for:404}  of Proposition \ref{po:9} still hold if we replace $z_i$ by $\mathfrak{z}_i$, $i=1,\,2$;

\smallskip

\item \label{for:405}   the following relations hold:
   \begin{gather*}
 \mathfrak{g}_{\textbf{a}_1}^{+}\subseteq \mathfrak{g}^+_{\mathfrak{z}_1},\quad [\mathfrak{g}^+_{\mathfrak{z}_1},\,\mathfrak{g}_{\textbf{a}_1}^{+}]\subseteq
\mathfrak{g}_{\textbf{a}_1}^{+},\quad \mathfrak{g}_{\textbf{a}_1}^{-}\subseteq \mathfrak{g}^-_{\mathfrak{z}_1}  \quad\text{and}\\
 [\mathfrak{g}^-_{\mathfrak{z}_1},\,\mathfrak{g}_{\textbf{a}_1}^{-}]\subseteq
 \mathfrak{g}_{\textbf{a}_1}^{-}
\end{gather*}
(see \eqref{for:273} of Section \ref{sec:3});

\smallskip

\item\label{for:426} for any $m\geq1$, there is $0<\eta_m(\textbf{a}_1)<\eta(\textbf{a}_1,\textbf{a}_2)$ such that if $w_1\in \exp(\mathfrak{g}_{1,A})\cap B_{\eta_m(\textbf{a}_1)}(e)$, then
    \begin{align*}
  \norm{\text{Ad}_{\mathfrak{z}_1^j}|_{\mathfrak{g}_{1,A}}}\leq Ce^{\frac{|j|\textbf{l}}{m}}(\abs{j}+1)^{\dim\mathfrak{G}},\qquad \forall\,j\in\ZZ.
\end{align*}

\end{enumerate}

\end{corollary}

The proof is left for Appendix \ref{sec:51}.

\subsection{The perturbation of $\alpha_A$ }\label{sec:1}
We fix a compact generating set $E$ of $A$ such that each element of $E$ is regular and both $\textbf{a}_i$, $i=1,2$ (see \eqref{for:62} of Section \ref{sec:3}) are in $E$.
 We use $\mathfrak{j}$ to denote the natural
projection from $\mathcal{M}\rightarrow \mathcal{X}$.

For any $\varepsilon>0$ we say that $E'$ is a \emph{$\varepsilon$-standard   perturbation of $E$}  if $E'=\{w_v\cdot \alpha_v:\,v\in E,\,w_v\in B_{\varepsilon}(e)\cap \exp(\mathfrak{g}_{1,A})\cap Z(L)\}$. For any $v\in E$, we use $z_v=w_v\cdot \alpha_v$ to denote the corresponding element in $E'$.
We emphasize that elements in $E'$ do not need to commute.
 By using \eqref{for:490} of Section \ref{sec:3}, we see that
 $\mathfrak{g}_{\mu,A}$, $\mu\in \Delta_A$ are invariant subspaces for $\text{Ad}_{z_v}$, $v\in E$.

Suppose
$\tilde{\alpha}_A$ is a $C^\infty$ small perturbation of $\alpha_A$ on $\mathcal{X}$ and $E'$ is a $\varepsilon$-standard  small perturbation of $E$.
Then  we can write
\begin{align}\label{for:410}
 \tilde{\alpha}_{v}(\mathfrak{j}(x))=\mathfrak{j}\big(\exp(\mathfrak{p}_v(x))\cdot
z_v\cdot x\big),\qquad \forall\,\,v\in E
\end{align}
for any $x\in \mathcal{M}$, where $\mathfrak{p}_v\in C^\infty(\mathcal{M},\,\mathcal{L}^\bot)$ (see \eqref{for:17} of Section \ref{sec:3}) are small maps.

Since $\tilde{\alpha}_A$ is $L$-fiber preserving, we have the following result, the proof of which  is left for Appendix \ref{sec:52}:
\begin{lemma}\label{ob:1} $\mathfrak{p}_v$ is adjoint $L$-invariant (see \eqref{for:117} of Section \ref{sec:3}) for any $v\in E$.
\end{lemma}
 For any $v,\,u\in E$ we define
\begin{gather*}
  r_v=\int_{\mathcal{M}}\mathfrak{p}_v,\quad \mathfrak{p}'_v=\mathfrak{p}_{v}-r_v, \qquad\text{and}\\
   \mathbb{D}_{v,u}=\mathfrak{p}'_{v}\circ z_{u}-\mathfrak{i}_{\mathcal{L}^\bot}\text{Ad}_{z_{u}}\mathfrak{p}'_{v}-\big(\mathfrak{p}'_{u}\circ z_{v}-\mathfrak{i}_{\mathcal{L}^\bot}\text{Ad}_{z_{v}}\mathfrak{p}'_{u}\big)
\end{gather*}
 ($\mathfrak{i}_{\mathcal{L}^\bot}$ is defined in \eqref{for:54} of Section \ref{sec:3}). Set
\begin{align*}
 \norm{\mathfrak{p}}_{C^m}=\max_{z\in E}\{\norm{\mathfrak{p}_z}_{C^m}\}\text{ and }\norm{\mathbb{D}}_{C^m}=\max_{u,v\in E}\{\norm{\mathbb{D}_{u,v}}_{C^m}\}\quad \forall\,m\geq0.
\end{align*}
Similar to Lemma \ref{le:1}, we have the following result:
\begin{lemma}\label{le:8} For any $m\in\NN$, there is $\delta_m>0$ such that if $\norm{\mathfrak{p}}_{C^{m}}\leq\delta_m$, then
\begin{align}
 \max_{u,v\in E}\norm{\log[z_{v},z_{u}]}&\leq C\norm{\mathfrak{p}}_{C^0}, \quad\text{ and }\\
 \norm{\mathbb{D}}_{C^{m}}
 &\leq C\norm{\mathfrak{p}}_{C^{m+1}}\norm{\mathfrak{p}}_{C^{m}}.
\end{align}
For any $1\neq\mu\in \Delta_A$, we have
\begin{align}\label{for:39}
  \big\|(r_v-\text{Ad}_{z_u}r_v-r_u+\text{Ad}_{z_v}r_u)\big|_{\mathfrak{g}_{\mu,A}}\big\|\leq C\norm{\mathfrak{p}}_{C^1}\norm{\mathfrak{p}}_{C^0}.
\end{align}

\end{lemma}
Proof of Lemma \ref{le:8} is left for Appendix \ref{sec:4}.

\subsection{The twisted cohomological equation}\label{sec:14} In this section, we assume $z_1=w\textbf{a}_1$ where $w\in \exp(\mathfrak{g}_{1,A})\cap B_{\eta(\textbf{a}_1,\textbf{a}_2)}(e)\cap Z(L)$ ($\eta(\textbf{a}_1,\textbf{a}_2)$ is from Proposition \ref{po:9}). We study the twisted equation:
\begin{align}\label{for:437}
\varphi\circ z_1-\mathfrak{i}_{\mathcal{L}^\bot}\text{Ad}_{z_1}\varphi=\omega,
\end{align}
where $\varphi,\,\omega\in \mathcal{L}^\bot(\mathcal{O})$ (see \eqref{for:54} of Section \ref{sec:3}).

We recall $q(z_1)$ and $\mathfrak{z}_1=qz_1q^{-1}$, as described in Corollary \ref{for:129}. To study equation \eqref{for:437}, it  is equivalent to study
the twisted equation:
\begin{align}\label{for:5}
\varphi_q\circ \mathfrak{z}_1-\mathfrak{i}_{\mathcal{L}^\bot}\text{Ad}_{z_1}\varphi_q=\omega_q,
\end{align}
where $\varphi_q=\varphi\circ q^{-1}$ and $\omega_q=\omega\circ q^{-1}$.

\smallskip
\emph{Note}. Since $q$ is uniformly bounded, passing to equation \eqref{for:5} will not affect the results in this part.

\smallskip

From \eqref{for:54} of Section \ref{sec:3}, we see that \eqref{for:5} decomposes into equations of the following forms:
\begin{align}\label{for:105}
\varphi_\mu\circ \mathfrak{z}_1-(\text{Ad}_{z_1}|_{J_\mu})\varphi_\mu=\omega_\mu,\qquad \mu\in \Delta_A.
\end{align}
where $\varphi_\mu=\varphi\circ q^{-1}|_{J_\mu}$ and $\omega_\mu=\omega\circ q^{-1}|_{J_\mu}$.

 The subsequent parts of this section are devoted to the study of the decomposed equation \eqref{for:105}.

\subsubsection{Notations}\label{sec:53} The positive constants $\eta, \sigma,\,\emph{\textbf{l}}$ (see Section \ref{for:432}) will appear in is section.
Suppose $\xi\in J_\mu(\mathcal{O}^m)$, $m\geq \sigma$  and $z_1=w\textbf{a}_1$ where $w\in \exp(\mathfrak{g}_{1,A})\cap B_{\eta(\textbf{a}_1,\textbf{a}_2)}(e)\cap Z(L)$.

\begin{enumerate}
  \item\label{for:412} $\mathcal{P}_{\mu}(\xi)\stackrel{\text{def}}{=} \delta(\mu)\sum_{n\in L_\mu} (\text{Ad}_{z_1}|_{J_\mu})^{n-1} \xi\circ \mathfrak{z}_1^{-n}$ (see \eqref{for:1} of Section \ref{sec:3}).

    \smallskip
  \item\label{for:413} For any $\textbf{u}\in \mathcal{U}_{\mu,\sigma}$ (see \eqref{for:1} of Section \ref{sec:3}) and any subset $\mathbb{S}\subseteq\ZZ$
\begin{align*}
 \mathfrak{B}_{\mu,\textbf{u},\mathbb{S}}(\xi)\stackrel{\text{def}}{=}\sum_{n\in \mathbb{S}} (\text{Ad}_{z_1}|_{J_\mu})^{n-1}\textbf{u}(\xi\circ \mathfrak{z}_1^{-n}).
\end{align*}

  \item\label{for:414} For any subset $\mathbb{S}\subseteq\ZZ$
  \begin{align*}
 \norm{\mathfrak{B}_{\mu,\mathbb{S}}(\xi)}&\stackrel{\text{def}}{=}\max_{\textbf{u}\in \mathcal{U}_{\mu,\sigma}}\sup_{\stackrel{\psi\in J_\mu(\mathcal{O}^{\sigma+1}),}{\norm{\psi}_{\sigma+1}=1}} \big\|\langle \mathfrak{B}_{\mu,\textbf{u},\mathbb{S}}(\xi),\,\psi\rangle\big\|
\end{align*}
  (recall \eqref{for:17} of Section \ref{sec:3} for the definition of $\langle \cdot \rangle$).

  \smallskip

  \item\label{for:415} For any subset $\mathbb{S}\subseteq\ZZ$, any $\tau\geq1$, any $\textbf{u}\in \mathcal{U}_{\mu,\sigma}$ and any $v\in\mathfrak{g}^{-\delta(\mu)}_{\mathfrak{z}_1}\cap\mathfrak{G}^1$ (see \eqref{for:273}, \eqref{for:1} of Section \ref{sec:3})
\begin{align*}
 \textbf{u}\mathfrak{D}_{\mu,v,\tau,\mathbb{S}}(\xi)\stackrel{\text{def}}{=}\sum_{n\in \mathbb{S}} (\text{Ad}_{z_1}|_{J_\mu})^{n-1}\textbf{u}\pi_{v}(f_1\circ \tau^{-1})(\xi\circ \mathfrak{z}_1^{-n})
\end{align*}
($f_1$ is defined in \eqref{sect:4} of Section \ref{sec:19}).

\smallskip

\item\label{for:465} For any subset $\mathbb{S}\subseteq\ZZ$, any $\tau\geq1$, any $\textbf{u}\in \mathcal{U}_{\mu,\sigma}$
\begin{align*}
 \norm{\textbf{u}\mathfrak{D}_{\mu,\tau,\mathbb{S}}(\xi)}&\stackrel{\text{def}}{=}\sup_{v\in\mathfrak{g}^{-\delta(\mu)}_{z_1}\cap\mathfrak{G}^1}\sup_{\stackrel{\psi\in J_\mu(\mathcal{O}^{\sigma+1}),}{\norm{\psi}_{\sigma+1}=1}} \big\|\langle \textbf{u}\mathfrak{D}_{\mu,v,\tau,\mathbb{S}}(\xi),\,\psi\rangle\big\|
\end{align*}
for any $\xi\in J_\mu(\mathcal{O}^{\sigma+1})$ (recall \eqref{for:17} of Section \ref{sec:3} for the definition of $\langle \cdot \rangle$).


\end{enumerate}

\subsubsection{Main results} Suppose $\xi,\,\omega_\mu\in J_\mu(\mathcal{O}^\infty)$ and $\mu\in \Delta_A$. Also suppose $\textbf{u}\in \mathcal{U}_{\mu,k}$ (see \eqref{for:1} of Section \ref{sec:3}), $k\geq\sigma$. Then:
\begin{enumerate}
  \item (\eqref{for:434} of Lemma \ref{le:12})  for any $n\in L_\mu^c$, $\textbf{u}(\text{Ad}_{z_1}|_{J_\mu})^{n-1}(\xi\circ \mathfrak{z}_1^{-n})\in J_\mu(\mathcal{O})$ with estimate
     \begin{align*}
\big\|(\text{Ad}_{z_1}|_{J_\mu})^{n-1}\textbf{u}(\xi\circ \mathfrak{z}_1^{-n})\big\|
  \leq Ce^{-10\emph{\textbf{l}}|n|}\norm{\xi}_{k};
\end{align*}

  \smallskip
  \item (\eqref{for:448} of Lemma \ref{le:12}) for any $n\in \ZZ$  we have
     \begin{align*}
\Big\|\big\langle (\text{Ad}_{z_1}|&_{J_\mu})^{n-1}(\xi\circ \mathfrak{z}_1^{-n}),\,\textbf{u}\psi\big\rangle\Big\|
  \leq Ce^{-3\emph{\textbf{l}}|n|}\norm{\xi}_{k}\norm{\psi}_{k+1};
\end{align*}
for any $\psi\in J_\mu(\mathcal{O}^{k+1})$;

\smallskip
   \item (\eqref{for:431} of Lemma \ref{le:7}) $\mathcal{P}_{\mu}(\omega_\mu)$ is a distributional solution of equation \eqref{for:105};

   \smallskip

   \item (\eqref{for:132} Lemma \ref{le:12}) if equation \eqref{for:105} has a solution $\varphi_\mu\in J_\mu(\mathcal{O}^\sigma)$,  then
  \begin{align*}
   \norm{\mathfrak{B}_{\mu,\ZZ}(\omega_\mu)}=0.
  \end{align*}

  \smallskip
  \item (\eqref{for:28} of Lemma \ref{le:7}) for any $m\geq \sigma+6$,  if $z_1=w\textbf{a}_1$ where $w\in \exp(\mathfrak{g}_{1,A})\cap B_{\eta_m(\textbf{a}_1)}(e)\cap Z(L)$ (see \eqref{for:426} of Corollary \ref{for:129}) and if $\mathfrak{B}_{\mu,\ZZ}(\omega_\mu)=0$, then
equation \eqref{for:105} has a solution
$\varphi_\mu\in J_{\mu}(\mathcal{O}^{m-6})$ with estimates
\begin{align*}
 \norm{\varphi_\mu}_{m-6}\leq C_{m}\norm{\omega_\mu}_{m}.
\end{align*}
\item (Lemma \ref{for:422})  Suppose $\tau\geq1$ and $\xi\in J_\mu(\mathcal{O}^m)$, $m\geq \sigma$. Also suppose  $\textbf{u}\in \mathcal{U}_{\mu,\sigma}$. Then for any $v\in\mathfrak{g}^{-\delta(\mu)}_{z_1}\cap\mathfrak{G}^1$, $\textbf{u}\mathfrak{D}_{\mu,v,\tau,L_\mu^c}(\xi)\in J_{\mu}(\mathcal{O})$ (see \eqref{for:415} of Section \ref{sec:53}) with estimates
\begin{align*}
 \big\| \textbf{u}\mathfrak{D}_{\mu,v,\tau,L_\mu^c}(\xi)\big\|\leq C_m\tau^{-m+\sigma}\norm{\xi}_m.
\end{align*}
\item (\eqref{for:64} of Proposition \ref{po:10})  If $m\geq \sigma+1$ and $\textbf{u}\in \mathcal{U}_{\mu,\sigma}$, then for any $\psi \in J_\mu(\mathcal{O}^{\sigma+1})$ we have
\begin{align*}
 \big\|\langle \textbf{u}\mathfrak{D}&_{\mu,v,\tau,\ZZ}(\xi),\,\psi\rangle\big\|\leq C\norm{\mathfrak{B}_{\mu,\ZZ}(\xi)}\norm{\psi}_{\sigma+1}.
\end{align*}
\end{enumerate}

\subsubsection{Useful facts for $z_1$ and $\mathfrak{z}_1$}\label{sec:54} We list the facts that will be used for $z_1$ and $\mathfrak{z}_1$.
\begin{enumerate}
  \item\label{for:442} (\eqref{for:78} of Corollary \ref{for:129}) $\mathfrak{z}_1$ is perfect  and
  $p(\mathfrak{z}_1)\in A_0$. Thus, $\mathfrak{z}_1^n=p(\mathfrak{z}_1)^ny_{\mathfrak{z}_1}^n$ for any $n\in\ZZ$.

  \smallskip

    \item\label{for:439} (\eqref{for:1} of Section \ref{sec:3}) For any $\mu\in\Delta_A$, $\mu(\textbf{a}_1)^{n}\leq1$ if $n\in L_\mu$ and $\mu(\textbf{a}_1)^{n}>1$ if $n\in L_\mu^c$.

If $\lambda_i\in\Delta_{\textbf{a}_1}^{-\delta(\mu)}$, $1\leq i\leq j$ then $(\mu,\lambda_1,\cdots, \lambda_j)$ is a poor row.
  \smallskip

  \item (\eqref{for:108} of Proposition \ref{po:9} and \eqref{for:435} of Corollary \ref{for:129}) For both $\text{Ad}_{z_1}$ and $\text{Ad}_{\mathfrak{z}_1}$:
  \begin{enumerate}
    \item $\mathfrak{g}_{\mu,A}$ are invariant subspaces for any $\mu\in\Delta_A$;

    \smallskip

    \item\label{for:441} \eqref{for:108} of Proposition \ref{po:9} holds true when estimating both $\norm{\text{Ad}_{z_1^n}|_{\mathfrak{g}_{\mu,A}}}$ and $\norm{\text{Ad}_{\mathfrak{z}_1^n}|_{\mathfrak{g}_{\mu,A}}}$.

  \end{enumerate}

  \smallskip

 \item\label{for:440} (\eqref{for:8} in Corollary \eqref{for:129} and \eqref{for:1} of Section \ref{sec:3}) Estimates of $\norm{\text{Ad}_{\mathfrak{z}_1^n}}$ on each Lyapunov space $\mathfrak{g}_{\chi,\mathfrak{z}_i}$, $\chi\in \Delta_{\mathfrak{z}_i}$ are obtained.

 We also note that $\chi(p(\mathfrak{z}_1))^{-n}\leq1$ if $n\in L_\mu$ for any $\chi\in \Delta_{\mathfrak{z}_1}^{\delta(\mu)}\cup\Delta_{\mathfrak{z}_1}^{0}$ and
 $\chi(p(\mathfrak{z}_1))^{n}>1$ if $n\in L_\mu^c$ for any $\chi\in \Delta_{\mathfrak{z}_1}^{-\delta(\mu)}$.

  \smallskip
  \item\label{for:438} (\eqref{for:403} of Proposition \ref{po:9} and \eqref{for:435} of Corollary \ref{for:129}) Decay of matrix coefficients of $\mathfrak{z}_1$ is provided.
\end{enumerate}

\subsubsection{A typical example}\label{sec:62} We use a typical example described in Section \ref{sec:32} to explain some frequently used nations in this part.
We use $u_{i,1,2}$ (resp. $u_{i,2,1}$) to denote the upper (lower) triangular $2\times 2$  unipotent matrix in the $i$-th copy of $\mathfrak{sl}(2,\RR)$ and use
$X_{i}$  to denote the diagonal $2\times 2$ matrix with entries $1$ and $-1$ in the $i$-th copy of $\mathfrak{sl}(2,\RR)$, $1\leq i\leq 3$. For example,
\begin{align*}
 u_{2,1,2}&=\begin{pmatrix}
0 & 0 \\
0 & 0
\end{pmatrix}\times \begin{pmatrix}
0 & 1 \\
0 & 0
\end{pmatrix}\times \begin{pmatrix}
0 & 0 \\
0 & 0
\end{pmatrix},\\
 X_{3}&=\begin{pmatrix}
0 & 0 \\
0 & 0
\end{pmatrix}\times \begin{pmatrix}
0 & 0 \\
0 & 0
\end{pmatrix}\times \begin{pmatrix}
1 & 0 \\
0 & -1
\end{pmatrix}.
\end{align*}
We see that $\Delta_A=\{\lambda_1,\lambda_2,\,\lambda_3=\lambda_1^{-1},\,\lambda_4=\lambda_2^{-1}, 1\}$ and
\begin{align*}
 \lambda_1(d\textbf{a}_1^md\textbf{a}_2^n)=(\text{\small$\frac{1}{2}$})^{m+n}\quad\text{and}\quad \lambda_2(d\textbf{a}_1^md\textbf{a}_2^n)=(\text{\small$\frac{1}{4}$})^{m-n}.
\end{align*}
$\Delta_{z_1}=\{\chi_1, \chi_2,\,\chi_3=\chi_1^{-1},\,\chi_4=\chi_2^{-1}, \chi_5,\,\chi_6=\chi_5^{-1},\,1\}$ and
\begin{align*}
 \chi_1(dz_1)=\text{\small$\frac{1}{2(1+\epsilon)^2}$},\quad \chi_2(dz_1)=\text{\small$\frac{1}{4(1+\epsilon)^2}$}\quad\text{and}\quad \chi_5(dz_1)=\text{\small$\frac{1}{(1+\epsilon)^2}$}.
\end{align*}
We also suppose $\xi\in \mathfrak{g}_\mu(\mathcal{O}^\infty)$. Then:
\begin{enumerate}
\item $A_0$ is the subgroup generated by $\exp(X_1)$, $\exp(X_2)$ and $\exp(X_3)$. $z_1$ is perfect and $z_1\in A_0$. In this case, $q=I$ and $z_1=\mathfrak{z}_1$ (see Corollary \ref{for:129}).

\smallskip
  \item $\Delta_{z_1}^{+}=\{\chi_3, \chi_4, \chi_6\}$, $\Delta_{z_1}^{-}=\{\chi_1,\chi_2, \chi_5\}$ and $\Delta_{z_1}^{0}=\{1\}$.
  \smallskip

  \item We note that
  \begin{gather*}
   \mathfrak{g}_{z_1}^{0}=\{X_1, X_2,X_3\}\subseteq\mathfrak{g}_{\textbf{a}_1}^{0}=\{X_1, X_2,X_3, u_{3,1,2}, u_{3,2,1}\}\\
   \mathfrak{g}_{\textbf{a}_1}^{+}=\{u_{1,1,2},u_{2,1,2}\}\subseteq \mathfrak{g}_{z_1}^{+}=\{u_{1,1,2},u_{2,1,2},u_{3,1,2}\}\\
   \mathfrak{g}_{\textbf{a}_1}^{-}=\{u_{1,2,1},u_{2,2,1}\}\subseteq \mathfrak{g}_{z_1}^{-}=\{u_{1,2,1},u_{2,2,1},u_{3,2,1}\}.
  \end{gather*}
 We point out $z_1$ is a regular element for  $\alpha_{A_0}$, the full Cartan action.  However, $\textbf{a}_1$ is NOT a regular element for the full Cartan action, even though it is regular for  $\alpha_{A}$. As a result, $z_1$ and $\textbf{a}_1$ have different Lyapunov exponents
and different Lyapunov subspaces.

 \smallskip

  \item If $\mu\in \{\lambda_3,\lambda_4,1\}$, $L_\mu=\ZZ^-\cup\{0\}$, $\delta(\mu)=-$ and $L_\mu^c=\ZZ^+$.

  \smallskip
  \item If $\mu\in \{\lambda_1,\lambda_2\}$, $L_\mu=\ZZ^+\cup\{0\}$, $\delta(\mu)=+$ and $L_\mu^c=\ZZ^-$.

  \smallskip

  \item If $\mu\in \{\lambda_3,\lambda_4,1\}$, then $(\mu,\,\lambda_1,\cdots,\lambda_j)$ is a poor row if $\lambda_k\in\{\lambda_3,\lambda_4\}$, $1\leq k\leq j$.

  \smallskip

   \item If $\mu\in \{\lambda_1,\lambda_2\}$, then $(\mu,\,\lambda_1,\cdots,\lambda_j)$ is a poor row if $\lambda_k\in\{\lambda_1,\lambda_2\}$, $1\leq k\leq j$.

   \smallskip

   \item If $\mu\in\{\lambda_3,\lambda_4,1\}$
   \begin{align*}
    \mathcal{P}_{\mu}(\xi)&= -\sum_{n\leq 0} (\text{Ad}_{z_1}|_{\mathfrak{g}_{\mu}})^{n-1} \xi\circ z_1^{-n}
   \end{align*}
   is a distribution. Further,
for any $v_i\in \mathfrak{g}^-_{z_1}\oplus\mathfrak{g}^0_{z_1}
=\mathfrak{g}^{\delta(\mu)}_{z_1}\oplus\mathfrak{g}^0_{z_1}$, we see that
\begin{align*}
 v_1v_2\cdots v_l\mathcal{P}_{\mu}(\xi)
\end{align*}
are also a distributions (see \eqref{for:421} and \eqref{for:106} of Lemma \ref{le:12}).

Then $\mathcal{P}_{\mu}(\xi)$ is $\infty$ differentiable along $\mathfrak{g}^{\delta(\mu)}_{z_1}\oplus\mathfrak{g}^0_{z_1}$.
However, $\mathcal{P}_{\mu}(\xi)$ probably loses differentiability along $\mathfrak{g}^{-\delta(\mu)}_{z_1}$.

   \smallskip

   \item If $\mu\in\{\lambda_1,\lambda_2\}$
   \begin{align*}
    \mathcal{P}_{\mu}(\xi)= \sum_{n\geq0} (\text{Ad}_{z_1}|_{\mathfrak{g}_{\mu}})^{n-1} \xi\circ z_1^{-n}\in \mathfrak{g}_{\mu}(\mathcal{O})
   \end{align*}
    as the twisted $0<\text{Ad}_{z_1}|_{\mathfrak{g}_{\mu}}<1$. Further,
for any $v_i\in \mathfrak{g}^+_{z_1}\oplus\mathfrak{g}^0_{z_1}
=\mathfrak{g}^{\delta(\mu)}_{z_1}\oplus\mathfrak{g}^0_{z_1}$,
\begin{align*}
 v_1v_2\cdots v_l\mathcal{P}_{\mu}(\xi)\in \mathfrak{g}_{\mu}(\mathcal{O}).
\end{align*}
Then we see that $\mathcal{P}_{\mu}(\xi)$ is $\infty$ differentiable along $\mathfrak{g}_{\lambda}$ if $(\mu,\lambda)$ is a NOT poor pair.

Then $\mathcal{P}_{\mu}(\xi)$ is $\infty$ differentiable along $\mathfrak{g}^{\delta(\mu)}_{z_1}\oplus\mathfrak{g}^0_{z_1}$.
However, $\mathcal{P}_{\mu}(\xi)$ probably loses differentiability along $\mathfrak{g}^{-\delta(\mu)}_{z_1}$.
\end{enumerate}

\subsubsection{Main ideas} Below, we use a typical example described in Section \ref{sec:32} to explain the main idea for the sake of clarity and transparency. Throughout Section \ref{sec:24}, it is recommended
to keep this example in mind.

We recall the twisted equation \ref{for:389} and recall
\begin{align*}
u_{1,1,2}=\begin{pmatrix}
0 & 1 \\
0 & 0
\end{pmatrix}\times \begin{pmatrix}
0 & 0 \\
0 & 0
\end{pmatrix}\times \begin{pmatrix}
0 & 0 \\
0 & 0
\end{pmatrix}  \quad\text{see \eqref{sec:62}. }
\end{align*}
Then $\text{Ad}_{z_1}u_{1,1,2}=2(1+\epsilon)^2u_{1,1,2}$.
 $u_{1,1,2}$ is a unstable direction for both $z_1$ and $\textbf{a}_1$. For the following
formal sum
\begin{align}\label{for:526}
 u_{1,1,2}^m\mathfrak{D}:=\sum_{n=-\infty}^\infty (2(1+\epsilon)^2)^{-(n+1)} u_{1,1,2}^m(\omega\circ z_1^{n})
\end{align}
we will determine $m\in\NN$ such that it is a distribution. We can write
\begin{align*}
 \mathfrak{D}=\mathfrak{D}_{+}+\mathfrak{D}_-
\end{align*}
where
\begin{align*}
    \mathfrak{D}_{+}=\sum_{n\geq 0} (2(1+\epsilon)^2)^{-(n+1)}\omega\circ z_1^n\quad\text{and}\quad \mathfrak{D}_{-}=\sum_{n<0} (2(1+\epsilon)^2)^{-(n+1)} \omega\circ z_1^n.
\end{align*}
Firstly, we show the convergence of $u_{1,1,2}^m\mathfrak{D}_{+}$: for any $\psi \in \mathcal{O}^\infty$, we have
\begin{align}\label{for:527}
 |\langle u_{1,1,2}^m\mathfrak{D}_{+},\,\psi\rangle|&=\big| \sum_{n\geq 0} (2(1+\epsilon)^2)^{-(n+1)} \langle\omega\circ z_1^n,\, u_{1,1,2}^m\psi \rangle\big|\notag\\
 &\overset{\text{(1)}}{\leq} \sum_{n\geq 0} (2(1+\epsilon)^2)^{-(n+1)} \norm{\omega\circ z_1^n}\cdot \norm{u_{1,1,2}^m\psi }\notag\\
 & \overset{\text{(2)}}{\leq} \sum_{n\geq 0} (2(1+\epsilon)^2)^{-(n+1)} \norm{\omega}\cdot \norm{\psi }_m\notag\\
 &\overset{\text{(3)}}{\leq} C\norm{\omega}\cdot \norm{\psi }_m.
\end{align}
Here in $(1)$ we use cauchy schwarz inequality; in $(2)$ we use the fact that
\begin{align*}
\norm{\varphi}=\norm{\varphi\circ z_1^n},\qquad \forall\,\varphi\in \mathcal{O};
\end{align*}
in $(3)$ we note that $\epsilon$ is sufficiency close to $0$.

\smallskip
\emph{Note}. When $n\geq 0$, the twist $(2(1+\epsilon)^2)^{-(n+1)}$ is not an obstacle for convergence; rather; it is helpful for convergence.

\smallskip

Secondly, we show the convergence of $u_{1,1,2}^m\mathfrak{D}_{-}$. Before that, we show how the twist is ``balanced" by derivatives along $v$, which guarantees the convergence of $u_{1,1,2}^m\mathfrak{D}_{-}$. We recall that $\text{Ad}_{z_1}v=2v$. Thus,
\begin{align*}
  u_{1,1,2}^m(\omega\circ z_1^{n})=(2(1+\epsilon)^2)^{nm}(u_{1,1,2}^m\omega)\circ z_1^n.
\end{align*}
It follows that for any $\psi \in \mathcal{O}^\infty$ and $n<0$, we have
\begin{align}\label{for:524}
 \big|\langle u_{1,1,2}^m&(\omega\circ z_1^n),\, \psi \rangle\big|=\big| \langle (2(1+\epsilon)^2)^{nm}(u_{1,1,2}^m\omega)\circ z_1^n,\, \psi \rangle\big|\notag\\
 &\overset{\text{(1)}}{\leq} (2(1+\epsilon)^2)^{nm}\norm{(u_{1,1,2}^m\omega)\circ z_1^n}\cdot \norm{\psi}\notag\\
&\overset{\text{(2)}}{=} (2(1+\epsilon)^2)^{nm}\norm{u_{1,1,2}^m\omega}\cdot \norm{\psi}\notag\\
&\leq (2(1+\epsilon)^2)^{nm}\norm{\omega}_m\cdot \norm{\psi}.
\end{align}
Here in $(1)$ we use cauchy schwarz inequality; in $(2)$ the fact that
\begin{align*}
\norm{\varphi}=\norm{\varphi\circ z_1^n},\qquad \forall\,\varphi\in \mathcal{O}.
\end{align*}
\emph{Note}. As $n<0$, $(2(1+\epsilon)^2)^{nm}=(\frac{1}{2(1+\epsilon)^2})^{|n|m}$. Noting that $\omega\in \mathcal{O}^\infty$, we can choose $m$ to be as large as we want.  This is how we use
directional derivatives to ``beat" the twist.
\smallskip

By choosing $m>1$, the convergence of $v^m\mathfrak{D}_{-}$ follows directly from \eqref{for:524}:
\begin{align}\label{for:528}
 |\langle u_{1,1,2}^m\mathfrak{D}_{-},\,\psi\rangle|&=\big| \sum_{n<0} (2(1+\epsilon)^2)^{-(n+1)} \langle u_{1,1,2}^m(\omega\circ z_1^n),\, \psi \rangle\big|\notag\\
 &\leq\sum_{n<0} (2(1+\epsilon)^2)^{-(n+1)}\big|  \langle u_{1,1,2}^m(\omega\circ z_1^n),\, \psi \rangle\big|\notag\\
 &\leq\sum_{n<0} (2(1+\epsilon)^2)^{-(n+1)}\cdot (2(1+\epsilon)^2)^{nm}\norm{\omega}_m\cdot \norm{\psi}\notag\\
 &\leq C_{m}\norm{\omega}_{m}\cdot \norm{\psi }.
\end{align}
Hence, we showed the convergence of $u_{1,1,2}^m\mathfrak{D}_{-}$.

It follows from \eqref{for:527} and \eqref{for:528} that
\begin{align}\label{for:529}
|\langle u_{1,1,2}^m\mathfrak{D},\,\psi\rangle|\leq C\norm{\omega}\cdot \norm{\psi }_m+C_{m}\norm{\omega}_{m}\cdot \norm{\psi }
\leq C_{m,1}\norm{\omega}_{m}\norm{\psi }_m.
\end{align}
Next, we show how to obtain uniform estimates of  $u_{1,1,2}^m\mathfrak{D}_{\mu,u_{3,1,2},\tau,\ZZ}(\omega)$ (recall \eqref{for:391}) if $m>1$.
We note that: for any $n\in\ZZ$
\begin{align}\label{for:525}
u_{1,1,2}^m\pi_{u_{3,1,2}}(f_1\circ \tau^{-1})\omega\circ z_1^n\overset{\text{(1)}}{=}\pi_{u_{3,1,2}}(f_1\circ \tau^{-1})u_{1,1,2}^m(\omega\circ z_1^n).
\end{align}
Here in $(1)$ we note that $[u_{3,1,2},u_{1,1,2}]=0$ and then use \eqref{for:10} of Section \ref{sec:27}.

It follows that
\begin{align}\label{for:530}
 &|\langle u_{1,1,2}^m \mathfrak{D}_{\mu,u_{3,1,2},\tau,\ZZ}(\omega),\,\psi\rangle|\notag\\
 &=\Big| \sum_{n\in\ZZ} 2^{-(n+1)} \big\langle  u_{1,1,2}^m\pi_{u_{3,1,2}}(f_1\circ \tau^{-1})(\omega\circ z_1^n),\, \psi \big\rangle\Big|\notag\\
 &\overset{\text{(1)}}{=}\Big| \sum_{n\in\ZZ} 2^{-(n+1)} \big\langle  \pi_{u_{3,1,2}}(f_1\circ \tau^{-1})u_{1,1,2}^m(\omega\circ z_1^n),\, \psi \big\rangle\Big|\notag\\
&\overset{\text{(2)}}{=}\Big| \sum_{n\geq 0} 2^{-(n+1)} \big\langle  u_{1,1,2}^m(\omega\circ z_1^n),\, \pi_{u_{3,1,2}}(f_1\circ \tau^{-1})\psi \big\rangle\Big|\notag\\
&\overset{\text{(3)}}{=}\Big| \big\langle  u_{1,1,2}^m\mathfrak{D},\, \pi_{u_{3,1,2}}(f_1\circ \tau^{-1})\psi \big\rangle\Big|\\
&\overset{\text{(4)}}{\leq}C_{m}\norm{\omega}_{m}\big\|\pi_{u_{3,1,2}}(f_1\circ \tau^{-1})\psi \big\|_m\notag\\
&\overset{\text{(5)}}{\leq}C_{m,1}\norm{\omega}_{m}\big\|\psi \big\|_m
\end{align}
Here in $(1)$ we use \eqref{for:525}; in $(2)$ we use \eqref{for:93} of Section \ref{sec:27}; in $(3)$ we recall \eqref{for:526};
in $(4)$ we use \eqref{for:529}; in $(5)$ we use \eqref{for:7} of Lemma \ref{cor:1}.

Hence we also get the uniform estimates of  $u_{1,1,2}^m\mathfrak{D}_{\mu,u_{3,1,2},\tau,\ZZ}(\omega)$.
\begin{remark}\label{re:9} Comments from the above discussion:

\begin{enumerate}
\item From \eqref{for:527} we see that $\mathfrak{D}_{+}$ is also distribution, which is in fact a solution of the twisted equation \ref{for:389}. As $\text{Ad}_{z_1}u_{1,1,2}=2(1+\epsilon)^2u_{1,1,2}$ and $u_{1,1,2}\in \mathfrak{g}_{\lambda_3,A}$,
$\mathfrak{D}_{+}$ is the sum over $\ZZ^+=L_{\lambda_3}$ (see Section \ref{sec:62}), where  $\lambda_3\in \Delta_{\textbf{a}_1}^+$.

If $|\lambda_3(\textbf{a}_1)|=1$, we need to use decay of matrix coefficients to show the convergence of $\mathfrak{D}_{+}$.

We point out that $\mathfrak{D}_{+}$ is exactly the $\mathcal{P}_{\lambda_3}(\omega)$ defined in \eqref{for:412} of Section \ref{sec:53}, where
$L_{\lambda_3}$ determines on which subset of $\ZZ$ we define $\mathcal{P}_{\lambda_3}(\omega_{\lambda_3})$.

\smallskip

\item $\mathfrak{D}_{-}$ is the sum over $\NN^c=L_{\lambda_3}^c$, which is not a distribution. To obtain uniform estimates for $v^m\mathfrak{D}_{-}$, we choose $v=u_{1,1,2}$ and $m>1$. From \eqref{for:524},  we see that if $v_i\in\mathfrak{g}_{\textbf{a}_1}^+$, $1\leq i\leq m$ and $m>c(z_1,+)$ (see \eqref{for:523}), then $v_1\cdots v_m\mathfrak{D}_{-}$ still has uniform estimates. It is clear that we will not choose any $v_i=u_{3,1,2}\in \mathfrak{g}_{\textbf{a}_1}^0\cap \mathfrak{g}_{z_1}^+$, as the expansion rate of
$\text{Ad}_{z_1}$ along $u_{3,1,2}$ is not uniformly bounded away from $1$. This is the reason we define $\mathcal{U}_{\lambda,k}$ (see \eqref{for:1} of Section \ref{sec:3}). Thus, for any $\textbf{u}\in \mathcal{U}_{\lambda_3, \lfloor c(z_1,+)\rfloor+1}$ where $\lfloor \cdot\rfloor$ denotes the floor function, $\textbf{u}\mathfrak{D}_{-}$ has uniform estimates. We point out that $\textbf{u}\mathfrak{D}_{-}$ is exactly the $\mathfrak{B}_{\mu,\textbf{u},\NN^c}$ defined in
\eqref{for:413} of Section \ref{sec:53}.

  \smallskip

\item We note that the notation $\textbf{u}\mathfrak{D}_{\mu,v,\tau,\mathbb{S}}$ in \eqref{for:413} of Section \ref{sec:53} is different from the $\mathfrak{D}_{\mu,\lambda,u,\tau,\mathbb{S}}$ in \eqref{for:285} of Section \ref{sec:39}. For the latter, $\lambda(dz_1)$, the deformation rate is crucial for the convergence. Thus, in the notation we put $\lambda$ in the subscript. However, for the former, it is $\textbf{u}$,  the directional derivatives aid the convergence (see \eqref{for:530}). Therefore,  we remove $\lambda$ from the subscript.

\item We note that, generally,  we need Proposition \ref{po:10} to handle the `` commutator relations" in \eqref{for:525}.

\smallskip

 \item We will see that $\mathfrak{B}_{\mu,\textbf{u},\ZZ}$ are also invariant distributions (see \eqref{for:132} of Lemma \ref{le:12}). However, we do not use them to classify the obstructions to solving \eqref{for:105}. One reason is that the existence of a low regularity solution implies vanishing of $\mathfrak{B}_{\mu,\textbf{u},\ZZ}$  (see \eqref{for:132} of Lemma \ref{le:12}). However, the existence of a low regularity solution does not necessarily  imply the existence of a high regularity solution (see \eqref{for:531} of Section \ref{sec:32}). Another reason is vanishing of $\mathfrak{B}_{\mu,z_1,\ZZ}(\omega_\mu)$ ensures smoothness of the solution along $\mathfrak{g}_{\textbf{a}_1}^{+}$ and $\mathfrak{g}_{\textbf{a}_1}^{-}$, but not necessarily along $\mathfrak{g}_{\textbf{a}_1}^{0}$ (see Remark \ref{re:4}). However, for non-accessible systems, this is not sufficient for global smoothness.

 We use $\mathfrak{B}_{\mu,\textbf{u},\ZZ}$ for two purposes. Firstly, they play auxiliary roles in obtaining uniform estimates of $\textbf{u}\mathfrak{D}_{\mu, \tau,\ZZ}$ (see \eqref{for:64} of Proposition \ref{po:10}), which are essential for the iteration scheme. Secondly, they are used to obtain uniform estimates for the solution (see \eqref{for:28} of Lemma \ref{le:7} and Corollary \ref{cor:10})).

\end{enumerate}

\end{remark}

\begin{lemma}\label{le:12} Suppose $\xi\in J_\mu(\mathcal{O}^m)$, $m\geq1$ then:
\begin{enumerate}
\item \label{for:9} if $n\in L_\mu$,  we have
\begin{align*}
  \big\|\langle (\text{Ad}_{z_1}|&_{J_\mu})^{n-1}(\xi\circ \mathfrak{z}_1^{-n}),\,\psi\rangle\big\|
  \leq
  Ce^{-3\textbf{l}|n|}\norm{\xi}_{1}\norm{\psi}_{1}
\end{align*}
(see \eqref{for:432} of Section \ref{sec:3}) for any $\psi\in J_\mu(\mathcal{O}^1)$;

\smallskip


\item\label{for:421} $\mathcal{P}_{\mu}(\xi)$ (see \eqref{for:412} of Section \ref{sec:53}) is a distribution, satisfying
\begin{align*}
 \|\langle \mathcal{P}_{\mu}(\xi),\,\psi\rangle\|\leq C\norm{\xi}_{1}\norm{\psi}_{1}
\end{align*}
for any $\psi\in J_{\mu}(\mathcal{O}^{1})$;

\smallskip
\item\label{for:106} for any $u_i\in (\mathfrak{g}_{\mathfrak{z}_1}^{\delta(\mu)}\cup\mathfrak{g}_{\mathfrak{z}_1}^{0}) \bigcap \mathfrak{G}^1$ (see \eqref{for:273} of Section \ref{sec:3}), $1\leq i\leq r$, $0\leq r\leq m-1$, we have
\begin{align*}
 \|\langle u_1\cdots u_r\mathcal{P}_{\mu}(\xi),\,\psi\rangle\|\leq C_r\norm{\xi}_{r+1}\norm{\psi}_{1}
\end{align*}
for any $\psi\in \mathfrak{g}_{\mu}(\mathcal{O}^1)$;

\smallskip

\item \label{for:229} for any $n\in L_\mu^c$ and any poor row $(\mu,\lambda_1,\cdots, \lambda_k )$ (see \eqref{for:1} of Section \ref{sec:3}), $k\geq\sigma$ (see \eqref{for:432} of Section \ref{sec:3}) we have
\begin{align*}
&(\Pi_{i=1}^k\big\|\text{Ad}_{\mathfrak{z}_1^{-n}}|_{\mathfrak{g}_{\lambda_i,A}}\big\|)\big\|\mu(\textbf{a}_1)^{n-1}e^{\textbf{l}|n-1|}\notag\\
 &\leq C\big(\max_{\lambda\in \Delta_A}\exp(\big|\log\lambda(\textbf{a}_1)\big|)\big)^{-2|n|};
\end{align*}

\item \label{for:231} for any $\textbf{u}\in \mathcal{U}_{\mu,k}$ (see \eqref{for:1} of Section \ref{sec:3}), $k\geq\sigma$ we have:
   \begin{enumerate}
     \item\label{for:433} for any $n\in L_\mu$ and any $\psi\in J_\mu(\mathcal{O}^{k+1})$ we have
\begin{align*}
  \Big\|\big\langle (\text{Ad}_{z_1}|&_{J_\mu})^{n-1}(\xi\circ \mathfrak{z}_1^{-n}),\,\textbf{u}\psi\big\rangle\Big\|
  \leq Ce^{-3\textbf{l}|n|}\norm{\xi}_{1}\norm{\psi}_{k+1};
\end{align*}

     \item\label{for:434} if $m\geq k$, for any $n\in L_\mu^c$, we see that $\textbf{u}(\text{Ad}_{z_1}|_{J_\mu})^{n-1}(\xi\circ \mathfrak{z}_1^{-n})\in J_\mu(\mathcal{O})$ with estimate
     \begin{align*}
\big\|(\text{Ad}_{z_1}|_{J_\mu})^{n-1}\textbf{u}(\xi\circ \mathfrak{z}_1^{-n})\big\|
  \leq Ce^{-10\textbf{l}|n|}\norm{\xi}_{k};
\end{align*}

   \item\label{for:448} if $m\geq k$, for any $n\in \ZZ$ and any $\psi\in J_\mu(\mathcal{O}^{k+1})$ we have
     \begin{align*}
\Big\|\big\langle (\text{Ad}_{z_1}|&_{J_\mu})^{n-1}(\xi\circ \mathfrak{z}_1^{-n}),\,\textbf{u}\psi\big\rangle\Big\|
  \leq Ce^{-3\textbf{l}|n|}\norm{\xi}_{k}\norm{\psi}_{k+1};
\end{align*}
   \end{enumerate}

\item\label{for:131} if $m\geq\sigma$, for any $\textbf{u}\in \mathcal{U}_{\mu,\sigma}$ and any subset $\mathbb{S}\subseteq\ZZ$ we have
\begin{enumerate}
  \item for any $\psi\in J_{\mu}(\mathcal{O}^{\sigma+1})$, we have
  \begin{align*}
   \|\langle \mathfrak{B}_{\mu,\textbf{u},\mathbb{S}}(\xi),\,\psi\rangle\|\leq C\norm{\xi}_{\sigma}\norm{\psi}_{\sigma+1}
  \end{align*}
  (see \eqref{for:413} of Section \ref{sec:53}).  Moreover, we have
\begin{align*}
 \norm{\mathfrak{B}_{\mu,\mathbb{S}}(\xi)}\leq C\norm{\xi}_{\sigma}
\end{align*}
(see \eqref{for:414} of Section \ref{sec:53});

  \smallskip
  \item if $\mathbb{S}\subseteq L_\mu^c$, then $\mathfrak{B}_{\mu,\textbf{u},\mathbb{S}}(\xi)\in J_\mu(\mathcal{O})$ with estimates
  \begin{align*}
 \norm{\mathfrak{B}_{\mu,\textbf{u},\mathbb{S}}(\xi)}\leq C\norm{\xi}_{\sigma};
\end{align*}
\end{enumerate}

\smallskip
  \item\label{for:132} if equation \eqref{for:105} has a solution $\varphi_\mu\in J_\mu(\mathcal{O}^\sigma)$,  then
  \begin{align*}
   \norm{\mathfrak{B}_{\mu,\ZZ}(\omega_\mu)}=0.
  \end{align*}

    \end{enumerate}
\end{lemma}
\begin{proof}
\eqref{for:9}: For any
$n\in L_\mu$,  we have
\begin{align*}
 &\big\|\langle (\text{Ad}_{z_1}|_{J_\mu})^{n-1} (\xi\circ \mathfrak{z}_1^{-n}),\,\psi\rangle\big\|\notag\\
&\overset{\text{(1)}}{\leq} C\mu(\textbf{a}_1)^{n-1}e^{\emph{\textbf{l}}|n-1|}
 \big\|\langle \xi\circ \mathfrak{z}_1^{-n},\,\psi\rangle\big\|\\
 &\overset{\text{(2)}}{\leq} C\mu(\textbf{a}_1)^{n-1}e^{\emph{\textbf{l}}|n-1|}
 e^{-4 \emph{\textbf{l}}|n|}\norm{\xi}_{1}\norm{\psi}_{1}\notag\\
 &\overset{\text{(3)}}{\leq} C_1
 e^{-3 \emph{\textbf{l}}|n|}\norm{\xi}_{1}\norm{\psi}_{1}
\end{align*}
for any $\psi\in J_\mu(\mathcal{O}^1)$. Here in $(1)$ we use \eqref{for:108} of Proposition \ref{po:9}; in $(2)$ we use \eqref{for:438} of Section \ref{sec:54};
in $(3)$ we use the fact that $\mu(\textbf{a}_1)^{n}\leq1$ if $n\in L_\mu$ (see \eqref{for:439} of Section \ref{sec:54}).

\smallskip

\eqref{for:421}: Follows immediately from \eqref{for:9}.

\smallskip


\eqref{for:106}: For any $u=u_1\cdots u_r$ where $u_i\in \mathfrak{G}$, $1\leq i\leq r$ and any $n\in\ZZ$
we note that
\begin{align}
&u_1u_2\cdots u_r (\xi\circ \mathfrak{z}_1^{-n})=\norm{\text{Ad}_{\mathfrak{z}_1^{-n}}u_1}\cdots \norm{\text{Ad}_{\mathfrak{z}_1^{-n}}u_r}(\textbf{u}_{n}\xi)\circ \mathfrak{z}_1^{-n},
\end{align}
where
\begin{align*}
\textbf{u}_{n}=(\norm{\text{Ad}_{\mathfrak{z}_1^{-n}}u_1}\cdots \norm{\text{Ad}_{\mathfrak{z}_1^{-n}}u_r})^{-1}(\text{Ad}_{\mathfrak{z}_1^{-n}}u_1)\cdots (\text{Ad}_{\mathfrak{z}_1^{-n}}u_r).
\end{align*}
Then for any $u_i\in (\mathfrak{g}_{\mathfrak{z}_1}^{\delta(\mu)}\cup\mathfrak{g}_{\mathfrak{z}_1}^{0}) \bigcap \mathfrak{G}^1$, $1\leq i\leq r$ where $0\leq r\leq m-1$ and
any
$n\in L_\mu$,  we have
\begin{align}
 &\big\|\langle (\text{Ad}_{z_1}|_{J_\mu})^{n-1}u_1u_2\cdots u_r (\xi\circ \mathfrak{z}_1^{-n}),\,\psi\rangle\big\|\notag\\
 &=\norm{\text{Ad}_{\mathfrak{z}_1^{-n}}u_1}\cdots \norm{\text{Ad}_{\mathfrak{z}_1^{-n}}u_r}\big\|\langle(\text{Ad}_{z_1}|_{J_\mu})^{n-1}(\textbf{u}_{n}\xi)\circ \mathfrak{z}_1^{-n},\,\psi\rangle\big\|\notag\\
 &\leq \norm{\text{Ad}_{\mathfrak{z}_1^{-n}}u_1}\cdots \norm{\text{Ad}_{\mathfrak{z}_1^{-n}}u_r}\norm{(\text{Ad}_{z_1}|_{J_\mu})^{n-1}}\big\|\langle(\textbf{u}_{n}\xi)\circ \mathfrak{z}_1^{-n},\,\psi\rangle\big\|\notag\\
 &\overset{\text{(1)}}{\leq} C\big\|\text{Ad}_{\mathfrak{z}_1^{-n}}|_{\mathfrak{g}_{\mathfrak{z}_1}^{\delta(\mu)}\cup\mathfrak{g}_{\mathfrak{z}_1}^{0}}\big\|^r\mu(\textbf{a}_1)^{n-1}e^{\emph{\textbf{l}}|n-1|}
 \big\|\langle(\textbf{u}_{n}\xi)\circ \mathfrak{z}_1^{-n},\,\psi\rangle\big\|\notag\\
  &\overset{\text{(2)}}{\leq} C_r(\max_{\chi\in \Delta_{\mathfrak{z}_1}^{\delta(\mu)}\cup\Delta_{\mathfrak{z}_1}^{0}}\big\|\text{Ad}_{\mathfrak{z}_1^{-n}}|_{\mathfrak{g}_{\chi,\mathfrak{z}_1}}\big\|)^re^{\emph{\textbf{l}}|n-1|}
 \big\|\langle(\textbf{u}_{n}\xi)\circ \mathfrak{z}_1^{-n},\,\psi\rangle\big\|\notag\\
 &\overset{\text{(3)}}{\leq} C_{r,1}(\abs{n}+1)^{r\dim\mathfrak{G}}e^{\emph{\textbf{l}}|n-1|}
 \big\|\langle(\textbf{u}_{n}\xi)\circ \mathfrak{z}_1^{-n},\,\psi\rangle\big\|\notag\\
&\overset{\text{(4)}}{\leq} C_{r,2}(\abs{n}+1)^{r\dim\mathfrak{G}}e^{\emph{\textbf{l}}|n-1|}
 e^{-4 \emph{\textbf{l}}|n|}\norm{\textbf{u}_{n}\xi}_{1}\norm{\psi}_{1}\notag\\
 &\leq C_{r,3}(\abs{n}+1)^{r\dim\mathfrak{G}}
 e^{-3 \emph{\textbf{l}}|n|}\norm{\xi}_{r+1}\norm{\psi}_{1}\label{for:109}
\end{align}
for any $\psi\in J_\mu(\mathcal{O}^1)$. Here in $(1)$ we use \eqref{for:108} of Proposition \ref{po:9}; in $(2)$
we use the fact that $\mu(\textbf{a}_1)^{n}\leq1$ if $n\in L_\mu$ (see \eqref{for:439} of Section \ref{sec:54});
in $(3)$ we use \eqref{for:8} of Corollary \eqref{for:129} and note that $\chi(\mathfrak{z}_1)^{-n}\leq1$ if $n\in L_\mu$ for any $\chi\in \Delta_{\mathfrak{z}_1}^{\delta(\mu)}\cup\Delta_{\mathfrak{z}_1}^{0}$ (see \eqref{for:440} of Section \ref{sec:54}); in $(4)$ we use \eqref{for:438} of Section \ref{sec:54}.

Then \eqref{for:106} follows directly from \eqref{for:109}.

\smallskip

\eqref{for:229}: We have
\begin{align*}
 &(\Pi_{i=1}^k\big\|\text{Ad}_{\mathfrak{z}_1^{-n}}|_{\mathfrak{g}_{\lambda_i,A}}\big\|)\big\|\mu(\textbf{a}_1)^{n-1}e^{\emph{\textbf{l}}|n-1|}\notag\\
 &\overset{\text{(1)}}{\leq} C\big(\Pi_{i=1}^k (\lambda_i(\textbf{a}_1)^{-n}e^{\emph{\textbf{l}}|n|})\big)
 \mu(\textbf{a}_1)^{n-1}e^{\emph{\textbf{l}}|n-1|}\notag\\
 &\overset{\text{(2)}}{\leq} C\big(\Pi_{i=1}^k (\lambda_i(\textbf{a}_1)^{-\frac{4}{5}n})\big)\mu(\textbf{a}_1)^{n-1}e^{\emph{\textbf{l}}|n-1|}\notag\\
 &\overset{\text{(3)}}{\leq} C\big(\max_{\lambda\in \Delta_A}\exp(\big|\log\lambda(\textbf{a}_1)\big|)\big)^{-4|n|}\mu(\textbf{a}_1)^{n-1}e^{\emph{\textbf{l}}|n-1|}\notag\\
 &\overset{\text{(4)}}{\leq} C_1\big(\max_{\lambda\in \Delta_A}\exp(\big|\log\lambda(\textbf{a}_1)\big|)\big)^{-2|n|}.
\end{align*}
Here in $(1)$ we use \eqref{for:441} of Section \ref{sec:54};
in $(2)$ we note that $(\lambda_i(\textbf{a}_1))^{n}>1$, $1\leq i\leq k$ if $n\in L_\mu^c$ (see \eqref{for:439} of Section \ref{sec:54}) and we use \eqref{for:409} of Proposition \ref{po:9}; in $(3)$ we note that:
\begin{align*}
 \Pi_{i=1}^k (\lambda_i(\textbf{a}_1)^{-\frac{4}{5}n})&\overset{\text{(a)}}{=}\Pi_{i=1}^k\Big(\exp(\big|\log\lambda_i(\textbf{a}_1)\big|)\Big)^{-\frac{4}{5}|n|}\\
 &\leq
 \min_{1\neq\lambda\in \Delta_A}\Big(\exp(\big|\log\lambda(\textbf{a}_1)\big|)\Big)^{-\frac{4}{5}|n|k}\\
 &\overset{\text{(b)}}{\leq}
 (e^{5\emph{\textbf{l}}})^{-\frac{4}{5}|n|k}=(e^{\emph{\textbf{l}}k})^{-4|n|}\\
 &\overset{\text{(c)}}{\leq} \big(\max_{\lambda\in \Delta_A}\exp(\big|\log\lambda(\textbf{a}_1)\big|)\big)^{-4|n|}.
\end{align*}
Here in $(a)$ we recall that $(\lambda_i(\textbf{a}_1))^{n}>1$, $1\leq i\leq \sigma$ if $n\in L_\mu^c$; in $(b)$ we use \eqref{for:409} of Proposition \ref{po:9}; in $(c)$ we use \eqref{for:27} of Section \ref{sec:3} and recall $k\geq\sigma$; in $(4)$
we use \eqref{for:409} of Proposition \ref{po:9}.

Then we get the result.

\smallskip

\eqref{for:231}:  For any $n\in L_\mu$ and any $\psi\in J_\mu(\mathcal{O}^{\sigma+1})$ by \eqref{for:9} we have
\begin{align*}
  \big\|\langle (\text{Ad}_{z_1}|&_{J_\mu})^{n-1}\xi\circ \mathfrak{z}_1^{-n},\,\textbf{u}\psi\rangle\big\|
  \leq
  Ce^{-3\emph{\textbf{l}}|n|}\norm{\xi}_{1}\norm{\textbf{u}\psi}_{1}\notag\\
  &\leq Ce^{-3\emph{\textbf{l}}|n|}\norm{\xi}_{1}\norm{\psi}_{k+1}.
\end{align*}
Thus we get \eqref{for:433}.

Next, we prove \eqref{for:434}. We write $\textbf{u}=u_1\cdots u_k$,  $u_i\in \mathfrak{g}_{\lambda_i,A}\cap\mathfrak{G}^1$, $1\leq i\leq k$,  where $(\mu,\lambda_1,\cdots, \lambda_k)$ is a poor row.  For any $n\in\ZZ$ we denote
\begin{align*}
\textbf{u}_{n}=(\norm{\text{Ad}_{\mathfrak{z}_1^{-n}}u_1}\cdots \norm{\text{Ad}_{\mathfrak{z}_1^{-n}}u_\sigma})^{-1}(\text{Ad}_{\mathfrak{z}_1^{-n}}u_1)\cdots (\text{Ad}_{\mathfrak{z}_1^{-n}}u_k).
\end{align*}
For any $\psi\in J_\mu(\mathcal{O})$ and $n\in L_\mu^c$ we have
\begin{align}
 &\big\|\langle \textbf{u}(\text{Ad}_{z_1}|_{J_\mu})^{n-1}\xi\circ \mathfrak{z}_1^{-n},\,\psi\rangle\big\|\notag\\
  &=\norm{\text{Ad}_{\mathfrak{z}_1^{-n}}u_1}\cdots \norm{\text{Ad}_{\mathfrak{z}_1^{-n}}u_k}\big\|\langle(\text{Ad}_{z_1}|_{J_\mu})^{n-1}(\textbf{u}_{n}\xi)\circ \mathfrak{z}_1^{-n},\,\psi\rangle\big\|\notag\\
 &\leq \norm{\text{Ad}_{\mathfrak{z}_1^{-n}}u_1}\cdots \norm{\text{Ad}_{\mathfrak{z}_1^{-n}}u_k}\norm{(\text{Ad}_{z_1}|_{J_\mu})^{n-1}}\Big\| \big\langle(\textbf{u}_{n}\xi)\circ \mathfrak{z}_1^{-n},\,\psi\big\rangle\Big\|\notag\\
 &\overset{\text{(1)}}{\leq} C(\Pi_{i=1}^k\big\|\text{Ad}_{\mathfrak{z}_1^{-n}}|_{\mathfrak{g}_{\lambda_i,A}}\big\|)\cdot \big(\mu(\textbf{a}_1)^{n-1}e^{\emph{\textbf{l}}|n-1|}\big)\cdot
 \Big\|\big\langle(\textbf{u}_{n}\xi)\circ \mathfrak{z}_1^{-n},\,\psi\big\rangle\Big\|\notag\\
 &\overset{\text{(2)}}{\leq} C_1\big(\max_{\lambda\in \Delta_A}\exp(\big|\log\lambda(\textbf{a}_1)\big|)\big)^{-2|n|}\cdot
 \Big\|\big\langle(\textbf{u}_{n}\xi)\circ \mathfrak{z}_1^{-n},\,\psi\big\rangle\Big\|\notag\\
 &\overset{\text{(3)}}{\leq} C_1(e^{5\emph{\textbf{l}}})^{-2|n|}\cdot
 \big\|(\textbf{u}_{n}\xi)\circ \mathfrak{z}_1^{-n}\big\|\cdot\norm{\psi}\notag\\
 &\leq C_2e^{-10\emph{\textbf{l}}|n|}\norm{\xi}_{k}\norm{\psi}.\label{for:25}
\end{align}
Here in $(1)$ we use \eqref{for:441} of Section \ref{sec:54}; in $(2)$ we use \eqref{for:229}; in $(3)$ we use \eqref{for:409} of Proposition \ref{po:9}
and cauchy schwarz inequality.

Then  \eqref{for:434} follows from \eqref{for:25}.

\eqref{for:448} follows directly from \eqref{for:433} and \eqref{for:434}.

\smallskip

\eqref{for:131}: Follows directly from \eqref{for:231}.

\smallskip
\eqref{for:132}: \eqref{for:131} shows that both $\mathfrak{B}_{\mu,\textbf{u},\ZZ}(\varphi_\mu)$ and
$\mathfrak{B}_{\mu,\textbf{u},\ZZ}(\omega_\mu)$ are well defined. Applying the operator $\mathfrak{B}_{\mu,\textbf{u},\ZZ}$ on both sides of equation \eqref{for:105}
we have
\begin{align*}
0=\mathfrak{B}_{\mu,\textbf{u},\ZZ}(\varphi_\mu)-\mathfrak{B}_{\mu,\textbf{u},\ZZ}(\varphi_\mu)=\mathfrak{B}_{\mu,\textbf{u},\ZZ}(\omega_\mu).
\end{align*}
Then we get the result.

\end{proof}

\begin{lemma}\label{for:422} Suppose $\tau\geq1$ and $\xi\in J_\mu(\mathcal{O}^m)$, $m\geq \sigma$.  Then for  any $\textbf{u}\in \mathcal{U}_{\mu,\sigma}$ and any $v\in\mathfrak{g}^{-\delta(\mu)}_{\mathfrak{z}_1}\cap\mathfrak{G}^1$, $\textbf{u}\mathfrak{D}_{\mu,v,\tau,L_\mu^c}(\xi)\in J_{\mu}(\mathcal{O})$ (see \eqref{for:415} of Section \ref{sec:53}) with estimates
\begin{align*}
 \big\| \textbf{u}\mathfrak{D}_{\mu,v,\tau,L_\mu^c}(\xi)\big\|\leq C_m\tau^{-m+\sigma}\norm{\xi}_m.
\end{align*}
\end{lemma}
\begin{proof} We write $\textbf{u}=u_1u_2\cdots u_\sigma$,  $u_i\in \mathfrak{g}_{\lambda_i,A}\cap\mathfrak{G}^1$, $1\leq i\leq \sigma$. For any $n\in\ZZ$, we set
\begin{align*}
\textbf{u}_{n}=(\norm{\text{Ad}_{p(\mathfrak{z}_1)^{-n}}u_1}\cdots \norm{\text{Ad}_{p(\mathfrak{z}_1)^{-n}}u_\sigma})^{-1}(\text{Ad}_{p(\mathfrak{z}_1)^{-n}}u_1)\cdots (\text{Ad}_{p(\mathfrak{z}_1)^{-n}}u_\sigma)
\end{align*}
and
 \begin{align*}
 c_n=\norm{\text{Ad}_{p(\mathfrak{z}_1)^{-n}}(v)}^{-1},\quad v_n=c_n\text{Ad}_{p(\mathfrak{z}_1)^{-n}}(v)\,\,\text{ and }\,\,\xi_n=\pi(y_{\mathfrak{z}_1}^n)\xi
\end{align*}
(see \eqref{sect:11} of Section \ref{sec:3}).

\smallskip

\emph{Step 1: We show that}: for any
$n\in L_\mu^c$
\begin{align}\label{for:33}
 &\big\|\textbf{u}\pi_{v}(f_1\circ \tau^{-1})(\xi\circ \mathfrak{z}_1^{-n})\big\|\notag\\
 &\leq(\Pi_{i=1}^\sigma\norm{\text{Ad}_{p(\mathfrak{z}_1)^{-n}}u_i})\big\|\pi_{v_n}(f_1\circ (\tau c_n)^{-1})\xi_n\big\|_\sigma.
\end{align}
We have
 \begin{align*}
 &\big\|u_1\cdots u_\sigma\pi_{v}(f_1\circ \tau^{-1})(\xi\circ \mathfrak{z}_1^{-n})\big\|\\
 &\overset{\text{(1)}}{=}\big\|u_1\cdots u_\sigma\pi_{v}(f_1\circ \tau^{-1})\pi(\mathfrak{z}_1^{n})\xi\big\|\\
 &\overset{\text{(2)}}{=}\Big\|u_1\cdots u_\sigma\big(\pi_{v}(f_1\circ \tau^{-1})\pi(p(\mathfrak{z}_1)^n)\big)\big(\pi(y_{\mathfrak{z}_1}^n)\xi\big)\Big\|\\
 &\overset{\text{(3)}}{=}\Big\|u_1\cdots u_\sigma\big(\pi(p(\mathfrak{z}_1)^{n})\pi_{v_n}(f_1\circ (\tau c_n)^{-1})\big)\xi_n\Big\|\\
 &\overset{\text{(4)}}{=}(\Pi_{i=1}^\sigma\norm{\text{Ad}_{p(\mathfrak{z}_1)^{-n}}u_i})\big\|\pi(p(\mathfrak{z}_1)^{n})\textbf{u}_n\pi_{v_n}(f_1\circ (\tau c_n)^{-1})\xi_n\big\|\\
 &\overset{\text{(1)}}{=}(\Pi_{i=1}^\sigma\norm{\text{Ad}_{p(\mathfrak{z}_1)^{-n}}u_i})\big\|\textbf{u}_n\pi_{v_n}(f_1\circ (\tau c_n)^{-1})\xi_n\big\|\\
 &\leq (\Pi_{i=1}^\sigma\norm{\text{Ad}_{p(\mathfrak{z}_1)^{-n}}u_i})\big\|\pi_{v_n}(f_1\circ (\tau c_n)^{-1})\xi_n\big\|_\sigma.
 \end{align*}
 Here in $(1)$ we note that $\pi$ is the extended regular representation (see \eqref{for:17} of Section \ref{sec:3});
 in $(2)$  we recall \eqref{for:442} of Section \ref{sec:54};  in $(3)$ we use Lemma \ref{le:4}; in $(4)$ we note that
 \begin{align*}
  u_1\cdots u_\sigma\big(\pi(p(\mathfrak{z}_1)^{n})=(\Pi_{i=1}^\sigma\norm{\text{Ad}_{p(\mathfrak{z}_1)^{-n}}u_i})\pi(p(\mathfrak{z}_1)^{n})\textbf{u}_{n}.
 \end{align*}
Hence, we get the result.

 \smallskip

\noindent\emph{Step 2: We show that}:
 \begin{align}\label{for:227}
   \Pi_{i=1}^\sigma\norm{(\text{Ad}_{p(\mathfrak{z}_1)^{-n}}u_i)}\leq
  C(\abs{n}+1)^{\sigma\dim\mathfrak{G}}(\Pi_{i=1}^\sigma\big\|\text{Ad}_{\mathfrak{z}_1^{-n}}|_{\mathfrak{g}_{\lambda_i,A}}\big\|).
 \end{align}
We have
 \begin{align*}
 \Pi_{i=1}^\sigma&\norm{(\text{Ad}_{p(\mathfrak{z}_1)^{-n}}u_i)}=\Pi_{i=1}^\sigma\norm{(\text{Ad}_{y_{\mathfrak{z}_1}^{n}\mathfrak{z}_1^{-n}}u_i)}\notag\\
 &\leq \norm{\text{Ad}_{y_{\mathfrak{z}_1}^{n}}}^\sigma(\Pi_{i=1}^\sigma\norm{(\text{Ad}_{\mathfrak{z}_1^{-n}}u_i)})\notag\\
 &\overset{\text{(1)}}{\leq} C(\abs{n}+1)^{\sigma\dim\mathfrak{G}}(\Pi_{i=1}^\sigma\big\|\text{Ad}_{\mathfrak{z}_1^{-n}}|_{\mathfrak{g}_{\lambda_i,A}}\big\|).
\end{align*}
Here in $(1)$ we use \eqref{for:57} of Corollary \eqref{for:129}.

  \smallskip
  \noindent\emph{Step 3: We show that}: for any
$n\in L_\mu^c$
  \begin{gather}
    \big\|\pi_{v_n}(f_1\circ (\tau c_n)^{-1})\xi_n\big\|\leq C\tau^{-m}c_n^{-m}\norm{\xi_n}_m \qquad\text{and}\label{for:423}\\
    \big\|\pi_{v_n}(f_1\circ (\tau c_n)^{-1})\xi_n\big\|_m\leq C_m\norm{\xi_n}_m. \label{for:424}
  \end{gather}
\eqref{for:423}: It follows directly from Lemma \ref{cor:4}.

\smallskip
\noindent\eqref{for:424}:  We note that
$c_n>1$ for any $n\in L_\mu^c$ (see \eqref{for:440} of Section \ref{sec:54}). Then \eqref{for:424} is from \eqref{for:7} of Lemma \ref{cor:1}.

\smallskip
\noindent\emph{Step 4: We show that}: for any
$n\in L_\mu^c$
 \begin{align}\label{for:228}
  \big\|\pi_{v_n}(f_1\circ (\tau c_n)^{-1})\xi_n\big\|_\sigma\leq C_{m}\tau^{-m+\sigma}(\abs{n}+1)^{m\dim\mathfrak{G}}\norm{\xi}_m
 \end{align}
for any $m\geq\sigma$.

 We have
\begin{align*}
&\big\|\pi_{v_n}(f_1\circ (\tau c_n)^{-1})\xi_n\big\|_\sigma\notag\\
&\overset{\text{(1)}}{\leq} C_m\big\|\pi_{v_n}(f_1\circ (\tau c_n)^{-1})\xi_n\big\|^{1-\frac{\sigma}{m}}
(\big\|\pi_{v_n}(f_1\circ (\tau c_n)^{-1})\xi_n\big\|_m)^{\frac{\sigma}{m}}\notag\\
&\overset{\text{(2)}}{\leq}C_{m,1}\big( \tau^{-m}c_n^{-m}\norm{\xi_n}_m\big)^{1-\frac{\sigma}{m}}(\norm{\xi_n}_m)^{\frac{\sigma}{m}}\notag\\
&\overset{\text{(3)}}{\leq}
C_{m,1}\tau^{-m+\sigma}\norm{\xi_n}_m\\
&\leq C_{m,2}\tau^{-m+\sigma}\norm{\text{Ad}_{y_{\mathfrak{z}_1}^{n}}}^m\norm{\xi}_m\notag\\
&\overset{\text{(4)}}{\leq}C_{m,3}\tau^{-m+\sigma}(\abs{n}+1)^{m\dim\mathfrak{G}}\norm{\xi}_m.\notag
\end{align*}
Here in $(1)$ we use interpolation inequalities; in $(2)$ we use \eqref{for:423} and \eqref{for:424}; in $(3)$ we recall $c_n>1$; in $(4)$ we use
\eqref{for:57} of Corollary \eqref{for:129}.

\smallskip

  \noindent\emph{Step 5: We show that}: for any
$n\in L_\mu^c$
\begin{align}\label{for:425}
 &\big\|(\text{Ad}_{z_1}|_{J_\mu})^{n-1}\textbf{u}\pi_{v}(f_1\circ \tau^{-1})(\xi\circ \mathfrak{z}_1^{-n})\big\|\notag\\
 &\leq C_{m}\big(\max_{\lambda\in \Delta_A}\exp(\big|\log\lambda(\textbf{a}_1)\big|)\big)^{-|n|}\tau^{-m+\sigma}\norm{\xi}_m
\end{align}
for any $m\geq\sigma$.

We have
 \begin{align*}
 &\big\|(\text{Ad}_{z_1}|_{J_\mu})^{n-1}\textbf{u}\pi_{v}(f_1\circ \tau^{-1})(\xi\circ \mathfrak{z}_1^{-n})\big\|\notag\\
 &\leq\norm{(\text{Ad}_{z_1}|_{J_\mu})^{n-1}}\big\|\textbf{u}\pi_{v}(f_1\circ \tau^{-1})(\xi\circ \mathfrak{z}_1^{-n})\big\|\notag\\
 &\overset{\text{(1)}}{\leq} C\mu(\textbf{a}_1)^{n-1}e^{\emph{\textbf{l}}(n-1)}\cdot (\Pi_{i=1}^\sigma\norm{\text{Ad}_{p(\mathfrak{z}_1)^{-n}}u_i})\big\|\pi_{v_n}(f_1\circ (\tau c_n)^{-1})\xi_n\big\|_\sigma\notag\\
 &\overset{\text{(2)}}{\leq} C\mu(\textbf{a}_1)^{n-1}e^{\emph{\textbf{l}}(n-1)}\cdot C(\abs{n}+1)^{\sigma\dim\mathfrak{G}}(\Pi_{i=1}^\sigma\big\|\text{Ad}_{\mathfrak{z}_1^{-n}}|_{\mathfrak{g}_{\lambda_i,A}}\big\|)\\
 &\quad\cdot C_{m}\tau^{-m+\sigma}(\abs{n}+1)^{m\dim\mathfrak{G}}\norm{\xi}_m\notag\\
 &\overset{\text{(3)}}{\leq}  C_1(\abs{n}+1)^{\sigma\dim\mathfrak{G}}\big(\max_{\lambda\in \Delta_A}\exp(\big|\log\lambda(\textbf{a}_1)\big|)\big)^{-2|n|}\\
 &\quad\cdot C_{m}\tau^{-m+\sigma}(\abs{n}+1)^{m\dim\mathfrak{G}}\norm{\xi}_m\\
 &\leq C_{m,1}\tau^{-m+\sigma}\big(\max_{\lambda\in \Delta_A}\exp(\big|\log|\lambda(\textbf{a}_1)|\big|)\big)^{-|n|}\norm{\xi}_m.
\end{align*}
Here in $(1)$ we use \eqref{for:108} of Proposition \ref{po:9} and \eqref{for:33}; in $(2)$ we use \eqref{for:227} and \eqref{for:228};
in $(3)$ we recall that $(\mu,\lambda_1,\cdots, \lambda_\sigma)$ is a poor row and $n\in L_\mu^c$. Then we use \eqref{for:229} of Lemma \ref{le:12}.

Then we finish the proof.

\smallskip

\emph{Conclusion}:  We point out that \eqref{for:425} implies
\begin{align*}
 \big\| \textbf{u}\mathfrak{D}_{\mu,v,\tau,L_\mu^c}(\xi)\big\|\leq C_m\tau^{-m+\sigma}\norm{\xi}_m.
\end{align*}
Hence we finish the proof of Lemma \ref{for:422}.
\end{proof}

The next result studies the commutator relation between $\pi_{v}$ and the directional derivatives, which will be used in the proof of Corollary \ref{cor:2}, as well as in the subsequent part.

\begin{proposition}\label{po:10} Suppose $\xi\in J_\mu(\mathcal{O}^m)$. For any $l\in\NN$, any $\tau\geq1$ and any $\textbf{u}\in \mathcal{U}_{\mu,\sigma}$, the following hold:

\begin{enumerate}

\item\label{re:2}  for any $g_i\in \tilde{\mathcal{S}}(\RR)$ (see \eqref{for:239} of Section \ref{sec:36}) and any $v_i\in\mathfrak{g}^{-\delta(\mu)}_{\mathfrak{z}_1}\cap\mathfrak{G}^1$, $1\leq i\leq l$, we have
\begin{gather*}
 \pi_{v_1}(g_1\circ \tau^{-1})\cdots\pi_{v_l}(g_l\circ \tau^{-1})\textbf{u}\\
 =\sum_i k_{i}\tau^{-\sum_{n=1}^lj_{n,i}}\textbf{u}_{i}\pi_{v_1}(q_{j_{1,i}}\circ \tau^{-1})\cdots \pi_{v_l}(q_{j_{l,i}}\circ \tau^{-1})
\end{gather*}
and
\begin{gather*}
 \textbf{u}\pi_{v_1}(g_1\circ \tau^{-1})\cdots\pi_{v_l}(g_l\circ \tau^{-1})\\
 =\sum_i k_{i}'\tau^{-\sum_{n=1}^lj'_{n,i}}\pi_{v_1}(p_{(j'_{1,i})}\circ \tau^{-1})\cdots \pi_{v_l}(p_{(j'_{l,i})}\circ \tau^{-1})\textbf{u}_{i}'
\end{gather*}
 where $\textbf{u}_{i},\,\textbf{u}_{i}'\in \mathcal{U}_{\mu,\sigma}$, $j_{n,i},\,j_{n,i}'
 \geq0$, $q_{j_{n,i}},\,p_{j'_{n,i}}\in \tilde{\mathcal{S}}(\RR)$  and $k_{i},\,k_{i}'\in\CC$;

\smallskip

  \item\label{for:64} if $m\geq \sigma+1$, then for any $v\in\mathfrak{g}^{-\delta(\mu)}_{\mathfrak{z}_1}\cap\mathfrak{G}^1$ and any $\psi \in J_\mu(\mathcal{O}^{\sigma+1})$ we have
\begin{align*}
 \big\|\langle \textbf{u}\mathfrak{D}&_{\mu,v,\tau,\ZZ}(\xi),\,\psi\rangle\big\|\leq C\norm{\mathfrak{B}_{\mu,\ZZ}(\xi)}\norm{\psi}_{\sigma+1}.
\end{align*}
\end{enumerate}

\end{proposition}
\begin{proof}

\smallskip

\eqref{re:2}: By using \eqref{for:405} of Corollary \ref{for:129}, for any $v\in\mathfrak{g}^{-\delta(\mu)}_{\mathfrak{z}_1}$ and any $u\in\mathfrak{g}^{-\delta(\mu)}_{\textbf{a}_1}$, we see that $(*)$:
\begin{align*}
 \text{ad}^{j}_{v}(u)\in \mathfrak{g}^{-\delta(\mu)}_{\textbf{a}_1},\qquad \forall\,j\in\NN.
\end{align*}
Keeping using $(*)$ and \eqref{for:234} (resp. \eqref{for:98}) of Section \ref{sec:27},  for any $\tau\geq1$, any $g\in \tilde{\mathcal{S}}(\RR)$,  any $\textbf{u}\in \mathcal{U}_{\mu,\sigma}$ and any $v\in\mathfrak{g}^{-\delta(\mu)}_{\mathfrak{z}_1}\cap\mathfrak{G}^1$, we have
\begin{align*}
 &\pi_{v}(g\circ \tau^{-1})\textbf{u}=\sum_i k_i\tau^{-j_{i}}\textbf{u}_{i}\pi_{v}(g_{j_{i}}\circ \tau^{-1})\\
 \big(\text{resp. }&\textbf{u}\pi_{v}(g\circ \tau^{-1})=\sum_i k_{i}'\tau^{-l_{i}}\pi_{v}(g_{l_{i}}\circ \tau^{-1})\textbf{u}_{i}'\big)
\end{align*}
where $\textbf{u}_{i},\,\textbf{u}_{i}'\in \mathcal{U}_{\mu,\sigma}$, $j_{i},\,l_i\geq0$, $g_{j_{i}},\,g_{l_{i}}\in \tilde{\mathcal{S}}(\RR)$ and $k_{i},\,k_{i}'\in\CC$.

Keeping using the above equation, we get the result.

\smallskip
\eqref{for:64}: We recall \eqref{for:415} of Section \ref{sec:53}. For any $v\in\mathfrak{g}^{-\delta(\mu)}_{\mathfrak{z}_1}\cap\mathfrak{G}^1$ and any $\psi \in J_\mu(\mathcal{O}^{\sigma+1})$ we have
\begin{align*}
 &\big\|\langle \textbf{u}\mathfrak{D}_{\mu,v,\tau,\ZZ}(\xi),\,\psi\rangle\big\|\\
 &\overset{\text{(1)}}{=}\Big\|
 \big\langle \sum_{n\in \ZZ} \textbf{u}\pi_{v}(f_1\circ \tau^{-1})(\text{Ad}_{z_1}|_{J_\mu})^{n-1}(\xi\circ \mathfrak{z}_1^{-n}),\,\psi\big\rangle\Big\|\\
 &\overset{\text{(2)}}{=}\Big\|
 \sum_i k_i\tau^{-j_i}\big\langle \sum_{n\in \ZZ} \pi_{v}(f_{j_i}\circ \tau^{-1})\textbf{u}_i(\text{Ad}_{z_1}|_{J_\mu})^{n-1}(\xi\circ \mathfrak{z}_1^{-n}),\,\psi\big\rangle\Big\|\\
 &\overset{\text{(3)}}{=}\Big\|
 \sum_i k_i\tau^{-j_i}\big\langle \mathfrak{B}_{\mu,\textbf{u}_i,\ZZ}(\xi),\,\pi_{v}(f_{j_i}\circ \tau^{-1})\psi\big\rangle\Big\|\\
 &\overset{\text{(4)}}{\leq} C\sum_i\norm{\mathfrak{B}_{\mu,\ZZ}(\xi)}\norm{\pi_{v}(f_{j_i}\circ \tau^{-1})\psi}_{\sigma+1}\\
   &\overset{\text{(5)}}{\leq}C_1\norm{\mathfrak{B}_{\mu,\ZZ}(\xi)}\norm{\psi}_{\sigma+1}.
\end{align*}
Here in $(1)$ we use \eqref{for:73} of Section \ref{sec:27}; in $(2)$ we use \eqref{re:2} where $\textbf{u}_{i}\in \mathcal{U}_{\mu,\sigma}$, $j_{i}
 \geq0$, $f_{j_{i}}\in \tilde{\mathcal{S}}(\RR)$; in $(3)$ we use \eqref{for:93} of Section \ref{sec:27} and recall \eqref{for:413} of
 Section \ref{sec:53}; in $(4)$ we recall \eqref{for:414} of
 Section \ref{sec:53}; in $(5)$ we use \eqref{for:7} of Lemma \ref{cor:1}.

 Hence, we get the result.

\end{proof}

Next, we show that vanishing of $\mathfrak{B}_{\mu,\ZZ}(\omega_\mu)$ implies solvability of \eqref{for:5} in some Sobolev spaces.  Tame and uniform Sobolev estimates of the solution are also formulated.
\begin{lemma}\label{le:7} Suppose $\omega_\mu\in J_\mu(\mathcal{O}^m)$, $m\geq1$. Then:
\begin{enumerate}
\item\label{for:431} $\mathcal{P}_{\mu}(\omega_\mu)$ is a distributional solution of equation \eqref{for:105};

\smallskip

\item\label{for:457} for any $u_i\in (\mathfrak{g}_{\mathfrak{z}_1}^{\delta(\mu)}\cup\mathfrak{g}_{\mathfrak{z}_1}^{0}) \bigcap \mathfrak{G}^1$, $1\leq i\leq r$, $0\leq r\leq m-1$, we have
\begin{align*}
 \|\langle u_1\cdots u_r\mathcal{P}_{\mu}(\xi),\,\psi\rangle\|\leq C_r\norm{\xi}_{r+1}\norm{\psi}_{1}
\end{align*}
for any $\psi\in \mathfrak{g}_{\mu}(\mathcal{O}^1)$

\smallskip

  \item \label{for:120} if equation \eqref{for:105} has a solution $\varphi_\mu\in J_\mu(\mathcal{O}^1)$,  then $\varphi_\mu$ is unique;

  \smallskip
  \item \label{for:28} for any $m\geq \sigma+6$,  if $z_1=w\textbf{a}_1$ where $w\in \exp(\mathfrak{g}_{1,A})\cap B_{\eta_m(\textbf{a}_1)}(e)\cap Z(L)$ (see \eqref{for:426} of Corollary \ref{for:129}) and if $\mathfrak{B}_{\mu,z_1,\ZZ}(\omega_\mu)=0$, then
equation \eqref{for:105} has a solution
$\varphi_\mu\in J_{\mu}(\mathcal{O}^{m-6})$ with estimates
\begin{align*}
 \norm{\varphi_\mu}_{m-6}\leq C_{m}\norm{\omega_\mu}_{m}.
\end{align*}

\end{enumerate}
\end{lemma}
\begin{remark}\label{re:4} From the proof, we see that the assumption $\mathfrak{B}_{\mu,z_1,\ZZ}(\omega_\mu)=0$ implies that $\varphi_\mu$ is $C^\infty$ along $\mathfrak{g}_{\textbf{a}_1}^{+}$ and $\mathfrak{g}_{\textbf{a}_1}^{-}$ if $\omega_\mu\in J_\mu(\mathcal{O}^\infty)$. However, it does not provide any information regarding the differentiability of $\varphi_\mu$ along $\mathfrak{g}_{\textbf{a}_1}^{0}$. Instead,  the condition that $w\in B_{\eta_m(\textbf{a}_1)}(e)$ is what truly assists. Of course, to
improve the regularity of $\varphi_\mu$, we have to make $w$ closer to $e$.

\end{remark}
\begin{proof}
\eqref{for:431}: It follows directly from \eqref{for:421} of Lemma \ref{le:12}.

\smallskip

\eqref{for:457}:   It follows directly from \eqref{for:106} of Lemma \ref{le:12}.

\eqref{for:120}: It suffices to show that: if $\xi\in J_\mu(\mathcal{O}^1)$ satisfies the equation
\begin{align}\label{for:261}
 \xi\circ \mathfrak{z}_1-(\text{Ad}_{z_1}|_{J_\mu})\xi=0,
\end{align}
then $\xi=0$. From \eqref{for:261} we have
\begin{align*}
 \xi\circ \mathfrak{z}_1^n=(\text{Ad}_{z_1}|_{J_\mu})^n\xi,\qquad \forall\,n\in\NN.
\end{align*}
For any $n\in L_\mu$ the above equation implies that
\begin{align*}
 \|\langle \xi,\,\xi\rangle\|=\big\|\langle (\text{Ad}_{z_1}|_{J_\mu})^{-n}(\xi\circ \mathfrak{z}_1^n),\,\xi\rangle\big\|
 \overset{\text{(1)}}{\leq}
  Ce^{-3\emph{\textbf{l}}|n|}\norm{\xi}_{1}^2
\end{align*}
Here in $(1)$ we use \eqref{for:9} of Lemma \ref{le:12}.

Letting $|n|\to\infty$, we have $\langle \xi,\,\xi\rangle=0$. Then we get $\xi=0$.

\smallskip

\eqref{for:28}: Firstly, we show differentiability of $\xi$ along $\mathfrak{g}_{\textbf{a}_1}^{\delta(\mu)}$ (see \eqref{for:273} of Section \ref{sec:3}).  \eqref{for:405} of Corollary \ref{for:129} shows that $\Delta_{\textbf{a}_1}^{+}\subset \mathfrak{g}^+_{\mathfrak{z}_1}$ and $\Delta_{\textbf{a}_1}^{-}\subset \mathfrak{g}^-_{\mathfrak{z}_1}$.
This means that $\mathfrak{g}_{\textbf{a}_1}^{\delta(\mu)}\subseteq \mathfrak{g}^{\delta(\mu)}_{\mathfrak{z}_1}$.
It follows from \eqref{for:106} of Lemma \ref{le:12} that
for any $u_i\in \mathfrak{g}_{\textbf{a}_1}^{\delta(\mu)}\bigcap \mathfrak{G}^1$, $1\leq i\leq m-1$,  we have
\begin{align}\label{for:267}
 \big\|\langle u_1u_2\cdots u_{m-1}(\mathcal{P}_{\mu}(\omega_\mu)),\,\psi\rangle \big\|\leq C_m\norm{\omega_\mu}_{m}\norm{\psi}_{1}
\end{align}
for any $\psi\in J_\mu(\mathcal{O}^1)$.

Next we show differentiability of $\xi$ along $\mathfrak{g}_{1,A}=\mathfrak{g}_{\textbf{a}_1}^{0}$ (see \eqref{for:11} of Section \ref{sec:3}). For any $u_i\in \mathfrak{g}_{1,A}\bigcap \mathfrak{G}^1$, $1\leq i\leq m-1$ and any $n\in L_\mu$ we have
\begin{align*}
 &\big\|\langle u_1u_2\cdots u_{m-1}(\text{Ad}_{z_1}|_{J_\mu})^{n-1} (\omega_\mu\circ \mathfrak{z}_1^{-n}),\,\psi\rangle\big\|\notag\\
 &\overset{\text{(1)}}{=}\norm{\text{Ad}_{\mathfrak{z}_1^{-n}}u_1}\cdots \norm{\text{Ad}_{\mathfrak{z}_1^{-n}}u_{m-1}}\big\|\langle(\text{Ad}_{z_1}|_{J_\mu})^{n-1}(\textbf{u}_{n}\omega_\mu)\circ \mathfrak{z}_1^{-n},\,\psi\rangle\big\|\notag\\
 &\leq \norm{\text{Ad}_{\mathfrak{z}_1^{-n}}u_1}\cdots \norm{\text{Ad}_{\mathfrak{z}_1^{-n}}u_r}\norm{(\text{Ad}_{z_1}|_{J_\mu})^{n-1}}\Big\| \big\langle(\textbf{u}_{n}\omega_\mu)\circ \mathfrak{z}_1^{-n},\,\psi\big\rangle\Big\|\notag\\
 &\overset{\text{(2)}}{\leq} C\big\|\text{Ad}_{\mathfrak{z}_1^{-n}}|_{\mathfrak{g}_{1,A}}\big\|^{m-1}\mu(\textbf{a}_1)^{n-1}e^{\emph{\textbf{l}}|n-1|}
 \big\|\langle(\textbf{u}_{n}\omega_\mu)\circ \mathfrak{z}_1^{-n},\,\psi\rangle\big\|\\
 &\overset{\text{(3)}}{\leq} Ce^{\frac{|n|(m-1)\emph{\textbf{l}}}{m}}(\abs{n}+1)^{(m-1)\dim\mathfrak{G}}e^{\emph{\textbf{l}}|n-1|}
 \big\|\langle(\textbf{u}_{n}\omega_\mu)\circ \mathfrak{z}_1^{-n},\,\psi\rangle\big\|\notag\\
&\overset{\text{(4)}}{\leq} C_1e^{|n|\emph{\textbf{l}}}(\abs{n}+1)^{(m-1)\dim\mathfrak{G}}e^{\emph{\textbf{l}}|n-1|}
 e^{-4 \emph{\textbf{l}}|n|}\norm{\textbf{u}_{n}\omega_\mu}_{1}\norm{\psi}_{1}\notag\\
 &\leq C_2(\abs{n}+1)^{(m-1)\dim\mathfrak{G}}
 e^{-2 \emph{\textbf{l}}|n|}\norm{\omega_\mu}_{m}\norm{\psi}_{1}
\end{align*}
for any $\psi\in J_\mu(\mathcal{O}^1)$.

Here in $(1)$
\begin{align*}
\textbf{u}_{n}=(\norm{\text{Ad}_{\mathfrak{z}_1^{-n}}u_1}\cdots \norm{\text{Ad}_{\mathfrak{z}_1^{-n}}u_{m-1}})^{-1}(\text{Ad}_{\mathfrak{z}_1^{-n}}u_1)\cdots (\text{Ad}_{\mathfrak{z}_1^{-n}}u_{m-1});
\end{align*}
in $(2)$ we use \eqref{for:108} of Proposition \ref{po:9}; in $(3)$ we use \eqref{for:426} of Corollary \ref{for:129} and note that $\mu(\textbf{a}_1)^{n}\leq1$ if $n\in L_\mu$ (see \eqref{for:439} of Section \ref{sec:54}); in $(4)$ we use \eqref{for:438} of Section \ref{sec:54}.

Then it follows that: for any $u_i\in \mathfrak{g}_{1,A}\bigcap \mathfrak{G}^1$, $1\leq i\leq m-1$
\begin{align}\label{for:268}
 \big\|\langle u_1u_2\cdots u_{m-1}(\mathcal{P}_{\mu}(\omega_\mu)),\,\psi\rangle\big\|\leq C_m\norm{\omega_\mu}_{m}\norm{\psi}_{1}
\end{align}
for any $\psi\in J_\mu(\mathcal{O}^1)$.

Finally,  we show differentiability of $\xi$ along $\mathfrak{g}_{\textbf{a}_1}^{-\delta(\mu)}$. For any $u_i\in \mathfrak{g}_{\lambda_i,A} \bigcap \mathfrak{G}^1$ where $\lambda_i\in\Delta_{\textbf{a}_1}^{-\delta(\mu)}$, $1\leq i\leq m$,  we have
\begin{align}\label{for:264}
&\big\|\langle u_1u_2\cdots u_{m}(\mathcal{P}_{\mu}(\omega_\mu)),\,\psi\rangle \big\|\notag\\
&=\big\|\langle u_{m-\sigma+1}\cdots u_{m}(\mathcal{P}_{\mu}(\omega_\mu)),\,u_{m-\sigma}\cdots u_1\psi\rangle\big\|\notag\\
&\overset{\text{(1)}}{=}\big\|\langle \mathfrak{B}_{\mu,z_1,\textbf{u},L_\mu}(\omega_\mu),\,u_{m-\sigma}\cdots u_1\psi\rangle\big\|\notag\\
&\overset{\text{(2)}}{=}\big\|\langle \mathfrak{B}_{\mu,z_1,\textbf{u},L_\mu^c}(\omega_\mu),\,u_{m-\sigma}\cdots u_1\psi\rangle\big\|\notag\\
&=\big\|\langle u_1\cdots u_{m-\sigma}\mathfrak{B}_{\mu,z_1,\textbf{u},L_\mu^c}(\omega_\mu),\, \psi\rangle\big\|\notag\\
&=\Big\|\big\langle \sum_{n\in L_\mu^c}u_1u_2\cdots u_{m}(\text{Ad}_{z_1}|_{J_\mu})^{n-1}(\omega_\mu\circ \mathfrak{z}_1^{-n}),\, \psi\big\rangle\Big\|
\end{align}
for any $\psi\in J_\mu(\mathcal{O}^\infty)$.

Here in $(1)$ $\textbf{u}=u_{m-\sigma+1}\cdots u_{m}$. Since $\lambda_i\in\Delta_{\textbf{a}_1}^{-\delta(\mu)}$, $1\leq i\leq m$, $\textbf{u}\in\mathcal{U}_{\mu,\sigma}$ (see \eqref{for:439} of Section \ref{sec:54}). In $(2)$ we use $\mathfrak{B}_{\mu,\textbf{u},\ZZ}(\omega_\mu)=0$, which is implied by the assumption $\mathfrak{B}_{\mu,\ZZ}(\omega_\mu)=0$.

Let $\textbf{v}=u_1u_2\cdots u_{m}$. Since $\lambda_i\in\Delta_{\textbf{a}_1}^{-\delta(\mu)}$, $1\leq i\leq m$, $\textbf{v}\in\mathcal{U}_{\mu,m}$ (see \eqref{for:439} of Section \ref{sec:54}). Hence, we have
\begin{align}\label{for:265}
 &\Big\|\big\langle \sum_{n\in L_\mu^c}\textbf{v}(\text{Ad}_{z_1}|_{J_\mu})^{n-1}(\omega_\mu\circ \mathfrak{z}_1^{-n}),\, \psi\big\rangle\Big\|\notag\\
 &\overset{\text{(1)}}{\leq} \sum_{n\in L_\mu^c}\big\|\textbf{v}(\text{Ad}_{z_1}|_{J_\mu})^{n-1}(\omega_\mu\circ \mathfrak{z}_1^{-n})\big\| \norm{\psi}\notag\\
 &\overset{\text{(2)}}{\leq} \sum_{n\in L_\mu^c}Ce^{-10\emph{\textbf{l}}|n|}\norm{\omega_\mu}_{m} \norm{\psi}\notag\\
 &\leq C_1\norm{\omega_\mu}_{m} \norm{\psi}.
\end{align}
Here in $(1)$ we use cauchy schwarz inequality; in $(2)$ we use \eqref{for:434} of Lemma \ref{le:12}.

It follows from \eqref{for:264} and \eqref{for:265} that
\begin{align*}
 \big\|\langle \textbf{v}(\mathcal{P}_{\mu}(\omega_\mu)),\,\psi\rangle \big\|\leq C\norm{\omega_\mu}_{m}\norm{\psi}
\end{align*}
for any $\psi\in J_\mu(\mathcal{O}^\infty)$, which gives $\textbf{v}(\mathcal{P}_{\mu}(\omega_\mu))\in J_\mu(\mathcal{O})$ with estimates
\begin{align}\label{for:266}
 \|\textbf{v}(\mathcal{P}_{\mu}(\omega_\mu))\|\leq C\norm{\omega_\mu}_{m}.
\end{align}
We note that $\mathfrak{g}_{\textbf{a}_1}^{\delta(\mu)}$, $\mathfrak{g}_{\textbf{a}_1}^{0}$ and $\mathfrak{g}_{\textbf{a}_1}^{-\delta(\mu)}$ span a full set of directions.
It follows from \eqref{for:267}, \eqref{for:268} and \eqref{for:266} and Lemma \ref{le:3} that $\xi\in J_\mu(\mathcal{O}^{m-6})$ with the estimate
\begin{align*}
 \norm{\mathcal{P}_{\mu}(\omega_\mu)}_{m-6}\leq C_{m}\norm{\omega_\mu}_{m}.
\end{align*}
Since $m-6\geq \sigma>1$, $\mathcal{P}_{z_1,\mu}(\omega_\mu)\in J_\mu(\mathcal{O}^1)$. By \eqref{for:120} we see that $\mathcal{P}_{\mu}(\omega_\mu)=\varphi_\mu$. Hence we finish the proof.
\end{proof}
The next result follows from \eqref{for:132} of Lemma \ref{le:12} and \eqref{for:28} of Lemma \ref{le:7}:
\begin{corollary}\label{cor:10}  Suppose  $\omega_\mu\in \mathcal{L}^\bot(\mathcal{O}^m)$, $m\geq \sigma+6$ and $z_1=w\textbf{a}_1$ where $w\in \exp(\mathfrak{g}_{1,A})\cap B_{\eta_m(\textbf{a}_1)}(e)\cap Z(L)$. If equation \eqref{for:5} has a solution $\varphi\in \mathcal{L}^\bot(\mathcal{O}^1)$, then $\varphi\in \mathcal{L}^\bot(\mathcal{O}^{m-6})$ with the estimate
\begin{align*}
 \norm{\varphi}_{m-6}\leq C_{m}\norm{\omega}_{m}.
\end{align*}
\end{corollary}

\subsection{The almost twisted cocycle equation}
In this section, we assume $z_i=w_i\textbf{a}_i$ where $w_i\in \exp(\mathfrak{g}_{1,A})\cap B_{\eta}(e)\cap Z(L)$, $i=1,2$. The positive constants $\eta, \sigma,\,c_0$ (see Section \ref{for:432}) will appear in is section.

 We study the twisted equation:
\begin{align}\label{for:443}
\varphi\circ z_1-\mathfrak{p}_{\mathcal{L}^\bot}\text{Ad}_{z_1}\varphi=\theta\circ z_2-\mathfrak{p}_{\mathcal{L}^\bot}\text{Ad}_{z_2}\theta+\omega,
\end{align}
where $\varphi,\,\,\theta,\,\omega\in \mathcal{L}^\bot(\mathcal{O})$.

We recall $q(z_1)$ and $\mathfrak{z}_i=qz_iq^{-1}$, $i=1,\,2$ as described in Corollary \ref{for:129}. To study equation \eqref{for:443}, it  is equivalent to study
the twisted equation:
\begin{align}\label{for:13}
\varphi_q\circ \mathfrak{z}_1-\mathfrak{p}_{\mathcal{L}^\bot}\text{Ad}_{z_1}\varphi_q=\theta_q\circ \mathfrak{z}_2-\mathfrak{p}_{\mathcal{L}^\bot}\text{Ad}_{z_2}\theta_q+\omega_q,
\end{align}
where $\varphi_q=\varphi\circ q^{-1}$, $\omega_q=\omega\circ q^{-1}$ and $\theta_q=\theta\circ q^{-1}$.

\smallskip
\emph{Note}. Since $q$ is uniformly bounded, passing to equation \eqref{for:13} will not affect the results in this part.

\smallskip

From \eqref{for:54} of Section \ref{sec:3}, we see that \eqref{for:13} decomposes into equations of the following forms:
\begin{align}\label{for:15}
\varphi_\mu\circ z_1-(\text{Ad}_{z_1}|_{J_\mu})\varphi_\mu=\theta_\mu\circ z_2-(\text{Ad}_{z_2}|_{J_\mu})\theta_\mu+\omega_\mu, \quad\mu\in\Delta_A
\end{align}
where $\varphi_\mu=\varphi\circ q^{-1}|_{J_\mu}$,  $\theta_\mu=\theta\circ q^{-1}|_{J_\mu}$ and $\omega_\mu=\omega\circ q^{-1}|_{J_\mu}$.

In Proposition \ref{po:1}  we estimate the norm of invariant distributions $\mathfrak{B}_{\mu,\textbf{u},\ZZ}$. However, the estimate
plays an auxiliary role: it will not be used in subsequent parts.
We use the norm of $\mathfrak{B}_{\mu,\textbf{u},\ZZ}$ to obtain uniform estimate for the norm of invariant distributions $\textbf{u}\mathfrak{D}_{\mu,v,\tau,\ZZ}$ in Corollary \ref{cor:2}. This is the only reason we introduce the auxiliary norm. As we explained in Section \ref{sec:32},
even though we know that the norm of $\mathfrak{D}_{\mu,v,\tau,\ZZ}$ is important in the construction of approximate solution, it seems unlikely to obtain a uniform estimate for all $z_1$ in a small neighbourhood of $\textbf{a}_1$ by following the arguments in Section \ref{sec:13}. So we instead estimate the norms of $\textbf{u}\mathfrak{D}_{\mu,v,\tau,\ZZ}$ using directional derivatives, which are uniformly bounded.

\begin{proposition}\label{po:1} Suppose $\varphi,\,\theta,\,\omega\in \mathcal{L}^\bot(\mathcal{O}^{\sigma+2})$. Then there is a constant $0<c_0\leq\frac{4}{7}$ only dependent on $\textbf{a}_1$ and $\textbf{a}_2$ such that for any $\mu\in\Delta_A$ we have:
\begin{align*}
    \norm{\mathfrak{B}_{\mu,\ZZ}(\theta_\mu)}\leq C\norm{\log[z_1,z_2]}^{c_0}\max\{\norm{\theta_\mu}_{\sigma+2},\,\norm{\varphi_\mu}_{\sigma+1}\}+C\norm{\omega_\mu}_{\sigma+1}
  \end{align*}
(see \eqref{for:414} of Section \ref{sec:53}).
\end{proposition}
The proof is left for Appendix \ref{sec:16}. We recall \eqref{for:465} of Section \ref{sec:53}. The following result can be deduced from \eqref{for:64} of Proposition \ref{po:10} and Proposition \ref{po:1}:

\begin{corollary}\label{cor:2} Suppose $\varphi,\,\theta,\,\omega\in \mathcal{L}^\bot(\mathcal{O}^{\sigma+1})$. For any $\tau\geq1$ and any $\textbf{u}\in \mathcal{U}_{\mu,\sigma}$ we have
\begin{align*}
 \norm{\textbf{u}\mathfrak{D}_{\mu,\tau,\ZZ}(\theta_\mu)}\leq C\norm{\log[z_1,z_2]}^{c_0}\max\{\norm{\theta_\mu}_{\sigma+2},\,\norm{\varphi_\mu}_{\sigma+1}\}+C\norm{\omega_\mu}_{\sigma+1}
\end{align*}
for any $\mu\in\Delta_A$.
\end{corollary}

\subsection{Construction of approximate $C^\infty$ solution}\label{sec:7} We recall notations in  Section \ref{sec:1}.
In this part we give an approximate $C^\infty$ solution $\Theta$ of the almost twisted cocycle equations (see \eqref{for:240}) involving every element in the compact generating set $E$. The entire section is dedicated to proving Theorem \ref{th:9}.

 \subsubsection{Notations}\label{sec:56} The following notations will be used in this part:
 \begin{enumerate}

\item The positive constants $\sigma,\,\eta$ and $\beta$ (see Section \ref{for:432}) will appear in this section.  Let
\begin{align}
 \varrho=\max\{&\beta+5,\,\sigma+6,\,\text{\small$\frac{5}{2}$}\sigma+5,\beta+6,\,\beta+11, \, 2\beta+17 \}.
\end{align}

\item\label{for:456} We choose $\delta_0$ such that $\delta_0<\eta_{\sigma+6+\beta}(\textbf{a}_1)$ (see \eqref{for:488} of Section \ref{for:432}). Note that $\eta_{\sigma+6+\beta}(\textbf{a}_1)<\eta$.

\smallskip

 \item\label{for:509}  For any $\delta_0$-standard perturbation $E'(z)$ of $E$ (see Section \ref{sec:1}) we can write $z_v=w_v\cdot
\alpha_v(x)$,  where $w_v\in \exp(\mathfrak{g}_{1,A})\cap B_{\delta_0}(e)\cap Z(L)$, $v\in E$. We emphasize that we do not require the set $\{z_v: v\in E\}$ to be abelian.

\smallskip

   \item For any $\varepsilon$-standard perturbation $E'(z)$ of $E$ (see Section \ref{sec:1}) and
any $\mathfrak{p}_v\in \mathcal{L}^\bot(\mathcal{O}^\infty)$, $v\in E$, we set
\begin{enumerate}
  \item$\norm{\log[z,z]}=\max_{u,v\in E}\norm{\log[z_{v},z_{u}]}$;

  \smallskip
  \item $\norm{\mathfrak{p}}_{C^m}=\max_{v\in E}\{\norm{\mathfrak{p}_{v}}_{C^m}\}$;
  \item\label{for:240} for any $u,v\in E$
  \begin{align*}
 \mathbb{D}_{v,u}=\mathfrak{p}_{v}\circ z_{u}-\mathfrak{i}_{\mathcal{L}^\bot}\text{Ad}_{z_{u}}\mathfrak{p}_{v}-\big( \mathfrak{p}_{u}\circ z_{v}-\mathfrak{i}_{\mathcal{L}^\bot}\text{Ad}_{z_{v}}\mathfrak{p}_{u}\big);
\end{align*}

\item $\norm{\mathbb{D}}_{C^m}=\max_{u,\,v\in E}\{\norm{\mathbb{D}_{v,u}}_{C^m}\}$, for any $m\geq0$. $\norm{\mathfrak{p}}_{m}$ and $\norm{\mathbb{D}}_{m}$ are defined accordingly.
\end{enumerate}

\smallskip

\item\label{for:455} We recall $q(z_{\textbf{a}_1})$ and $\mathfrak{z}_1=qz_{\textbf{a}_i}q^{-1}$, $i=1,\,2$ as described in Corollary \ref{for:129}. We also set
$\mathfrak{z}_v=qz_{v}q^{-1}$, $v\in E$.

 \end{enumerate}

\subsubsection{Basic facts}\label{sec:55}
Since $\mathfrak{z}_1$ is perfect  and $p(\mathfrak{z}_1)\in A_0$  (see \eqref{for:78} of Corollary \ref{for:129}), each $\mathfrak{g}_{\chi, \mathfrak{z}_1}$, $\chi\in \Delta_{\mathfrak{z}_1}$ (see \eqref{sec:48} of Section \ref{sec:3}) is a set of unions of $\mathfrak{g}_{\phi,A_0}$, $\phi\in \Delta_{A_0}$. More precisely, there exists a subset $\mathfrak{s}_\chi\subset \Delta_{A_0}$  such that
\begin{align*}
  \mathfrak{g}_{\chi, \mathfrak{z}_1}=\bigcup_{\lambda\in \mathfrak{s}_\chi}\mathfrak{g}_{\lambda,A_0},\qquad \forall\, \chi\in \Delta_{\mathfrak{z}_1}.
\end{align*}
Then we define:
 \begin{enumerate}
   \item For each $\phi\in \Delta_{A_0}$, we fix a basis
$\{u_{(\phi,1)},\cdots, u_{(\phi,x_\phi)}\}$ of $\mathfrak{g}_{\phi,A_0}$, where $x_\phi=\dim \mathfrak{g}_{\phi,A_0}$.

\smallskip
  \item Let $\mathfrak{s}^+=\bigcup_{\chi\in \Delta_{\mathfrak{z}_1}^+}\mathfrak{s}_\chi$ and $\mathfrak{s}^-=\bigcup_{\chi\in \Delta_{\mathfrak{z}_1}^-}\mathfrak{s}_\chi$. Set
\begin{align*}
 \Lambda^+&=\{u_{\phi,i}:\phi\in \mathfrak{s}^+,\,\,1\leq i\leq x_\phi\}\qquad\text{and}\\
 \Lambda^-&=\{u_{\phi,i}:\phi\in \mathfrak{s}^-,\,\,1\leq i\leq x_\phi\}.
\end{align*}
Then $\Lambda^+$ is a basis of $\mathfrak{g}_{\mathfrak{z}_1}^+$ and $\Lambda^-$ is a basis of $\mathfrak{g}_{\mathfrak{z}_1}^-$.

\smallskip

\item\label{for:458} We can write
$\Lambda^+=\{\nu_1,\cdots, \nu_{p_+}\}$ (resp. $\Lambda^-=\{\nu_1,\cdots, \nu_{p_-}\}$), where $\nu_n=u_{(\phi_{i_n},j_n)}$, $1\leq n\leq p_+$ (resp. $1\leq n\leq p_-$) in such a way that if $n\leq m$ then
\begin{align*}
 \phi_{i_n}(p(\mathfrak{z}_1))\leq \phi_{i_m}(p(\mathfrak{z}_1))\quad\big(\text{resp.}\quad\phi_{i_n}(p(\mathfrak{z}_1))\geq \phi_{i_m}(p(\mathfrak{z}_1))\big).
\end{align*}
It is clear that
\begin{align}\label{for:270}
 [\nu_n,\nu_m] &\text{ is inside the subspace
linearly spanned by }\notag\\
&\nu_k, \,\, k>\max\{n,\,m\}.
\end{align}

\item\label{for:454} We define:
\begin{align*}
 \mathcal{E}_{+,\tau}&\stackrel{\text{def}}{=}\pi_{\nu_1}(f\circ \tau^{-1})\pi_{\nu_2}(f\circ \tau^{-1})\cdots \pi_{\nu_{p_{+}}}(f\circ \tau^{-1})\\
 (resp. \quad \mathcal{E}_{-,\tau}&\stackrel{\text{def}}{=}\pi_{\nu_1}(f\circ \tau^{-1})\pi_{\nu_2}(f\circ \tau^{-1})\cdots \pi_{\nu_{p_{-}}}(f\circ \tau^{-1})).
\end{align*}
(see Section \ref{sec:27} and \eqref{sect:4} of Section \ref{sec:19}).

We also define:
\begin{align*}
 \mathcal{E}_{+,\tau}'&\stackrel{\text{def}}{=}\pi_{\nu_{p_{+}}}(f\circ \tau^{-1})\cdots\pi_{\nu_2}(f\circ \tau^{-1})\pi_{\nu_1}(f\circ \tau^{-1}), \\
 \mathcal{E}_{+,\tau,k}&\stackrel{\text{def}}{=}\pi_{\nu_k}(f\circ \tau^{-1})\pi_{\nu_2}(f\circ \tau^{-1})\cdots \pi_{\nu_{p_{+}}}(f\circ \tau^{-1}),
\end{align*}
for any $1\leq k\leq p_{+}$.

$\mathcal{E}_{-,\tau}'$ and $\mathcal{E}_{-,\tau,k}$ are defined accordingly.

\smallskip

\item For any $\mu\in\Delta_A$,
\begin{align*}
 \mathcal{E}_{\mu,\tau}\stackrel{\text{def}}{=} \mathcal{E}_{-\delta(\mu),\tau},\quad \mathcal{E}'_{\mu,\tau}\stackrel{\text{def}}{=} \mathcal{E}'_{-\delta(\mu),\tau}
 \,\,\text{ and }\,\,\mathcal{E}_{\mu,\tau,k}\stackrel{\text{def}}{=} \mathcal{E}_{-\delta(\mu),\tau,k}
\end{align*}
(see \eqref{for:1} of Section \ref{sec:3}).
    Set
$p_\mu=p_{-\delta(\mu)}$.
 \end{enumerate}

\subsection{Main result} Theorem \ref{th:9} is the central part of the proof of Theorem \ref{th:7}. We recall notations in Section \ref{sec:56}.

\begin{theorem}\label{th:9} There is $\delta_0>0$ such that: for any $\delta_0$ standard perturbation $E'(z)$ of $E$, any $\tau>1$ and
any $\mathfrak{p}_v\in \mathcal{L}^\bot(\mathcal{O}^\infty)$, $v\in E$,
there exist $\Theta\in \mathcal{L}^\bot(\mathcal{O}^\infty)$ and $\mathcal{R}_v\in \mathcal{L}^\bot(\mathcal{O}^\infty)$ for each $v\in E$ such that
\begin{align}\label{for:77}
 \mathfrak{p}_v=\Theta\circ z_v-\mathfrak{i}_{\mathcal{L}^\bot}\text{Ad}_{z_{v}}\Theta+\mathcal{R}_v
\end{align}
with estimates: for any $v\in E$
\begin{align}\label{for:114}
 \max\{\norm{\Theta}_{C^m},\,\norm{\mathcal{R}_v}_{C^m}\}
 \leq C_{m}(\tau^{m+\varrho}\norm{\mathfrak{p}}_{C^\varrho}+\norm{\mathfrak{p}}_{C^{m+\varrho}})
\end{align}
for any $m\geq0$;  and
\begin{align}\label{for:121}
 \norm{\mathcal{R}_v}&_{C^0}\leq C\tau^{\varrho}\norm{\log[z,z]}^{\frac{c_0}{8}}\norm{\mathfrak{p}}_{C^\varrho}
 +C\tau^{\varrho}(\norm{\mathbb{D}}_{C^\varrho})^{\frac{1}{8}}(\norm{\mathfrak{p}}_{C^\varrho})^{\frac{7}{8}}\notag\\
 &+C\tau^{\varrho}(C_{m}\tau^{-m+\varrho}\norm{\mathfrak{p}}_{C^m})^{\frac{1}{8}}(\norm{\mathfrak{p}}_{C^\varrho})^{\frac{7}{8}}+C
\norm{\mathbb{D}}_{C^{\varrho}}
\end{align}
for any $m\geq \varrho$.

\end{theorem}

\subsection{Proof of Theorem \ref{th:9}} We follow the proof line of Theorem \ref{th:6}. Typically distinctions between nilmanifold examples and
(twisted) symmetric space examples are minor, and appear mostly at
 \emph{Step 4, 5 and 6}, where we need to take care of directional derivatives.

\smallskip
\emph{Step 1: Construction of $\Theta$ and $R_v$, $v\in E$. }

Let
\begin{align}
 \mathfrak{p}_{v,\mu}=\mathfrak{p}_{v}\circ q^{-1}|_{J_{\mu}},\quad
 v\in E\text{ and }\mu\in \Delta_A
\end{align}
(see \eqref{for:455} of Section \ref{sec:55}).

We denote by $z_i=z_{\textbf{a}_i}$, $i=1,\,2$,  to simplify notations.  Set
\begin{align}\label{for:461}
 \Phi_{\mu}\stackrel{\text{def}}{=} \mathcal{E}_{\mu,\tau}'\big(\mathcal{P}_{z_{1},\mu}(\mathfrak{p}_{\textbf{a}_1,\mu})\big), \quad \mu\in \Delta_A
\end{align}
(see \eqref{for:412} of Section \ref{sec:53} and \eqref{for:454} of Section \ref{sec:55}).

Set $\Theta$ to be the  map with coordinate maps $\Phi_{\mu}\circ q$. Let
\begin{align}
\mathcal{R}_v&\stackrel{\text{def}}{=}\mathfrak{p}_v-(\Theta\circ z_v-\mathfrak{i}_{\mathcal{L}^\bot}\text{Ad}_{z_{v}}\Theta)\quad \quad\forall\,v\in E\quad\text{and}\label{for:134}\\
 \mathcal{R}_{\textbf{a}_1,\mu}&\stackrel{\text{def}}{=}\mathcal{R}_{\textbf{a}_1}\circ q^{-1}|_{J_{\mu}},\quad\quad \mu\in\Delta.\notag
\end{align}
It follows from \eqref{for:134} that
\begin{align}\label{for:460}
\mathcal{R}_{\textbf{a}_1,\mu}&\overset{\text{(1)}}{=}\mathfrak{p}_{\textbf{a}_1,\mu}-\big(\Phi_\mu\circ \mathfrak{z}_{1}-(\text{Ad}_{z_{1}}|_{J_\mu})\Phi_\mu\big),\qquad \mu\in\Delta.
\end{align}
Here in $(1)$ we use the fact that each $J_{\mu}$ is invariant under $\text{Ad}_{z_1}$ (see \eqref{for:407}  of Proposition \ref{po:9}).

\smallskip

\noindent\emph{Step 2: Estimates for $\Theta$ and $\mathcal{R}_v$, $v\in E$. }

We recall that $z_{i}=w_{\textbf{a}_i}\cdot
\textbf{a}_i$,  where $w_{\textbf{a}_i}\in \exp(\mathfrak{g}_{1,A})\cap B_{\delta_0}(e)\cap Z(L)$ with $\delta_0<\eta$, $i=1,\,2$ (see \eqref{for:456} and \eqref{for:509} of Section \ref{sec:56}). It follows from \eqref{for:431} and \eqref{for:457} of Lemma \ref{le:7} that $\mathcal{P}_{z_{\textbf{a}_1},\mu}(\mathfrak{p}_{\textbf{a}_1,\mu})$ is a distributional solution of the equation
\begin{align}\label{for:237}
\mathcal{P}_{z_{1},\mu}(\mathfrak{p}_{\textbf{a}_1,\mu})\circ \mathfrak{z}_1-(\text{Ad}_{z_{1}}|_{J_\mu})\mathcal{P}_{z_{1},\mu}(\mathfrak{p}_{\textbf{a}_1,\mu})=\mathfrak{p}_{\textbf{a}_1,\mu}
\end{align}
(see \eqref{for:455} of Section \ref{sec:55}) satisfying
\begin{align}\label{for:82}
   \big\|u_1u_2\cdots u_m\mathcal{P}_{z_{1},\mu}(\mathfrak{p}_{\textbf{a}_1,\mu})\|_{-1}\leq C_m \norm{\mathfrak{p}_{\textbf{a}_1,\mu}}_{m+1}
  \end{align}
for any $m\geq0$ and any $u_i\in (\mathfrak{g}_{\mathfrak{z}_{1}}^{\delta(\mu)}\cup\mathfrak{g}_{\mathfrak{z}_{1}}^{0}) \bigcap \mathfrak{G}^1$, $1\leq i\leq m$.

We note that for any $\psi\in \mathfrak{g}_\mu(\mathcal{O}^\infty)$, keeping using \eqref{for:93} of Section \ref{sec:27} we see that
 \begin{align}\label{for:459}
  \langle \Phi_{\mu}, \, \psi\rangle=\langle\mathcal{P}_{z_1,\mu}(\mathfrak{p}_{\textbf{a}_1,\mu}), \, \mathcal{E}_{\mu,\tau}\psi\rangle.
 \end{align}
We will use Proposition \ref{cor:6} to estimate $\Phi_{\mu}$. To applying Proposition \ref{cor:6},
firstly, we  define the following sets:
 \begin{enumerate}
   \item  Let $Q_i=\{\nu_i\}$, $1\leq i\leq p_{-\delta(\mu)}$ (see \eqref{for:458} of Section \ref{sec:55}).

   \smallskip
   \item Let $Q'(\mu)$ be the subgroup with Lie algebra $\mathfrak{g}_{\mathfrak{z}_{1}}^{-\delta(\mu)}$.

   \smallskip
   \item Let
$Q_0(\mu)$ be the subgroup with Lie algebra $\mathfrak{g}_{\mathfrak{z}_{1}}^{\delta(\mu)}\oplus\mathfrak{g}_{\mathfrak{z}_{1}}^{0}$.
 \end{enumerate}
 Secondly, we need to show the following conditions are satisfied:
 \begin{enumerate}
  \item by definition, we see that $\text{Lie}(Q')$ is linearly spanned by  $\nu_i$, $1\leq i\leq p_{\mu}$;

 \smallskip
  \item it is clear that $\mathfrak{G}=\mathfrak{g}_{\mathfrak{z}_{1}}^{-\delta(\mu)}\oplus\mathfrak{g}_{\mathfrak{z}_{1}}^{\delta(\mu)}\oplus\mathfrak{g}_{\mathfrak{z}_{1}}^{0}$;

  \smallskip
  \item by \eqref{for:270} of Section \ref{sec:55} we have
\begin{align*}
 [\nu_n,\nu_m] &\text{ is inside the subspace
linearly spanned by }\notag\\
&\nu_k, \,\, k>\max\{n,\,m\};
\end{align*}
 \item  estimates of derivatives of $\mathcal{P}_{z_1,\mu}(\mathfrak{p}_{\textbf{a}_1,\mu})$ along $\mathfrak{g}_{\mathfrak{z}_{1}}^{\delta(\mu)}\cup\mathfrak{g}_{\mathfrak{z}_{1}}^{0}$ are available (see \eqref{for:82}).
\end{enumerate}
\eqref{for:459} shows that
 \begin{align*}
 \Phi_{\mu}=\mathcal{P}_{f,\tau}\big(\mathcal{P}_{z_1,\mu}(\mathfrak{p}_{\textbf{a}_1,\mu})\big)
 \end{align*}
 in the sense described  in Proposition \ref{cor:6}. The above discussion justifies the application of Proposition \ref{cor:6} to $\Phi_{\mu}$.  As a consequence, we have
\begin{align}\label{for:269}
  \norm{\Phi_{\mu}}_m&\leq  C_{m}(\tau^{m+5}\norm{\mathfrak{p}_{\textbf{a}_1,\mu}}_1+\norm{\mathfrak{p}_{\textbf{a}_1,\mu}}_{m+5})
\end{align}
for any $m\geq0$.

We note that $q$ is uniformly bounded.  As an immediate consequence of \eqref{for:269} we have
\begin{align*}
  \norm{\Theta}_m\leq  C_{m}(\tau^{m+5}\norm{\mathfrak{p}}_1+\norm{\mathfrak{p}}_{m+5}),\qquad \forall\,m\geq0.
\end{align*}
It follows that: for any $v\in E$ and any $m\geq0$
\begin{align}\label{for:243}
  \norm{\mathcal{R}_v}_m&\leq \norm{\mathfrak{p}_v}_m+\norm{\Theta\circ z_v-\mathfrak{i}_{\mathcal{L}^\bot}\text{Ad}_{z_{v}}\Theta}_m\notag\\
  &\leq  C_{m}(\tau^{m+5}\norm{\mathfrak{p}}_1+\norm{\mathfrak{p}}_{m+5}).
\end{align}
Then we have
\begin{align}\label{for:90}
 \max\{\norm{&\Theta}_{C^m},\,\norm{\mathcal{R}_v}_{C^m}\}\overset{\text{(1)}}{\leq} C_{m}\max\{\norm{\Theta}_{m+\beta},\,\norm{\mathcal{R}_v}_{m+\beta}\}\notag\\
 &\leq C_{m}(\tau^{m+5+\beta}\norm{\mathfrak{p}}_1+\norm{\mathfrak{p}}_{m+5+\beta})\notag\\
 &\overset{\text{(1,2)}}{\leq} C_{m}(\tau^{m+\varrho}\norm{\mathfrak{p}}_{C^1}+\norm{\mathfrak{p}}_{C^{m+\varrho}}).
\end{align}
for any $v\in E$ and any $m\geq0$.

Here in $(1)$ we use Sobolev embedding inequality \eqref{for:453} of Section \ref{sec:3}; in $(2)$ we use the definition of $\varrho$ in \eqref{for:236} of \eqref{for:432} of Section \ref{sec:3}.

Hence, \eqref{for:114} follows from \eqref{for:90}.

\smallskip

\smallskip
\noindent\emph{Step 3:  In this part we show that:}
\begin{align}\label{for:238}
 \langle \mathcal{R}_{\textbf{a}_1,\mu}, \, \psi\rangle=\big\langle (I-\mathcal{E}'_{\mu,\tau})\mathcal{P}_{z_{1},\mu}(\mathfrak{p}_{\textbf{a}_1,\mu}), \, \psi\circ \mathfrak{z}_{1}^{-1}-(\text{Ad}_{z_{1}}|_{J_\mu})^*\psi\big\rangle
\end{align}
for any $\mu\in\Delta$ and any $\psi\in J_\mu(\mathcal{O})$. Here $(\text{Ad}_{z_{1}}|_{J_\mu})^*$ is the transpose of $\text{Ad}_{z_{1}}|_{J_\mu}$.

\medskip

It follows from \eqref{for:460} that
\begin{align*}
 \mathcal{R}_{\textbf{a}_1,\mu}&=\mathfrak{p}_{\textbf{a}_1,\mu}-\big(\Phi_\mu\circ \mathfrak{z}_{1}-(\text{Ad}_{z_{1}}|_{J_\mu})\Phi_\mu\big)\notag\\
 &\overset{\text{(1)}}{=}\big(\mathcal{P}_{z_{1},\mu}(\mathfrak{p}_{\textbf{a}_1,\mu})-\Phi_\mu\big)\circ \mathfrak{z}_{1}-(\text{Ad}_{z_{1}}|_{J_\mu})
 \big(\mathcal{P}_{z_1,\mu}(\mathfrak{p}_{\textbf{a}_1,\mu})-\Phi_\mu\big)\notag\\
 &\overset{\text{(2)}}{=}\big((I-\mathcal{E}'_{\mu,\tau})\mathcal{P}_{z_{1},\mu}(\mathfrak{p}_{\textbf{a}_1,\mu})\big)\circ \mathfrak{z}_{1}\\
 &-(\text{Ad}_{z_{1}}|_{J_\mu})
 \big((I-\mathcal{E}'_{\mu,\tau})\mathcal{P}_{z_{1},\mu}(\mathfrak{p}_{\textbf{a}_1,\mu})\big)
\end{align*}
for any $\mu\in\Delta_A$.

Here in $(1)$ we use \eqref{for:237}; in $(2)$ we recall \eqref{for:461}.

Then for any $\psi\in J_\mu(\mathcal{O})$ we have
\begin{align*}
 \langle \mathcal{R}_{\textbf{a}_1,\mu}, \, \psi\rangle&=\big\langle\big((I-\mathcal{E}'_{\mu,\tau})\mathcal{P}_{\mu}(\mathfrak{p}_{\textbf{a}_1,\mu})\big)\circ \mathfrak{z}_1, \, \psi\big\rangle\\
 &-\big\langle(\text{Ad}_{z_{1}}|_{J_\mu})
 \big((I-\mathcal{E}'_{\mu,\tau})\mathcal{P}_{z_1,\mu}(\mathfrak{p}_{\textbf{a}_1,\mu})\big), \, \psi\big\rangle\\
 &=\big\langle(I-\mathcal{E}'_{\mu,\tau})\mathcal{P}_{z_1,\mu}(\mathfrak{p}_{\textbf{a}_1,\mu}), \, \psi\circ \mathfrak{z}_1^{-1}\big\rangle\\
 &-\big\langle(I-\mathcal{E}'_{\mu,\tau})\mathcal{P}_{z_1,\mu}(\mathfrak{p}_{\textbf{a}_1,\mu}), \, (\text{Ad}_{z_{1}}|_{J_\mu})^*\psi\big\rangle\\
 &=\big\langle (I-\mathcal{E}_{\mu,\tau}')\mathcal{P}_{z_1,\mu}(\mathfrak{p}_{\textbf{a}_1,\mu}), \, \psi\circ \mathfrak{z}_1^{-1}-(\text{Ad}_{z_{1}}|_{J_\mu})^*\psi\big\rangle.
\end{align*}
We thus finish the proof.

\smallskip
\noindent\emph{Step 4:  In this part we show that:}

for any $\textbf{u}\in \mathcal{U}_{\mu,\sigma}$ (see \eqref{for:1} of Section \ref{sec:3}) and any $\psi\in J_\mu(\mathcal{O}^\sigma)$
\begin{align}\label{for:462}
   (I-\mathcal{E}_{\mu,\tau})\textbf{u}\psi=\sum_{l=1}^{p_{\mu}}\pi_{\nu_l}(f_1\circ \tau^{-1})\psi_l,\qquad \mu\in\Delta,
  \end{align}
where
\begin{gather*}
 \psi_{l}=\sum_i \tau_{l,i}\tau^{-j_{l,i}}\textbf{u}_{l,i}\psi_{l,i},
\end{gather*}
where $\textbf{u}_{l,i}\in \mathcal{U}_{\mu,\sigma}$, $j_{l,i}\geq0$ and $\tau_{l,i}\in\CC$,   and
\begin{align*}
 \psi_{l,i}=\pi_{\nu_{l+1}}(f^{(j_{l,l+1,i})}\circ \tau^{-1})\pi_{\nu_{l+2}}(f^{(j_{l,l+2,i})}\circ \tau^{-1})\cdots \pi_{\nu_{p_{\mu}}}(f^{(j_{l,p_{\mu},i})}\circ \tau^{-1})\psi,
\end{align*}
$1\leq i\leq l$, where $f^{(j_{l,l,i})}\in \mathcal{S}(\RR)$, $l+1\leq l\leq p_{\mu}$.

\medskip
\emph{Note}. Keeping using \eqref{for:7} of Lemma \ref{cor:1}, we have
\begin{align}\label{for:466}
\|\psi_{l,i}\|_{\sigma+1}\leq C\|\psi\|_{\sigma+1},\qquad \forall\,\psi \in J_\mu(\mathcal{O}^{\sigma+1})
\end{align}
for any $l,\,i$.

\smallskip
\eqref{for:238} shows that to estimate $\mathcal{R}_{\textbf{a}_1,\mu}$, we need to estimate $(I-\mathcal{E}'_{\mu,\tau})\mathcal{P}_{z_{1},\mu}(\mathfrak{p}_{\textbf{a}_1,\mu})$ first.
We also note that
\begin{align}\label{for:463}
 \big\langle (I-\mathcal{E}'_{\mu,\tau})&\mathcal{P}_{z_1,\mu}(\mathfrak{p}_{\textbf{a}_1,\mu}), \, \textbf{u}\psi\big\rangle=\big\langle \mathcal{P}_{z_1,\mu}(\mathfrak{p}_{\textbf{a}_1,\mu}), \, \textbf{u}\psi\big\rangle-\big\langle \mathcal{E}'_{\mu,\tau}\mathcal{P}_{z_1,\mu}(\mathfrak{p}_{\textbf{a}_1,\mu}), \, \textbf{u}\psi\big\rangle\notag\\
 &\overset{\text{(1)}}{=}\big\langle \mathcal{P}_{z_1,\mu}(\mathfrak{p}_{\textbf{a}_1,\mu}), \, \textbf{u}\psi\big\rangle-\big\langle \mathcal{P}_{z_1,\mu}(\mathfrak{p}_{\textbf{a}_1,\mu}), \, \mathcal{E}_{\mu,\tau}\textbf{u}\psi\big\rangle\notag\\
 &=\big\langle \mathcal{P}_{z_1,\mu}(\mathfrak{p}_{\textbf{a}_1,\mu}), \, (I-\mathcal{E}_{\mu,\tau})\textbf{u}\psi\big\rangle.
\end{align}
Here in $(1)$ we keep using \eqref{for:93} of Section \ref{sec:27}.

Then it is natural for us to rewrite $(I-\mathcal{E}_{\mu,\tau})\textbf{u}\psi$ in this step.

Using the same arguments as in \emph{Step 4} of Section \ref{sec:57}, similar to \eqref{for:205}, we have
\begin{align*}
 I-\mathcal{E}_{\mu,\tau}=\sum_{l=1}^{p_{\mu}}\pi_{\nu_l}(f_1\circ \tau^{-1})\pi_{\nu_{l+1}}(f\circ \tau^{-1})\cdots\pi_{\nu_{p_{\mu}}}(f\circ \tau^{-1}).
\end{align*}
Then the result follows from \eqref{re:2} of Proposition \ref{po:10}.

\smallskip
\noindent\emph{Step 5:  In this part we show that:}
\begin{align}\label{for:242}
 &\Big\|\big\langle (I-\mathcal{E}'_{\mu,\tau})\mathcal{P}_{z_{1},\mu}(\mathfrak{p}_{\textbf{a}_1,\mu}), \, \textbf{u}\psi\big\rangle\Big\|\notag\\
 &\leq C\big(\norm{\log[z,z]}^{c_0}\norm{\mathfrak{p}}_{\sigma+2}
 +\norm{\mathbb{D}}_{\sigma+1}+C_{m}\tau^{-m+\sigma}\norm{\mathfrak{p}}_m\big)
 \|\psi\|_{\sigma+1}
 \end{align}
for any $m\geq\sigma$,  any $\psi\in J_\mu(\mathcal{O}^{\sigma+1})$ and any $\textbf{u}\in \mathcal{U}_{\mu,\sigma}$.

\medskip

For any $\psi\in J_\mu(\mathcal{O}^{\sigma+1})$ and any $\textbf{u}\in \mathcal{U}_{\mu,\sigma}$ we have
\begin{align*}
 &\big\langle (I-\mathcal{E}'_{\mu,\tau})\mathcal{P}_{z_1,\mu}(\mathfrak{p}_{\textbf{a}_1,\mu}), \, \textbf{u}\psi\big\rangle\notag\\
 &\overset{\text{(1)}}{=}\big\langle \mathcal{P}_{z_1,\mu}(\mathfrak{p}_{\textbf{a}_1,\mu}), \, (I-\mathcal{E}_{\mu,\tau})\textbf{u}\psi\big\rangle\notag\\
 &\overset{\text{(2)}}{=}\sum_{l=1}^{p_{\mu}}\big\langle
 \mathcal{P}_{z_1,\mu}(\mathfrak{p}_{\textbf{a}_1,\mu}), \, \pi_{\nu_l}(f_1\circ \tau^{-1})\psi_l\big\rangle\notag\\
 &\overset{\text{(3)}}{=}\sum_{l=1}^{p_\mu}\big\langle
 \pi_{\nu_l}(f_1\circ \tau^{-1})\mathcal{P}_{z_1,\mu}(\mathfrak{p}_{\textbf{a}_1,\mu}), \, \psi_l\big\rangle\notag\\
 &\overset{\text{(2)}}{=}\sum_{l=1}^{p_\mu}\sum_i \overline{k_{l,i}}\tau^{-j_{l,i}}\big\langle
 \mathfrak{D}_{\mu,\nu_l,\tau,L_\mu}(\mathfrak{p}_{\textbf{a}_1,\mu}), \, \textbf{u}_{l,i}\psi_{l,i}\big\rangle\notag\\
 &\overset{\text{(4)}}{=}P_1-P_2.
\end{align*}
Here in $(1)$ we use \eqref{for:463}; in $(2)$ we use \eqref{for:462}; in $(3)$ we use \eqref{for:93} of Section \ref{sec:27};
in $(4)$ we set
\begin{align*}
 P_1&=\sum_{l=1}^{p_\mu}\sum_i \overline{k_{l,i}}\tau^{-j_{l,i}}\big\langle
 \mathfrak{D}_{\mu,\nu_l,\tau,\ZZ}(\mathfrak{p}_{\textbf{a}_1,\mu}), \, \textbf{u}_{l,i}\psi_{l,i}\big\rangle,\qquad\text{and}\\
 P_2&=\sum_{l=1}^{p_\mu}\sum_i \overline{k_{l,i}}\tau^{-j_{l,i}}\big\langle
 \mathfrak{D}_{\mu,\nu_l,\tau,L_\mu^c}(\mathfrak{p}_{\textbf{a}_1,\mu}), \, \textbf{u}_{l,i}\psi_{l,i}\big\rangle,
\end{align*}
where $\textbf{u}_{l,i}\in \mathcal{U}_{\mu,\sigma}$ and $j_{l,i}\geq0$.

It follows that
\begin{align}\label{for:468}
 \Big\|\big\langle (I-\mathcal{E}'_{\mu,\tau})\mathcal{P}_{z_1,\mu}(\mathfrak{p}_{\textbf{a}_1,\mu}), \, \textbf{u}\psi\big\rangle\Big\|\leq \norm{P_1}+\norm{P_2}.
\end{align}
Next, we estimate $\norm{P_1}$ and $\norm{P_2}$, respectively. We will use Corollary \ref{cor:2} to estimate $\norm{P_1}$. Thus we turn to the almost twisted cocycle equation equation for $\mathfrak{z}_1$ and $\mathfrak{z}_2$:
\begin{align*}
 \mathfrak{p}_{\textbf{a}_1}\circ z_{2}-\mathfrak{i}_{\mathcal{L}^\bot}\text{Ad}_{z_{2}}\mathfrak{p}_{\textbf{a}_1}&=\mathfrak{p}_{\textbf{a}_2}\circ z_{1}-\mathfrak{i}_{\mathcal{L}^\bot}\text{Ad}_{z_{2}}\mathfrak{p}_{\textbf{a}_2}+ \mathbb{D}_{\textbf{a}_1,\textbf{a}_2}
\end{align*}
(see \eqref{for:240} of Section \ref{sec:56}).

Since $J_{\mu}$ is invariant under $\text{Ad}_{z_1}$ (see \eqref{for:407}  of Proposition \ref{po:9}), we
have the decomposition:
\begin{align}\label{for:464}
 \mathfrak{p}_{\textbf{a}_1,\mu}\circ \mathfrak{z}_2-(\text{Ad}_{z_{2}}|_{J_\mu})\mathfrak{p}_{\textbf{a}_1,\mu}&=\mathfrak{p}_{\textbf{a}_2,\mu}\circ \mathfrak{z}_1-(\text{Ad}_{z_{1}}|_{J_\mu})\mathfrak{p}_{\textbf{a}_2,\mu}+\mathbb{D}_{\textbf{a}_2,\textbf{a}_1,\mu},
\end{align}
where $\mathbb{D}_{\textbf{a}_2,\textbf{a}_1,\mu}=\mathbb{D}_{\textbf{a}_2,\textbf{a}_1}\circ q^{-1}|_{J_{\mu}}$ .

Applying Corollary \ref{cor:2} to \eqref{for:464}, we have
\begin{align}\label{for:467}
 \norm{\textbf{u}\mathfrak{D}_{\mu,\tau,\ZZ}(\theta_\mu)}&\leq C\norm{\log[z_{1},z_{2}]}^{c_0}\max\{\norm{\mathfrak{p}_{\textbf{a}_1,\mu}}_{\sigma+2},\,\norm{\mathfrak{p}_{\textbf{a}_2,\mu}}_{\sigma+1}\}\notag\\
 &+C\norm{\mathbb{D}_{\textbf{a}_2,\textbf{a}_1,\mu}}_{\sigma+1}\notag\\
 &\leq C\norm{\log[z,z]}^{c_0}\norm{\mathfrak{p}}_{\sigma+2}+C\norm{\mathbb{D}}_{\sigma+1}
\end{align}
for any $\textbf{u}\in \mathcal{U}_{\mu,\sigma}$.

It follows that
\begin{align}\label{for:276}
 P_1&\overset{\text{(1)}}{\leq}\sum_{l=1}^{p_\mu}\sum_i C\Big\| \big\langle
 \textbf{u}_{l,i}\mathfrak{D}_{\mu,\nu_l,\tau,\ZZ}(\mathfrak{p}_{\textbf{a}_1,\mu}), \, \psi_{l,i}\big\rangle \Big\|\notag\\
 &\overset{\text{(2)}}{\leq}\sum_{l=1}^{p_\mu}\sum_i C\norm{\textbf{u}_{l,i}\mathfrak{D}_{\mu,\tau,\ZZ}(\mathfrak{p}_{\textbf{a}_1,\mu})}
 \|\psi_{l,i}\|_{\sigma+1}\notag\\
 &\overset{\text{(3)}}{\leq}\sum_{l=1}^{p_\mu}\sum_iC_1\big(\norm{\log[z,z]}^{c_0}\norm{\mathfrak{p}}_{\sigma+2}+\norm{\mathbb{D}}_{\sigma+1}\big)
 \|\psi_{l,i}\|_{\sigma+1}\notag\\
 &\overset{\text{(4)}}{\leq} C_2\big(\norm{\log[z,z]}^{c_0}\norm{\mathfrak{p}}_{\sigma+2}+\norm{\mathbb{D}}_{\sigma+1}\big)
 \|\psi\|_{\sigma+1}.
\end{align}
Here in $(1)$ we use the fact that $j_{l,i}\geq0$ and $\tau\geq1$; in $(2)$
we recall the definition of $\norm{\textbf{u}_{l,i}\mathfrak{D}_{\mu,z_1,\tau,\ZZ}(\mathfrak{p}_{\textbf{a}_1,\mu})}$ in  \eqref{for:465} of Section \ref{sec:53}  and
use the fact that $\textbf{u}_{l,i}\in \mathcal{U}_{\mu,\sigma}$; in $(3)$ we use \eqref{for:467}; in $(4)$ we use \eqref{for:466}.

Next, we estimate $P_2$.
\begin{align}\label{for:277}
 \norm{P_2}&\overset{\text{(1)}}{\leq} \sum_{l=1}^{p_\mu}\sum_i C\big\| \textbf{u}_{l,i}\mathfrak{D}_{\mu,\nu_l,\tau,L_\mu^c}(\mathfrak{p}_{\textbf{a}_1,\mu})\big\|
 \norm{\psi_{l,i}}\notag\\
 &\overset{\text{(2)}}{\leq} \sum_{l=1}^{p_\mu}\sum_iC_m\tau^{-m+\sigma}\norm{\mathfrak{p}_{\textbf{a}_1,\mu}}_m\norm{\psi_{l,i}}\notag\\
 &\overset{\text{(3)}}{\leq} \sum_{k=1}^{p_\mu}C_{m,1}\tau^{-m}\norm{\mathfrak{p}}_m
 \|\psi\|
\end{align}
for any $m\geq\sigma$.

Here in $(1)$ we use cauchy schwarz inequality and the fact that $j_{l,i}\geq0$ and $\tau\geq1$; in $(2)$ we recall that $\nu_l\in\mathfrak{g}^{-\delta(\mu)}_{\mathfrak{z}_1}\cap\mathfrak{G}^1$ (see \eqref{for:454} of Section \ref{sec:55})   and $\textbf{u}_{l,i}\in \mathcal{U}_{\mu,\sigma}$. Then we use Lemma \ref{for:422}; in $(3)$ we use \eqref{for:466}.

Hence, \eqref{for:242} follows from \eqref{for:468}, \eqref{for:276} and \eqref{for:277}.

\smallskip
\noindent\emph{Step 5:  In this part, we prove that}
\begin{align}\label{for:244}
 \norm{v^{\frac{\sigma}{2}}\mathcal{R}_{\textbf{a}_1,\mu}}\leq& C\tau^{\frac{\varrho}{2}}\big(\norm{\log[z,z]}^{c_0}\norm{\mathfrak{p}}_{\varrho}
 \notag\\
 &+\norm{\mathbb{D}}_{\varrho}+C_{m}\tau^{-m+\sigma}\norm{\mathfrak{p}}_m\big)^{\frac{1}{2}}(\norm{\mathfrak{p}}_{\varrho})^{\frac{1}{2}}
\end{align}
for any $\mu\in \Delta_A$ and any $m\geq\sigma$.

For any $\mu\in \Delta_A$, choose $\lambda\in \Delta_{A}$ such that $(\mu,\lambda)$ is a poor row (see \eqref{for:1} of Section \ref{sec:3}) and $\mathfrak{g}_{\lambda,A}\cap \text{Lie}(G)\neq\emptyset$. Fix $v\in \mathfrak{g}_{\lambda,A}\cap \text{Lie}(G)\cap\mathfrak{G}^1$. It is clear that $(\mu,\underbrace{\lambda,\cdots, \lambda}_{\sigma\text{ copies of }\lambda})$ is also a a poor row and thus $v^\sigma\in \mathcal{U}_{\mu,\sigma}$.

By using \eqref{for:238} we have
\begin{align*}
  &\langle \mathcal{R}_{\textbf{a}_1,\mu}, \, v^{\sigma}\mathcal{R}_{\textbf{a}_1,\mu}\rangle\\
 &=\big\langle (I-\mathcal{E}'_{\mu,\tau})\mathcal{P}_{z_{1},\mu}(\mathfrak{p}_{\textbf{a}_1,\mu}), \, (v^{\sigma}\mathcal{R}_{\textbf{a}_1,\mu})\circ \mathfrak{z}_{1}^{-1}-(\text{Ad}_{z_{1}}|_{J_\mu})^* v^{\sigma}\mathcal{R}_{\textbf{a}_1,\mu}\big\rangle\\
 &=\big\langle (I-\mathcal{E}_{\mu,\tau})\mathcal{P}_{z_{1},\mu}(\mathfrak{p}_{\textbf{a}_1,\mu}), \, cv_1^{\sigma}(\mathcal{R}_{\textbf{a}_1,\mu}\circ \mathfrak{z}_{1}^{-1})\big\rangle\\
 &-\big\langle (I-\mathcal{E}_{\mu,\tau})\mathcal{P}_{z_{1},\mu}(\mathfrak{p}_{\textbf{a}_1,\mu}), \, v^{\sigma}(\text{Ad}_{z_{1}}|_{J_\mu})^* \mathcal{R}_{\textbf{a}_1,\mu}\big\rangle.
\end{align*}
where $v_1=\norm{\text{Ad}_{\mathfrak{z}_{1}}(v)}^{-1}\text{Ad}_{\mathfrak{z}_{1}}(v)$ and $c=\norm{\text{Ad}_{\mathfrak{z}_{1}}(v)}^{\sigma}$.

Since $v^\sigma\in \mathcal{U}_{\mu,\sigma}$, $v_1^\sigma\in \mathcal{U}_{\mu,\sigma}$ as $\text{Ad}_{\mathfrak{z}_{1}}$ preserves $\mathfrak{g}_{\lambda,A}$ (see \eqref{for:404} of Proposition \ref{po:9} and \eqref{for:435} of Corollary \ref{for:129}).
It is harmless to assume that $\sigma$
is even. Then we have
\begin{align*}
&\norm{v^{\frac{\sigma}{2}}\mathcal{R}_{\textbf{a}_1,\mu}}^2=\big\|\langle \mathcal{R}_{\textbf{a}_1,\mu}, \, v^{\sigma}\mathcal{R}_{\textbf{a}_1,\mu}\rangle\big\|\notag\\
&\leq \Big\|\big\langle (I-\mathcal{E}_{\mu,\tau})\mathcal{P}_{z_{1},\mu}(\mathfrak{p}_{\textbf{a}_1,\mu}), \, cv_1^{\sigma}(\mathcal{R}_{\textbf{a}_1,\mu}\circ \mathfrak{z}_{1}^{-1})\big\rangle\Big\|\notag\\
 &+\Big\|\big\langle (I-\mathcal{E}_{\mu,\tau})\mathcal{P}_{z_{1},\mu}(\mathfrak{p}_{\textbf{a}_1,\mu}), \, v^{\sigma}(\text{Ad}_{z_{1}}|_{J_\mu})^* \mathcal{R}_{\textbf{a}_1,\mu}\big\rangle\Big\|\notag\\
&\overset{\text{(1)}}{\leq} C\big(\norm{\log[z,z]}^{c_0}\norm{\mathfrak{p}}_{\sigma+2}
 +\norm{\mathbb{D}}_{\sigma+1}+C_{m}\tau^{-m+\sigma}\norm{\mathfrak{p}}_m\big)\notag\\
 &\cdot (\|\mathcal{R}_{\textbf{a}_1,\mu}\circ \mathfrak{z}_{1}^{-1}\|_{\sigma+1}+\big\|(\text{Ad}_{z_{1}}|_{J_\mu})^* \mathcal{R}_{\textbf{a}_1,\mu}\big\|_{\sigma+1})\notag\\
&\leq C_1\big(\norm{\log[z,z]}^{c_0}\norm{\mathfrak{p}}_{\sigma+2}
 +\norm{\mathbb{D}}_{\sigma+1}+C_{m}\tau^{-m+\sigma}\norm{\mathfrak{p}}_m\big)\big\|\mathcal{R}_{\textbf{a}_1,\mu}\big\|_{\sigma+1}\notag\\
  &\overset{\text{(2)}}{\leq} C_2\tau^{\varrho}\big(\norm{\log[z,z]}^{c_0}\norm{\mathfrak{p}}_{\varrho}
 +\norm{\mathbb{D}}_{\varrho}+C_{m}\tau^{-m+\sigma}\norm{\mathfrak{p}}_m\big)\norm{\mathfrak{p}}_{\varrho}
\end{align*}
for any $m\geq\sigma$.

Here in $(1)$ we use \eqref{for:242}; in $(2)$ we use \eqref{for:243} and the definition of $\varrho$ in \eqref{for:236} of \eqref{for:432} of Section \ref{sec:3}.

Hence we get the result.

\smallskip
\noindent\emph{Step 6:  In this part, we prove that}.
\begin{align}\label{for:245}
 \norm{\mathcal{R}_{\textbf{a}_1}}\leq &C\tau^{\varrho}\big(\norm{\log[z,z]}^{c_0}\norm{\mathfrak{p}}_{\varrho}
 +\norm{\mathbb{D}}_\varrho+C_{m}\tau^{-m+\varrho}\norm{\mathfrak{p}}_m\big)^{\frac{1}{4}}(\norm{\mathfrak{p}}_{\varrho})^{\frac{3}{4}}
\end{align}
for any $m\geq\varrho$.

The Jacobson-Morosov theorem asserts that $v$ imbeds in a subalgebra $\mathfrak{s}\subset \text{Lie}(G)$ which is isomorphic to $\mathfrak{sl}(2,\RR)$. Let $\exp(\mathfrak{s})$ be the subgroup with Lie algebra $\mathfrak{s}$. Theorem \ref{th:3}
shows that $\pi|_{\exp(\mathfrak{s})}$ has a spectral gap. Let $\xi=v^{\frac{\sigma}{2}}\mathcal{R}_{\textbf{a}_1,\mu}$.  Keeping applying Theorem \ref{th:8} to the equation
\begin{align*}
 v^{\frac{\sigma}{2}}\mathcal{R}_{\textbf{a}_1,\mu}=\xi,
\end{align*}
we have
\begin{align*}
 \norm{\mathcal{R}&_{\textbf{a}_1,\mu}}\leq C\norm{\xi}_{\exp(\mathfrak{s}),\sigma}\leq C\norm{\xi}_{\sigma}=C\norm{v^{\frac{\sigma}{2}}\mathcal{R}_{\textbf{a}_1,\mu}}_{\sigma}\\
 &\overset{\text{(1)}}{\leq} C_1(\norm{v^{\frac{\sigma}{2}}\mathcal{R}_{\textbf{a}_1,\mu}})^{\frac{1}{2}}(\norm{v^{\frac{\sigma}{2}}\mathcal{R}_{\textbf{a}_1,\mu}}_{2\sigma})^{\frac{1}{2}}\\
 &\overset{\text{(2)}}{\leq}C_2\tau^{\frac{\varrho}{4}}\big(\norm{\log[z,z]}^{c_0}\norm{\mathfrak{p}}_{\sigma+2}
 +\norm{\mathbb{D}}_{\sigma+1}+C_{m}\tau^{-m+\sigma}\norm{\mathfrak{p}}_m\big)^{\frac{1}{4}}\\
 &\cdot (\norm{\mathfrak{p}}_{\varrho})^{\frac{1}{4}}(\norm{\mathcal{R}_{\textbf{a}_1,\mu}}_{\frac{5}{2}\sigma})^{\frac{1}{2}}\\
 &\overset{\text{(3)}}{\leq}C_3\tau^{\frac{3}{4}\varrho}\big(\norm{\log[z,z]}^{c_0}\norm{\mathfrak{p}}_{\varrho}
 +\norm{\mathbb{D}}_\varrho+C_{m}\tau^{-m+\sigma}\norm{\mathfrak{p}}_m\big)^{\frac{1}{4}}(\norm{\mathfrak{p}}_{\varrho})^{\frac{3}{4}}
\end{align*}
for any $m\geq\sigma$.

Here in $(1)$ we use interpolation inequalities; in $(2)$ we use \eqref{for:244} to estimate $\norm{v^{\frac{\sigma}{2}}\mathcal{R}_{\textbf{a}_1,\mu}}$; in $(3)$ we use \eqref{for:243} to estimate $\norm{\mathcal{R}_{\textbf{a}_1,\mu}}_{\frac{5}{2}\sigma}$:
\begin{align*}
  \norm{\mathcal{R}_{\textbf{a}_1,\mu}}_{\frac{5}{2}\sigma}\leq C(\tau^{\frac{5}{2}\sigma+5}\norm{\mathfrak{p}}_1+\norm{\mathfrak{p}}_{\frac{5}{2}\sigma+5})
  \overset{\text{(a)}}{\leq}  2C\tau^{\varrho}\norm{\mathfrak{p}}_{\varrho}
\end{align*}
Here in $(a)$ we use the definition of $\varrho$ in \eqref{for:236} of \eqref{for:432} of Section \ref{sec:3}.

The above estimate implies \eqref{for:245} by noting that $\tau\geq1$.

\smallskip
\noindent\emph{Step 7:  Estimates for $\mathcal{R}_{v}$, $v\in E$}.

We use the almost twisted cocycle equation over $v$ and $\textbf{a}_1$ from \eqref{for:240} of Section \ref{sec:56}:
\begin{align*}
 \mathfrak{p}_{v}\circ z_{1}-\mathfrak{i}_{\mathcal{L}^\bot}\text{Ad}_{z_{1}}\mathfrak{p}_{v}=\mathfrak{p}_{\textbf{a}_1}\circ z_{v}-\mathfrak{i}_{\mathcal{L}^\bot}\text{Ad}_{z_{v}}\mathfrak{p}_{\textbf{a}_1}+ \mathbb{D}_{v,\textbf{a}_1},\quad \forall \,v\in E.
\end{align*}
We substitute  the expressions for $\mathfrak{p}_{\textbf{a}_1}$ and $\mathfrak{p}_{v}$
from \eqref{for:134} respectively into the above equation. Then we get: for each $v\in E$
\begin{align*}
\big(\mathcal{R}_v+(\Theta\circ z_v-\mathfrak{i}&_{\mathcal{L}^\bot}\text{Ad}_{z_{v}}\Theta)\big)\circ z_{1}-\mathfrak{i}_{\mathcal{L}^\bot}\text{Ad}_{z_{1}}\big(\mathcal{R}_v+(\Theta\circ z_v-\mathfrak{i}_{\mathcal{L}^\bot}\text{Ad}_{z_{v}}\Theta)\big)\\
&=\big(\mathcal{R}_{\textbf{a}_1}+(\Theta\circ z_{1}-\mathfrak{i}_{\mathcal{L}^\bot}\text{Ad}_{z_{1}}\Theta)\big)\circ z_v\\
&-\mathfrak{i}_{\mathcal{L}^\bot}\text{Ad}_{z_{v}}\big(\mathcal{R}_{\textbf{a}_1}+(\Theta\circ z_{1}-\mathfrak{i}_{\mathcal{L}^\bot}\text{Ad}_{z_{1}}\Theta)\big)+\mathbb{D}_{\textbf{a}_1,v}.
\end{align*}
 The above equation can be simplified as: for each $v\in E$
\begin{align}\label{for:469}
&\mathcal{R}_v\circ z_{1}-\mathfrak{i}_{\mathcal{L}^\bot}\text{Ad}_{z_{1}}\mathcal{R}_v\notag\\
&=\mathcal{R}_{\textbf{a}_1}\circ z_v-\mathfrak{i}_{\mathcal{L}^\bot}\text{Ad}_{z_{v}}\mathcal{R}_{\textbf{a}_1}+\Theta\circ (z_{1}z_v)-\Theta\circ (z_vz_{1})\notag\\
&+\big(\mathfrak{i}_{\mathcal{L}^\bot}\text{Ad}_{z_vz_{1}}-\mathfrak{i}_{\mathcal{L}^\bot}\text{Ad}_{z_{1}z_{v}}\big)\Theta+\mathbb{D}_{\textbf{a}_1,v}.
\end{align}
Here in $(1)$ we use the fact that
\begin{align*}
(\mathfrak{i}_{\mathcal{L}^\bot}\text{Ad}_{z_{1}})(\mathfrak{i}_{\mathcal{L}^\bot}\text{Ad}_{z_{v}})&
=\mathfrak{i}_{\mathcal{L}^\bot}\text{Ad}_{z_{1}z_{v}}\qquad\text{and}\\
(\mathfrak{i}_{\mathcal{L}^\bot}\text{Ad}_{z_{v}})(\mathfrak{i}_{\mathcal{L}^\bot}\text{Ad}_{z_{1}})&=\mathfrak{i}_{\mathcal{L}^\bot}\text{Ad}_{z_{v}z_{1}}
\end{align*}
(see \eqref{for:54} of Section \ref{sec:3}).

Then we have
 \begin{align}\label{for:186}
  \norm{\mathcal{R}_v}&_{C^0}\overset{\text{\tiny$(1)$}}{\leq} C\norm{\mathcal{R}_v}_{\beta}\notag\\
  &\overset{\text{\tiny$(2)$}}{\leq} C_1\big\|\mathcal{R}_{\textbf{a}_1}\circ z_v-\mathfrak{i}_{\mathcal{L}^\bot}\text{Ad}_{z_{v}}\mathcal{R}_{\textbf{a}_1}+\Theta\circ (z_{1}z_v)-\Theta\circ (z_vz_{1})\notag\\
&+\big(\mathfrak{i}_{\mathcal{L}^\bot}\text{Ad}_{z_vz_{1}}-\mathfrak{i}_{\mathcal{L}^\bot}\text{Ad}_{z_{1}z_{v}}\big)\Theta+\mathbb{D}_{\textbf{a}_1,v}\big\|_{\beta+6}\notag\\
&\leq C_1\norm{\mathcal{R}_{\textbf{a}_1}}_{\beta+6}+C_1\emph{\textbf{d}}(z_{1}z_v,z_vz_{1})\norm{\Theta}_{\beta+6}+C
\norm{\mathbb{D}_{\textbf{a}_1,v}}_{\beta+6}\notag\\
&\overset{\text{\tiny$(3)$}}{\leq} C_1\norm{\mathcal{R}_{\textbf{a}_1}}_{\beta+6}+C_2\norm{\log[z_{1},z_v] }\norm{\Theta}_{\beta+6}+C
\norm{\mathbb{D}}_{\beta+6}\notag\\
&\overset{\text{\tiny$(4,1)$}}{\leq} C_1\norm{\mathcal{R}_{\textbf{a}_1}}_{\beta+6}+C_2\tau^\varrho\norm{\log[z,z]}\norm{\mathfrak{p}}_{C^\varrho}+C
\norm{\mathbb{D}}_{C^\varrho}.
\end{align}
 Here in $(1)$ by Sobolev embedding inequality \eqref{for:453} of Section \ref{sec:3}; in $(2)$ we note that $\delta_0<\eta_{\sigma+6+\beta}(\textbf{a}_1)$ (see \eqref{for:456} of Section \ref{sec:56}).
Then we use Corollary \ref{cor:10} to equation \ref{for:469}; in $(4)$
we use
\begin{align}\label{for:485}
 \emph{\textbf{d}}(z_{1}z_v,z_vz_{1})\overset{\text{(a)}}{=}\textbf{\emph{d}}([z_{1},z_v],e)\overset{\text{(b)}}{\leq} C\norm{\log[z_{1},z_v] }
\end{align}
here in $(a)$ we use right invariance of $\emph{\textbf{d}}$ (see \eqref{for:67} of Section \ref{sec:3}); in $(b)$ we use \eqref{for:35} of Section \ref{sec:3} as $\delta_0$ is sufficiently small; in $(4)$ we use \eqref{for:90} and the definition of $\varrho$ in \eqref{for:236} of \eqref{for:432} of Section \ref{sec:3}.

\smallskip
Next, we estimate $\norm{\mathcal{R}_{\textbf{a}_1}}_{\beta+6}$. By interpolation inequalities we have
\begin{align}\label{for:154}
 &\norm{\mathcal{R}_{\textbf{a}_1}}_{\beta+6}\leq C(\norm{\mathcal{R}_{\textbf{a}_1}})^{\frac{1}{2}}(\norm{\mathcal{R}_{\textbf{a}_1}}_{2\beta+12})^{\frac{1}{2}}\notag\\
 &\overset{\text{\tiny$(1)$}}{\leq} C_1\tau^{\frac{1}{2}\varrho}\big(\norm{\log[z,z]}^{c_0}\norm{\mathfrak{p}}_{\varrho}
 +\norm{\mathbb{D}}_\varrho+C_{m}\tau^{-m+\varrho}\norm{\mathfrak{p}}_m\big)^{\frac{1}{8}}\notag\\
 &\cdot(\norm{\mathfrak{p}}_{\varrho})^{\frac{3}{8}}(\norm{\mathcal{R}_{\textbf{a}_1}}_{2\beta+12})^{\frac{1}{2}}\notag\\
 &\overset{\text{\tiny$(2)$}}{\leq} C_2\tau^{\varrho}\big(\norm{\log[z,z]}^{c_0}\norm{\mathfrak{p}}_{\varrho}
 +\norm{\mathbb{D}}_\varrho+C_{m}\tau^{-m+\varrho}\norm{\mathfrak{p}}_m\big)^{\frac{1}{8}}(\norm{\mathfrak{p}}_{\varrho})^{\frac{7}{8}}\notag\\
 &\overset{\text{\tiny$(3)$}}{\leq} C_2\tau^{\varrho}\norm{\log[z,z]}^{\frac{c_0}{8}}\norm{\mathfrak{p}}_{\varrho}
 +C_2\tau^{\varrho}(\norm{\mathbb{D}}_\varrho)^{\frac{1}{8}}(\norm{\mathfrak{p}}_{\varrho})^{\frac{7}{8}}\notag\\
 &+C_2\tau^{\varrho}(C_{m}\tau^{-m+\varrho}\norm{\mathfrak{p}}_m)^{\frac{1}{8}}(\norm{\mathfrak{p}}_{\varrho})^{\frac{7}{8}}\notag\\
 &\overset{\text{\tiny$(4)$}}{\leq} C_2\tau^{\varrho}\norm{\log[z,z]}^{\frac{c_0}{8}}\norm{\mathfrak{p}}_{C^\varrho}
 +C_2\tau^{\varrho}(\norm{\mathbb{D}}_{C^\varrho})^{\frac{1}{8}}(\norm{\mathfrak{p}}_{C^\varrho})^{\frac{7}{8}}\notag\\
 &+C_2\tau^{\varrho}(C_{m}\tau^{-m+\varrho}\norm{\mathfrak{p}}_{C^m})^{\frac{1}{8}}(\norm{\mathfrak{p}}_{C^\varrho})^{\frac{7}{8}}
\end{align}
for any $m\geq\varrho$.

Here in $(1)$ we use \eqref{for:245}; in $(2)$ by \eqref{for:243} and the definition of $\varrho$ in \eqref{for:236} of \eqref{for:432} of Section \ref{sec:3} we have
\begin{align*}
 \norm{\mathcal{R}_{\textbf{a}_1}}_{2\beta+12}\leq \tau^{\varrho}\norm{\mathfrak{q}}_{\varrho};
\end{align*}
in $(3)$ we use the inequality:
\begin{align}\label{for:217}
 (x+y+w)^c\leq x^c+y^c+w^c,\qquad \forall\, x,\,y,\,w>0,\,\,\,0<c<1;
\end{align}
in $(4)$ we use Sobolev embedding inequality \eqref{for:179}.

It follows from \eqref{for:186} and \eqref{for:154} that
\begin{align*}
 \norm{\mathcal{R}_v}&_{C^0}\leq C\tau^{\varrho}\norm{\log[z,z]}^{\frac{c_0}{8}}\norm{\mathfrak{p}}_{C^\varrho}
 +C\tau^{\varrho}(\norm{\mathbb{D}}_{C^\varrho})^{\frac{1}{8}}(\norm{\mathfrak{p}}_{C^\varrho})^{\frac{7}{8}}\notag\\
 &+C\tau^{\varrho}(C_{m}\tau^{-m+\varrho}\norm{\mathfrak{q}}_{C^m})^{\frac{1}{8}}(\norm{\mathfrak{p}}_{C^\varrho})^{\frac{7}{8}}\notag\\
 &+C\tau^\varrho\norm{\log[z,z]}\norm{\mathfrak{p}}_{C^\varrho}+C
\norm{\mathbb{D}}_{C^{\varrho}}\\
&\overset{\text{\tiny$(1)$}}{\leq} 2C\tau^{\varrho}\norm{\log[z,z]}^{\frac{c_0}{8}}\norm{\mathfrak{p}}_{C^\varrho}
 +C\tau^{\varrho}(\norm{\mathbb{D}}_{C^\varrho})^{\frac{1}{8}}(\norm{\mathfrak{p}}_{C^\varrho})^{\frac{7}{8}}\notag\\
 &+C\tau^{\varrho}(C_{m}\tau^{-m+\varrho}\norm{\mathfrak{p}}_{C^m})^{\frac{1}{8}}(\norm{\mathfrak{p}}_{C^\varrho})^{\frac{7}{8}}+C
\norm{\mathbb{D}}_{C^{\varrho}}.
\end{align*}
for any $m\geq\varrho$.

Here in $(1)$ we note that $\norm{\log[z,z]}<1$ if $\delta_0$ is sufficiently small and note that $c_0<1$ (see \eqref{for:470}  of \eqref{for:432} of Section \ref{sec:3}).

Then we get \eqref{for:121}. Thus we finish the proof.

\subsubsection{Revised version of Theorem \ref{th:9}}

 In this part, similarly to Corollary \ref{cor:8}, we use the smoothing operators $\mathfrak{s}_{b}$ (see \eqref{for:191} of Section \ref{sec:3}) to solve the fixed loss of regularity in solving the twisted almost cocycle equations in Theorem \ref{th:9}. As a result, we also introduce two more parameters $b$ and $\ell$. We also obtain  adjoint $L$-invariant (see \eqref{for:117} of Section \ref{sec:3}) approximation for adjoint $L$-invariant almost coboundaries that is needed in the iteration process.
\begin{corollary}\label{cor:7}
There is $\delta_0>0$ such that: for any $\delta_0$ standard perturbation $E'(z)$ of $E$, any $b,\,\tau>1$ and
any adjoint $L$-invariant (see \eqref{for:117} of Section \ref{sec:3}) $\mathfrak{p}_v\in \mathcal{L}^\bot(\mathcal{O}^\infty)$, $v\in E$,
there exist  adjoint $L$-invariant $\Theta\in C^\infty(\mathcal{M},\,\mathcal{L}^\bot)$ (see \eqref{for:17} of Section \ref{sec:3}) and  adjoint $L$-invariant $\mathcal{R}_v\in C^\infty(\mathcal{M},\,\mathcal{L}^\bot)$, $v\in E$ such that
\begin{align*}
 \mathfrak{p}_v=\Theta\circ z_v-\mathfrak{i}_{\mathcal{L}^\bot}\text{Ad}_{z_{v}}\Theta+\mathcal{R}_v,\quad \forall\,v\in E
\end{align*}
with estimates: for any $v\in E$
\begin{align}\label{for:81}
 \max\{\norm{\Theta}_{C^m},\,\norm{\mathcal{R}_v}_{C^m}\}
 \leq C_{m}b^\varrho(\tau^{m}\norm{\mathfrak{p}}_{C^\varrho}+\norm{\mathfrak{p}}_{C^{m}})
\end{align}
for any $m\geq\varrho$;  and
\begin{align}\label{for:83}
 \norm{\mathcal{R}_v}&_{C^0}\leq C\tau^{\varrho}\norm{\log[z,z]}^{\frac{c_0}{8}}\norm{\mathfrak{p}}_{C^\varrho}
 +C\tau^{\varrho}(\norm{\mathbb{D}}_{C^\varrho})^{\frac{1}{8}}(\norm{\mathfrak{p}}_{C^\varrho})^{\frac{7}{8}}\notag\\
 &+C\tau^{\varrho}(C_{m}\tau^{-m+\varrho}\norm{\mathfrak{p}}_{C^m})^{\frac{1}{8}}(\norm{\mathfrak{p}}_{C^\varrho})^{\frac{7}{8}}+C
\norm{\mathbb{D}}_{C^{\varrho}}\notag\\
&+ C_{\ell} b^{-\ell+\varrho}(\tau^{\ell}\norm{\mathfrak{p}}_{C^\varrho}+\norm{\mathfrak{p}}_{C^{\ell}})
\end{align}
for any $m,\,\ell\geq \varrho$,  where $\mathbb{D}_{v,u}$, $\norm{\log[z,z]}$, $\norm{\mathfrak{p}}_{C^m}$, $\norm{\mathbb{D}}_{C^m}$ and
$\norm{\mathbb{D}}_{C^m}$ are as defined in Section \ref{sec:56}.
\end{corollary}

\begin{proof} By Theorem \ref{th:9}, there exist $\Theta'\in \mathcal{L}^\bot(\mathcal{O}^\infty)$ and $\mathcal{R}'_v\in \mathcal{L}^\bot(\mathcal{O}^\infty)$, $v\in E$ such that
\begin{align}\label{for:85}
 \mathfrak{q}_v=\Theta'\circ z_v-\mathfrak{i}_{\mathcal{L}^\bot}\text{Ad}_{z_{v}}\Theta'+\mathcal{R}'_v, \quad\forall\,v\in E
\end{align}
and the estimates \eqref{for:114} and \eqref{for:121} are satisfied for $\Theta'$ and $\mathcal{R}'_v$, $v\in E$.  Let
\begin{align*}
  \Theta&=\int_{L}\text{Ad}_{l}(\mathfrak{s}_{b}\Theta')\circ l^{-1}dl,\qquad \text{and}\\
  \mathcal{R}_v&=\mathfrak{q}_v-(\Theta\circ z_v-\mathfrak{i}_{\mathcal{L}^\bot}\text{Ad}_{z_{v}}\Theta ), \quad v\in E
\end{align*}
where $dl$ is the Haar measure on $L$.

\smallskip

\noindent\emph{Note}. Even though $\Theta'\in \mathcal{L}^\bot(\mathcal{O}^\infty)$, it is not necessary that $\mathfrak{s}_{b}\Theta'\in \mathcal{L}^\bot(\mathcal{O}^\infty)$ since $\int_{\mathcal{M}} \mathfrak{s}_{b}\Theta'$ may not vanish. As a consequence, $\Theta$, as well as  $\mathcal{R}_v$, $v\in E$ are not necessarily inside $\mathcal{L}^\bot(\mathcal{O}^\infty)$ either.
But it is clear that they are all $C^\infty$ maps.

\smallskip

Then we have
\begin{align*}
 \norm{\Theta}_{C^m}&\leq C_m\norm{\mathfrak{s}_{b}\Theta'}_{C^m}\overset{\text{\tiny$(1)$}}{\leq} C_{m,1}b^\varrho\norm{\Theta'}_{C^{m-\varrho}}\notag\\
 &\overset{\text{\tiny$(2)$}}{\leq} C_{m,2}b^\varrho(\tau^{m}\norm{\mathfrak{q}}_{C^\varrho}+\norm{\mathfrak{q}}_{C^{m}})
 \end{align*}
for any $m\geq\varrho$.  Here in $(1)$ we use \eqref{for:45} of Section \ref{sec:3}; in $(2)$ we use \eqref{for:114}.

It follows that: for any $v\in E$ and any $m\geq\varrho$
\begin{align*}
  \norm{\mathcal{R}_v}_{C^m}&\leq \norm{\mathfrak{p}_v}_{C^m}+\norm{\Theta\circ z_v-\mathfrak{i}_{\mathcal{L}^\bot}\text{Ad}_{z_{v}}\Theta}_m\notag\leq
  \norm{\mathfrak{q}_v}_{C^m}+C\norm{\Theta}_{C^m}\\
  &\leq  C_{m}b^\varrho(\tau^{m}\norm{\mathfrak{p}}_{C^\varrho}+\norm{\mathfrak{p}}_{C^{m}}).
\end{align*}
Hence we get \eqref{for:81}.

Next, we show how to estimate $\norm{\mathcal{R}_v}_{C^0}$, $v\in E$.  From \eqref{for:85}, we have
\begin{align*}
 \mathfrak{p}_v=(\mathfrak{s}_{b}\Theta')\circ z_v-\mathfrak{i}_{\mathcal{L}^\bot}\text{Ad}_{z_{v}}\mathfrak{s}_{b}\Theta'+T_v, \quad\forall\,v\in E
\end{align*}
where
\begin{align*}
 T_v=\mathcal{R}'_v+\big((I-\mathfrak{s}_{b})\Theta'\big)\circ z_v-\mathfrak{i}_{\mathcal{L}^\bot}\text{Ad}_{z_{v}}\big((I-\mathfrak{s}_{b})\Theta'\big),\quad \forall\,v\in E.
\end{align*}
As $\mathfrak{q}_v$ is adjoint $L$-invariant and $z_v\in Z(L)$ (we recall $E'(z)$ is a standard perturbation of $E$, see Section \ref{sec:1}) for each $v\in E$, we see that
\begin{align*}
 \mathcal{R}_v=\int_{L}\text{Ad}_{l}T_v\circ l^{-1}dl,\qquad \forall\,v\in E.
\end{align*}
It follows that
\begin{align}\label{for:86}
 \norm{\mathcal{R}_v}_{C^0}&\leq C\norm{T_v}_{C^0}\leq C\norm{\mathcal{R}'_v}_{C^0}+C\norm{(I-\mathfrak{s}_{b})\Theta'}_{C^0}\notag\\
 &\overset{\text{\tiny$(1)$}}{\leq} C\norm{\mathcal{R}'_v}_{C^0}+C_\ell b^{-\ell+\varrho}\norm{\Theta'}_{C^{\ell-\varrho}}\notag\\
 &\overset{\text{\tiny$(2)$}}{\leq} C\norm{\mathcal{R}'_v}_{C^0}+C_{\ell,1} b^{-\ell+\varrho}(\tau^{\ell}\norm{\mathfrak{p}}_{C^\varrho}+\norm{\mathfrak{p}}_{C^{\ell}}).
\end{align}
Here in $(1)$ we use \eqref{for:471} of Section \ref{sec:3}; in $(2)$ we use \eqref{for:114}.

 Hence \eqref{for:83} follows from \eqref{for:86} and \eqref{for:121}.
\end{proof}

\subsection{ Iterative step and error estimate}\label{sec:29} In this part for any given  sufficiently small $C^\infty$ perturbation $\tilde{\alpha}_A$
of the action $\alpha_A$, which is also close to $\alpha_{E'(z)}$ on $E$ where $E'(z)$ is a standard perturbation of $E$,
we construct a $C^\infty$ diffeomorphism $h$ and a new standard perturbation $E'(z^{(1)})$ of $E$;  and then estimate the closeness between the new action $\tilde{\alpha}'_A=h^{-1}\tilde{\alpha}h$ and $\alpha_{E'(z^{(1)})}$ on $E$.

We recall notations in Section \ref{sec:1}.  We point out that the positive constants $\varrho$, $c_0$, $\delta_0$  and $\delta_\varrho$ (see Section \ref{for:432}) will be used in this part.

\begin{proposition}\label{po:5} There exists $0<\bar{c}<1$ such that the following holds: if
$\tilde{\alpha}_A$ is a $C^\infty$ small perturbation of $\alpha_A$ on $\mathcal{X}$  which can be written as
\begin{align}\label{for:123}
 \tilde{\alpha}_{v}(\mathfrak{j}(x))=\mathfrak{j}\big(\exp(\mathfrak{p}_v(x))\cdot
z_v\cdot x\big),\qquad \forall\,\,v\in E,\quad \forall\, x\in \mathcal{M}
\end{align}
where $z_v=w_v\cdot \alpha_v$, $w_v\in \exp(\mathfrak{g}_{1,A})\cap Z(L)$ and $\mathfrak{p}_v\in C^\infty(\mathcal{M},\,\mathcal{L}^\bot)$
satisfying
\begin{align*}
 \norm{\mathfrak{p}}_{C^{\varrho+1}}:=\max_{z\in E}\norm{\mathfrak{p}_z}_{C^{\varrho+1}}\leq \bar{c}
 \quad\text{and}\quad\textbf{d}(w,e)=\max_{v\in E}\textbf{d}(w_v,e)\leq \bar{c},
\end{align*}
then for any $b,\,\tau>1$,  there is $h\in C^\infty(\mathcal{X})$
such that:
\begin{enumerate}
  \item\label{for:101}   for any $m\geq\varrho$
  \begin{align*}
   d(h, I)_{C^m}\leq C_{m}b^\varrho(\tau^{m}\norm{\mathfrak{p}}_{C^\varrho}+\norm{\mathfrak{p}}_{C^{m}});
  \end{align*}
if $b,\,\tau$ are well chosen such that
 \begin{align*}
   b^\varrho\tau^{\varrho}\norm{\mathfrak{p}}_{C^{\varrho}}< \bar{c},
 \end{align*}
then $h$ is invertible;

\smallskip
\item\label{for:112} for $\tilde{\alpha}^{(1)}=h^{-1}\tilde{\alpha}h$, there exist $\mathfrak{p}_v^{(1)}\in C^\infty(\mathcal{M},\,\mathcal{L}^\bot)$ and $w_v^{(1)}\in \exp(\mathfrak{g}_{1,A})\cap Z(L)$, $v\in E$ such that
\begin{align*}
 \tilde{\alpha}_{v}^{(1)}(\mathfrak{j}(x))=\mathfrak{j}\big(\exp(\mathfrak{p}^{(1)}_v(x))\cdot
z_v^{(1)}\cdot x\big),\qquad \forall\,v\in E,\,x\in \mathcal{M}
\end{align*}
where $z_v^{(1)}=w_v^{(1)}\cdot z_v$;  and
\begin{align*}
 \max_{z\in E} \textbf{d}(z_v^{(1)},\,z_v)<C\norm{\mathfrak{p}}_{C^0};
\end{align*}

 \item\label{for:119} for any $m\geq\varrho$ we have:
 \begin{align*}
  \norm{\mathfrak{p}^{(1)}}_{C^m}\leq C_{m}(\tau^{m}b^\varrho\norm{\mathfrak{p}}_{C^\varrho}+b^\varrho\norm{\mathfrak{p}}_{C^{m}}+1);
 \end{align*}

\item\label{for:69} for $\mathfrak{p}^{(1)}$ the following estimates hold:
 \begin{align*}
\norm{\mathfrak{p}^{(1)}}_{C^0}&\leq Cb^{2\varrho}\tau^{2\varrho}(\norm{\mathfrak{p}}_{C^{\varrho+1}})^{1+y_0}\notag\\
&+C\tau^{\varrho}(C_{m}\tau^{-m+\varrho}\norm{\mathfrak{p}}_{C^m})^{\frac{1}{8}}
(\norm{\mathfrak{p}}_{C^\varrho})^{\frac{7}{8}}\notag\\
&+C_{\ell} b^{-\ell+\varrho}(\tau^{\ell}\norm{\mathfrak{p}}_{C^\varrho}+\norm{\mathfrak{p}}_{C^{\ell}});
\end{align*}
for any $m,\,\ell\geq\varrho$, where $y_0=\frac{c_0}{8}$  and
\begin{align*}
\norm{\mathfrak{p}}_{C^m}=\max_{v\in E}\{\norm{\mathfrak{p}_{v}}_{C^m}\}\quad\text{and}\quad
\norm{\mathfrak{p}^{(1)}}_{C^m}:=\max_{v\in E}\norm{\mathfrak{p}^{(1)}_v}_{C^m};
\end{align*}

 \end{enumerate}
\end{proposition}
\begin{remark}\label{re:5} Comparing  Proposition \ref{po:7} and Proposition \ref{po:5}, we see that the same estimates hold for $ d(h, I)_{C^m}$,
$\norm{\mathfrak{p}^{(1)}}_{C^0}$ and $\norm{\mathfrak{p}^{(1)}}_{C^m}$, except that $y_0$ takes different values.
\end{remark}

\begin{proof} We follow the proof line of Proposition \ref{po:7}. Typically distinctions between nilmanifold examples and
(twisted) symmetric space examples are minor, and appear mostly at notations.

\smallskip
\emph{Basic notations}:

Let
\begin{gather}
  r_v=\int_{\mathcal{M}}\mathfrak{p}_v,\quad \mathfrak{p}'_v=\mathfrak{p}_{v}-r_v, \qquad\text{and}\notag\\
   \mathbb{D}_{v,u}(x)=\mathfrak{p}'_{v}\circ z_{u}-\mathfrak{i}_{\mathcal{L}^\bot}\text{Ad}_{z_{u}}\mathfrak{p}'_{v}-\big(\mathfrak{p}'_{u}\circ z_{v}-\mathfrak{i}_{\mathcal{L}^\bot}\text{Ad}_{z_{v}}\mathfrak{p}'_{u}\big)\label{for:38}
\end{gather}
for any $v,\,u\in E$; and denote
\begin{align*}
 \norm{\log[z,z]}=\max_{u,v\in E}\norm{\log[z_{v},z_{u}]}\quad\text{and}\quad \norm{\mathbb{D}}_{C^m}=\max_{u,\,v\in E}\{\norm{\mathbb{D}_{v,u}}_{C^m}\}.
\end{align*}
We can assume $\bar{c}<\delta_\varrho$ (see Lemma \ref{le:8}). By Lemma \ref{le:8} we have
\begin{align}
\norm{\log[z,z]}\leq C\norm{\mathfrak{p}}_{C^0}, \quad\text{and}\quad
 \norm{\mathbb{D}}_{C^\varrho}\leq C\norm{\mathfrak{p}}^2_{C^{\varrho+1}}. \label{for:253}
\end{align}
Let
\begin{align}\label{for:478}
 c_{\mu}=\big((I-\text{Ad}_{z_{\textbf{a}_1}})|_{\mathfrak{g}_{\mu,A}}\big)^{-1}(r_{\textbf{a}_1}|_{\mathfrak{g}_{\mu,A}});\quad 1\neq\mu\in \Delta_A, \,\,\text{and}\,\,
 c_{1}=0;
\end{align}
and let $c_{\textbf{a}_1}\in \mathfrak{G}$ be the vector with coordinates $c_{\mu}$.

\smallskip

\emph{Construction of $\Theta$ and $\mathcal{R}_v$, $v\in E$}:

 We note that $\mathfrak{p}_v$ is adjoint $L$-invariant for each $v\in E$ (see Lemma \ref{ob:1}). Obviously,
 $\mathfrak{p}'_v$ is also adjoint $L$-invariant for each $v\in E$. By  applying Corollary \ref{cor:7} to almost twisted cocycle equation \eqref{for:38}, we have:
for any $b,\,\tau>1$ there exist adjoint $L$-invariant $\Theta\in C^\infty(\mathcal{M},\,\mathcal{L}^\bot)$ and $\mathcal{R}_v\in C^\infty(\mathcal{M},\,\mathcal{L}^\bot)$ for each $v\in E$ such that
\begin{align}\label{for:92}
 \mathfrak{p}'_v=\Theta\circ z_v-\mathfrak{i}_{\mathcal{L}^\bot}\text{Ad}_{z_{v}}\Theta+\mathcal{R}_v
\end{align}
with estimates: for any $v\in E$
\begin{align}\label{for:44}
 \max\{\norm{\Theta}_{C^m},\,\norm{\mathcal{R}_v}_{C^m}\}
 \leq C_{m}b^\varrho(\tau^{m}\norm{\mathfrak{p}}_{C^\varrho}+\norm{\mathfrak{p}}_{C^{m}})
\end{align}
for any $m\geq\varrho$;  and
\begin{align}\label{for:96}
 \norm{\mathcal{R}_v}&_{C^0}\leq C\tau^{\varrho}\norm{\log[z,z]}^{\frac{c_0}{8}}\norm{\mathfrak{p}}_{C^\varrho}
 +C\tau^{\varrho}(\norm{\mathbb{D}}_{C^\varrho})^{\frac{1}{8}}(\norm{\mathfrak{p}}_{C^\varrho})^{\frac{7}{8}}\notag\\
 &+C\tau^{\varrho}(C_{m}\tau^{-m+\varrho}\norm{\mathfrak{p}}_{C^m})^{\frac{1}{8}}(\norm{\mathfrak{p}}_{C^\varrho})^{\frac{7}{8}}+C
\norm{\mathbb{D}}_{C^{\varrho}}\notag\\
&+C_{\ell} b^{-\ell+\varrho}(\tau^{\ell}\norm{\mathfrak{p}}_{C^\varrho}+\norm{\mathfrak{p}}_{C^{\ell}})\notag\\
&\overset{\text{\tiny$(1)$}}{\leq}C_1\tau^{\varrho}(\norm{\mathfrak{p}}_{C^{0}})^{\frac{c_0}{8}}\norm{\mathfrak{p}}_{C^\varrho}
 +C_1\tau^{\varrho}(\norm{\mathfrak{p}}_{C^{\varrho+1}})^{\frac{2}{8}}(\norm{\mathfrak{p}}_{C^\varrho})^{\frac{7}{8}}\notag\\
 &+C\tau^{\varrho}(C_{m}\tau^{-m+\varrho}\norm{\mathfrak{p}}_{C^m})^{\frac{1}{8}}(\norm{\mathfrak{p}}_{C^\varrho})^{\frac{7}{8}}+C
\norm{\mathfrak{p}}^2_{C^{\varrho+1}}\notag\\
&+C_{\ell} b^{-\ell+\varrho}(\tau^{\ell}\norm{\mathfrak{p}}_{C^\varrho}+\norm{\mathfrak{p}}_{C^{\ell}})\notag\\
&\overset{\text{\tiny$(2)$}}{\leq}3C\tau^{\varrho}(\norm{\mathfrak{p}}_{C^{\varrho+1}})^{1+\frac{c_0}{8}}
+C\tau^{\varrho}(C_{m}\tau^{-m+\varrho}\norm{\mathfrak{p}}_{C^m})^{\frac{1}{8}}(\norm{\mathfrak{p}}_{C^\varrho})^{\frac{7}{8}}\notag\\
&+C_{\ell} b^{-\ell+\varrho}(\tau^{\ell}\norm{\mathfrak{p}}_{C^\varrho}+\norm{\mathfrak{p}}_{C^{\ell}})
\end{align}
for any $m\geq \varrho$ and $\ell>0$. Here in $(1)$ we use \eqref{for:253}; in $(2)$ we use the assumption $\norm{\mathfrak{p}}_{C^{\varrho+1}}\leq \bar{c}<1$ and the fact $0<c_0<1$.

\smallskip

\emph{Construction of $h$ and $z^{(1)}_v$, $v\in E$}:

Let $\mathfrak{h}=\Theta+c_{\textbf{a}_1}$
and let
\begin{gather*}
 \tilde{h}(x)=\exp(\mathfrak{h}(x))x,\quad \forall\,x\in \mathcal{M},\\
 w_v^{(1)}=\exp(r_v|_{\mathfrak{g}_{1,A}}),\quad z_v^{(1)}=w_v^{(1)}z_v,\,\,v\in E.
\end{gather*}
Since $\mathfrak{h}(x)$ is adjoint $L$-invariant, $\tilde{h}(x)$ descends to a map $h$ on $\mathcal{X}$:
\begin{align}
 h(\mathfrak{j}(x))=\mathfrak{j}(\tilde{h}(x)),\qquad \forall\,x\in x\in \mathcal{M}.
\end{align}

\emph{Note}.  $w_v^{(1)}$, $v\in E$ is the coordinate change, which is determined by $r_v$, $v\in E$.

\smallskip
\emph{Estimates of $c_{\textbf{a}_1}$ and $w^{(1)}_v$, $v\in E$}:

From \eqref{for:472} of Proposition \ref{po:9} we have
\begin{align}\label{for:43}
  \Big\|\big((I-\text{Ad}_{z_{\textbf{a}_1}})|_{\mathfrak{g}_{\mu,A}}\big)^{-1}\Big\|\leq C,\qquad \forall\,1\neq\mu\in \Delta_A.
\end{align}
 \eqref{for:43} shows that
\begin{align}\label{for:66}
 \norm{c_{\textbf{a}_1}}\leq C\norm{\mathfrak{p}}_{C^0}.
\end{align}
We note that
\begin{align}\label{for:122}
  \max_{v\in E}\emph{\textbf{d}}&(z_v^{(1)}, v)=\max_{v\in E}\emph{\textbf{d}}(w_v^{(1)}z_v, v)\notag\\
  &\leq \max_{v\in E}\emph{\textbf{d}}(w_v^{(1)}z_v, z_v)+\max_{v\in E}\emph{\textbf{d}}(z_v, v)\notag\\
  &\overset{\text{\tiny$(1)$}}{=} \max_{v\in E}\emph{\textbf{d}}(w_v^{(1)}, e)+\max_{v\in E}\emph{\textbf{d}}(z_v, v)\notag\\
   &\overset{\text{\tiny$(2)$}}{\leq} \max_{v\in E}C\norm{r_v|_{\mathfrak{g}_{1,A}}}+\max_{v\in E}\emph{\textbf{d}}(z_v, v)\notag\\
  &\leq C\norm{\mathfrak{p}}_{C^0}+\varepsilon.
\end{align}
Here in $(1)$ we use right invariance of $\emph{\textbf{d}}$ (see \eqref{for:67} of Section \ref{sec:3}); here in $(2)$ we use
\eqref{for:35} of Section \ref{sec:3}.

From \eqref{for:122}, we also have
\begin{align}\label{for:476}
 \max_{z\in E} \textbf{\emph{d}}(w_v^{(1)},\,e)=\max_{v\in E}\emph{\textbf{d}}(z_v^{(1)}, z_v)\leq C\norm{\mathfrak{p}}_{C^0}.
\end{align}
Next, we show that
\begin{align}\label{for:473}
 w_v^{(1)}\in \exp(\mathfrak{g}_{1,A})\cap Z(L).
\end{align}
It is clear that $w_v^{(1)}\in \exp(\mathfrak{g}_{1,A})$. Next, we show that $w_v^{(1)}\in  Z(L)$.  Since $\mathfrak{p}_v$ is adjoint $L$-invariant, we have
\begin{align*}
 \text{Ad}_{l}r_v=r_v, \text{ for }\forall\,l\in L\overset{\text{\tiny$(1)$}}{\Longrightarrow} \text{Ad}_{l}(r_v|_{\mathfrak{g}_{1,A}})=r_v|_{\mathfrak{g}_{1,A}}, \text{ for }\forall\,l\in L.
\end{align*}
Here in $(1)$ we note that $\mathfrak{g}_{\mu,A}$ is invariant under $\text{Ad}_{L}$ for each $\mu\in \Delta_A$ (see \eqref{for:54} of Section \ref{sec:3}).
This implies that $w_v\in \exp(\mathfrak{g}_{1,A})\cap Z(L)$.

Next, we claim that
\begin{align}\label{for:474}
 r_v|_{\mathfrak{g}_{1,A}}\in \mathcal{L}^\bot,
\end{align}
which will be used later. We note that
\begin{align*}
 r_v|_{\mathfrak{g}_{1,A}}=r_v-\sum_{1\neq\mu\in \Delta_A}r_v|_{\mathfrak{g}_{\mu,A}}
\end{align*}
Then \eqref{for:474} follows from that facts that
$r_v\in \mathcal{L}^\bot$ and $\mathfrak{g}_{\mu,A}\in \mathcal{L}^\bot$, $1\neq\mu\in \Delta_A$ (see \eqref{for:54} of Section \ref{sec:3}).

\smallskip
\eqref{for:101}: We have
\begin{align}\label{for:79}
 d_\mathcal{X}&(h, I)_{C^m}\leq d_\mathcal{M}(\tilde{h}, I)_{C^m}\leq C_m\norm{\mathfrak{h}}_{C^m}\notag\\
 &\leq C_{m,1}\norm{\Theta}_{C^m}+C_{m,1}\norm{c_{\textbf{a}_1}}\notag\\
 &\overset{\text{\tiny$(1)$}}{\leq} C_{m}b^\varrho(\tau^{m}\norm{\mathfrak{p}}_{C^\varrho}+\norm{\mathfrak{p}}_{C^{m}})
 \end{align}
for any $m\geq\varrho$.  Here in $(1)$ we use \eqref{for:44} and \eqref{for:66}.

In \eqref{for:79} letting $m=\varrho$, we have
\begin{align}\label{for:477}
  d_\mathcal{X}(h, I)_{C^\varrho}\leq  d_\mathcal{M}(\tilde{h}, I)_{C^\varrho}\leq Cb^\varrho\tau^{\varrho}\norm{\mathfrak{p}}_{C^\varrho}<C\bar{c},
\end{align}
We can assume that $\bar{c}$ is sufficiently small such that $h$ and $\tilde{h}$ are both invertible; moreover,
\begin{align}
  d_\mathcal{X}(h^{-1}, I)_{C^m}&\leq C_md_\mathcal{X}(h, I)_{C^m} \qquad\text{and}\label{for:80}\\
  d_\mathcal{M}(\tilde{h}^{-1}, I)_{C^m}&\leq C_md_\mathcal{M}(\tilde{h}, I)_{C^m}
\end{align}
for any $m\geq 0$, see \cite[Lemma AII.26]{Lazutkin}.

\smallskip
\eqref{for:112}: \eqref{for:473} together with \eqref{for:122} show that $E'(z^{(1)})$ is a standard $(\varepsilon+C\norm{\mathfrak{p}}_{C^0})$-perturbation of $E$.

\eqref{for:80} justifies the definition of $\tilde{\alpha}_A^{(1)}$. For any $v\in E$ we have
\begin{align}\label{for:272}
  d_\mathcal{X}(\tilde{\alpha}^{(1)}_v,& \,\alpha_{z_v^{(1)}})_{C^0}=d_\mathcal{X}(h^{-1}\circ \tilde{\alpha}_v\circ h, \,\alpha_{z_v^{(1)}})_{C^0} \notag\\
  &\leq d_\mathcal{X}(h^{-1}\circ \tilde{\alpha}_v\circ h, \tilde{\alpha}_v\circ h)_{C^0}+d_\mathcal{X}(\tilde{\alpha}_v\circ h, \alpha_{z_v}\circ h)_{C^0}\notag\\
  &+d_\mathcal{X}(\alpha_{z_v}\circ h, \alpha_{z_v})_{C^0}+d_\mathcal{X}(\alpha_{z_v}, \alpha_{z_v^{(1)}})_{C^0}\notag\\
  &\leq d_\mathcal{X}(h^{-1}, I)_{C^0}+d_\mathcal{X}(\tilde{\alpha}_v, \alpha_{z_v})_{C^0}+\norm{\alpha_{z_v}}_{C^1(\mathcal{X})}d_\mathcal{X}(h, I)_{C^0}\notag\\
  &+C\textbf{\emph{d}}(z_v,\,z_v^{(1)})\notag\\
  &\overset{\text{\tiny$(1)$}}{\leq} Cd_\mathcal{X}(h, I)_{C^0}+C\norm{\mathfrak{p}}_{C^0}+Cd_\mathcal{X}(h, I)_{C^0}+C\norm{\mathfrak{p}}_{C^0}\notag\\
  &\overset{\text{\tiny$(2)$}}{\leq}C_1b^\varrho\tau^{\varrho}\norm{\mathfrak{p}}_{C^\varrho}\leq C_2\bar{c}.
\end{align}
Here in $(1)$ we use \eqref{for:80} and \eqref{for:476};  $(2)$ we use \eqref{for:477}.

Hence, if $\bar{c}$ is sufficiently small there exists $\mathfrak{p}_v^{(1)}\in C^\infty(\mathcal{M},\,\mathcal{L}^\bot)$  for each $v\in E$ such that
\begin{align}\label{for:124}
  \tilde{\alpha}_{v}^{(1)}(\mathfrak{j}(x))=\mathfrak{j}\big(\exp(\mathfrak{p}^{(1)}_v(x))\cdot
z_v^{(1)}\cdot x\big),\qquad \forall\,\,v\in E,\quad \forall\, x\in \mathcal{M}.
\end{align}
The above discussion with \eqref{for:476} complete the proof of \eqref{for:112}.

\medskip
\eqref{for:119}: From \eqref{for:124} we have
\begin{align}
 \norm{\mathfrak{p}^{(1)}}_{C^m}&\leq C_m\max_{v\in E}d_\mathcal{X}(\tilde{\alpha}^{(1)}_v, \,\alpha_{z_v^{(1)}})_{C^m}\leq C_m(\max_{v\in E}\norm{h^{-1}\tilde{\alpha}_vh}_{C^m}+1)\notag\\
 &\overset{\text{\tiny$(1)$}}{\leq} C_{m,1}(\norm{h^{-1}}_{C^m}+\norm{h}_{C^m}+\max_{v\in E}\norm{\tilde{\alpha}_v}_{C^m}+1)\notag\\
 &\overset{\text{\tiny$(2)$}}{\leq} C_{m,2}(\norm{h}_{C^m}+\norm{\mathfrak{p}}_{C^m}+1)\notag\\
 &\overset{\text{\tiny$(3)$}}{\leq}C_{m,3}(\tau^{m}b^\varrho\norm{\mathfrak{p}}_{C^\varrho}+b^\varrho\norm{\mathfrak{p}}_{C^{m}}+1)
\end{align}
for any $m\geq\varrho$.

Here $(1)$ is a standard ``linear" estimate, see for example, \cite[Propositions 5.5]{Llave}; in $(2)$
we use \eqref{for:80}; in $(3)$ we use \eqref{for:79}.

\medskip

\eqref{for:69}: Let
\begin{align}\label{for:88}
  \beta_{v}(x)=\exp(\mathfrak{p}^{(1)}_v(x))\cdot
z_v^{(1)}\cdot x\quad\text{and}\quad \gamma_{v}(x)=\exp(\mathfrak{p}_v(x))\cdot
z_v\cdot x
\end{align}
for any $v\in E$ and any $x\in \mathcal{M}$.

Since  $\tilde{h}$ is a lift of $h$ on $\mathcal{M}$, $\beta_v$ is a lift of $\tilde{\alpha}_v^{(1)}$ and $\gamma_v$ is a lift of $\tilde{\alpha}_v$, $v\in E$, by using $h\circ \tilde{\alpha}^{(1)}_v=\tilde{\alpha}_v\circ h$, we have
\begin{align*}
 \tilde{h}\circ \beta_v(x)=l_v(x)\cdot\tilde{\alpha}_v\circ \tilde{h}(x),\qquad \forall\,x\in \mathcal{M},\,\,\forall\,v\in E,
\end{align*}
where $l_v:\mathcal{M}\to L$,  $v\in E$ are small maps.

It follows that
\begin{align*}
&\exp\big(\mathfrak{h}(\beta_{v}(x))\big)\cdot\exp(\mathfrak{p}^{(1)}_v(x))\cdot
\exp(r_v|_{\mathfrak{g}_{1,A}})z_v\cdot x\\
&=l_v(x)\exp\big(\mathfrak{p}_v(\tilde{h}(x))\big)\cdot
z_v\cdot\big(\exp(\mathfrak{h}(x))\cdot x\big)
\end{align*}
for any $v\in E$ and any $x\in \mathcal{M}$.

We recall $z_v^{(1)}=w_v^{(1)}\cdot z_v$ and $\log w_v^{(1)}=r_v|_{\mathfrak{g}_{1,A}}$. Then the above equation  further gives
\begin{align}
&\exp\big(\mathfrak{h}(\beta_{v}(x))\big)\cdot\exp(\mathfrak{p}^{(1)}_v(x))\cdot \exp(r_v|_{\mathfrak{g}_{1,A}})\notag\\
&=l_v(x)\exp\big(\mathfrak{p}_v(\tilde{h}(x))\big)\cdot
\exp(\text{Ad}_{z_v}\mathfrak{h}(x))
\end{align}
for any $v\in E$ and any $x\in \mathcal{M}$.

By  using Baker-Campbell-Hausdorff formula we have
\begin{align}\label{for:91}
 \log l_v&=\mathfrak{h}\circ \beta_{v}+\mathfrak{p}^{(1)}_v+r_v|_{\mathfrak{g}_{1,A}}-\mathfrak{p}_v\circ \tilde{h}-\text{Ad}_{z_v}\mathfrak{h}+\Psi_v
\end{align}
where
\begin{align}\label{for:111}
 \max_{v\in E}\big\|\Psi_v\big\|_{C^0}&\leq C\max_{v\in E}\{\norm{\mathfrak{h}}_{C^0}, \norm{\mathfrak{p}}_{C^0},\norm{\mathfrak{p}^{(1)}}_{C^0},\norm{r_v}\}^2\notag\\
 &\overset{\text{\tiny$(1)$}}{\leq}Cb^{2\varrho}\tau^{2\varrho}(\norm{\mathfrak{p}}_{C^\varrho})^2.
\end{align}
Here in $(1)$ we use \eqref{for:79} and \eqref{for:272}.

By taking projection to $\mathcal{L}^\bot$ on each side of \eqref{for:91},  we have
\begin{align*}
 0=\mathfrak{h}\circ \beta_{v}+\mathfrak{p}^{(1)}_v+r_v|_{\mathfrak{g}_{1,A}}-\mathfrak{p}_v\circ \tilde{h}-\mathfrak{i}_{\mathcal{L}^\bot}\text{Ad}_{z_v}
 \mathfrak{h}+\mathfrak{i}_{\mathcal{L}^\bot}\Psi_v,
\end{align*}
by noting that $\mathfrak{h}$, $\mathfrak{p}_v$ and $\mathfrak{p}^{(1)}_v$ take values inside $\mathcal{L}^\bot$ and $r_v|_{\mathfrak{g}_{1,A}}\in \mathcal{L}^\bot$ (see \eqref{for:474}).

It follows that
\begin{align}\label{for:479}
 \mathfrak{p}^{(1)}_v&=-\mathfrak{h}\circ \beta_{v}-r_v|_{\mathfrak{g}_{1,A}}+\mathfrak{p}_v\circ \tilde{h}+\mathfrak{i}_{\mathcal{L}^\bot}\text{Ad}_{z_v}
 \mathfrak{h}-\mathfrak{i}_{\mathcal{L}^\bot}\Psi_v\notag\\
 &=\big(-\mathfrak{h}\circ \beta_{v}+\mathfrak{h}\circ z_v^{(1)}\big)+\big(-\mathfrak{h}\circ z_v^{(1)}+\mathfrak{h}\circ z_v \big)\notag\\
 &+\big(-\mathfrak{h}\circ z_v-r_v|_{\mathfrak{g}_{1,A}}+\mathfrak{p}_v+\mathfrak{i}_{\mathcal{L}^\bot}\text{Ad}_{z_v}
 \mathfrak{h} \big)\notag\\
 &+\big(\mathfrak{p}_v\circ \tilde{h}-\mathfrak{p}_v\big)-\mathfrak{i}_{\mathcal{L}^\bot}\Psi_v.
\end{align}
We note that
\begin{align}\label{for:480}
&-\mathfrak{h}\circ z_v-r_v|_{\mathfrak{g}_{1,A}}+\mathfrak{p}_v+\mathfrak{i}_{\mathcal{L}^\bot}\text{Ad}_{z_v}
 \mathfrak{h}\notag\\
 &=-(\Theta+c_{\textbf{a}_1})\circ z_v-r_v|_{\mathfrak{g}_{1,A}}+\mathfrak{p}_v+\mathfrak{i}_{\mathcal{L}^\bot}\text{Ad}_{z_v}
 (\Theta+c_{\textbf{a}_1})\notag\\
 &\overset{(1)}{=}\mathcal{R}_v+\big(r_v-r_v|_{\mathfrak{g}_{1,A}}+(\text{Ad}_{z_v}-I)c_{\textbf{a}_1}\big).
\end{align}
Here in $(1)$ we recall
\begin{align*}
 \mathcal{R}_v&=\mathfrak{p}_v-r_v-\Theta\circ z_v+\mathfrak{i}_{\mathcal{L}^\bot}\text{Ad}_{z_{v}}\Theta
\end{align*}
(see  \eqref{for:92}) and recall $c_{\textbf{a}_1}|_{\mathfrak{g}_{1,A}}=0$ and the fact that $\mathfrak{g}_{\mu,A}$, $\mu\in \Delta_A$ are invariant subspaces for $\text{Ad}_{z_v}$, $v\in E$ (see \eqref{for:490} of Section \ref{sec:3}).

Consequently, we see that: if $1\neq \mu\in \Delta_A$
\begin{align}\label{for:483}
 \big(r_v-r_v|_{\mathfrak{g}_{1,A}}+(\text{Ad}_{z_v}-I)c_{\textbf{a}_1}\big)|_{\mathfrak{g}_{\mu,A}}
 =\big(r_v+(\text{Ad}_{z_v}-I)c_{\textbf{a}_1}\big)|_{\mathfrak{g}_{\mu,A}}
\end{align}
and
\begin{align}\label{for:484}
 \big(r_v-r_v|_{\mathfrak{g}_{1,A}}+(\text{Ad}_{z_v}-I)c_{\textbf{a}_1}\big)|_{\mathfrak{g}_{1,A}}=0.
\end{align}
It follows from \eqref{for:479} and \eqref{for:480} that
\begin{align}\label{for:481}
 \norm{\mathfrak{p}^{(1)}}_{C^0}\leq \mathcal{E}_1+\mathcal{E}_2+\mathcal{E}_3+\mathcal{E}_4,
\end{align}
where
\begin{gather*}
  \mathcal{E}_1=\max_{v\in E}\|-\mathfrak{h}\circ \beta_{v}+\mathfrak{h}\circ z_v^{(1)}\|_{C^0},\\
  \mathcal{E}_2=\max_{v\in E}\|-\mathfrak{h}\circ z_v^{(1)}+\mathfrak{h}\circ z_v \|_{C^0}
  +\max_{v\in E}\|\mathfrak{p}_v\circ \tilde{h}-\mathfrak{p}_v \|_{C^0},\\
  \mathcal{E}_3=\max_{v\in E}\|\mathcal{R}_v\|_{C^0}
  +\max_{v\in E}\|\mathfrak{i}_{\mathcal{L}^\bot}\Psi_v \|_{C^0}\\
  \mathcal{E}_4=\max_{v\in E}\|r_v-r_v|_{\mathfrak{g}_{1,A}}+(\text{Ad}_{z_v}-I)c_{\textbf{a}_1}\|.
\end{gather*}
We note that
\begin{align}\label{for:482}
\mathcal{E}_1&\leq C\norm{\mathfrak{h}}_{C^1}\max_{v\in E}d(\beta_{v},\alpha_{z_v^{(1)}})_{C^0}\leq   C_1\norm{\mathfrak{h}}_{C^1} \norm{\mathfrak{p}^{(1)}}_{C^0}\notag\\
&\overset{\text{\tiny$(1)$}}{\leq} C_2b^\varrho\tau^{\varrho}\norm{\mathfrak{p}}_{C^{\varrho}}\norm{\mathfrak{p}^{(1)}}_{C^0}\overset{\text{\tiny$(2)$}}{<}\text{\small$\frac{1}{2}$}\norm{\mathfrak{p}^{(1)}}_{C^0}.
\end{align}
Here in $(1)$ we use \eqref{for:79}; in $(2)$ by assumption $b^\varrho\tau^{\varrho}\norm{\mathfrak{p}}_{C^{\varrho}}<\bar{c}$. Thus
$C_2b^\varrho\tau^{\varrho}\norm{\mathfrak{p}}_{C^{\varrho}}<\frac{1}{2}$ providing $\bar{c}$ is sufficiently small.

It follows from \eqref{for:481} and \eqref{for:482} that
\begin{align}
 \norm{\mathfrak{p}^{(1)}}_{C^0}\leq 2\mathcal{E}_2+2\mathcal{E}_3+2\mathcal{E}_4.
\end{align}
Firstly, we estimate $\mathcal{E}_2$:
\begin{align}\label{for:486}
  \mathcal{E}_2&\leq \norm{\mathfrak{h}}_{C^1}\max_{v\in E}\emph{\textbf{d}}(z_v^{(1)},z_v)+\norm{\mathfrak{p}}_{C^1}d(\tilde{h},I)_{C^0}\notag\\
  &\overset{(1)}{\leq} Cb^\varrho\tau^{\varrho}\norm{\mathfrak{p}}_{C^{\varrho}}
 \norm{\mathfrak{p}}_{C^{0}}+Cb^\varrho\tau^{\varrho}\norm{\mathfrak{p}}_{C^{\varrho}}
 \norm{\mathfrak{p}}_{C^{1}}
\end{align}
Here in $(1)$ we use \eqref{for:79} and \eqref{for:476}.

Secondly, we estimate $\mathcal{E}_3$:
\begin{align}\label{for:487}
  \mathcal{E}_3&\overset{(1)}{\leq} C\tau^{\varrho}(\norm{\mathfrak{p}}_{C^{\varrho+1}})^{1+\frac{c_0}{8}}
+C\tau^{\varrho}(C_{m}\tau^{-m+\varrho}\norm{\mathfrak{p}}_{C^m})^{\frac{1}{8}}(\norm{\mathfrak{p}}_{C^\varrho})^{\frac{7}{8}}\notag\\
&+C_{\ell} b^{-\ell+\varrho}(\tau^{\ell}\norm{\mathfrak{p}}_{C^\varrho}+\norm{\mathfrak{p}}_{C^{\ell}})+Cb^{2\varrho}\tau^{2\varrho}(\norm{\mathfrak{p}}_{C^\varrho})^2.
\end{align}
for any $m,\,\ell\geq\varrho$.

Here in $(1)$ we use \eqref{for:96} and \eqref{for:111}.

Finally, we estimate $\mathcal{E}_4$: we show that
\begin{align}\label{for:42}
 \mathcal{E}_4\leq C(\norm{\mathfrak{p}}_{C^1})^2.
\end{align}
From \eqref{for:483} and  \eqref{for:484}, it suffices to prove: for any $1\neq\mu\in \Delta_A$
\begin{align}\label{for:100}
 \big\|r_v|_{\mathfrak{g}_{\mu,A}}-(I-\text{Ad}_{z_v})c_{\mu}\big\|\leq C(\norm{\mathfrak{p}}_{C^1})^2,\quad \forall\,v\in E.
\end{align}
 We note that
\begin{align}\label{for:271}
  &(I-\text{Ad}_{z_{\textbf{a}_1}})\big(r_v|_{\mathfrak{g}_{\mu,A}}-(I-\text{Ad}_{z_v})c_{\mu}\big)\notag\\
  &=(I-\text{Ad}_{z_{\textbf{a}_1}})(r_v|_{\mathfrak{g}_{\mu,A}})-(I-\text{Ad}_{z_{v}})(I-\text{Ad}_{z_{\textbf{a}_1}})c_{\mu}\notag\\
  &+(\text{Ad}_{z_{v}z_{\textbf{a}_1}}-\text{Ad}_{z_{\textbf{a}_1}z_{v}})c_{\mu}\notag\\
  &\overset{\text{\tiny$(1)$}}{=}(I-\text{Ad}_{z_{\textbf{a}_1}})(r_v|_{\mathfrak{g}_{\mu,A}})-(I-\text{Ad}_{z_{v}})(r_{\textbf{a}_1}|_{\mathfrak{g}_{\mu,A}})\notag\\
  &+(\text{Ad}_{z_{v}z_{\textbf{a}_1}}-\text{Ad}_{z_{\textbf{a}_1}z_{v}})c_{\mu}\notag\\
  &=(r_v-\text{Ad}_{z_{\textbf{a}_1}}r_v-r_{\textbf{a}_1}+\text{Ad}_{z_v}r_{\textbf{a}_1})\big|_{\mathfrak{g}_{\mu,A}}\notag\\
  &+(\text{Ad}_{z_{v}z_{\textbf{a}_1}}-\text{Ad}_{z_{\textbf{a}_1}z_{v}})c_{\mu}.
\end{align}
Here in $(1)$ we recall definition of $c_{\mu}$; in $(2)$ we recall the fact that $\mathfrak{g}_{\mu,A}$, $\mu\in \Delta_A$ are invariant subspaces for $\text{Ad}_{z_v}$, $v\in E$ (see Section \ref{sec:1}).

It follows that: for any $1\neq\mu\in \Delta_A$
\begin{align*}
  &\big\|r_v|_{\mathfrak{g}_{\mu,A}}-(I-\text{Ad}_{z_v})c_{\mu}\big\|\notag\\
  &\overset{\text{\tiny$(1)$}}{\leq}C\Big\|(I-\text{Ad}_{z_{\textbf{a}_1}})\big(r_v|_{\mathfrak{g}_{\mu,A}}-(I-\text{Ad}_{z_v})c_{\mu}\big)\Big\|\notag\\
  &\overset{\text{\tiny$(2)$}}{\leq} \big\|(r_v-\text{Ad}_{z_{\textbf{a}_1}}r_v-r_{\textbf{a}_1}+\text{Ad}_{z_v}r_{\textbf{a}_1})\big|_{\mathfrak{g}_{\mu,A}} \big\|\notag\\
  &+\big\|(\text{Ad}_{z_{\textbf{a}_1}z_v}-\text{Ad}_{z_vz_{\textbf{a}_1}})c_{\mu}\big\|\notag\\
  &\overset{\text{\tiny$(3)$}}{\leq} C(\norm{\mathfrak{p}}_{C^1})^2+C\norm{\log[z_{\textbf{a}_1},z_v]}\norm{c_{\textbf{a}_1}}\notag\\
 &\overset{\text{\tiny$(4)$}}{\leq}  C(\norm{\mathfrak{p}}_{C^1})^2+C_1(\norm{\mathfrak{p}}_{C^0})^{2}.
\end{align*}
Here in $(1)$ we use \eqref{for:43}; in $(2)$ we use \eqref{for:271};
in $(3)$ we use \eqref{for:39} of Lemma \ref{le:8} and note that
\begin{align*}
  \norm{\text{Ad}_{z_{\textbf{a}_1}z_v}-\text{Ad}_{z_vz_{\textbf{a}_1}}}&\leq C\emph{\textbf{d}}(z_{\textbf{a}_1}z_v,\,z_vz_{\textbf{a}_1})
  \overset{\text{\tiny$(a)$}}{\leq} C_1\norm{\log[z_{\textbf{a}_1},z_v]}.
\end{align*}
Here in $(a)$ we use \eqref{for:485}; in $(4)$ we use \eqref{for:253} and \eqref{for:66}.

Hence we get \eqref{for:100}.

\medskip
Finally, it follows from \eqref{for:486}, \eqref{for:487} and \eqref{for:42} that
\begin{align}\label{for:113}
\norm{\mathfrak{p}^{(1)}}_{C^0}&\leq Cb^\varrho\tau^{\varrho}\norm{\mathfrak{p}}_{C^{\varrho}}
 \norm{\mathfrak{p}}_{C^{0}}+Cb^\varrho\tau^{\varrho}\norm{\mathfrak{p}}_{C^{\varrho}}
 \norm{\mathfrak{p}}_{C^{1}}\notag\\
 &+C\tau^{\varrho}(\norm{\mathfrak{p}}_{C^{\varrho+1}})^{1+\frac{c_0}{8}}
+C\tau^{\varrho}(C_{m}\tau^{-m+\varrho}\norm{\mathfrak{p}}_{C^m})^{\frac{1}{8}}(\norm{\mathfrak{p}}_{C^\varrho})^{\frac{7}{8}}\notag\\
&+C_{\ell} b^{-\ell+\varrho}(\tau^{\ell}\norm{\mathfrak{p}}_{C^\varrho}+\norm{\mathfrak{p}}_{C^{\ell}})+Cb^{2\varrho}\tau^{2\varrho}(\norm{\mathfrak{p}}_{C^\varrho})^2\notag\\
&+C(\norm{\mathfrak{p}}_{C^1})^2
\end{align}
for any $m,\,\ell\geq\varrho$.

Noting that $\norm{\mathfrak{p}}_{C^{\varrho+1}}<1$, $b,\,\tau>1$ and $0<c_0<1$ (\eqref{for:470} of \eqref{for:432} of Section \ref{sec:3}),  we can simplify \eqref{for:113} as
\begin{align*}
\norm{\mathfrak{p}^{(1)}}_{C^0}
&\leq 5Cb^{2\varrho}\tau^{2\varrho}(\norm{\mathfrak{p}}_{C^{\varrho+1}})^{1+\frac{c_0}{8}}+C\tau^{\varrho}(C_{m}\tau^{-m+\varrho}\norm{\mathfrak{p}}_{C^m})^{\frac{1}{8}}
(\norm{\mathfrak{p}}_{C^\varrho})^{\frac{7}{8}}\notag\\
&+C_{\ell} b^{-\ell+\varrho}(\tau^{\ell}\norm{\mathfrak{p}}_{C^\varrho}+\norm{\mathfrak{p}}_{C^{\ell}}).
\end{align*}
Hence we finish the proof.

\end{proof}
\subsection{The iterative process }\label{sec:6} In this part we obtain a smooth
conjugacy between the initial perturbation $\tilde{\alpha}_A$ and $\alpha_A$ by using the KAM iteration. We follow the line
of arguments in Section \ref{sec:21}.

We recall the following positive constants defined in previous sections:
\begin{enumerate}
  \item $\varrho$ is defined in Section \ref{sec:56};

  \smallskip
  \item $\bar{c}$ is from Proposition \ref{po:5};

  \smallskip
  \item $y_0=\frac{c_0}{8}$, where $c_0$ is from Proposition \ref{po:1}.
\end{enumerate}

Let $y_0=\frac{c_0}{8}$ and let $\gamma$, $\gamma_1$, $x_0$ and $\nu$ be as defined in Section \ref{sec:21}. Let $\ell>\varrho$  be sufficiently large such that \eqref{for:220} to \eqref{for:254} are satisfied.
We fix an increasing sequence $\beta_n\to \infty$ with $\beta_1>2\ell$. We construct $\mathfrak{p}_v^{(n)}$, $E'(z^{(n)})$, $v\in E$ and $h_n$ inductively as follows. We denote
\begin{align*}
 \tilde{\alpha}_v(\mathfrak{j}(x))=\mathfrak{j}\big(\exp(\mathfrak{p}_v(x))\cdot
z_v\cdot x\big),\qquad \forall\,v\in E,\,x\in \mathcal{X}.
\end{align*}
Set
\begin{gather*}
 \tilde{\alpha}^{(0)}=\tilde{\alpha},\quad \mathfrak{p}_v^{(0)}=\mathfrak{p}_v,\quad h_{0}=I, \\
 E'(z^{(0)})=E'(z),\quad \epsilon_n=\epsilon^{\text{\tiny$\gamma^n$}}
\end{gather*}
where $0<\epsilon^{\text{\tiny$\frac{1}{2}$}}<\bar{c}$ is sufficiently small so that the following hold:
\begin{align*}
\norm{\mathfrak{p}^{(0)}}_{C^0}\leq \epsilon_0=\epsilon,\quad \norm{\mathfrak{p}^{(0)}}_{C^{\ell+\varrho}}\leq \epsilon_0^{-\gamma_1}, \quad d(h_0, I)_{C^1}<\epsilon_0^{\text{\tiny$\frac{1}{2}$}}
\end{align*}
and $E'(z)$ is a $\frac{\delta_0}{2}$-standard perturbation of $E$, where $\norm{\mathfrak{p}^{(0)}}_{C^r}=\max_{v\in E}\norm{\mathfrak{p}^{(0)}_z}_{C^r}$ for any $r\geq0$.

Suppose  inductively that
\begin{align*}
 \tilde{\alpha}^{(n)}_v(\mathfrak{j}(x))=\mathfrak{j}\big(\exp(\mathfrak{p}^{(n)}_v(x))\cdot
z^{(n)}_v\cdot x\big),\qquad \forall\,v\in E,\,x\in \mathcal{X}
\end{align*}
where $\tilde{\alpha}^{(n)}=h_{n}^{-1}\circ \tilde{\alpha}^{(n-1)}\circ h_{n}$, $E'(z^{(n)})$ is a $(\frac{\delta_0}{2}+\sum_{i=0}^{n-1}\epsilon_i^{\frac{1}{2}})$-standard perturbation of $E$ and
\begin{gather*}
 \norm{\mathfrak{p}^{(n)}}_{C^0}\leq \epsilon_n, \quad \norm{\mathfrak{p}^{(n)}}_{C^{\ell+\varrho}}\leq \epsilon^{-\gamma_1}_n,\\
 \max_{z\in E} d(z_v^{(n)},\,z_v^{(n-1)})<\epsilon_n^{\frac{1}{2}},\quad d(h_{n},I)_{C^1}<\epsilon_n^{\frac{1}{2}},\\
\norm{\mathfrak{p}^{(n)}}_{C^{\beta_l}}<r_l^n\epsilon^{\text{\tiny$-8\gamma$}}_n(\norm{\mathfrak{p}^{(l-1)}}_{C^{\beta_l}}+1)
 \end{gather*}
for any $l\leq n$; where $r_l$ is a constant dependent only on $l$ and $\norm{\mathfrak{p}^{(n)}}_{C^m}=\max_{v\in E}\norm{\mathfrak{p}^{(n)}_v}_{C^m}$ for any $m\geq0$.

By interpolation inequalities \eqref{for:183} and \eqref{for:76} hold as well. We note that
\begin{align}\label{for:125}
 \text{\small$\frac{\delta_0}{2}$}+\sum_{i=0}^{n-1}\epsilon_i^{\frac{1}{2}}+C\norm{\mathfrak{p}^{(n)}}_{C^{0}}\leq \text{\small$\frac{\delta_0}{2}$}+\sum_{i=0}^{n-1}\epsilon_i^{\frac{1}{2}}+C\epsilon_{n}
 < \text{\small$\frac{\delta_0}{2}$}+\sum_{i=0}^{n}\epsilon_i^{\frac{1}{2}}<\delta_0.
\end{align}
\eqref{for:76} and \eqref{for:125}
allow us to apply Proposition \ref{po:5} to obtain  the new iterates $\mathfrak{p}^{(n+1)}$, $h_{n+1}$ and $E'(z^{(n+1)})$.

\subsection{Convergence}\label{sec:30} Let $b=\epsilon^{-\text{\small$\frac{2\nu}{\ell}$}}_n$. Set $\zeta_l$,  $1\leq l\leq n$ and $\zeta$ as in \eqref{for:214}.

\textbf{Firstly, we consider the case of $\zeta>b^{\frac{1}{2}}$.} Set $\tau=b^{\frac{1}{2}}=\epsilon^{-\text{\small$\frac{\nu}{\ell}$}}_n$ and $m=\ell+\varrho$.
By Proposition \ref{po:5} there are $h_{n+1}\in C^\infty(\mathcal{X})$, $\mathfrak{p}_v^{(n+1)}\in C^\infty(\mathcal{M},\,\mathcal{L}^\bot)$, $v\in E$
and $E'(z^{(n+1)})$ is a standard perturbation of $E$, such that $\tilde{\alpha}^{(n+1)}=h_{n+1}^{-1}\circ \tilde{\alpha}^{(n)}\circ h_{n+1}$ can be written as
\begin{align*}
 \tilde{\alpha}^{(n+1)}_v(\mathfrak{j}(x))=\mathfrak{j}\big(\exp(\mathfrak{p}^{(n+1)}_v(x))\cdot
z^{(n+1)}_v\cdot x\big),\qquad \forall\,v\in E,\,x\in \mathcal{X};
\end{align*}
and the following estimates hold:

\smallskip

$(\textrm{I})$ For $E'(z^{(n+1)})$, by \eqref{for:112} and the assumption, $E'(z^{(n+1)})$ is a standard $(\frac{\delta_0}{2}+\sum_{i=0}^{n-1}\epsilon_i^{\frac{1}{2}}+C\norm{\mathfrak{p}^{(n)}}_{C^{0}})$-perturbation $E'(z^{(1)})$ of $E$ and
\begin{align}\label{for:248}
  \max_{z\in E} d(z_v^{(n+1)},\,z_v^{(n)})\leq C\norm{\mathfrak{p}^{(n+1)}}_{C^0}<\epsilon_n^{\frac{1}{2}}.
\end{align}
$(\textrm{II})$ For $d(h_{n+1},I)_{C^{1}}$, $\norm{\mathfrak{p}^{(n+1)}}_{C^{0}}$, $\norm{\mathfrak{p}^{(n+1)}}_{C^{\ell+\varrho}}$, $\norm{\mathfrak{p}^{(n+1)}}_{C^{\beta_l}}$ and $d(h_{n+1},I)_{C^{\beta_l}}$, $1\leq l\leq n+1$:  Remark \ref{re:5} shows that following exactly the same arguments as in Section \ref{sec:21},  \eqref{for:219}, \eqref{for:256}, \eqref{for:115}, \eqref{for:225}, \eqref{for:192} and
\eqref{for:218} hold as well.

\smallskip

\textbf{Next we consider the case of $\zeta\leq b^{\frac{1}{2}}$.} Choose $1\leq p\leq n$ such that $\zeta_p=\zeta$.  Set $\tau=\zeta$ and $m=\beta_p$.
Since $\tau=\zeta\leq b^{\frac{1}{2}}= \epsilon^{-\text{\small$\frac{\nu}{\ell}$}}_n$, the estimates for $d(h_{n+1},I)_{C^1}$, $E'(z^{(n+1)})$, $\norm{\mathfrak{p}^{(n+1)}}_{C^{\ell+\varrho}}$,
 $d(h_{n+1},I)_{C^{\beta_{n+1}}}$ and $\norm{\mathfrak{p}^{(n+1)}}_{C^{\beta_{n+1}}}$
 still hold. For $\norm{\mathfrak{p}^{(n+1)}}_{C^{\beta_l}}$, $1\leq l\leq n$, since $\tau=\zeta\leq \zeta_l$ for each $1\leq l\leq n$, \eqref{for:226} still holds.
 As a result, the estimates for $\norm{\mathfrak{p}^{(n+1)}}_{C^{\beta_l}}$, $1\leq l\leq n$ hold as well.
 So we just need to get the estimate for $\norm{\mathfrak{p}^{(n+1)}}_{C^{0}}$. Following exactly the same arguments as in Section \ref{sec:21}, we see that \eqref{for:23} holds as well.


\medskip
\textbf{Conclusion}: Hence we can obtain infinite sequences $h_n$, $E'(z^{(n)})$, $\mathfrak{p}_v^{(n)}$, $v\in E$ inductively. Set
\begin{align*}
 H_n=h_0\circ \cdots \circ h_n,\quad\text{and}\quad a_v=\lim_nz_v^{(n)},\quad v\in E.
\end{align*}
As $E'(z^{(n)})$ is a standard perturbation of $E$ for each $n$ \eqref{for:248} shows that $z_v^{(n)}$ converges to an element inside $Z(L)$ for each $v\in E$.
\eqref{for:219} shows that $H_n$ converges in $C^1$ topology to a $C^1$ conjugacy $h$ between $\tilde{\alpha}_A$ and $\alpha_A$; moreover, \eqref{for:218} shows that
 the convergence of the sequence $H_n$ in $C^{\frac{\beta_l}{8}}$ for any $l\in\NN$. Hence we see that
 $h$ is of class $C^{\infty}$. The convergence step shows that:
\begin{align*}
  h^{-1}\circ\tilde{\alpha}_v(hx)=\alpha_{a_v}(x), \qquad \forall x\in \mathcal{X},\,\forall\,v\in E.
\end{align*}
Since $\tilde{\alpha}$ is abelian, $\{a_v:\,v\in E\}$ is an abelian set. As $E$ is a generating set of the group $A$, the map $v\to a_v$, $v\in E$ extends to a close to identity group homeomorphism $i_0$ from $A$ to an abelian group $A'$; and thus we have
\begin{align*}
  h^{-1}\circ\tilde{\alpha}_a(hx)=\alpha_{i_0(a)}(x),\quad \text{for all }a\in A,\,x\in \mathcal{X}.
\end{align*}
Then we complete the proof of Theorem \ref{th:7}.

\appendix

\section{Proof of Theorem \ref{th:3}}\label{sec:8}

It suffices to prove the following: there exist constants $\sigma,\,C>0$, dependent only on $G$, $\Gamma$ and $\rho$ such that for any
$f_i\in \mathcal{O}^1$, $i=1,2$ and  any $g\in G$, we have
\begin{align}\label{for:59}
\abs{\langle f_1(g),f_2\rangle}&\leq
C\norm{f_1}_{1}\norm{f_2}_{1}e^{-\sigma d(e,g)}.
\end{align}
For the case of $\GG=G$, the result follows from Theorem \ref{th:10} and \ref{th:2} immediately.

For the twisted case, firstly we show that: $(*)$ Any noncompact one parameter subgroups of $G$ acts ergodictly on $\GG/\Upsilon$.
 Brezin and Moore  \cite{Brezin}
show that a one parameter subgroup of $G$ acts ergodically on
$\GG/\Upsilon$ if and only if the quotient flows
on the maximal Euclidean quotient and the maximal semisimple
quotient are ergodic. Suppose $\GG/E$ is the Euclidean quotient. Since $\rho(G)$-fixed points in $V/[V,V]$ are $\{0\}$, we have $[\GG,\GG]=\GG$ using \eqref{for:12}.
Then by Brezin and Moore $\GG\subset E$. Hence the maximal Euclidean quotient is a point. The latter quotient is
isomorphic to $G/\Gamma$. Since $\Gamma$ is irreducible, all
noncompact one parameter subgroups of $G$ are ergodic on $G/\Gamma$ (see Theorem \ref{th:10}).
Then we get $(*)$. $(*)$ implies that for each simple factor $G_i$ of $G$, $\pi|_{G_i}$ has no $G_i$-fixed vectors, where $\pi$ denotes the regular representation on $\mathcal{O}$. Then we have

$(**)$: if $G_i$ is not isomorphic to
$SO(n,1)$ or $SU(n,1)$, $\pi|_{G_i}$ has a spectral gap \cite{Zimmer}.

Now we suppose $G_i$ is isomorphic to
$SO(n,1)$ or $SU(n,1)$. We recall the definition of  relative Property $(T)$ for triples \cite[Rem. 0.2.2, p.3]{jaudon}.
Let $H$ and $M$ be closed subgroups of
a locally compact group $S$. We say that the triple $(S, H, M)$ has relative Property $(T)$, if for any unitary representation $\pi$ of $S$, such that the restriction $\pi|_H$
has almost-invariant vectors, then there exist nonzero $\pi(M)$-invariant vectors.
By Corollary 1.9 of \cite{Indira} the triple $(G_i\ltimes V, G_i, V)$ has relative Property $(T)$ if and only if
the triple $(G_i\ltimes (V/[V,V]), G_i, V/[V,V])$ has relative Property $(T)$. Theorem 5.5 of \cite{Shalom1} shows   that the triple $(G_i\ltimes (V/[V,V]), G_i, V/[V,V])$ has relative Property $(T)$ if $\rho$ is excellent. Hence we see that the triple $(G_i\ltimes V, G_i, V)$ has relative Property $(T)$, which shows that:

$(***)$ if $G_i$ is isomorphic to
$SO(n,1)$ or $SU(n,1)$, then for any unitary representation of $\GG$ without $V$-fixed vectors, the restriction to $G_i$ is isolated from the trivial representation.

Next, we consider the regular representation. For any $f\in \mathcal{O}$, let
\begin{align*}
  \tilde{f}(g,v)=f(g,v)-\int_{V/\Lambda}f(g,v)dv\quad\text{and}\quad \bar{f}(g,v)=\int_{V/\Lambda}f(g,v)dv.
\end{align*}
where $dv$ is the Haar measure on $V/\Lambda$. Let $\tilde{\mathcal{O}}=\{\tilde{f}: f\in \mathcal{O}\}$ and $\bar{\mathcal{O}}=\{\bar{f}: f\in \mathcal{O}\}$.
We note that $\mathcal{O}=\tilde{\mathcal{O}}\oplus\bar{\mathcal{O}}$. Then there exist constants $\alpha_i,\,C_i>0$, $i=1,2$ dependent only on $G$, $\Gamma$ and $\rho$ such that
for any $f_1,\,f_2\in \mathcal{O}^1$ and any $g\in G_i$ we have
\begin{align}\label{for:127}
\abs{\langle\pi(g)f_1,&f_2\rangle}\overset{\text{(1)}}{\leq} \abs{\langle\pi(g)\tilde{f_1},\tilde{f_2}\rangle}+\abs{\langle\pi(g)\bar{f_1},\bar{f_1}\rangle}\notag\\
&\overset{\text{(2)}}{\leq} C_1\norm{\tilde{f_1}}_1\norm{\tilde{f_2}}_1e^{-\alpha_1 d(e,g)}+C_2\norm{\bar{f_1}}_1\norm{\bar{f_2}}_1e^{-\alpha_2 d(e,g)}\notag\\
&\leq 2C_1\norm{f_1}_1\norm{f_2}_1e^{-\alpha_1 d(e,g)}+C_2\norm{f_1}_1\norm{f_2}_1e^{-\alpha_2 d(e,g)}.
\end{align}
Here in $(1)$ we use the fact the subspaces $\tilde{\mathcal{O}}$ and $\bar{\mathcal{O}}$ are orthogonal to each other; in $(2)$ we note that
$\tilde{\mathcal{O}}$ contains no $V$-fixed vectors and $\bar{\mathcal{O}}$ is a subspace of $L^2_0(G/\Gamma)$. Then by $(**)$, $(***)$ we see that
$\pi|_{G_i}$ has a spectral gap restricted on $\tilde{\mathcal{O}}$ for each $i$. It follows from
Theorem \ref{th:2} we have
\begin{align*}
 \abs{\langle\pi(g)\tilde{f_1},\tilde{f_2}\rangle}\leq C_1\norm{\tilde{f_1}}_1\norm{\tilde{f_2}}_1e^{-\alpha_1 d(e,g)}.
\end{align*}
By arguments for the case of $\GG=G$ we have
\begin{align*}
 \abs{\langle\pi(g)\bar{f_1},\bar{f_2}\rangle}\leq C_2\norm{\bar{f_1}}_1\norm{\bar{f_2}}_1e^{-\alpha_2 d(e,g)}.
\end{align*}
\eqref{for:127} shows that $\pi|_{G_i}$ has a spectral gap on $\mathcal{O}$. Then \eqref{for:59} follows from Theorem \ref{th:2}.

\section{Proof of \eqref{for:93}}\label{sec:35}

We have
\begin{align}
&\bigl\langle \pi_{\mathfrak{u}}(f)\xi,\eta\bigl\rangle=\text{\tiny$\frac{1}{\sqrt{2\pi}}$}\int_{\RR^m}\hat{f}(t)\langle \pi(t)\xi,\eta\rangle dt\label{for:14}\\
&=\text{\tiny$\frac{1}{\sqrt{2\pi}}$}\int_{\RR}\int_{\widehat{\RR}}\hat{f}(t)\chi(t)\langle \xi_\chi,\,\eta_\chi \rangle d\mu(\chi) dt\notag\\
&=\int_{\widehat{\RR}}\langle \xi_\chi,\,\eta_\chi \rangle\big(\text{\tiny$\frac{1}{\sqrt{2\pi}}$}\int_{\RR}\hat{f}(t)\chi(t) dt\big)d\mu(\chi)\notag\\
&=\int_{\widehat{\RR}}f(\chi)\langle \xi_\chi,\,\eta_\chi \rangle d\mu(\chi).\notag
\end{align}
Since
\begin{align*}
 \norm{\pi_{\mathfrak{u}}(f)}\leq \norm{f}_\infty,\qquad \forall\,f\in \mathcal{S}(\RR)
\end{align*}
we can extend $\pi_{\mathfrak{u}}$ from $\mathcal{S}(\RR)$ to $L^\infty(\RR)$ by taking strong limits of operators and pointwise monotone increasing limits of non-negative functions (see \cite{Lang} for a detailed treatment). Hence $\pi_{\mathfrak{u}}$ is a homomorphism of $L^\infty(\RR^m)$ to bounded operators on $\mathcal{H}$. Moreover, for any $f\in L^\infty(\RR^m)$ we see that
\begin{align}\label{for:247}
 \pi_{\mathfrak{u}}(f)(\xi)=\int_{\widehat{\RR}}f(\chi)\xi_\chi d\mu(\chi),\qquad \forall\,\xi\in \mathcal{H}.
\end{align}
Then the result follows from \eqref{for:247}.

\section{Proof of \eqref{for:234} and \eqref{for:98}}\label{sec:31}
Since $\mathfrak{u}$ is nilpotent, $(\text{ad}_\mathfrak{u})^{\dim\mathfrak{S}}=0$, we have
\begin{align}\label{for:235}
 \pi(t)v=\text{Ad}_{\exp(t\mathfrak{u})}(v)\pi(t)=\sum_{j=0}^{\dim\mathfrak{S}-1}\text{\small$\frac{t^j}{j!}$}\text{ad}^{j}_\mathfrak{u}(v)\pi(t),\quad \forall\,v\in \mathfrak{S}.
\end{align}
We suppose $f\in \mathcal{S}(\RR)$, $\vartheta\in \mathcal{H}^1$ and $\eta\in \mathcal{H}$. Then we have
\begin{align*}
&\bigl\langle \pi_{\mathfrak{u}}(f)(v\vartheta),\eta\bigl\rangle=\text{\tiny$\frac{1}{\sqrt{2\pi}}$}\int_{\RR}\bigl\langle \hat{f}(t)\pi(t)(v\vartheta),\eta\rangle dt\notag\\
&\overset{\text{(1)}}{=}\text{\tiny$\frac{1}{\sqrt{2\pi}}$}\sum_{j=0}^{\dim\mathfrak{S}-1}\text{\small$\frac{1}{j!}$}\int_{\RR}\bigl\langle \hat{f}(t)t^j\text{ad}^{j}_\mathfrak{u}(v)\pi(t)\vartheta,\,\eta\bigl\rangle dt\\
&=\text{\tiny$\frac{1}{\sqrt{2\pi}}$}\sum_{j=0}^{\dim\mathfrak{S}-1}\text{\small$\frac{\textrm{i}^{-j}}{j!}$}\int_{\RR}\bigl\langle \widehat{f^{(j)}}(t)\text{ad}^{j}_\mathfrak{u}(v)\pi(t)\vartheta,\,\eta\bigl\rangle dt\\
&=\sum_{j=0}^{\dim\mathfrak{S}-1}\text{\small$\frac{\textrm{i}^{-j}}{j!}$}
\bigl\langle \text{ad}^{j}_\mathfrak{u}(v)\pi_{\mathfrak{u}}(f^{(j)})\vartheta,\eta\bigl\rangle
\end{align*}
Here in $(1)$ we use \eqref{for:235}.

This shows that for any $v\in \mathfrak{S}$ and any $f\in \mathcal{S}(\RR)$ we have
\begin{align*}
  \pi_{\mathfrak{u}}(f)v=\sum_{j=0}^{\dim\mathfrak{S}-1}\text{\small$\frac{\textrm{i}^{-j}}{j!}$}\text{ad}^{j}_\mathfrak{u}(v)\pi_{\mathfrak{u}}(f^{(j)}).
\end{align*}
By arguments at the end of Appendix \ref{sec:35}, we can extend the above equation from $\mathcal{S}(\RR)$ to $\tilde{\mathcal{S}}(\RR)$.
Then we get \eqref{for:234}.

We also note that
\begin{align*}
v\pi(t)=\pi(t)\text{Ad}_{\exp(-t\mathfrak{u})}(v)=\sum_{j=0}^{\dim\mathfrak{S}-1}\text{\small$\frac{(-1)^jt^j}{j!}$}\pi(t)\text{ad}^{j}_\mathfrak{u}(v),\quad \forall\,v\in \mathfrak{S}
\end{align*}
The arguments follow in a similar way.

\section{Proof of Lemma \eqref{cor:1}}\label{sec:33}

\eqref{for:128}: It follows from \eqref{for:41} that
\begin{align*}
 \big\|\mathfrak{u}^{l}\big(\pi_{\mathfrak{u}}(f\circ \tau^{-1})\xi\big)\big\|=
 \tau^{l}\big\|\pi_{\mathfrak{u}}(f_{l}\circ \tau^{-1})\xi\big\|\overset{\text{(1)}}{\leq} \tau^l\norm{f_{l}}_\infty\norm{\xi}.
\end{align*}
Here in $(1)$ we use \eqref{for:215}. This implies the result.

\smallskip

\eqref{for:292}: Keeping using \eqref{for:234}, for any vectors $v_i\in \mathfrak{S}$, $1\leq i\leq l$ we have
\begin{align}
  &v_l\cdots v_2v_1\pi_{\mathfrak{u}}(f)\notag\\
  &=\sum_{\substack{ 0\leq j_i\leq\dim\mathfrak{S}-1,\\ 1\leq i\leq l}}c_{j_1,\cdots,j_l}\pi_{\mathfrak{u}}(f^{(\sum_{i=1}^lj_i)})
  \big((\text{ad}_\mathfrak{u})^{j_l}v_l)\cdots (\text{ad}_\mathfrak{u})^{j_1}v_1)\big).
\end{align}
It follows that
\begin{align*}
&\norm{v_l\cdots v_2v_1(\pi_{\mathfrak{u}}(f)\xi)}\\
  &\leq\sum_{\substack{ 0\leq j_i\leq\dim\mathfrak{S}-1,\\ 1\leq i\leq l}}|c_{j_1,\cdots,j_l}|
  \Big\|\pi_{\mathfrak{u}}(f^{(\sum_{i=1}^lj_i)})
  \big((\text{ad}_\mathfrak{u})^{j_l}v_l)\cdots (\text{ad}_\mathfrak{u})^{j_1}v_1)\xi\big)\Big\|\\
  &\overset{\text{(1)}}{\leq}\sum_{\substack{ 0\leq j_i\leq\dim\mathfrak{S}-1,\\ 1\leq i\leq l}}|c_{j_1,\cdots,j_l}|\cdot\norm{f^{(\sum_{i=1}^lj_i)}}_{\text{\small$L^\infty(\RR)$}}\\
  &\cdot\big\|(\text{ad}_\mathfrak{u})^{j_l}v_l)\cdots (\text{ad}_\mathfrak{u})^{j_1}v_1)\xi\big\|\\
  &\leq C_{l}\norm{f}_{\tilde{\mathcal{S}}(\RR^m),l \dim\mathfrak{S}}\norm{\xi}_l.
\end{align*}
Here in $(1)$ we use \eqref{for:215}. This  implies that
\begin{align}\label{for:241}
 \norm{\pi_{\mathfrak{u}}(f)\xi}_l\leq C_{l}\norm{f}_{\tilde{\mathcal{S}}(\RR),l\dim\mathfrak{S}}\norm{\xi}_l, \qquad \forall\, 0\leq l\leq s.
\end{align}
\eqref{for:489}: In \eqref{for:241} replacing $f$ by $f\circ \tau^{-1}$ we get \eqref{re:7}.
In \eqref{re:7}, noting that $\tau\geq1$, we get \eqref{for:7}.

\section{Proof of Lemma \ref{cor:4}}\label{sec:34}

For set $X\subset \RR$ we use $I_X$ to denote the characteristic function of $X$. From \eqref{for:84} for any $r>0$ we see that
\begin{align*}
 \norm{\mathfrak{u}^{s}\xi}^2&=\int_{\widehat{\RR}}|\chi|^{2s}\norm{\xi_\chi}^2 d\mu(\chi)
 \geq \int_{\widehat{\RR}}|\chi|^{2s}I_{\norm{\chi}\geq r}^2\norm{\xi_\chi}^2 d\mu(\chi)\\
 &\geq r^{2s}\int_{\widehat{\RR}}I_{\norm{\chi}\geq r}^2\norm{\xi_\chi}^2 d\mu(\chi)\overset{\text{(1)}}{=}r^{2s}\bigl\langle \pi_{\mathfrak{u}}(I_{\norm{t}\geq r}^2)\xi,\xi\bigl\rangle\\
 &\overset{\text{(2)}}{=}r^{2s}\bigl\langle \pi_{\mathfrak{u}}(I_{\norm{t}\geq r})\xi,\,\pi_{\mathfrak{u}}(I_{\norm{t}\geq r})\xi\bigl\rangle=r^{2s}\big\|\pi_{\mathfrak{u}}(I_{\norm{t}\geq r})\xi\big\|^2.
\end{align*}
Here in $(1)$ we use \eqref{for:93}; in $(2)$ we use \eqref{for:2}.

This shows that for any $r>0$,
\begin{align}\label{for:246}
\norm{\xi}_{\{\mathfrak{u}\},s}\geq|r^s|\norm{\pi_{\mathfrak{u}}(I_{\norm{x}\geq r})\xi}.
\end{align}
Hence  we have
\begin{align*}
  \norm{\xi-\pi_{\mathfrak{u}}(f\circ \tau^{-1})\xi}&\overset{\text{(1)}}{=}\norm{\pi_{\mathfrak{u}}(1-f\circ \tau^{-1})\xi}\overset{\text{(2)}}{\leq } (\norm{f}_{C^0}+1)\norm{\pi_{\mathfrak{u}}(I_{\norm{x}\geq \tau})\xi}\\
  &\overset{\text{(3)}}{\leq }  C_{f} \tau^{-s}\norm{\xi}_{\{\mathfrak{u}\},s}.
\end{align*}
Here in $(1)$ we use \eqref{for:93} since $\xi=\pi_{\mathfrak{u}}(1)\xi$; in $(2)$ we use the fact that $1-f\circ \tau^{-1}=0$ if $\norm{x}\leq \tau$; in $(3)$ we use \eqref{for:246}.

\section{Proof of Proposition \ref{cor:6}}\label{sec:18}

\subsection{Notations} We define:

\begin{enumerate}
  \item For any $n\in\NN$ set
\begin{align*}
 \mathfrak{M}^{n_0}\stackrel{\text{def}}{=}\{\vec{h}=(h_1,\cdots,h_{n_0}): h_i\text{ is supported on $[-2,2]$ for each }i\}.
\end{align*}
  \item For $f\in C^\infty(\RR)$, $\vec{f}$ denotes $(f,\cdots,f)$ in $\mathfrak{M}^{n_0}$. We write $\vec{h}\hookrightarrow \vec{g}$, if $h_i(t)=g_i^{(l_i)}(t)(t\textrm{i})^{k_i}$, where $l_i,\,k_i\geq0$ for each $i$ and $\textrm{i}^2=-1$.
  \item For any $\vec{h}\in \mathfrak{M}^{n_0}$ we define
\begin{align}
 \mathcal{E}_{\vec{h},\tau}\stackrel{\text{def}}{=}\pi_{\mathfrak{u}_{1}}(h_1\circ \tau^{-1})\pi_{\mathfrak{u}_{2}}(h_2\circ \tau^{-1})\cdots \pi_{\mathfrak{u}_{n_0}}(h_{n_0}\circ \tau^{-1});
\end{align}
and for any $\mathfrak{D}\in \mathcal{E}(\mathcal{H})$ and $\omega\in \mathcal{H}^\infty$, we define
\begin{align*}
 (\mathcal{E}_{\vec{h},\tau}\mathfrak{D})(\omega)\stackrel{\text{def}}{=}\mathfrak{D}\big(\pi_{\mathfrak{u}_{1}}(h_{1}\circ \tau^{-1}) \pi_{\mathfrak{u}_{2}}(h_2\circ \tau^{-1})\cdots\pi_{\mathfrak{u}_{n_0}}(h_{n_0}\circ \tau^{-1})\omega \big).
\end{align*}

\item For any $\textbf{v}\in \mathfrak{S}^1$ and any $1\leq i\leq n_0$
\begin{gather*}
 \big((\textbf{v}\dag_{i})\mathcal{E}_{\vec{h},\tau}\mathfrak{D}\big)(\omega)\\
 \stackrel{\text{def}}{=}\mathfrak{D}\Big(\pi_{\mathfrak{u}_{1}}(h_1\circ \tau^{-1})
 \cdots\big(\pi_{\mathfrak{u}_{i}}(h_{i}\circ \tau^{-1})\textbf{v}\pi_{\mathfrak{u}_{i+1}}(h_{i+1}\circ \tau^{-1})\big)
 \cdots \pi_{\mathfrak{u}_{n_0}}(h_{n_0}\circ \tau^{-1})\omega\Big).
\end{gather*}
Then we see that $\textbf{v}\dag_{i}$ means that $\textbf{v}$ is inserted between $\pi_{\mathfrak{u}_{i}}$ and $\pi_{\mathfrak{u}_{i+1}}$.

\item  For any $\textbf{v}\in \mathfrak{S}^1$ we define
\begin{gather*}
\big((\textbf{v}\dag_{0})\mathcal{E}_{\vec{h},\tau}\mathfrak{D}\big)(\omega)\\
  \mathfrak{D}\big(\textbf{v}\pi_{\mathfrak{u}_{1}}(h_{1}\circ \tau^{-1}) \pi_{\mathfrak{u}_{2}}(h_2\circ \tau^{-1})\cdots\pi_{\mathfrak{u}_{n_0}}(h_{n_0}\circ \tau^{-1})\omega \big).
\end{gather*}
\item Let $S_k$ be the subgroup generated by one parameter unipotent subgroups $\{\mathfrak{u}_{i}\}$, where $k\leq i\leq n_0$. We note that $S_1=Q'$.  Assumption \eqref{for:71} shows that
$\text{Lie}(S_k)$ is the linear space spanned by $\mathfrak{u}_{i}$, $i\geq k$.
\end{enumerate}

\subsection{Proof of Proposition \ref{cor:6}}

By \eqref{for:7} of Lemma \ref{cor:1} we see that $\mathcal{E}_{\vec{h},\tau}\mathfrak{D}\in \mathcal{O}^{-r_0}$ with estimate
\begin{align}\label{for:65}
  \norm{\mathcal{E}_{\vec{h},\tau}\mathfrak{D}}_{-r_0}\leq C_{\vec{h}}\norm{\mathfrak{D}}_{-r_0}.
\end{align}

\emph{Step 1}: For each $0\leq k\leq n_0$,  we show that: for any $v\in \text{Lie}(S_k)\cap \mathfrak{S}^1$ and any $i$ we have
\begin{align}\label{for:68}
 (v\dag_{i})\mathcal{E}_{\vec{h},\tau}\mathfrak{D}=\sum_{l}a_l\tau^{k_l}\mathcal{E}_{\overrightarrow{h_l},\tau}\mathfrak{D}
\end{align}
where $a_l\in\CC$ is uniformly bounded, $k_l\leq 1$ and $\overrightarrow{h_l}\hookrightarrow \vec{h}$ for each $l$.

\medskip
Assumption \eqref{for:71} shows that $S_{n_0}$ is inside the center of $S_1$. By \eqref{for:10} of Section \ref{sec:5} we have
\begin{align*}
 \mathfrak{u}_{n_0}\pi_{\mathfrak{u}_i}(h_i\circ \tau^{-1})=\pi_{\mathfrak{u}_i}(h_i\circ \tau^{-1})\mathfrak{u}_{n_0},\quad \forall\,i.
 \end{align*}
Then we have
 \begin{align*}
  (\mathfrak{u}_{n_0}\dag_{i})\mathcal{E}_{\vec{h},\tau}&=\pi_{\mathfrak{u}_1}(h_1\circ \tau^{-1})\cdots
 \big(\mathfrak{u}_{n_0}\pi_{\mathfrak{u}_{n_0}}(h_{n_0}\circ \tau^{-1})\big)\\
 &\overset{\text{(1)}}{=}\tau\pi_{\mathfrak{u}_1}(h_1\circ \tau^{-1})\cdots
 \pi_{\mathfrak{u}_{n_0}}(h'_{n_0}\circ \tau^{-1})      \qquad \forall\,i,
 \end{align*}
 where $h'_{n_0}=h_{n_0}(t)(t\textrm{i})^1$. Here in $(1)$ we use \eqref{for:41} of Section \ref{sec:5}. Hence, we have
\begin{align}\label{for:19}
(\mathfrak{u}_{n_0}\dag_{i})\mathcal{E}_{\vec{h},\tau}\mathfrak{D}=\tau\mathcal{E}_{\vec{g},\tau}\mathfrak{D},\qquad \forall\,i
\end{align}
 where $(\vec{g})_{n_0}=h_{n_0}(t)(t\textrm{i})^1$ and $(\vec{g})_{j}=h_j$ otherwise. This implies that \eqref{for:68} holds for any $v\in \text{Lie}(S_{n_0})\cap \mathfrak{S}^1$ and any $i$.
 We proceed by induction. Suppose
\eqref{for:68} holds for any $v\in \text{Lie}(S_{j})\cap \mathfrak{S}^1$ and any $i$ if $j\geq k+1$.

To make the induction work, it suffices to show that \eqref{for:68} holds for $\mathfrak{u}_{k}$ and any $i$.
We  consider $(\mathfrak{u}_{k}\dag_{i})\mathcal{E}_{\vec{h},\tau}\mathfrak{D}$. If $i=k$, then \eqref{for:19} holds if we replace $n_0$ by $k$.
If $i<k$, \eqref{for:98} of Section \ref{sec:5} shows that
\begin{align*}
(\mathfrak{u}_{k}\dag_{i})\mathcal{E}_{\vec{h},\tau}\mathfrak{D}=(\mathfrak{u}_{k}\dag_{i+1})\mathcal{E}_{\vec{h},\tau}\mathfrak{D}
+\sum_{j\geq1}c_j\tau^{-j}(v_j\dag_{i+1})\mathcal{E}_{\overrightarrow{g_j},\tau}\mathfrak{D}
\end{align*}
where $v_j\in \text{Lie}(S_{k+1})$ for each $j$; and $(\overrightarrow{g_j})_{i}=h_i^{(j)}$ and $(\overrightarrow{g_j})_{l}=h_{l}$ otherwise for each $j$. Since $v_j\in \text{Lie}(S_{k+1})$ for each $j$,  by using inductive assumption we have
\begin{align*}
(\mathfrak{u}_{k}\dag_{i})\mathcal{E}_{\vec{h},\tau}\mathfrak{D}&=(\mathfrak{u}_{k}\dag_{i+1})\mathcal{E}_{\vec{h},\tau}\mathfrak{D}
+\sum_{r}b_r\tau^{-k_r}\mathcal{E}_{\overrightarrow{p_r},\tau}\mathfrak{D}
\end{align*}
where $b_r\in\CC$,  $k_r\leq0$ and $\overrightarrow{p_r}\hookrightarrow \vec{h}$ for each $r$.

Keeping repeating the above arguments, we see that
\begin{align*}
 (\mathfrak{u}_{k}\dag_{i})\mathcal{E}_{\vec{h},\tau}\mathfrak{D}&=(\mathfrak{u}_{k}\dag_{k})\mathcal{E}_{\vec{h},\tau}\mathfrak{D}
 +\sum_{r}b'_r\tau^{-k'_r}\mathcal{E}_{\overrightarrow{p_r},\tau}\mathfrak{D}\\
 &=\tau\mathcal{E}_{\vec{g},\tau}\mathfrak{D}+\sum_{r}b'_r\tau^{-k'_r}\mathcal{E}_{\overrightarrow{p_r},\tau}\mathfrak{D},
\end{align*}
where $b_r'\in\CC$,  $k'_r\leq0$ and $\overrightarrow{p_r}\hookrightarrow \vec{h}$ for each $r$; $g$ is as defined in \eqref{for:19} by replacing $n_0$ by $k$.

Hence we finish the proof for $i<k$. Similar arguments hold for $i>k$ if we use \eqref{for:234} of Section \ref{sec:5}. The above discussion shows that the result holds for any $v\in \text{Lie}(S_{k})\cap \mathfrak{S}^1$ and any $i$.  Hence we finish the proof.

\smallskip
\emph{Step 2}: Keeping using Step 1,  we have: for any $v_j\in \text{Lie}(S_1)\cap \mathfrak{S}^1$, $1\leq j\leq m$ and any $i_j$ we have
\begin{align*}
 v_1\cdots v_m\mathcal{E}_{\vec{h},\tau}\mathfrak{D}=\sum_{l}a_l\tau^{k_l}\mathcal{E}_{\overrightarrow{g_l},\tau}\mathfrak{D}
\end{align*}
where $a_l\in\CC$ is uniformly bounded, $k_l\leq m$ and  $\overrightarrow{g_l}\hookrightarrow \vec{h}$ for each $l$.  It follows from \eqref{for:65} that
\begin{align}\label{for:32}
 \norm{v_1\cdots v_m\mathcal{E}_{\vec{h},\tau}\mathfrak{D}}_{-r_0}\leq C_{\vec{h}, m}\tau ^m E_0,\quad \forall\,m\geq0.
\end{align}

\emph{Step 3}:  We show that: for any $v\in\text{Lie}(Q_0)\cap \mathfrak{S}^1$
and any $i$ we have
\begin{align*}
 (v\dag_{i})\mathcal{E}_{\vec{h},\tau}\mathfrak{D}=\sum_n b_n\tau^{j_n}\mathcal{E}_{\overrightarrow{\vartheta_n},\tau}(v_n\mathfrak{D})+\sum_{l}a_l\tau^{k_l}\mathcal{E}_{\overrightarrow{g_l},\tau}\mathfrak{D}
\end{align*}
where $b_n,a_l\in\CC$ are uniformly bounded, $j_n,k_l\leq 0$, $v_n\in \text{Lie}(Q_0)$, $\overrightarrow{\vartheta_n}\hookrightarrow \vec{h}$ and $\overrightarrow{g_l}\hookrightarrow \vec{h}$ for each $l,\,n$.

\medskip
We note that
\begin{align*}
 (v\dag_{0})\mathcal{E}_{\vec{h},\tau}D=-\mathcal{E}_{\vec{h},\tau}(vD)
\end{align*}
Then the argument holds for $i=0$. We proceed by induction. Suppose
the argument holds for $(v\dag_{j})\mathcal{E}_{\vec{h},\tau}D$, if $j\leq k-1$. We consider $(v\dag_{k})\mathcal{E}_{\vec{h},\tau}D$.
\eqref{for:234} of Section \ref{sec:27} shows that
\begin{align}\label{for:74}
 (v\dag_{k})\mathcal{E}_{\vec{h},\tau}D=(v\dag_{k-1})\mathcal{E}_{\vec{h},\tau}D+\sum_{j\geq1}d_j\tau^{-j}(v_j\dag_{k-1})\mathcal{E}_{\overrightarrow{g_j},\tau}D
\end{align}
where $v_j\in \mathfrak{S}$ and $\overrightarrow{g_j}\hookrightarrow \vec{h}$ for each $j$.

By assumption \eqref{for:72} we can rewrite \eqref{for:74} as
\begin{align*}
 (v\dag_{k})\mathcal{E}_{\vec{h},\tau}D&=(v\dag_{k-1})\mathcal{E}_{\vec{h},\tau}D+\sum_{j\geq1}d_j
 \tau^{-j}(u_j\dag_{k-1})\mathcal{E}_{\overrightarrow{g_j},\tau}D\\
 &+\sum_{j\geq1}d_j\tau^{-j}(w_j\dag_{k-1})\mathcal{E}_{\overrightarrow{g_j},\tau}D
\end{align*}
where $u_j\in \text{Lie}(Q_0)$ and $w_j\in \text{Lie}(Q')$ for each $j$. By using inductive assumption and \emph{Step 2}, we have
\begin{align*}
 (v\dag_{k})\mathcal{E}_{\vec{h},\tau}D =\sum_no_n\tau^{j_n}\mathcal{E}_{\overrightarrow{q_n},\tau}(v_nD)+\sum_{l}e_l\tau^{k_l}\mathcal{E}_{\overrightarrow{y_l},\tau}D
\end{align*}
where $o_n,e_l\in\CC$, $j_n,k_l\leq 0$, $v_n\in \text{Lie}(Q_0)$, $\overrightarrow{q_n}\hookrightarrow \vec{h}$ and $\overrightarrow{y_l}\hookrightarrow \vec{h}$ for each $l,\,n$.

  Then the result holds for $i=k$.  Hence we finish the proof.

\smallskip
\emph{Step 4}:   Keeping using Step 3,  we have: for any $v_j\in \text{Lie}(Q_0)\cap \mathfrak{S}^1$, $1\leq j\leq m$ and any $i_j$ we have
\begin{align*}
 v_1\cdots v_m\mathcal{E}_{\vec{h},\tau}D=\sum_nb_n\tau^{j_n}\mathcal{E}_{\overrightarrow{\vartheta_n},\tau}(\textbf{v}_nD)+\sum_{l}a_l\tau^{k_l}\mathcal{E}_{\overrightarrow{g_l},\tau}D
\end{align*}
where $b_n,a_l\in\CC$ are uniformly bounded, $j_n,k_l\leq 0$, $\textbf{v}_n\in U_{m}(\text{Lie}(Q_0))$ (see \eqref{for:67} of Section \ref{sec:3}), $\overrightarrow{\vartheta_n}\hookrightarrow \vec{h}$ and $\overrightarrow{g_l}\hookrightarrow \vec{h}$ for each $l,\,n$. It follows from \eqref{for:65} and \eqref{for:75} that
\begin{align}\label{for:34}
 \norm{v_1\cdots v_m\mathcal{E}_{\vec{h},\tau}D}_{-r_0}\leq C_{m}E_0+C_mE_m,\quad \forall\,m\leq s.
\end{align}
As a result, \eqref{for:70} is from  assumption \eqref{for:72}, \eqref{for:32}, \eqref{for:34} and
Lemma \ref{le:3}.

\section{Proof of \eqref{for:138} of Section \ref{sec:19}}\label{sec:43}
It is harmless to assume that $v=\textbf{a}_1$. To simplify notations, we write $z=w\cdot \textbf{a}_1$, $w\in V^o$. We note that
$dz=\text{Ad}_{w}\cdot d\textbf{a}_1$. Since $\log w\in \mathfrak{g}_1$, from \eqref{for:168}, we see that
each $\mathfrak{g}_{\lambda}$, $\lambda\in \Delta$ is invariant under $\text{Ad}_{w}$. As each $\mathfrak{g}_{\lambda}$, $\lambda\in \Delta$ is also invariant under $ d\textbf{a}_1$, they are invariant subspaces for $dz$.

Next, we show that: for any $n\in\ZZ$, we have $z^n=w_{n}\textbf{a}_1^n$,  where $w_n=w_{1,n}w_{2,n}\cdots w_{n,n}$, $w_{j,n}\in V^o$ with estimates
\begin{align}\label{for:495}
 \norm{\log w_{j,n}}\leq C(1+|j|)^{\dim\mathcal{V}}\norm{\log w},\quad   1\leq j\leq n.
\end{align}
and
\begin{align}\label{for:352}
 \norm{\text{Ad}_{w_{n}}}&\leq C_{\norm{w}}
 (|n|+1)^{\dim \mathcal{V}+(\dim \mathcal{V})^2}.
\end{align}
We note that
\begin{align}\label{for:335}
 z^n=w_{n}\textbf{a}_1^n,\quad \text{where }w_{n}=w\textbf{a}_1(w)\cdots \textbf{a}^{n-1}_1(w)
\end{align}
for any $n>0$.  We note that
\begin{align*}
 z^{-1}=\textbf{a}_1^{-1} (w^{-1})\cdot \textbf{a}_1^{-1}=w'\textbf{a}_i^{-1},\quad \text{where }\log w'=-d\textbf{a}_1^{-1} (\log w).
\end{align*}
Thus, if $n<0$ we have
\begin{align*}
 z^n=w_{n}\textbf{a}_1^n,\quad \text{where }w_{n}=w'\textbf{a}^{-1}_1(w')\cdots \textbf{a}^{n+1}_1(w').
\end{align*}
By Engel's theorem $\text{Ad}_{V}$ can all be simultaneously brought to an upper triangular form. Hence, we have
\begin{align}\label{for:493}
 \norm{\text{Ad}_{w_{n}}}\leq C(|n|+1)^{\dim \mathcal{V}}\max_{0\leq j\leq n-1}\{\norm{\text{Ad}_{\textbf{a}^{j}_1(w)}},\,\norm{\text{Ad}_{\textbf{a}^{-j}_1(w')}}\}^{\dim \mathcal{V}}.
\end{align}
We note that
\begin{align*}
 \log(\textbf{a}^{j}_1(w))=d\textbf{a}_1^j(\log w)\quad\text{and}\quad \log(\textbf{a}^{-j}_1(w'))=-d\textbf{a}_1^{-j-1}(\log w),
\end{align*}
which gives
\begin{align}\label{for:494}
&\max\{\norm{\text{Ad}_{\textbf{a}^{j}_1(w)}},\,\norm{\text{Ad}_{\textbf{a}^{-j}_1(w')}}\}\notag\\
&\leq C\max\{\norm{\log(\textbf{a}^{j}_1(w))},\,\norm{\log(\textbf{a}^{-j}_1(w'))}\}^{\dim \mathcal{V}}\notag\\
&\leq C\max\{\norm{d\textbf{a}_1^j(\log w)},\,\norm{d\textbf{a}_1^{-j-1}(\log w)}\}^{\dim \mathcal{V}}\notag\\
&\overset{\text{(1)}}{\leq}C_1\big((1+|j|)^{\dim\mathcal{V}}\norm{\log w}\big)^{\dim \mathcal{V}}.
\end{align}
Here in $(1)$ we use  \eqref{for:334} of Section \ref{sec:19} and the fact $w\in V^o$.

The above discussion shows that \eqref{for:495} holds. It follows from \eqref{for:493} and \eqref{for:494} that
\begin{align*}
 \norm{\text{Ad}_{w_{n}}}&\leq C(|n|+1)^{\dim \mathcal{V}}\Big\{\big(C(1+|n|)^{\dim\mathcal{V}}\norm{\log w}\big)^{\dim \mathcal{V}}\Big\}^{\dim \mathcal{V}}\notag\\
 &\leq C_{\norm{w}}
 (|n|+1)^{\dim \mathcal{V}+(\dim \mathcal{V})^2}.
\end{align*}
Thus, we get \eqref{for:352}.

By \eqref{for:352} and \eqref{for:334} of Section \ref{sec:19}, for any $n\in\ZZ$, we have
\begin{align*}
 \norm{dz^n|_{\mathfrak{g}_{\lambda}}}&\leq \norm{\text{Ad}_{w_{n}}}\norm{d\textbf{a}_1^n|_{\mathfrak{g}_{\lambda}}}\notag\\
 &\leq C_{\norm{w}}(|n|+1)^{\dim \mathcal{V}+(\dim \mathcal{V})^2}\cdot C|\lambda(d\textbf{a}_1)|^n(1+|n|)^{\dim \mathcal{V}}\\
 &\leq C_{\norm{w},1}|\lambda(d\textbf{a}_1)|^n(|n|+1)^{2\dim \mathcal{V}+(\dim \mathcal{V})^2}.
\end{align*}
Hence, we finish the proof.

\section{Proof of \eqref{sec:45} of Section \ref{sec:19}}\label{sec:58}
\begin{lemma}\label{for:353} Suppose $z_v=w_v\cdot v$, where $w_v\in V^o$ and $v\in E$. Then for any $m,\,n\in\ZZ$ we can write
\begin{align*}
  z_v^mz_u^n=w_{v,u,m,n}v^mu^n,\qquad\text{where }w_{v,u,m,n}\in V^o
\end{align*}
with estimates
\begin{align*}
&\norm{\text{Ad}_{w_{v,u,m,n}}}\\
 &\leq C_{\norm{w_v},\norm{w_u}}(|n|+1)^{\dim \mathcal{V}+\dim \mathcal{V}^2+\dim \mathcal{V}^3}(|m|+1)^{\dim \mathcal{V}+\dim \mathcal{V}^3}.
\end{align*}

\end{lemma}
\begin{proof}
It follows from \eqref{for:352} of Section \ref{sec:43}  that we can write
\begin{align*}
 z_v^n=w_{1,n}v^n\quad\text{ and }\quad z_u^m=w_{2,m}u^m,\quad\forall\,m,\,n\in\ZZ
\end{align*}
where $w_{1,n},\,w_{2,m}\in V^o$ with estimates
\begin{align}\label{for:497}
 \norm{\text{Ad}_{w_{1,n}}}&\leq C_{\norm{w_v}}
 (|n|+1)^{\dim \mathcal{V}+(\dim \mathcal{V})^2}
\end{align}
and $w_{2,m}=w_{2,m,1}w_{2,m,2}\cdots w_{2,m,m}$, $w_{2,m,j}\in V^o$ with estimates
\begin{align}\label{for:496}
 \norm{\log w_{2,m,j}}\leq C(1+|j|)^{\dim\mathcal{V}}\norm{\log w_u},\quad   1\leq j\leq m.
\end{align}
We note that
\begin{align*}
 z_v^nz_u^m=w_{1,n}v^n(w_{2,m}u^m)=w_{1,n}v^n(w_{2,m})\cdot v^nu^m,
\end{align*}
where
\begin{align*}
 v^n(w_{2,m})=(v^nw_{2,m,1})(v^n(w_{2,m,2}))\cdots (v^n(w_{2,m,m})).
\end{align*}
Let $w_{v,u,m,n}=w_{1,n}v^n(w_{2,m})$. It is clear that $w_{v,u,m,n}\in V^o$.

By Engel's theorem $\text{Ad}_{V}$ can all be simultaneously brought to an upper triangular form. Hence, we have
\begin{align}\label{for:498}
 \norm{\text{Ad}&_{v^n(w_{2,m})}}\leq C(|m|+1)^{\dim \mathcal{V}}\max_{1\leq j\leq m}\{\norm{\text{Ad}_{v^n(w_{2,m,j})}}\}^{\dim \mathcal{V}}\notag\\
 &\leq C_1(|m|+1)^{\dim \mathcal{V}}\max_{1\leq j\leq m}\Big\{\big(\norm{ dv^n\log w_{2,m,j}}\big)^{\dim \mathcal{V}}\Big\}^{\dim \mathcal{V}}\notag\\
 &\overset{\text{(1)}}{\leq} C_1(|m|+1)^{\dim \mathcal{V}}\notag\\
 &\cdot\max_{1\leq j\leq m}\Big\{\big(C_1(1+|n|)^{\dim\mathcal{V}}C(1+|m|)^{\dim\mathcal{V}}\norm{\log w_2}\big)^{\dim \mathcal{V}}\Big\}^{\dim \mathcal{V}}\notag\\
 &\leq C_{\norm{w_u}}(|m|+1)^{\dim \mathcal{V}+(\dim \mathcal{V})^3}(|n|+1)^{(\dim \mathcal{V})^3}.
\end{align}
Here in $(1)$ we use \eqref{for:334} of Section \ref{sec:19} and the fact $w_{2,m,j}\in V^o$:
\begin{align*}
 \norm{ dv^n\log w_{2,m,j}}&\leq C(1+|n|)^{\dim\mathcal{V}}\norm{\log w_{2,m,j}}\\
 &\overset{\text{(a)}}{\leq} C_1(1+|n|)^{\dim\mathcal{V}}C(1+|m|)^{\dim\mathcal{V}}\norm{\log w_u}.
\end{align*}
Here in $(a)$ we use \eqref{for:496}.

Hence, we have
\begin{align*}
&\norm{\text{Ad}_{w_{1,n}v^n(w_{2,m})}}\leq \norm{\text{Ad}_{w_{1,n}}}\cdot \norm{\text{Ad}_{v^n(w_{2,m})}}\notag\\
&\overset{\text{(1)}}{\leq}C_{\norm{w_v}}
 (|n|+1)^{\dim \mathcal{V}+\dim \mathcal{V}^2}\notag\\
 &\quad\cdot C_{\norm{w_u}}(|m|+1)^{\dim \mathcal{V}+\dim \mathcal{V}^3}(|n|+1)^{\dim \mathcal{V}^3}\notag\\
 &\leq C_{\norm{w_v},\norm{w_u}}(|n|+1)^{\dim \mathcal{V}+\dim \mathcal{V}^2+\dim \mathcal{V}^3}(|m|+1)^{\dim \mathcal{V}+\dim \mathcal{V}^3}.
\end{align*}
Here in $(1)$ we use \eqref{for:497} and \eqref{for:498}.

Hence we finish the proof.
\end{proof}
\begin{remark}\label{re:12} \eqref{for:352} of Section \ref{sec:43} shows that both the Lyapunov
exponents and the Lyapunov subspace of $dz_v$ are determined by $dv$. Lemma \ref{for:353} shows that if we denote by $\tilde{A}$ the subgroup of affine actions on $V$ generated by $z_v$, $v\in E$, even though $\tilde{A}$ is probably not abelian, \eqref{for:349} of Section \ref{sec:19} is also the Lyapunov subspace decomposition for $\tilde{A}$:
\begin{align*}
 \mathcal{V}= \sum_{\lambda\in \Delta}\mathfrak{g}_{\lambda}
\end{align*}
such that for any $\lambda\in \Delta$ and any $\tilde{a}\in \tilde{A}$, $\lambda(da)$ is the absolute value of the eigenvalue of $d\tilde{a}$ on $\mathfrak{g}_\lambda$, where $a$ is the automorphism part of $\tilde{a}$.
\end{remark}
It follows that Lemma \eqref{for:353} that we can write
\begin{align*}
  z_1^mz_1^n=w_{m,n}\textbf{a}_1^n\textbf{a}_2^m,\qquad\text{where }w_{m,n}\in V^o
\end{align*}
with estimates
\begin{align}\label{for:503}
&\norm{\text{Ad}_{w_{v,u,m,n}}}\notag\\
 &\leq C_{\norm{w_1},\norm{w_2}}(|n|+1)^{\dim \mathcal{V}+\dim \mathcal{V}^2+\dim \mathcal{V}^3}(|m|+1)^{\dim \mathcal{V}+\dim \mathcal{V}^3}.
\end{align}
Then we have
\begin{align*}
\abs{\langle f\circ (z_1^nz_2^m), g\rangle}&=\big|\langle f(w_{m,n}\cdot \textbf{a}_1^n\textbf{a}_2^m x), g(x)\rangle\big|\\
&=\big|\langle f(w_{m,n}x), g(\textbf{a}_2^{-m}\textbf{a}_1^{-n}x)\rangle\big|\\
&\overset{\text{(1)}}{\leq}
E\big\|f\circ (w_{m,n}))\big\|_{C^1}\norm{g}_{C^1}e^{-\gamma(|n|+|m|)}\\
&\leq E\norm{\text{Ad}_{w_{m,n})}}\norm{f}_{C^1}\norm{g}_{C^1}e^{-\gamma(|n|+|m|)}\\
&\overset{\text{(2)}}{\leq} C_{\norm{w_1},\norm{w_2}}(|n|+1)^{\dim \mathcal{V}+(\dim \mathcal{V})^2+(\dim \mathcal{V})^3}\\
&\cdot(|m|+1)^{\dim \mathcal{V}+(\dim \mathcal{V})^3} \\
&\cdot\norm{f}_{C^1}\norm{g}_{C^1}e^{-\gamma(|n|+|m|)}.
\end{align*}
Here in $(1)$ we use Theorem \ref{th:1}; in $(2)$ we use \eqref{for:503}.

Let $\gamma_1=\frac{\gamma}{2}$. Then we finish the proof.
\section{Proof of Lemma \ref{le:1}}\label{sec:12}

By commutativity of $\tilde{\alpha}_A$  we have
\begin{align}\label{for:340}
&\exp\bigl(\mathfrak{p}_{u}\bigl(\exp(\mathfrak{p}_{v}(x))\cdot z_vx\bigl)\bigl)\cdot
z_u\big(\exp(\mathfrak{p}_{v}(x))\cdot z_vx\big)\notag\\
&=\exp\bigl(\mathfrak{p}_{v}\bigl(\exp(\mathfrak{p}_{u}(x))\cdot z_ux\bigl)\bigl)\cdot z_v\big(\exp(\mathfrak{p}_{u}(x))\cdot z_ux\big)
\end{align}
for any $x\in \mathcal{X}$. Then we have
\begin{align}\label{for:339}
&\exp\bigl(\mathfrak{p}_{u}\bigl(\exp(\mathfrak{p}_{v}(x))\cdot z_vx\bigl)\bigl)\cdot
\exp\big(dz_u\mathfrak{p}_{v}(x)\big)\cdot z_uz_v(x)\notag\\
&=\exp\bigl(\mathfrak{p}_{v}\bigl(\exp(\mathfrak{p}_{u}(x))\cdot z_ux\bigl)\bigl)\cdot \exp\big(dz_v\mathfrak{p}_{u}(x)\big)\cdot z_vz_u(x).
\end{align}
Firstly, we prove \eqref{for:338} of Remark \ref{re:11}, which is the special case of
$z_v=v$ for each $v\in E$. If $\norm{\mathfrak{p}}_{C^m}$ is sufficiently small, from \eqref{for:339}, by the  Baker-Campbell-Hausdorff formula there are
$c_{u,v}:\mathcal{X}\to \mathcal{V}$ with
\begin{align*}
 \norm{c_{u,v}}_{C^m}\leq C_m\norm{\mathfrak{p}}_{C^m}^2,\qquad i=1,2
\end{align*}
such that
\begin{align*}
 \mathfrak{p}_{u}\bigl(\exp(\mathfrak{p}_{v}(x))\cdot vx\bigl)+du\mathfrak{p}_{v}(x)=\mathfrak{p}_{v}\bigl(\exp(\mathfrak{p}_{u}(x))\cdot ux\bigl)+dv\mathfrak{p}_{u}(x)+c_{u,v}(x).
\end{align*}
Then we have
\begin{align*}
\mathbb{D}_{u,v}(x)&=\mathfrak{p}_{u}(vx)-\mathfrak{p}_{u}\bigl(\exp(\mathfrak{p}_{v}(x))\cdot vx\bigl)\\
 &-\big(\mathfrak{p}_{v}(ux)-\mathfrak{p}_{v}\bigl(\exp(\mathfrak{p}_{u}(x))\cdot ux\bigl)\big)+c_{u,v}(x).
\end{align*}
We note that (see \cite[Appendix II]{Lazutkin})
\begin{align*}
\big\| \mathfrak{p}_{u}(vx)-\mathfrak{p}_{u}\bigl(\exp(\mathfrak{p}_{v}(x))\cdot vx\bigl)\big\|_{C^m}&\leq C_m\norm{\mathfrak{p}_u}_{C^{m+1}}\norm{\mathfrak{p}_v}_{C^0}\\
\big\| \mathfrak{p}_{v}(ux)-\mathfrak{p}_{v}\bigl(\exp(\mathfrak{p}_{u}(x))\cdot ux\bigl)\big\|_{C^m}&\leq C_m\norm{\mathfrak{p}_v}_{C^{m+1}}\norm{\mathfrak{p}_u}_{C^0}
\end{align*}
It follows that
\begin{align*}
 \|\mathbb{E}_{u,v}(x)\|_{C^m}&\leq C_m\norm{\mathfrak{p}_u}_{C^{m+1}}\norm{\mathfrak{p}_v}_{C^0}+C_m\norm{\mathfrak{p}_v}_{C^{m+1}}\norm{\mathfrak{p}_u}_{C^0}+\norm{c_{u,v}}_{C^m}\\
 &\leq C_{m,1}\norm{\mathfrak{p}}_{C^m}\norm{\mathfrak{p}}_{C^{m+1}}.
\end{align*}
Hence we finish the proof.

Next,  we proceed to the proof of Lemma \ref{le:1}. \eqref{for:347} is  a direct consequence of \eqref{for:339}.

\smallskip

\eqref{for:337}: \emph{Step 1}: We show that there are
$c_v:\mathcal{X}\to \mathcal{V}$, $v\in E$ with estimates
\begin{align}\label{for:341}
  \norm{c_v}_{C^m}\leq C_m\norm{\mathfrak{p}}_{C^{m+1}}\norm{\mathfrak{p}}_{C^{m}},\qquad v\in E
\end{align}
such that
\begin{gather*}
\exp\bigl(\mathfrak{p}'_{u}(z_v\cdot
x)+c_u(x)\bigl)\\
=\exp\bigl(\mathfrak{p}_{u}\bigl(\exp(\mathfrak{p}_{2}(x))\cdot z_v\cdot
x\bigl)\bigl)\exp(-r_u),\\
\exp\bigl(\mathfrak{p}'_{v}(z_u\cdot
x)+c_v(x)\bigl)\\
=\exp\bigl(\mathfrak{p}_{v}\bigl(\exp(\mathfrak{p}_{u}(x))\cdot z_u\cdot
x\bigl)\bigl)\exp(-r_v).
\end{gather*}
By the  Baker-Campbell-Hausdorff formula there are
$c_v':\mathcal{X}\to \mathcal{V}$ with
\begin{align*}
 \norm{c'_v}_{C^m}\leq C_m\norm{\mathfrak{p}}_{C^m}^2,\qquad v\in E
\end{align*}
such that
\begin{gather*}
\exp\bigl(\mathfrak{p}_{u}\bigl(\exp(\mathfrak{p}_{v}(x))\cdot z_v\cdot
x\bigl)\bigl)\exp(-r_u)\\
=\exp\bigl(\mathfrak{p}_{u}\bigl(\exp(\mathfrak{p}_{v}(x))\cdot z_v\cdot
x\bigl)-r_u+c_u'(x)\bigl)
\end{gather*}
and
\begin{gather*}
\exp\bigl(\mathfrak{p}_{v}\bigl(\exp(\mathfrak{p}_{u}(x))\cdot z_u\cdot
x\bigl)\bigl)\exp(-r_v)\\
=\exp\bigl(\mathfrak{p}_{v}\bigl(\exp(\mathfrak{p}_{u}(x))\cdot z_u\cdot
x\bigl)-r_v+c_v'(x)\bigl).
\end{gather*}
We also note that (see \cite[Appendix II]{Lazutkin})
\begin{align*}
\big\| \mathfrak{p}_{u}\bigl(\exp(\mathfrak{p}_{v}(x))\cdot z_v\cdot
x\bigl)-\mathfrak{p}_{u}(z_v\cdot
x)\big\|_{C^m}&\leq C_m\norm{\mathfrak{p}}_{C^{m+1}}\norm{\mathfrak{p}}_{C^0}\\
\big\| \mathfrak{p}_{v}\bigl(\exp(\mathfrak{p}_{u}(x))\cdot z_u\cdot
x\bigl)-\mathfrak{p}_{v}(z_u\cdot
x)\big\|_{C^m}&\leq C_m\norm{\mathfrak{p}}_{C^{m+1}}\norm{\mathfrak{p}}_{C^0}
\end{align*}
The above discussion implies the result.

\smallskip
\noindent\emph{Step 2}: We see that there are
$d_v:\mathcal{X}\to \mathcal{V}$, $v\in E$ with estimates
\begin{align}\label{for:342}
  \norm{d_v}_{C^m}\leq C_m\norm{\mathfrak{p}}_{C^{m}}^2,\qquad v\in E
\end{align}
such that
\begin{gather*}
\exp\bigl(\mathfrak{p}'_{u}(
x)+d_u(x)\bigl)=\exp(\mathfrak{p}_{u}(x))\exp(-r_u),\\
\exp\bigl(\mathfrak{p}'_{v}(
x)+d_v(x)\bigl)=\exp(\mathfrak{p}_{v}(x))\exp(-r_v).
\end{gather*}
The result is a direct consequence of the  Baker-Campbell-Hausdorff formula.

\smallskip

\noindent\emph{Step 3}: We show that
\begin{align}\label{for:345}
   \norm{\mathbb{D}_{u,v}}_{C^m}\leq C_m\norm{\mathfrak{p}}_{C^{m+1}}\norm{\mathfrak{p}}_{C^{m}},\qquad \forall\,u,\,v\in E
\end{align}
By \emph{Step 1}, \emph{Step 2} and \eqref{for:340} we have
\begin{align*}
&\exp\bigl(\mathfrak{p}'_{u}(z_v\cdot
x)+c_u(x)\bigl)\cdot
z_u'\big(\exp\bigl(\mathfrak{p}'_{v}(
x)+d_v(x)\bigl)\cdot z_v'x\big)\notag\\
&=\exp\bigl(\mathfrak{p}'_{v}(z_u\cdot
x)+c_v(x)\bigl)\cdot z_v'\big(\exp\bigl(\mathfrak{p}'_{u}(
x)+d_u(x)\bigl)\cdot z_u'x\big).
\end{align*}
where $z_v'=\exp(r_v)z_v$.

It follows that
\begin{align*}
&\exp\bigl(\mathfrak{p}'_{u}(z_v\cdot
x)+c_u(x)\bigl)\cdot
\exp\Bigl(dz_u'\big(\mathfrak{p}'_{v}(
x)+d_v(x)\big)\Bigl)z_u'z_v'\notag\\
&=\exp\bigl(\mathfrak{p}'_{v}(z_u\cdot
x)+c_v(x)\bigl)\cdot\exp\Bigl(dz_v'\big(\mathfrak{p}'_{u}(
x)+d_u(x)\big)\Bigl)z_v' z_u'.
\end{align*}
By \eqref{for:341}, \eqref{for:342} and the  Baker-Campbell-Hausdorff formula, there are
$q_{u,v}:\mathcal{X}\to \mathcal{V}$, $u,\,v\in E$ with
\begin{align}\label{for:343}
 \norm{q_{u,v}}_{C^m}\leq C_m\norm{\mathfrak{p}}_{C^{m+1}}\norm{\mathfrak{p}}_{C^{m}}
\end{align}
such that
\begin{align}\label{for:499}
 \log[z_v',z_u']=\mathbb{L}_{u,v}(x)+q_{u,v}(x)
\end{align}
(see \eqref{sect:9} of Section \ref{sec:19}) where
\begin{align*}
 \mathbb{L}_{u,v}(x)&=\mathfrak{p}'_{u}(z_v\cdot
x)+d z_u'(\mathfrak{p}'_{v}
x)-\mathfrak{p}'_{v}(z_u\cdot
x)-dz_v'(\mathfrak{p}'_{u}
x).
\end{align*}
It follows that
\begin{align*}
 &\log[z_v',z_u']\text{Vol}(\mathcal{X})=\int_{\mathcal{X}}\mathbb{L}_{u,v}+\int_{\mathcal{X}}q_{u,v}\|\overset{\text{(1)}}{=}\int_{\mathcal{X}}q_{u,v},
\end{align*}
which gives
\begin{align}\label{for:344}
 \norm{\log[z_v',z_u']}\leq C\big\|\int_{\mathcal{M}}q_{u,v}\big\|\overset{\text{(2)}}{\leq} C\norm{\mathfrak{p}}_{C^1}\norm{\mathfrak{p}}_{C^0}.
\end{align}
Here in $(1)$ we recall that both $\mathfrak{p}'_{u}$ and $\mathfrak{p}'_{v}$ have zero averages; in $(2)$ we use \eqref{for:343}.

By using equation \eqref{for:499}, the estimates in \eqref{for:343} and \eqref{for:344} show that
\begin{align}\label{for:346}
 \norm{\mathbb{L}_{u,v}}_{C^m}\leq C_m\norm{\mathfrak{p}}_{C^{m+1}}\norm{\mathfrak{p}}_{C^{m}},\qquad \forall\,u,\,v\in E.
\end{align}
Hence, we have
\begin{align}\label{for:500}
 \norm{\mathbb{D}_{u,v}&-\mathbb{L}_{u,v}}_{C^m}\leq \big\|d z_u'(\mathfrak{p}'_{v}
x)-d z_u(\mathfrak{p}'_{v}x)\big\|_{C^m}+\big\|dz_v'(\mathfrak{p}'_{u}
x)-dz_v(\mathfrak{p}'_{u}
x)\big\|_{C^m}\notag\\
&\overset{\text{(1)}}{\leq}\norm{\text{Ad}_{\exp(r_u)}-I}\big\|d z_u(\mathfrak{p}'_{v}x)\big\|_{C^m}+\norm{\text{Ad}_{\exp(r_v)}-I}\big\|dz_v(\mathfrak{p}'_{u}
x)\big\|_{C^m}\notag\\
&\overset{\text{(2)}}{\leq} C\norm{r_u}\norm{\mathfrak{p}}_{C^m}+C\norm{r_v}\norm{\mathfrak{p}}_{C^m}
\big\|_{C^m}\notag\\
&\leq C\norm{\mathfrak{p}}_{C^0}\norm{\mathfrak{p}}_{C^m}.
\end{align}
Here in $(1)$ we use \eqref{sect:8} of Section \ref{sec:19}; in $(2)$ we use \eqref{sect:9} of Section \ref{sec:19}.

Then \eqref{for:345} follows from the above estimate and \eqref{for:346}. Hence, we finish the proof.

\smallskip

\eqref{for:348}:
By definition we have
\begin{align*}
 \quad[z_u',z_v']\exp(r_v)\cdot z_v\big(\exp(r_u)\cdot z_u\big)&=\exp(r_u)\cdot z_u\big(\exp(r_v)\cdot z_v\big)\\
 \Rightarrow [z_u',z_v']\exp(r_v)\exp(dz_vr_u)\cdot z_vz_u&=\exp(r_u)\exp(dz_ur_v)\cdot z_uz_v\\
 \Rightarrow [z_u',z_v']\exp(r_v)\exp(dz_vr_u)&=\exp(r_u)\exp(dz_ur_v)[z_v,z_u]
\end{align*}
By the  Baker-Campbell-Hausdorff formula, \eqref{for:344} and \eqref{for:347} we have
\begin{align*}
  \norm{r_u-dz_vr_u-r_v+dz_ur_v+\log [z_u,z_v]}\leq C\norm{\mathfrak{p}}_{C^1}\norm{\mathfrak{p}}_{C^0}.
\end{align*}
As $\log[z_u,z_v]\in \mathfrak{g}_{1}$ and each $\mathfrak{g}_{\lambda}$, $\lambda\in \Delta$ is invariant under $dz_v$  and $dz_u$ (see \eqref{sect:10} of Section \ref{sec:19}), the above inequality implies the result.

\section{Proof of Proposition \ref{po:2}}\label{sec:44}

The proof is similar to the case  where $z_1$ and $z_2$ commute, but new arguments are also needed to handle the non-abelianness.
We fix $l>0$ ($l$ is determined in \emph{Step 6}, see \eqref{for:363}).  To simplify notations, we denote
\begin{align*}
B_i=dz_i|_{\mathfrak{g}_{\mu}},\,\,i=1,2\quad\text{and}\quad \Phi_{m}=\sum_{n=-l}^l B_2^{m}B_1^{n-1}\theta_\mu\circ (z_1^{-n}z_2^{-m}).
\end{align*}
Suppose $(\mu,\,\lambda)$ is a poor pair and $u\in \mathfrak{g}_{\lambda}\bigcap \mathcal{V}^1$. Also suppose $\psi\in \mathfrak{g}_\mu(\mathcal{O}^\kappa)$.

\smallskip
\noindent\emph{Step 1: We show that}
\begin{align}\label{for:367}
 \big\|\langle\mathfrak{D}_{\mu,\lambda,u,\tau,\ZZ},\,\psi\rangle\big\|\leq \big\|\langle \Phi_{0},\, \pi_{u}(f_1\circ \tau^{-1})\rangle\big\|+Ce^{-\gamma_2l}\norm{\theta_\mu}_\kappa\norm{\psi}_{\kappa}.
\end{align}
We have
\begin{align*}
    &\big\|\langle\mathfrak{D}_{\mu,\lambda,u,\tau,\ZZ},\,\psi\rangle\big\|\notag\\
    &=\Big\|\sum_{n\in\ZZ}\big\langle (dz_1|_{\mathfrak{g}_{\mu}})^{n-1} \xi\circ z_1^{-n},\, \pi_{u}(f_1\circ \tau^{-1})\psi\big\rangle\Big\|\notag\\
    &\leq \big\|\langle \Phi_{0},\, \pi_{u}(f_1\circ \tau^{-1})\psi\rangle\big\|\\
    &+\big\|\sum_{|n|>l}\langle (dz_1|_{\mathfrak{g}_{\mu}})^{n-1} \pi_{u}(f_1\circ \tau^{-1})(\xi\circ z_1^{-n}),\, \psi\rangle\big\|\notag\\
    &\overset{\text{(1)}}{\leq} \big\|\langle \Phi_{0},\, \pi_{u}(f_1\circ \tau^{-1})\psi\rangle\big\|+\sum_{|n|>l} Ce^{-\gamma_2 |n|}\norm{\theta_\mu}_\kappa\norm{\psi}_{\kappa}\\
    &\leq \big\|\langle \Phi_{0},\, \pi_{u}(f_1\circ \tau^{-1})\psi\rangle\big\|+C_1e^{-\gamma_2l}\norm{\theta_\mu}_\kappa\norm{\psi}_{\kappa}.
  \end{align*}
Here in $(1)$ we use \eqref{for:354} of Proposition \ref{po:6}. Then we get the result.

\smallskip
\noindent\emph{Step 2: We show that}
\begin{align}\label{for:357}
\big\|\langle \Phi_{0},\, \pi_{u}(f_1\circ \tau^{-1})\psi\rangle\big\|\leq\sum_{m\in L_{\mu}}\Big\|\big\langle \Phi_{m}-\Phi_{m+1}, \pi_{u}(f_1\circ \tau^{-1})\psi\big\rangle\Big\|.
\end{align}
For any $m\in L_{\mu}$,  it follows from Lemma \ref{le:2} that
\begin{align}
 \big\|\langle \Phi_{m}, \,\pi_{u}(f_1\circ \tau^{-1})\psi\rangle\big\|&\leq\sum_{n} Ce^{-\gamma_3(|n|+|m|)}\norm{\theta_\mu}_{\kappa}\norm{\psi}_\kappa\notag\\
  &\leq C_1e^{-\gamma_3|m|}\norm{\theta_\mu}_{\kappa}\norm{\psi}_\kappa.
\end{align}
Letting $|m|\to\infty$, we get
\begin{align*}
 \big\|\langle \Phi_{m}, \,\pi_{u}(f_1\circ \tau^{-1})\psi\rangle\big\|\to 0,\qquad\text{as }|m|\to\infty.
\end{align*}
This implies \eqref{for:357}.

\smallskip
\noindent\emph{Step 3: We show that}
\begin{align}\label{for:369}
 \sum_{m\in L_{\mu}}\Big\|\big\langle \Phi_{m}-\Phi_{m+1}, \pi_{u}(f_1\circ \tau^{-1})\psi\big\rangle\Big\|\leq \mathcal{E}_1+\mathcal{E}_2,
\end{align}
where
\begin{gather*}
 \mathcal{E}_1=\sum_{m\in L_{\mu}}\Big\|\langle B_2^{m}\Psi_{l}\circ z_2^{-(m+1)}, \,\pi_{u}(f_1\circ \tau^{-1})\psi\rangle\Big\| \qquad \text{and}\\
 \mathcal{E}_2=\sum_{m\in L_{\mu}}\Big\|\big\langle B_2^{m}\Psi_{l}\circ z_2^{-(m+1)}-(\Phi_{m}-\Phi_{m+1}), \,\pi_{u}(f_1\circ \tau^{-1})\psi\big\rangle\Big\|.
\end{gather*}
Here $\Psi_{l}$ is defined in \eqref{for:358}.

 From equation \ref{for:164}, we have
\begin{align}\label{for:359}
 &\sum_{n=-l}^l B_1^{n-1}\theta_\mu\circ(z_2z_1^{-n})-\sum_{n=-l}^l B_1^{n-1}B_2\theta_\mu\circ z_1^{-n}=\Psi_{l}
\end{align}
where
\begin{align}\label{for:358}
 \Psi_{l}=B_1^{-(l+1)}\varphi_\mu\circ z_1^{l+1} -B_1^{l}\varphi_\mu\circ z_1^{-l}-\sum_{n=-l}^l B_1^{n-1}\omega_\mu \circ z_1^{-n}.
\end{align}
We note that
\begin{align*}
 \Phi_{m}-\Phi_{m+1}&=B_2^{m}\Psi_{l}\circ z_2^{-(m+1)}\notag\\
 &-\big(B_2^{m}\Psi_{l}\circ z_2^{-(m+1)}-(\Phi_{m}-\Phi_{m+1})\big),\qquad \forall\,m\in\ZZ.
\end{align*}
This implies the result.

\smallskip
\noindent\emph{Step 4: We show that}
\begin{align}\label{for:370}
 \mathcal{E}_1\leq C(e^{-\gamma_3l}\norm{\varphi_\mu}_\kappa+\norm{\omega_\mu}_\kappa)\norm{\psi}_{\kappa}.
\end{align}
By using the form of $\Psi_{l}$ in \eqref{for:358}, we see that
\begin{align*}
 \mathcal{E}_1&\leq \sum_{m\in L_{\mu}}\Big\|\langle B_2^{m}B_1^{-(l+1)}\varphi_\mu\circ (z_1^{l+1}z_2^{-(m+1)}), \,\pi_{u}(f_1\circ \tau^{-1})\psi\rangle\Big\|\\
 &+\sum_{m\in L_{\mu}}\Big\|\langle B_2^{m}B_1^{l}\varphi_\mu\circ (z_1^{-l}z_2^{-(m+1)}), \,\pi_{u}(f_1\circ \tau^{-1})\psi\rangle\Big\|\\
 &+\sum_{m\in L_{\mu}}\sum_{n=-l}^l\Big\|\langle B_2^{m}B_1^{n-1}\omega_\mu \circ (z_1^{-n}z_2^{-(m+1)}), \pi_{u}(f_1\circ \tau^{-1})\psi\rangle\Big\|
\end{align*}
It follows from Lemma \ref{le:2} that
\begin{align*}
 \mathcal{E}_1&\leq \sum_{m\in L_{\mu}}Ce^{-\gamma_3((l+1)+|m|+1)}\norm{\varphi_\mu}_{\kappa}\norm{\psi}_\kappa\\
 &+\sum_{m\in L_{\mu}}Ce^{-\gamma_3(l+|m+1|)}\norm{\varphi_\mu}_{\kappa}\norm{\psi}_\kappa\\
 &+\sum_{m\in L_{\mu}}\sum_{n=-l}^lCe^{-\gamma_3(|n|+|m+1|)}\norm{\omega_\mu}_{\kappa}\norm{\psi}_\kappa\\
 &\leq C_1e^{-\gamma_3l}\norm{\varphi_\mu}_{\kappa}\norm{\psi}_\kappa
 +C_1\norm{\omega_\mu}_{\kappa}\norm{\psi}_\kappa.
\end{align*}
\smallskip
\noindent\emph{Step 5: We show that}
\begin{align}\label{for:368}
\mathcal{E}_2\leq \sum_{m\in L_{\mu}}\big\|\langle Y_{1,m}, \,\pi_{u}(f_1\circ \tau^{-1})\psi\rangle\big\|
+\sum_{m\in L_{\mu}}\big\|\langle Y_{2,m},\, \pi_{u}(f_1\circ \tau^{-1})\psi\rangle\big\|
\end{align}
where
\begin{gather}
 Y_{1,m}=\sum_{n=-l}^l B_2^{m}B_1^{n-1}(\theta_\mu\circ s_n-\theta_\mu)\circ(z_1^{-n}z_2^{-m}),\label{for:361}\\
 Y_{2,m}=\sum_{n=-l}^l B_2^{m+1}B_1^{n-1}\big((dk_n|_{J_\mu})-I\big)\theta_\mu\circ (z_1^{-n}z_2^{-(m+1)}), \qquad\text{and}\label{for:362}\\
 s_n=z_2z_1^{-n}z_2^{-1}z_1^n,\quad k_n=z_1^{-(n-1)}z_2^{-1}z_1^{n-1}z_2.\notag
\end{gather}
By using the form of $\Psi_{l}$ in \eqref{for:359}, we see that
\begin{align*}
  B_2^{m}\Psi_{l}\circ z_2^{-(m+1)}&-(\Phi_{m}-\Phi_{m+1})=Y_{1,m}-Y_{2,m}.
\end{align*}
This gives the result.

\smallskip
\noindent\emph{Step 6:}
We recall nations in \eqref{sect:9} of Section \ref{sec:19}. Before we estimate $Y_{i,m}$, $i=1,\,2$,  it is natural for us to estimate $\textbf{\emph{d}}(s_n,e)$ and $\textbf{\emph{d}}(k_n,e)$ with $|n|\leq l$ at first.

 Denote by $\norm{\log[z_1,z_2]}=x^2$, $x>0$. We can assume $x$ is sufficiently small such that $x\leq \delta$ (see \eqref{sect:9} of Section \ref{sec:19}). Choose $l>0$ such that
 \begin{align}\label{for:363}
  e^{l }=x^{-1}
 \end{align}
(see \eqref{sect:10} of Section \ref{sec:19}).
In this part we study the commutator relations between
$z_2^{\pm}$ and $z_1^{n}$. We show that
for any $|n|\leq l$
\begin{align}\label{for:366}
\max\{\textbf{\emph{d}}(s_n,e),\,\textbf{\emph{d}}(k_n,e)\}\leq Cx^{\frac{8}{7}}.
\end{align}
We can write $z_2z_1=\exp(\mathfrak{n})z_1z_2$. Then $\norm{\mathfrak{n}}=x^2$. We note that
\begin{align*}
  z_2z_1z_2^{-1}&=\exp(\mathfrak{n})z_1\Rightarrow z_2z_1^{|n|}z_2^{-1}=(\exp(\mathfrak{n})z_1)^{|n|}\\
  &\Rightarrow z_2z_1^{|n|}z_2^{-1}z_1^{-|n|}=\exp(\mathfrak{n})\exp(dz_1\mathfrak{n})\cdots \exp(dz_1^{|n|-1}\mathfrak{n}).
\end{align*}
From \eqref{for:138} of Section \ref{sec:19}, for any $0\leq j\leq |n|-1\leq l-1$ we have
\begin{align}\label{for:452}
 \norm{dz_i^j\mathfrak{n}}\leq C(1+j)^{c_1}\norm{\mathfrak{n}}\leq Cl^{c_1}x^2\overset{\text{(1)}}{\leq}x.
\end{align}
Here in $(1)$ we recall that $x$ is sufficiently small.

Hence we have
\begin{align}\label{for:360}
 \textbf{\emph{d}}(z_2&z_1^{|n|}z_2^{-1}z_1^{-|n|},\,e)\overset{\text{(1)}}{\leq} \sum_{j=0}^{|n|-1} \textbf{\emph{d}}(\exp(dz_1^{j}\mathfrak{n}),e)
 \overset{\text{(2)}}{\leq} \sum_{j=0}^{|n|-1} C\|dz_1^{j}\mathfrak{n}\|\notag\\
 &\overset{\text{(3)}}{\leq} C_1l \cdot l^{c_1}x^2\overset{\text{(4)}}{\leq}C_2x^{\frac{8}{7}}.
\end{align}
Here in $(1)$ we use right invariance of $\textbf{\emph{d}}$ (see \eqref{sect:9} of Section \ref{sec:19});  in $(2)$ we use
\eqref{sect:9} of Section \ref{sec:19}; in $(3)$ we use \eqref{for:452}; in $(4)$ we recall that $x$ is sufficiently small.

Since $z_2z_1^{|n|}z_2^{-1}z_1^{-|n|}=(z_1^{|n|}z_2z_1^{-|n|}z_2^{-1})^{-1}$ and $x$ is sufficiently small, by using
\eqref{sect:9} of Section \ref{sec:19}, we see that \eqref{for:360} also implies that
\begin{align*}
 \textbf{\emph{d}}(z_1^{|n|}z_2z_1^{-|n|}z_2^{-1},\,e)\leq C_3x^{\frac{8}{7}}.
\end{align*}
We also note that
\begin{align*}
  z_2z_1z_2^{-1}&=\exp(\mathfrak{n})z_1\Rightarrow z_2z_1^{-|n|}z_2^{-1}=(\exp(\mathfrak{n})z_1)^{-|n|}\\
  &\Rightarrow z_2z_1^{-|n|}z_2^{-1}\mathfrak{n}z_1^{|n|}=\exp(-dz_1^{-1}\mathfrak{n})\exp(-dz_1^{-2}\mathfrak{n})\cdots \exp(dz_1^{-|n|}\mathfrak{n}).
\end{align*}
By using the same reasoning, we can also show that
\begin{align*}
 \textbf{\emph{d}}(z_2z_1^{-|n|}z_2^{-1}z_1^{|n|},\,e)\leq Cx^{\frac{8}{7}}\quad\text{and}\quad\textbf{\emph{d}}(z_1^{-|n|}z_2z_1^{|n|}z_2^{-1},\,e)\leq Cx^{\frac{8}{7}}.
\end{align*}
Hence we finish the proof.

\smallskip
\noindent\emph{Step 7: We show that }
\begin{align}\label{for:371}
  \mathcal{E}_2\leq Cx^{\frac{8}{7}}\|\theta_\mu\|_{\kappa+1}\norm{\psi}_{\kappa}.
\end{align}
We have
\begin{align*}
\mathcal{E}_2&\overset{\text{(1)}}{\leq}\sum_{m\in L_{\mu}}\big\|\langle Y_{1,m}, \,\pi_{u}(f_1\circ \tau^{-1})\psi\rangle\big\|\\
&+\sum_{m\in L_{\mu}}\big\|\langle Y_{2,m},\, \pi_{u}(f_1\circ \tau^{-1})\psi\rangle\big\|\\
&\overset{\text{(2)}}{\leq} \sum_{n=-l}^l Ce^{-\gamma_3(|n|+|m|)}\norm{\theta_\mu\circ s_n-\theta_\mu}_{\kappa}\norm{\psi}_{\kappa}\notag\\
 &+\sum_{n=-l}^l Ce^{-\gamma_3(|n|+|m|)}\big\|\big((dk_n|_{J_\mu})-I\big)\theta_\mu\big\|_{\kappa}\norm{\psi}_{\kappa}\notag\\
 &\overset{\text{(3)}}{\leq}   \sum_{n=-l}^l C_1e^{-\gamma_3(|n|+|m|)}\textbf{\emph{d}}(s_n,e)\norm{\theta_\mu}_{\kappa+1}\norm{\psi}_{\kappa}\notag\\
 &+\sum_{n=-l}^l C_1e^{-\gamma_3(|n|+|m|)}\textbf{\emph{d}}(k_n,e)\|\theta_\mu\|_{\kappa}\norm{\psi}_{\kappa}\notag\\
 &\overset{\text{(4)}}{\leq} C_2x^{\frac{8}{7}}\|\theta_\mu\|_{\kappa+1}\norm{\psi}_{\kappa}.
\end{align*}
Here in $(1)$ we use \eqref{for:368}; in $(2)$ we use Lemma \ref{le:2}, \eqref{for:361} and \eqref{for:362};
in $(3)$ we note that
\begin{align*}
 \norm{\theta_\mu\circ s_n-\theta_\mu}_{\kappa}\leq C\norm{\text{Ad}_{s_n}-I}\|\theta_\mu\big\|_{\kappa+1}\overset{\text{(a)}}{\leq} C_1\textbf{\emph{d}}(s_n,e)\|\theta_\mu\|_{\kappa+1}
\end{align*}
and
\begin{align*}
\norm{dk_n|_{J_\mu}-I}\overset{\text{(b)}}{=}\norm{\text{Ad}_{k_n}|_{J_\mu}-I}\overset{\text{(a)}}{\leq}  C\textbf{\emph{d}}(k_n,e).
\end{align*}
Here in $(a)$ we use \eqref{for:364} of \eqref{sect:9} of Section \ref{sec:19} and \eqref{for:366}; in $(b)$ we note that
$k_n\in V^o$; in $(4)$ we use \eqref{for:366}. Then we finish the proof.

\medskip

\emph{Conclusion}: It follows from \eqref{for:367}, \eqref{for:357}, \eqref{for:369}, \eqref{for:370} and \eqref{for:371} that
\begin{align*}
 &\big\|\langle\mathfrak{D}_{\mu,\lambda,u,\tau,\ZZ},\,\psi\rangle\big\|\\
 &\leq Ce^{-\gamma_2l}\norm{\theta_\mu}_\kappa\norm{\psi}_{\kappa}+C(e^{-\gamma_3l}\norm{\varphi_\mu}_\kappa+\norm{\omega_\mu}_\kappa)\norm{\psi}_{\kappa}\\
 &+Cx^{\frac{8}{7}}\|\theta_\mu\|_{\kappa+1}\norm{\psi}_{\kappa}\\
 &\overset{\text{(1)}}{\leq} C_1(x^{\gamma_2}\norm{\theta_\mu}_\kappa+x^{\gamma_3}\norm{\varphi_\mu}_\kappa+x^{\frac{8}{7}}\|\theta_\mu\|_{\kappa+1})\norm{\psi}_{\kappa}.
\end{align*}
Here in $(1)$ we use \eqref{for:363}.

We recall the definition of $c_0$ in \eqref{for:304} of Section \ref{sec:40}. Hence we finish the proof.

\section{Proof of \eqref{for:18} of Section \ref{sec:3}}\label{sec:47}
If $\GG=G$, then $y_{c'}^j=k_{c'}^jn_{c'}^j$. Then we have
\begin{align*}
  \norm{\text{Ad}_{(y_{c'})^j}}&\overset{\text{(1)}}{\leq} C\norm{\text{Ad}_{(n_{c'})^j}}\overset{\text{(2)}}{\leq} C_{c} (\abs{j}+1)^{\dim \mathfrak{g}}
\end{align*}
for any $j\in \ZZ$. Here in $(1)$ we use the fact that $k_{c'}$ is  compact; in $(2)$ we use Engel's theorem and the fact that $n_{c'}$ is close to $n_{c}$ if $\delta(c)$ is small.

If $\GG=G\ltimes_\rho V$, \eqref{for:12} of Section \ref{sec:11} shows that
\begin{align*}
 y_{c'}^j=\Big(k_{g_{c'}}^j,\big(n_{g_{c'}},\rho(k_{g_{c'}})^{-(j-1)}v_{c'}\big)\big(n_{g_{c'}},\rho(k_{g_{c'}})^{-(j-2)}v_{c'}\big)\cdots (n_{g_{c'}},v_{c'})\Big)
\end{align*}
if $j\geq0$.

Since $k_{g_{c'}}$ is compact, we see that
\begin{align*}
  \max_{0\leq i\leq j-1}\{\big\|(n_{c'},\,\rho(k_{g_{c'}})^{-i}v_{c'})\big\|,\,\big\|(n_{c'},\,\rho(k_{g_{c'}})^{-i}v_{c'})^{-1}\big\|\}\leq C_c
\end{align*}
if $c'$ is close to $c$.

We also note that $n_{g_{c'}}$ generate a $\ZZ$ subgroup of $G$. It is clear that the subgroup $\ZZ\ltimes V\subseteq\GG$ is nilpotent.
Hence, it follows from Engel's theorem that
\begin{align*}
  \norm{\text{Ad}_{(y_{c'})^j}}&\leq C_{c} (\abs{j}+1)^{\dim V+1}\leq C_{c} (\abs{j}+1)^{\dim\mathfrak{G}}
\end{align*}
for any $j\in \ZZ$. Hence, we finish the proof.

\section{Proof of Proposition \ref{po:3} }\label{sec:9}

\eqref{for:393}: Proposition 4.2 of \cite{vw} shows that if $\GG=G$, then \eqref{for:393} holds. If $\GG=G_\rho$, \eqref{for:12} of Section \ref{sec:11} shows that $p$ only depends on the semisimple part of $A$. Thus, \eqref{for:393} holds as well.

\smallskip


\eqref{for:394}: If $\GG=G$, let $u$ be trivial. The conclusion follows from \eqref{for:393}.

  If $\GG=G_\rho$ and $V$ is abelian,  the conclusion follows from Proposition 4.4 of \cite{vw}. So we just need to show that \eqref{for:394} holds when $V$ is non-abelian.

  We prove by induction. Suppose the argument holds for any simply connected $V$ having a central series of length $1$, i.e., $V$ is abelian. Then the result follows from previous discussion.
Assume the result holds for any simply connected $V$ of length $k$. Next, we consider simply connected $V$ of length $k+1$.

\smallskip
\noindent\emph{Step 1: We show that: }there is $u\in V$ such that for any $c=(g_c, v_c)\in A$, $ucu^{-1}=\big(g_c, \,\exp(n_{c,0})\exp(n_{c,1})\big)$, where
\begin{align*}
 n_{c,0}\in \mathfrak{g}_{1,A}\quad\text{ and }\quad n_{c,1}\in \text{Lie}(Z(V))\bigcap \sum_{1\neq\mu\in \Delta_A}\mathfrak{g}_{\mu,A}.
\end{align*}
We note that $V/Z(V)$ is of length $k$. Since $V$ is simply connected, the exponential map is a bijection, which gives that $V/Z(V)$ is also simply connected.  This allows us to use the inductive assumption on $G\ltimes (V/Z(V))$:
 there is $u\in V$ such that for any $c=(g_c, v_c)\in A$, $ucu^{-1}=(g_c, v_{c,1})$, where $\mathfrak{t}(v_{c,1})$ commutes with any $p(a)$, $a\in A$. Here $\mathfrak{t}$ denotes the projection from $V$ to $V/Z(V)$. This means that
 \begin{align*}
 (p(a),0)(0,v_{c,1})(p(a)^{-1},0)(0,v_{c,1}^{-1})\in Z(V),\qquad \forall\,a\in A.
 \end{align*}
It follows from \eqref{for:12} of Section \ref{sec:11} that
 \begin{gather}
 \Rightarrow  \rho(p(a))(v_{c,1})v_{c,1}^{-1}=z_a, \quad\forall\,a\in A\text{ and }z_a\in Z(V);\notag\\
  \Rightarrow  d\rho(p(a))(\log v_{c,1})=\log v_{c,1}+\log z_a, \quad\forall\,a\in A\text{ and }z_a\in Z(V),\label{for:396}
\end{gather}
where $d\rho$ is the induced representation on $\text{Lie}(V)$.

Choose a regular $a\in A$ (see \eqref{for:62}). Since  $d\rho(p(a))-I$ is invertible on each $\mathfrak{g}_{\mu,A}$, $\mu\neq 1$, it follows from \eqref{for:396} that
\begin{align*}
 \log v_{c,1}=n_{c,0}+n_{c,1}
\end{align*}
where $n_{c,0}\in \mathfrak{g}_{1,A}$ and $n_{c,1}\in \text{Lie}(Z(V))\bigcap \sum_{1\neq\mu\in \Delta_A}\mathfrak{g}_{\mu,A}$.

Hence we have
\begin{align*}
 v_{c,1}=\exp(n_{c,0})\exp(n_{c,1})
\end{align*}
Then we finish the proof.


\smallskip
\noindent\emph{Step 2: We show that: }$A'$ is an abelian subgroup, where
\begin{align*}
 A'=\{(g_c, \exp(n_{c,1})): c\in A\}.
\end{align*}
 For any $c,\,a\in A$, we have
\begin{gather*}
 (g_c, v_{c,1})(g_a, v_{a,1})=(g_a, v_{a,1})(g_c, v_{c,1})\\
 \Leftrightarrow \rho(g_a^{-1})(v_{c,1})v_{a,1}=\rho(g_c^{-1})(v_{a,1})v_{c,1}\\
 \overset{\text{(1)}}{\Leftrightarrow} \rho(g_a^{-1})\big(\exp(n_{c,0}) \big)\exp(n_{a,0})\rho(g_a^{-1})\big( \exp(n_{c,1})\big)\exp(n_{a,1})\\
 =\rho(g_c^{-1})\big(\exp(n_{a,0}) \big)\exp(n_{c,0})\rho(g_c^{-1})\big( \exp(n_{a,1})\big)\exp(n_{c,1})\\
 \overset{\text{(2)}}{\Leftrightarrow} \exp(-n_{c,0})\rho(g_c^{-1})\big(\exp(-n_{a,0}) \big)\rho(g_a^{-1})\big(\exp(n_{c,0}) \big)\exp(n_{a,0})\\
 =\rho(g_c^{-1})\big( \exp(n_{a,1})\big)\exp(n_{c,1})\exp(-n_{a,1})\rho(g_a^{-1})\big( \exp(-n_{c,1})\big)\\
\overset{\text{(3)}}{\Leftrightarrow}\exp(-n_{c,0})\rho(g_c^{-1})\big(\exp(-n_{a,0}) \big)\rho(g_a^{-1})\big(\exp(n_{c,0}) \big)\exp(n_{a,0})\\
 =\exp\big( d\rho(g_c^{-1})n_{a,1}+n_{c,1}-n_{a,1}-d\rho(g_a^{-1})n_{c,1}\big).
\end{gather*}
Here in $(1)$ we note that $\exp(n_{c,1})\in Z(V)$ and $\exp(n_{a,1})\in Z(V)$ and $Z(V)$ in invariant under $\rho(G)$; in $(2)$ we note that
\begin{align*}
 \big(\rho(g)\exp(v) \big)^{-1}=\big(\exp(d\rho(g)v)\big)^{-1}=\exp(-d\rho(g)v)=\rho(g)\exp(-v)
\end{align*}
for any $g\in G$ and $v\in \text{Lie}(V)$;

 in $(3)$ using \emph{Step 1}, we see that the left side is inside the group $\exp(\mathfrak{g}_{1,A})$ (see \eqref{for:52} of Section \ref{sec:3}) and the right side is inside $\exp\big(\text{Lie}(Z(V))\bigcap \sum_{1\neq\mu\in \Delta_A}\mathfrak{g}_{\mu,A}\big)$. Hence we see that
\begin{align*}
d\rho(g_c^{-1})n_{a,1}+n_{c,1}-n_{a,1}-d\rho(g_a^{-1})n_{c,1}=0.
\end{align*}
This implies that $A'$ is abelian and inside $G\ltimes Z(V)$.

\smallskip
\noindent\emph{Step 3: We show that: }there is $u_1\in Z(V)\cap V$ such that
\begin{align*}
  u_1(g_c, \exp(n_{c,1}))u_1^{-1}=(g_c,z_c),\qquad \forall\,\,c\in A
\end{align*}
where $z_c\in Z(p(A))\cap Z(V)$.

As $Z(V)$ is simply connected and abelian, the claim follows from the induction assumption, i.e., Proposition 4.4 of \cite{vw}.

\smallskip
\noindent\emph{Conclusion: } for any $c\in A$
\begin{align*}
 u_1ucu^{-1}u_1^{-1}&\overset{\text{(1)}}{=}u_1\big(g_c, \,\exp(n_{c,0})\exp(n_{c,1})\big)u_1^{-1}\\
 &\overset{\text{(2)}}{=}u_1(g_c, \exp(n_{c,1})) (e, \exp(n_{c,0}))u_1^{-1}\\
 &\overset{\text{(3)}}{=}u_1(g_c, \exp(n_{c,1}))u_1^{-1} (e, \exp(n_{c,0}))\\
 &\overset{\text{(4)}}{=}(g_c,z_c)(e, \exp(n_{c,0}))\\
 &=(g_c, \,z_c\exp(n_{c,0})).
\end{align*}
Here in $(1)$ we use \emph{Step 1}, $(2)$ we note that $\exp(n_{c,1})\in Z(V)$ and we use \eqref{for:12} of Section \ref{sec:11}; in $(3)$ we recall
 $u_1\in Z(V)\cap V$; in $(4)$ we use \emph{Step 3}.

Since $\exp(n_{c,0})z_c\in Z(p(A))$ (see \emph{Step 1} and \emph{Step 3}), we proved the case of $k+1$. Hence we finish the proof.

\section{Proof of \eqref{sec:48} of Section \ref{sec:3}}\label{ob:5}

We recall notations in \eqref{sect:11} of Section \ref{sec:3}. Any
$c\in\GG$ generates a $\ZZ$ subgroup. Applying  Proposition \ref{po:3} to this $\ZZ$ subgroup.
Then there is $u\in V$ such that $ucu^{-1}\in Z(s_c)$. \eqref{for:12} of Section \ref{sec:11} shows that $g_{ucu^{-1}}=g_c$, which implies that
\begin{align*}
 y_{ucu^{-1}}s_{ucu^{-1}}&=y_{ucu^{-1}}s_{c}=(s_{ucu^{-1}}^{-1}ucu^{-1})s_{c}=(s_{c}^{-1}ucu^{-1})s_{c}\\
 &=s_{c}^{-1}(s_{c} ucu^{-1})=ucu^{-1}=s_{ucu^{-1}}y_{ucu^{-1}}=s_{c}y_{ucu^{-1}}.
\end{align*}
Hence, we have
\begin{align*}
 uc^ju^{-1}=s_{c}^jy_{ucu^{-1}}^j,\qquad \forall\,j\in\ZZ.
\end{align*}
This and \eqref{for:18} of Section \ref{sec:3} justifies the decomposition in \eqref{sec:48}.

\section{Proof of \eqref{sec:49} of Section \ref{sec:3}}\label{cor:9}
At first, we prove the following result, which is crucial for the proof of \eqref{sec:49}:
\begin{corollary}\label{cor:11} For any $c\in\GG$, if there is $u(c)\in V$ satisfying $u(c)cu(c)^{-1}\in Z(s_c)$, then for any $\epsilon>0$ there is $\delta(c,\epsilon)>0$ and  such that:   and for any $c'\in B_{\delta}(c)$ there is $u(c')\in V\bigcap B_{\epsilon}(u(c))$ such that $u(c')c'u(c')^{-1}\in Z(s_{c'})$.

\end{corollary}
\begin{proof} By the argument at the beginning of the proof of Proposition \ref{po:3}, it suffices to assume $\GG=G_\rho$. If $V$ is abelian,  the conclusion follows from the proof of Proposition 4.4 of \cite{vw}. We proceed by induction. Suppose the argument holds for any simply connected $V$ having a central series of length $k$. Next, we consider simply connected $V$ of length $k+1$.

\smallskip
\noindent\emph{Step 1: We show that: }for any $c\in\GG$, if there is $u(c)\in V$ satisfying $u(c)cu(c)^{-1}\in Z(s_c)$, then for any $\epsilon>0$ there is $\delta_1(c,\epsilon)>0$ such that for any $c'\in B_{\delta_1}(c)$ there is $q(c')\in V\bigcap B_{\epsilon}(u(c))$ such that
\begin{align*}
 qc'q^{-1}=\big(g_{c'}, \,\exp(n_{c',0})\exp(n_{c',1})\big)
\end{align*}
where
\begin{align*}
 n_{c',0}\in \mathfrak{g}^0_{c'}=Z(s_{c'})\quad\text{and}\quad n_{c',1}\in (\mathfrak{g}^-_{c'}\oplus\mathfrak{g}^+_{c'})\cap\text{Lie}(Z(V))
\end{align*}
(see \eqref{for:273} of Section \ref{sec:3}).

Using the inductive assumption on $G\ltimes (V/Z(V))$, we have: for any $\epsilon>0$ there there is $\delta(c,\epsilon)>0$ such that for any $c'\in B_{\delta}(c)$ there is $q(c')\in V\bigcap B_{\epsilon}(u(c))$ such that $qc'q^{-1}=(g_{c'}, v_{c',1})$, where $\mathfrak{t}(v_{c',1})$ commutes with
$s_{c'}$. Here $\mathfrak{t}$ denotes the projection from $V$ to $V/Z(V)$. Thus, we have
\begin{align*}
(s_{c'},0)(0,v_{c',1})(s_{c'}^{-1},0)(0,v_{c',1}^{-1})\in Z(V).
\end{align*}
It follows from \eqref{for:12} of Section \ref{sec:11} that
\begin{gather*}
 \Rightarrow  d\rho(s_{c'})(\log v_{c',1})=\log v_{c',1}+\log z_{c'}, \qquad\text{where }z_{c'}\in Z(V)\\
 \overset{\text{(1)}}{\Rightarrow}  v_{c',1}=\exp(n_{c',0})\exp(n_{c',1})
\end{gather*}
where $n_{c',0}\in \mathfrak{g}^0_{c'}$ and $n_{c',1}\in (\mathfrak{g}^-_{c'}\oplus\mathfrak{g}^+_{c'})\cap\text{Lie}(Z(V))$. Here in $(1)$ we note that $d\rho(s_{c'})-I$ is invertible on $\mathfrak{g}^-_{c'}\oplus\mathfrak{g}^+_{c'}$.

Then we finish the proof.

\smallskip

\noindent\emph{Step 2: we show that: }if $\delta(c,\epsilon)$ is sufficiently small, then there is $q_1\in Z(V)\cap B_{\epsilon}(e)$ such that
\begin{align}\label{for:397}
 q_1(g_{c'},\exp(n_{c',1}))q_1^{-1}\in Z(s_{c'}).
\end{align}
Since $ucu^{-1}$ is perfect by assumption, $n_{c,1}=0$. If $\delta$ is sufficiently small, then $n_{c',1}$ is small. Thus, $\exp(n_{c',1})$ is sufficiently close to $e$.
Since $(g_{c},0)\in G\ltimes Z(V)$ is perfect and $(g_{c'},\exp(n_{c',1}))\in G\ltimes Z(V)$ is sufficiently close to $(g_{c},0)$,  \eqref{for:397} follows from the
the inductive assumption on $G\ltimes Z(V)$.

\smallskip

\noindent\emph{Step 3: we show that: }if we let $u(c')=q_1q$, then
\begin{align*}
u(c')c'(u(c'))^{-1}\in Z(s_{c'}).
\end{align*}
 We have
  \begin{align*}
 q_1qc'q^{-1}q_1^{-1}&\overset{\text{(1)}}{=}q_1\big(g_{c'}, \,\exp(n_{c',0})\exp(n_{c',1})\big)q_1^{-1}\notag\\
 &\overset{\text{(2)}}{=}q_1(g_{c'}, \exp(n_{c',1})) (e, \exp(n_{c',0}))q_1^{-1}\notag\\
 &\overset{\text{(3)}}{=}q_1(g_{c'}, \exp(n_{c',1}))q_1^{-1} (e, \exp(n_{c',0})).
 \end{align*}
 Here in $(1)$ we use \emph{Step 1}, $(2)$ we note that $\exp(n_{c',1})\in Z(V)$ and we use \eqref{for:12} of Section \ref{sec:11}; in $(3)$ we recall
 $q\in Z(V)$.

\eqref{for:397} and \emph{Step 1} show that
\begin{align*}
 q(g_{c'}, \exp(n_{c',1}))q^{-1} (e, \exp(n_{c',0}))\in Z(s_{c'}).
\end{align*}
Then we finish the proof.

\smallskip

\noindent\emph{Step 4: we show that: }
\begin{align*}
\textbf{d}\big(u(c'), \,u(c)\big)<2\epsilon.
\end{align*}
We have
\begin{align*}
  &\textbf{d}\big(u(c'), \,u(c)\big)=\textbf{d}\big(q_1q, \,u(c)\big)<\textbf{d}(q_1q, q)+\textbf{d}\big(q, u(c)\big)\\
  &\overset{\text{(1)}}{=}\textbf{d}(q_1, e)+\textbf{d}\big(q, u(c)\big)<\epsilon+\epsilon.
\end{align*}
Here in $(1)$ we use the right invariance of $\textbf{d}$ (see \eqref{for:35} of Section \ref{sec:3}).

Hence we proved the case of $k+1$. Thus we finish the proof.

\end{proof}
Now we proceed to the proof of \eqref{sec:49} of Section \ref{sec:3}. It follows directly from \eqref{for:18} of Section \ref{sec:3} and Corollary \ref{cor:11}.

\section{Proof of Proposition \ref{po:9}}\label{sec:10}
\subsection{Preparative notations and results}
For any regular $a_1,\,a_2\in A$ (see \eqref{for:62} of Section \ref{sec:3}) we set
\begin{align}
 l_1(a_1,a_2)&=\min_{\lambda\neq\mu\in \Delta_A}\{\Big|\log\frac{\mu(p(a_1))}{\lambda(p(a_1))}\Big|,\,\Big|\log\frac{\mu(p(a_2))}{\lambda(p(a_2))}\Big|\}\label{for:110}\\
  l_2(a_1,a_2)&=\min_{\mu\neq\lambda\in \Delta_A\backslash\{1\}}\{\frac{|\log\mu(p(a_1))|}{|\log\lambda(p(a_1))|},\,\frac{|\log\mu(p(a_2))|}{|\log\lambda(p(a_2))|}\}\label{for:24}\\
  l_3(a_1,a_2)&=\min_{\substack{\abs{r_1}+\abs{r_2}=1\\(r_1,r_2)\in\RR^2}}\sum_{\mu\in\Delta_A}
  \abs{r_1\log\mu(p(a_1))+r_2\log\mu(p(a_2))}\label{for:58}\\
  \emph{\textbf{l}}(a_1,a_2)&=\min\{\text{\tiny$\frac{1}{5}$}\gamma l_1(a_1,a_2),\,\text{\tiny$\frac{1}{5}$} \gamma l_3(a_1,a_2),\,\text{\tiny$\frac{1}{5}$}\gamma l_1(a_1,a_2)l_2(a_1,a_2)\}\label{for:63}\\
   \emph{\textbf{l}}_1(a_1,a_2)&=\emph{\textbf{l}}(a_1,a_2)^2<\emph{\textbf{l}}(a_1,a_2)\label{for:446}
\end{align}
($0<\gamma<1$ is defined in Theorem \ref{th:3} corresponding to the split Cartan subgroup $A_0$ described in Proposition \ref{po:3}.).

\begin{definition}\label{de:1}
For any  elements $a_1,\,a_2\in A$, we say that $a_1$ and $a_2$ are \emph{genuinely  linearly independent} if $p(a_1)$ and $p(a_2)$ linearly independent.

We say that $(a_1,a_2)$ is a \emph{good} pair if the following are satisfied:
\begin{enumerate}
  \item\label{for:398} $a_1$ and $a_2$ are regular, genuinely  linearly independent;


  \item\label{for:408} $\mu(a_i)\neq \lambda(a_i)$, $i=1,2$,  for any $\mu,\,\lambda\in \Delta_A$ with $\mu\neq \lambda$ (see \eqref{for:62} of Section \ref{sec:3});

  \smallskip

  \item \label{for:56} $e^{\emph{\textbf{l}}_1(a_1,a_2)}>\max_{i=1,2}\,\norm{\text{Ad}_{y_{a_i}^\pm}}$ (see \eqref{sect:11} of Section \ref{sec:3});

      \smallskip
  \item\label{for:399} $y_{a_i}\in Z(p(A))$, $i=1,2$.
\end{enumerate}

\end{definition}

\begin{lemma}\label{re:1}
We can always find a good pair in $A$ if $\alpha_A$ is higher rank.
\end{lemma}
\begin{proof} We choose regular and genuinely  linearly independent $a_1$ and $a_2$.  Set
\begin{align*}
 a_{i,j}=a_i^j,\quad i=1,\,2,\quad j\in \ZZ.
\end{align*}
Hence \eqref{for:398} and \eqref{for:408} of Definition \ref{de:1}  are also satisfied for any of the pairs $(a_{1,j},\,a_{2,j})$.

Since $A$ satisfies property  \hyperlink{o.1}{(P)}, $y_{a_i}\in Z(p(A))$, $i=1,2$. Then \eqref{for:399} of Definition \ref{de:1}  is satisfied for any of the pairs $(a_{1,j},\,a_{2,j})$.

We note that
\begin{align*}
\max_{i=1,2}\,\norm{\text{Ad}_{(y_{a_{i,j}})^\pm}}=\max_{i=1,2}\,\norm{\text{Ad}_{(y_{a_i}^\pm)^j}},\qquad \forall\,j\in\ZZ.
\end{align*}
The above equation together with \eqref{for:18} of Section \ref{sec:3} show that $\max_{i=1,2}\,\norm{\text{Ad}_{(y_{a_{i,m}})^\pm}}$ increases polynomially with respect to $m$, in contrast with the
exponential increase of  $e^{\emph{\textbf{l}}(a_1^m,a_2^m)}$ with respect to $m$.  Then we can always choose $m\in\NN$ big enough such that \eqref{for:56} of Definition \ref{de:1} is satisfied if we replace $a_i$ by $a_i^m$, $i=1,\,2$. Hence,  $(a_1^m,a_2^m)$ is a good pair.
\end{proof}

\begin{lemma}\label{ob:3} For any good pair $(a_1,a_2)$ of $A$, there is $\eta(a_1,a_2)>0$ such that for any $w_i\in B_\eta(e)\cap Z(p(A))$, letting $z_i=w_ia_i$,  $i=1,2$, we have: for any $(m,n)\in\ZZ^2$
  \begin{align}\label{for:401}
 z_1^nz_2^m= (w_1y_{a_1})^n(w_2y_{a_2})^mp(a_1)^np(a_2)^m
\end{align}
and
\begin{align}\label{for:400}
  \norm{\text{Ad}&_{(w_1y_{a_1})^n(w_2y_{a_2})^m}}\leq e^{\textbf{l}_1(a_1,a_2)(|n|+|m|)}.
\end{align}

\end{lemma}
\begin{proof}
For any good pair $(a_1,a_2)$ of $A$, there is $\eta(a_1,a_2)>0$ such that for any $w\in B_\eta(e)$ we have
\begin{align}\label{for:103}
 \norm{\text{Ad}_{w^\pm}}\leq e^{\emph{\textbf{l}}_1(a_1,a_2)}.
\end{align}
For any $w_i\in Z(p(A))$, let $z_i=w_ia_i$,  $i=1,2$. By \eqref{for:399} of Definition \ref{de:1}  we have
\begin{align*}
 z_1^nz_2^m= (w_1y_{a_1})^n(w_2y_{a_2})^mp(a_1)^np(a_2)^m.
\end{align*}
Moreover, if $w_i\in B_\eta(e)\cap Z(p(A))$, $i=1,2$, we have
\begin{align*}
  \norm{\text{Ad}&_{(w_1y_{a_1})^n(w_2y_{a_2})^m}}\leq (\norm{\text{Ad}_{y_{a_1}^\pm}}\norm{\text{Ad}_{w_{1}^\pm}})^{|n|}(\norm{\text{Ad}_{y_{a_2}^\pm}}\norm{\text{Ad}_{w_{2}^\pm}})^{\abs{m}}\notag\\
  &\overset{\text{(1)}}{\leq} (\max\{\norm{\text{Ad}_{y_{a_1}^\pm}},\,\norm{\text{Ad}_{y_{a_2}^\pm}},\,e^{\emph{\textbf{l}}_1(a_1,a_2)}\})^{|n|+|m|}\notag\\
  &\overset{\text{(2)}}{\leq} e^{\emph{\textbf{l}}_1(a_1,a_2)(|n|+|m|)}
\end{align*}
for any $(m,n)\in\ZZ^2$. Here in $(1)$ we use \eqref{for:103}; in $(2)$ we use \eqref{for:56} of Definition \ref{de:1} .

\end{proof}

By Proposition \ref{po:3} there is a split Cartan subgroup $A_0$ such that $p(A)\subseteq A_0$. Let $\gamma$ and $E$ be as defined in Theorem \ref{th:3} for $A_0$. Then we have:
\begin{proposition}\label{po:4}
For any good pair $(a_1,a_2)$ of $A$, any $w_i\in B_{\eta(a_1,a_2)}(e)\cap Z(p(A))$ (see Lemma \ref{ob:3}) and any $\xi_i\in\mathcal{O}^1$, $i=1,2$, letting $z_i=w_ia_i$, $i=1,2$ we have: for any $(m,n)\in\ZZ^2$
\begin{align*}
\abs{\langle \xi_1\circ (z_1^nz_2^m),\xi_2\rangle}&\leq
E\norm{\xi_1}_{1}\norm{\xi_2}_{1}e^{-4
\textbf{l}(a_1,a_2)(|n|+|m|)}
\end{align*}


\end{proposition}
\begin{proof}
Let $\theta_{n,m}=\xi_1\circ\big((w_1y_{a_1})^n(w_2y_{a_2})^m\big)$. Then we have
\begin{align}\label{for:104}
 \norm{\theta_{n,m}}_1\leq C\norm{\text{Ad}_{(w_1y_{a_1})^n(w_2y_{a_2})^m}}\norm{\xi_1}_1\overset{\text{(1)}}{\leq} Ce^{\emph{\textbf{l}}_1(a_1,a_2)(|n|+|m|)}\norm{\xi_1}_1.
\end{align}
Here in $(1)$ we use \eqref{for:400} of Lemma \ref{ob:3}.

Then we have
\begin{align}
\abs{\langle \xi_1&\circ(z_1^nz_2^m),\xi_2\rangle}\overset{\text{(1)}}{\leq}\abs{\langle \theta_{n,m}\circ\big(p(a_1)^np(a_2)^m\big),\xi_2\rangle}\notag\\
&\overset{\text{(2)}}{\leq}
E\norm{\theta_{n,m}}_{1}\norm{\xi_2}_{1}e^{-\gamma
l_3(a_1,a_2)(|n|+|m|)}\notag\\
&\overset{\text{(3)}}{\leq} CE\norm{\xi_1}_{1}\norm{\xi_2}_{1}e^{\emph{\textbf{l}}_1(a_1,a_2)(|n|+|m|)}e^{-\gamma
l_3(a_1,a_2)(|n|+|m|)}\notag\\
&\overset{\text{(4)}}{\leq} CE\norm{\xi_1}_{1}\norm{\xi_2}_{1}e^{\emph{\textbf{l}}(a_1,a_2)(|n|+|m|)}e^{-\gamma
l_3(a_1,a_2)(|n|+|m|)}\notag\\
&\overset{\text{(5)}}{\leq} CE\norm{\xi_1}_{1}\norm{\xi_2}_{1}e^{-4
\emph{\textbf{l}}(a_1,a_2)(|n|+|m|)}\notag
\end{align}
Here in $(1)$ we use \eqref{for:401} of Lemma \ref{ob:3}; $(2)$ we use Theorem \ref{th:3} and \eqref{for:58}; in $(3)$ we use \eqref{for:104}; in $(4)$ we use \eqref{for:446}; in $(5)$ we use \eqref{for:63}. Hence we finish the proof.

\end{proof}

\subsection{Proof of Proposition \ref{po:9}} Lemma \ref{re:1} shows that we can choose $\textbf{a}_1$ and $\textbf{a}_2$ of $A$ such that
$(\textbf{a}_1,\textbf{a}_2)$ is a good pair.

\eqref{for:402}: Follows from \eqref{for:398} of Definition \ref{de:1}.

\smallskip

\eqref{for:409}: Follows from \eqref{for:110} and \eqref{for:63}.

\smallskip

\eqref{for:50} holds clearly.

\smallskip

\eqref{for:403}: Follows from Proposition \ref{po:4}.

\smallskip

\eqref{for:404}: \eqref{for:490} of Section \ref{sec:3} imply that $\mathfrak{g}_{\mu,A}$ are invariant subspaces for $\text{Ad}_{z_i}$, $i=1,2$. Moreover,
it follows from Lemma \ref{ob:3} that
\begin{align*}
  \norm{\text{Ad}_{z_i^n}|_{\mathfrak{g}_{\mu,A}}}&\leq  \norm{\text{Ad}_{p(\textbf{a}_i)^n}|_{\mathfrak{g}_{\mu,A}}}\cdot\norm{\text{Ad}_{(w_iy_{\textbf{a}_i})^n}}\\
  &\leq \mu(\textbf{a}_i)^ne^{\emph{\textbf{l}}_1|n|}\leq \mu(\textbf{a}_i)^ne^{\emph{\textbf{l}}|n|}
\end{align*}
for any $n\in\ZZ,\,\,i=1,2$.

\smallskip
\eqref{for:472}: It follows from \eqref{for:404}, \eqref{for:110} and \eqref{for:63}: if $\mu(\textbf{a}_i)>1$ (resp. $\mu(\textbf{a}_i)<1$), then for any $n<0$ (resp. $n>0$), we have
\begin{align*}
  \norm{\text{Ad}_{(w_i\textbf{a}_i)^n}|_{\mathfrak{g}_{\mu,A}}}< \mu(\textbf{a}_i)^ne^{\emph{\textbf{l}}(\textbf{a}_1,\,\textbf{a}_2)|n|}<\mu(\textbf{a}_i)^{\frac{4}{5}n}.
\end{align*}
This implies the result.


\smallskip

\eqref{for:407}: see \eqref{for:54} of Section \ref{sec:3}.

\smallskip

Hence, we finish the proof.

\section{Proof of Corollary \ref{for:129}}\label{sec:51}
\eqref{for:78}: By Corollary \ref{cor:11}, noting that $\textbf{a}_1$ is perfect, there is $q_1\in V$ sufficiently close to $e$ such that
$z_1'=q_1z_1q_1^{-1}$ is perfect. We recall the fact that
the split Cartan subgroup is unique up to inner automorphisms. This implies that there is $q_2\in G$ sufficiently close to $e$ such that
$q_2s_{z_1'}q_2^{-1}\in A_0$. Let $q=q_2q_1$ and $\mathfrak{z}_i=qz_iq^{-1}$, $i=1,\,2$. Then it is clear that $\emph{\textbf{d}}(q,e)+\emph{\textbf{d}}(q^{-1},e)<1$ and
$\mathfrak{z}_1$ is perfect and
  $p(\mathfrak{z}_1)\in A_0$.

 \smallskip

\eqref{for:57} is clear from \eqref{for:18} of Section \ref{sec:3}.

\smallskip
\eqref{for:8} is a direct consequence of \eqref{for:57}.

\smallskip

\eqref{for:435}: Since $\mathfrak{z}_1$ is perfect and
  $p(\mathfrak{z}_1)\in A_0$, $\mathfrak{g}_{\mu,A}$ are invariant subspaces for $\text{Ad}_{p(\mathfrak{z}_1)}$ for any $\mu\in\Delta_A$.
  Since $p(\mathfrak{z}_1)$ is sufficiently close to $\textbf{a}_1$, $p(\mathfrak{z}_1)$ also satisfies \eqref{for:408} of Definition \ref{de:1}.
This implies that $\mathfrak{g}_{\mu,A}$ are invariant subspaces for $\text{Ad}_{\mathfrak{z}_1}$ for any $\mu\in\Delta_A$. Moreover, as $\text{Ad}_{\mathfrak{z}_1}|_{\mathfrak{g}_{\mu,A}}$ has  Lyapunov
exponent sufficiently close to that of $\text{Ad}_{\textbf{a}_1}|_{\mathfrak{g}_{\mu,A}}$, $\text{Ad}_{q}$ preserves $\mathfrak{g}_{\mu,A}$ for any $\mu\in\Delta_A$.
We recall that $\mathfrak{g}_{\mu,A}$ are invariant subspaces for $\text{Ad}_{z_2}$ for any $\mu\in\Delta_A$. Thus, $\mathfrak{g}_{\mu,A}$ are also invariant subspaces for $\text{Ad}_{\mathfrak{z}_2}$ for any $\mu\in\Delta_A$. The above discussion shows that if we replace $z_i$ by $\mathfrak{z}_i$, $i=1,\,2$ in \eqref{for:403} and \eqref{for:404} of Proposition \ref{po:9}, the results still hold.

\smallskip

\eqref{for:405}: By \eqref{for:78}, we see that $\mathfrak{z}_1$ is perfect and $p(\mathfrak{z}_1)\in A_0$. This together with \eqref{for:57} imply that
the Lyapunov exponents of $\Ad_{\mathfrak{z}_1}$ are the same as the eigenvalues of $\Ad_{p(\mathfrak{z}_1)}$.
Additionally, the Lyapunov space decomposition of $\Ad_{\mathfrak{z}_1}$  matches the eigenspace decomposition of $\Ad_{p(\mathfrak{z}_1)}$.
If $\eta$ is sufficiently small, then $\Ad_{p(\mathfrak{z}_1)}$ is sufficiently close to $\Ad_{\textbf{a}_1}$. Thus, we have
   \begin{gather*}
 \mathfrak{g}_{\textbf{a}_1}^{+}\subseteq \mathfrak{g}^+_{\mathfrak{z}_1}\subseteq \mathfrak{g}_{\textbf{a}_1}^{+}\cup \mathfrak{g}_{\textbf{a}_1}^{0}\qquad \mathfrak{g}_{\textbf{a}_1}^{-}\subseteq \mathfrak{g}^-_{\mathfrak{z}_1}\subseteq \mathfrak{g}_{\textbf{a}_1}^{-}\cup \mathfrak{g}_{\textbf{a}_1}^{0}.
\end{gather*}
It follows that
\begin{align*}
 [\mathfrak{g}^+_{\mathfrak{z}_1},\,\mathfrak{g}_{\textbf{a}_1}^{+}]&\subseteq [\mathfrak{g}_{\textbf{a}_1}^{+}\cup \mathfrak{g}_{\textbf{a}_1}^{0},\,\mathfrak{g}_{\textbf{a}_1}^{+}]
 \subseteq \mathfrak{g}_{\textbf{a}_1}^{+},\qquad\text{and}\\
 [\mathfrak{g}^-_{\mathfrak{z}_1},\,\mathfrak{g}_{\textbf{a}_1}^{-}]&\subseteq [\mathfrak{g}_{\textbf{a}_1}^{-}\cup \mathfrak{g}_{\textbf{a}_1}^{0},\,\mathfrak{g}_{\textbf{a}_1}^{-}]
 \subseteq \mathfrak{g}_{\textbf{a}_1}^{-}.
\end{align*}
The above discussion implies \eqref{for:405}.

\smallskip

\eqref{for:426}: We recall that $\mathfrak{g}_{1,A}$ are invariant subspaces for $\text{Ad}_{\mathfrak{z}_1}$.
 As $\text{Ad}_{\textbf{a}_1}|_{\mathfrak{g}_{1,A}}$ has Lyapunov
exponent $0$, if $w_1$ is sufficiently close to $e$,   $\text{Ad}_{\mathfrak{z}_1}|_{\mathfrak{g}_{1,A}}$ has sufficiently small Lyapunov
exponents. Thus, we have
\begin{align*}
\big |\chi\big(p(\mathfrak{z}_1^\pm)\big)\big|<e^{\frac{\emph{\textbf{l}}(\textbf{a}_1,\textbf{a}_2)}{m}},
\end{align*}
for any $\chi\in \Delta_{\mathfrak{z}_1}$, if $\mathfrak{g}_{\chi,\mathfrak{z}_1}\subseteq\mathfrak{g}_{1,A}$.

This together with \eqref{for:8} give the result.

\section{Proof of Lemma \ref{ob:1}}\label{sec:52}
Since $\tilde{\alpha}_A$ is $L$-fiber preserving, from \eqref{for:410} we have
\begin{align*}
\exp(\mathfrak{p}_v(l\cdot x))\cdot
z_v\cdot l\cdot x=L_{v,x}(l)\cdot\exp(\mathfrak{p}_v(x))\cdot
z_v\cdot x\qquad
\forall\,l\in L
\end{align*}
where $L_{v,x}$ is map on $L$.

Since $z_v\in Z(L)$ we have
\begin{align*}
\exp(\mathfrak{p}_v(l\cdot x))\cdot l
&=L_{v,x}(l)\cdot\exp(\mathfrak{p}_v(x))\qquad
\forall\,l\in L \\
\Rightarrow \quad l\cdot\exp(\text{Ad}_{l^{-1}}\mathfrak{p}_v(l\cdot x))
&=L_{v,x}(l)\cdot\exp(\mathfrak{p}_v(x))\qquad
\forall\,l\in L.
\end{align*}
We recall that $\mathfrak{p}_v$ takes values in $\mathcal{L}^\bot$ and $\mathcal{L}^\bot$ is $\text{Ad}_L$-invariant  (see \eqref{for:54} of Section \ref{sec:3}).
If $\norm{\mathfrak{p}}_{C^0}$ is small
enough by unique decomposition property near identity (see \eqref{for:52} of Section \ref{sec:3}) it
follows that
\begin{align*}
\text{Ad}_{l^{-1}}\mathfrak{p}_v(l\cdot x)=\mathfrak{p}_v(x)\qquad \forall\,l\in L.
\end{align*}
Then we get the result.

\section{Proof of Lemma \ref{le:8}}\label{sec:4}


To simplify notations, we denote $z_{v}$ and $z_{u}$ by $z_1$ and $z_2$;  and denote $\mathfrak{p}_{v}$ and $\mathfrak{p}_{u}$ by $\mathfrak{p}_1$ and $\mathfrak{p}_2$.

   By commutativity of $\tilde{\alpha}_A$  we have
\begin{align}\label{for:30}
&\exp\bigl(\mathfrak{p}_{1}\bigl(\exp(\mathfrak{p}_{2}(x))\cdot z_2\cdot
x\bigl)\bigl)\cdot
z_1\cdot\exp(\mathfrak{p}_{2}(x))\cdot z_2\cdot x\notag\\
&=l(x)\exp\bigl(\mathfrak{p}_{2}\bigl(\exp(\mathfrak{p}_{1}(x))\cdot z_1\cdot
x\bigl)\bigl)\cdot z_2\cdot\exp(\mathfrak{p}_{1}(x))\cdot z_1\cdot x
\end{align}
for any $x\in \mathcal{M}$, where $l:\mathcal{M}\rightarrow L$. Let
\begin{align*}
 z_i'=\exp(r_i)z_i,\qquad \text{where }r_i=\int_{\mathcal{M}}\mathfrak{p}_i,\quad i=1,\,2.
\end{align*}
\emph{Step 0}: We show that
\begin{enumerate}
  \item\label{for:107} $\norm{l}_{C^{m}}\leq C_m$;

  \smallskip
  \item\label{for:419} $r_i\in \mathcal{L}^\bot,\,\,z_i'\in Z(L),\,\,i=1,2$;

  \smallskip
  \item\label{for:418} Let $[Z(L),Z(L)]$ denote the commutator group  of $Z(L)$. Then
  \begin{align*}
 \text{Lie}([Z(L),Z(L)])\cap \mathcal{L}=\{0\}.
\end{align*}
\end{enumerate}
\eqref{for:107} follows directly from \eqref{for:30} and the fact that $z_i$, $i=1,2$ are uniformly bounded.

\eqref{for:419}: It is clear that $r_i\in \mathcal{L}^\bot$,  as $\mathfrak{p}_i$ takes values in $\mathcal{L}^\bot$, $i=1,2$.  Since $\mathfrak{p}_i$, $i=1,2$ are adjoint $L$-invariant (see Lemma \ref{ob:1}), we have
\begin{align}
 \text{Ad}_{l}r_i=r_i,\qquad \forall\,\,l\in L,\,\,i=1,\,2.
\end{align}
This shows that $\exp(r_i)\in Z(L)$. Then $z_i'\in Z(L)$, $i=1,2$.

\smallskip

\eqref{for:418}: Let $\mathfrak{k}$ be the maximal compact subalgebra of $ \mathfrak{G}$ containing $\mathcal{L}$. Let $\mathfrak{T}$ be the maximal abelian subalgebra of $\mathfrak{k}$ that contains $\mathcal{D}$, a maximal abelian subalgebra of $\mathcal{L}$. For any subalgebra $s$, we use $s_\mathbb{C}$ to denote its complexification.
Then $\mathfrak{T}_\mathbb{C}$ is a Cartan subalgebra of $\mathfrak{k}_\mathbb{C} \cong  \mathfrak{G}_\mathbb{C}$. We consider the root space decomposition of $ \mathfrak{G}_\mathbb{C}$ with respect to $\mathfrak{T}_\mathbb{C}$:
\begin{align*}
 \mathfrak{G}_{\CC}=\big(\bigoplus_{\lambda\in\Delta} \mathcal{R}_\lambda \big) \bigoplus \mathfrak{T}_\mathbb{C}.
\end{align*}
We denote  $\text{Lie}(Z(L))_{\CC}$ by $\mathcal{Z}$. Then $\mathcal{Z}$ is contained in the centralizer of $\mathcal{D}_{\CC}$, i.e.,
\begin{align*}
 \mathcal{Z}\subseteq\big(\bigoplus_{\lambda\in\Delta,\,\lambda|_{\mathcal{D}_{\CC}}=0} \mathcal{R}_\lambda \big)\bigoplus \mathfrak{T}_\mathbb{C}.
\end{align*}
It follows that
\begin{align*}
 [\mathcal{Z},\mathcal{Z}]\subseteq\big(\bigoplus_{\lambda\in\Delta,\,\lambda|_{\mathcal{D}_\mathbb{C}}=0} \mathcal{R}_\lambda \big)\bigoplus \mathcal{S}
\end{align*}
where $\mathcal{S}\subseteq \mathfrak{T}_\mathbb{C}$ with $\mathcal{S}\cap \mathcal{D}_\mathbb{C}=\{0\}$.

Hence we see that $[\mathcal{Z},\mathcal{Z}]\bigcap \mathcal{L}_{\CC}=\{0\}$. Otherwise it contradicts that fact that $\mathcal{D}_\mathbb{C}$ is a maximal abelian
subalgebra of $\mathcal{L}_{\CC}$.

\smallskip

\emph{Step 1}: We show that there are
$c_i:\mathcal{M}\to \mathfrak{G}$, $i=1,2$ with estimates
\begin{align}\label{for:51}
  \norm{c_i}_{C^m}\leq C_m\norm{\mathfrak{p}}_{C^{m+1}}\norm{\mathfrak{p}}_{C^{m}},\qquad i=1,2
\end{align}
such that
\begin{gather*}
\exp\bigl(\mathfrak{p}'_{1}(z_2\cdot
x)+c_1(x)\bigl)\\
=\exp\bigl(\mathfrak{p}_{1}\bigl(\exp(\mathfrak{p}_{2}(x))\cdot z_2\cdot
x\bigl)\bigl)\exp(-r_1),\\
\exp\bigl(\mathfrak{p}'_{2}(z_1\cdot
x)+c_2(x)\bigl)\\
=\exp\bigl(\mathfrak{p}_{2}\bigl(\exp(\mathfrak{p}_{1}(x))\cdot z_1\cdot
x\bigl)\bigl)\exp(-r_2).
\end{gather*}
The proof is almost the same as \emph{Step 1} in the proof of Lemma \ref{le:1}  (see Appendix \ref{sec:12}). Distinctions appear at
the level of notations. We omit the proof here.

\smallskip
\noindent\emph{Step 2}: We see that there are
$d_i:\mathcal{M}\to \mathfrak{G}$ with estimates
\begin{align}\label{for:89}
  \norm{d_i}_{C^m}\leq C_m\norm{\mathfrak{p}}_{C^{m}}^2,\qquad i=1,2
\end{align}
such that
\begin{gather*}
\exp\bigl(\mathfrak{p}'_{1}(
x)+d_1(x)\bigl)=\exp(\mathfrak{p}_{1}(x))\exp(-r_1),\\
\exp\bigl(\mathfrak{p}'_{2}(
x)+d_2(x)\bigl)=\exp(\mathfrak{p}_{2}(x))\exp(-r_2).
\end{gather*}
The result is a direct consequence of the  Baker-Campbell-Hausdorff formula.

\smallskip
\noindent\emph{Step 3}: We show that
\begin{align*}
  \norm{\log[z_1',z_2']}\leq C\norm{\mathfrak{p}}_{C^1}^2,\quad \norm{\log l(x)}_{C^{m}}\leq C_m\norm{\mathfrak{p}}_{C^{m}}.
\end{align*}
By \emph{Step 1}, \emph{Step 2} and \eqref{for:30} we have
\begin{align*}
&\exp\bigl(\mathfrak{p}'_{1}(z_2\cdot
x)+c_1(x)\bigl)\cdot
z_1'\cdot\exp\bigl(\mathfrak{p}'_{2}(
x)+d_2(x)\bigl)\cdot z_2'\cdot x\notag\\
&=l(x)\exp\bigl(\mathfrak{p}'_{2}(z_1\cdot
x)+c_2(x)\bigl)\cdot z_2'\cdot\exp\bigl(\mathfrak{p}'_{1}(
x)+d_1(x)\bigl)\cdot z_1'\cdot x.
\end{align*}
It follows that
\begin{align*}
&\exp\bigl(\mathfrak{p}'_{1}(z_2\cdot
x)+c_1(x)\bigl)\cdot
\exp\Bigl(\text{Ad}_{z_1'}\big(\mathfrak{p}'_{2}(
x)+d_2(x)\big)\Bigl)z_1'z_2'\notag\\
&=\exp\Bigl(\text{Ad}_{l(x)}\big(\mathfrak{p}'_{2}(z_1\cdot
x)+c_2(x)\big)\Bigl)\\
&\,\,\,\cdot\exp\Bigl(\text{Ad}_{l(x)z_2'}\big(\mathfrak{p}'_{1}(
x)+d_1(x)\big)\Bigl)l(x)z_2' z_1'.
\end{align*}
By the  Baker-Campbell-Hausdorff formula, \eqref{for:107} of \emph{Step 0}, \eqref{for:51} and \eqref{for:89}, there is
$q_1:\mathcal{M}\to \mathfrak{G}$, with
\begin{align}\label{for:94}
 \norm{q_1}_{C^m}\leq C_m\norm{\mathfrak{p}}_{C^{m+1}}\norm{\mathfrak{p}}_{C^{m}}
\end{align}
such that
\begin{align}\label{for:99}
 l(x)[z_2',z_1']=\exp(\mathbb{L}(x)+q_1(x))
\end{align}
where
\begin{align*}
 \mathbb{L}(x)&=\mathfrak{p}'_{1}(z_2\cdot
x)+\text{Ad}_{z_1'}\big(\mathfrak{p}'_{2}(
x)\big)-\text{Ad}_{l(x)}\big(\mathfrak{p}'_{2}(z_1\cdot
x)\big)-\text{Ad}_{l(x)z_2'}\big(\mathfrak{p}'_{1}(
x)\big).
\end{align*}
 We note that $l(x)\in L$  and $[z_1',z_2']\in [Z(L),Z(L)]$ (see \eqref{for:419} of \emph{Step 0}). \eqref{for:418} of of \emph{Step 0} allows us to use
the uniqueness decomposition
property (see \eqref{for:52} of Section \ref{sec:3}) to \eqref{for:99}. Thus, we have
\begin{align}\label{for:95}
 \norm{\log l(x)}_{C^m}+\norm{\log [z_2',z_1']}\leq C_m\norm{\mathbb{L}(x)+q_1(x)}_{C^{m}}\overset{\text{(1)}}{\leq} C_m\norm{\mathfrak{p}}_{C^{m}}.
\end{align}
Here in $(1)$ we use \eqref{for:107} of \emph{Step 0} and \eqref{for:94}.

Set
\begin{align*}
 \mathbb{L}_1(x)&=\mathfrak{p}'_{1}(z_2\cdot
x)+\text{Ad}_{z_1}\big(\mathfrak{p}'_{2}(
x)\big)-\mathfrak{p}'_{2}(z_1\cdot
x)-\text{Ad}_{z_2}\big(\mathfrak{p}'_{1}(
x)\big).
\end{align*}
Then we have
\begin{align}\label{for:48}
 \norm{\mathbb{L}-\mathbb{L}_1}_{C^m}&\overset{\text{(1)}}{\leq}C\norm{r_1}\norm{\mathfrak{p}}_{C^m}+C_m\norm{\log l}_{C^m}\norm{\mathfrak{p}}_{C^m}\notag\\
 &+C_m\Big\|\log \big(l(x)\exp(r_2)\Big)\big\|_{C^m}\norm{\mathfrak{p}}_{C^m}\notag\\
 &\overset{\text{(2,3)}}{\leq}C_1\norm{\mathfrak{p}}_{C^0}\norm{\mathfrak{p}}_{C^m}+C_{m,1}\norm{\mathfrak{p}}_{C^m}^2\notag\\
 &+C_{m,1}(\|\log l(x)\big\|_{C^m}+\|r_2\|)\norm{\mathfrak{p}}_{C^m}\notag\\
 &\overset{\text{(2)}}{\leq} C_{m,2}\norm{\mathfrak{p}}_{C^{m}}^2.
\end{align}
Here in $(1)$ we follow the arguments similar to those in \eqref{for:500}; in $(2)$ we use
\eqref{for:95}; in $(3)$  we use the  Baker-Campbell-Hausdorff formula.

From \eqref{for:94}, \eqref{for:99} and \eqref{for:48}, there is
$q_2:\mathcal{M}\to \mathfrak{G}$, with
\begin{align}\label{for:47}
 \norm{q_2}_{C^m}\leq C_m\norm{\mathfrak{p}}_{C^{m+1}}\norm{\mathfrak{p}}_{C^{m}}
\end{align}
such that
\begin{align}\label{for:36}
 \log[z_1',z_2']+\log l(x)&=\mathbb{L}(x)+q_1(x)=\mathbb{L}_1(x)+q_2(x).
\end{align}
Then
\begin{align*}
 &\log[z_1',z_2']\text{Vol}(\mathcal{M})+\int_{\mathcal{M}}\log l=\int_{\mathcal{M}}\mathbb{L}_1+\int_{\mathcal{M}}q_2\overset{\text{(1)}}{=}\int_{\mathcal{M}}q_2.
\end{align*}
Here in $(1)$ we recall that both $\mathfrak{p}'_{1}$ and $\mathfrak{p}'_{2}$ have zero averages.

We note that $\int_{\mathcal{M}}\log l\in\mathcal{L}$  and $\log[z_1',z_2']\in$ $ \text{Lie}([Z(L),Z(L)])$ (see \eqref{for:419} of \emph{Step 0}). By \eqref{for:418} of \emph{Step 0} we see that
\begin{align}\label{for:49}
 \norm{\log[z_1',z_2']}\leq C\big\|\int_{\mathcal{M}}q_2\big\|\overset{\text{(1)}}{\leq} C\norm{\mathfrak{p}}_{C^1}\norm{\mathfrak{p}}_{C^0}.
\end{align}
Here in $(1)$ we use \eqref{for:47}.

\smallskip
\noindent\emph{Step 4}: We show that
\begin{align}\label{for:97}
  \norm{\log[z_1,z_2]}\leq C\norm{\mathfrak{p}}_{C^0}.
\end{align}
From \eqref{for:30} we have
\begin{align*}
&\exp\bigl(\mathfrak{p}_{1}\bigl(\exp(\mathfrak{p}_{2}(x))\cdot z_2\cdot
x\bigl)\bigl)\cdot
\exp\big(\text{Ad}_{z_1}\mathfrak{p}_{2}(x)\big)\cdot [z_1,z_2]\notag\\
&=l(x)\exp\bigl(\mathfrak{p}_{2}\bigl(\exp(\mathfrak{p}_{1}(x))\cdot z_1\cdot
x\bigl)\bigl)\cdot\exp\big(\text{Ad}_{z_2}\mathfrak{p}_{1}(x)\big).
\end{align*}
By \eqref{for:95} we see that $\emph{\textbf{d}}(l(x),e)<C\norm{\mathfrak{p}}_{C^0}$. This together with the above equation gives
\begin{align*}
 \emph{\textbf{d}}([z_1,z_2],e)<C\norm{\mathfrak{p}}_{C^0}.
\end{align*}
This together with \eqref{for:35} of Section \ref{sec:3} give \eqref{for:97}.

\smallskip
\noindent\emph{Step 5}: We show that
\begin{align*}
  \norm{\mathbb{D}_{1,2}}_{C^{m}}
 \leq C_m\norm{\mathfrak{p}}_{C^{m+1}}\norm{\mathfrak{p}}_{C^{m}}.
\end{align*}
We note that $\mathfrak{i}_{\mathcal{L}^\bot}\mathbb{L}_1=\mathbb{D}_{1,2}$. \eqref{for:36} shows that
\begin{align*}
  \mathfrak{i}_{\mathcal{L}^\bot}\mathbb{L}_1(x)=\mathfrak{i}_{\mathcal{L}^\bot}\big(\log[z_1',z_2']-q_2(x)\big).
\end{align*}
\eqref{for:47} and \eqref{for:49} show that
\begin{align}\label{for:87}
\norm{\mathfrak{i}_{\mathcal{L}^\bot}\mathbb{L}_1}_{C^m}\leq C_m\norm{\mathfrak{p}}_{C^{m+1}}\norm{\mathfrak{p}}_{C^{m}}.
\end{align}

\noindent\emph{Step 6}: We show that: for any $1\neq\mu\in \Delta_A$ we have
\begin{align*}
  \big\|r_1|_{\mathfrak{g}_{\mu,A}}-\text{Ad}_{z_2}(r_1|_{\mathfrak{g}_{\mu,A}})-r_2|_{\mathfrak{g}_{\mu,A}}+\text{Ad}_{z_1}(r_2|_{\mathfrak{g}_{\mu,A}})\big\|\leq C\norm{\mathfrak{p}}_{C^1}\norm{\mathfrak{p}}_{C^0}
\end{align*}
By definition we have
\begin{align*}
 &\quad[z_1',z_2']\exp(r_2)z_2\exp(r_1)z_1=\exp(r_1)z_1\exp(r_2)z_2\\
 \Rightarrow &[z_1',z_2']\exp(r_2)\exp(\text{Ad}_{z_2}r_1)z_2z_1=\exp(r_1)\exp(\text{Ad}_{z_1}r_2)z_1z_2\\
 \Rightarrow &[z_1',z_2']\exp(r_2)\exp(\text{Ad}_{z_2}r_1)=\exp(r_1)\exp(\text{Ad}_{z_1}r_2)[z_1,z_2]
\end{align*}
By the  Baker-Campbell-Hausdorff formula, \eqref{for:97} and \eqref{for:49} we have
\begin{align*}
  \norm{r_1-\text{Ad}_{z_2}r_1-r_2+\text{Ad}_{z_1}r_2+\log [z_1,z_2]}\leq C\norm{\mathfrak{p}}_{C^1}\norm{\mathfrak{p}}_{C^0}.
\end{align*}
As $\log[z_1,z_2]\in \mathfrak{g}_{1,A}$ and each $\mathfrak{g}_{\lambda}$, $\lambda\in \Delta_A$ is invariant under $\text{Ad}_{z_1}$  and $\text{Ad}_{z_2}$  (see \eqref{for:490} of Section \ref{sec:3})  the above inequality implies the result.

\section{Proof of Proposition \ref{po:1}}\label{sec:16}
The proof is similar to that of Proposition \ref{po:2}, but new arguments are also needed to handle the non-abelianness.  Firstly, we introduce some notations which will be used for the subsequent part. For any poor row $(\mu,\lambda_1,\cdots,\lambda_\sigma)$ (see \eqref{for:1} of Section \ref{sec:3}), we define $L_{(\mu,\lambda_1,\cdots,\lambda_\sigma)}$,  a subset of $\ZZ$ as follows:

1. Suppose $\mu=1$. Set $L_{(\mu,\lambda_1,\cdots,\lambda_\sigma)}=\ZZ^-\bigcup \{0\}$.

\smallskip
2. Suppose  $\mu\neq1$. From the definition of poor row, there is $c(\mu,\lambda_1,\cdots,\lambda_\sigma)>0$ such that
 \begin{align*}
  \lambda_1(\textbf{a}_1)\lambda_2(\textbf{a}_1)\cdots \lambda_\sigma(\textbf{a}_1)=\mu(\textbf{a}_1)^c.
  \end{align*}
 It follows that
 \begin{align}\label{for:445}
  c=\sum_{i=1}^\sigma \frac{\log\lambda_i(\textbf{a}_1)}{\log\mu(\textbf{a}_1)}\overset{\text{(1)}}{\geq} \sigma \min_{\mu_1\neq\mu_2\in \Delta_A\backslash\{1\}}\frac{|\log\mu_1(\textbf{a}_1)|}{|\log\mu_2(\textbf{a}_1)|}\overset{\text{(2)}}{>} \sigma\cdot 5\emph{\textbf{l}}.
 \end{align}
 Here in $(1)$ we note that $\frac{\log\lambda_i(\textbf{a}_1)}{\log\mu(\textbf{a}_1)}=\frac{|\log\mu_1(\textbf{a}_1)|}{|\log\mu_2(\textbf{a}_1)|}$, $1\leq i\leq \sigma$, as
 $(\mu,\lambda_1,\cdots,\lambda_\sigma)$ is a poor row; in $(2)$ we recall the definition of $\emph{\textbf{l}}(\textbf{a}_1,\textbf{a}_2)$, see \eqref{for:63}

  If $\lambda_1(\textbf{a}_2)\lambda_2(\textbf{a}_2)\cdots \lambda_\sigma(\textbf{a}_2)\geq \mu(\textbf{a}_2)^c$ (resp. $\lambda_1(\textbf{a}_2)\lambda_2(\textbf{a}_2)\cdots \lambda_\sigma(\textbf{a}_2)< \mu(\textbf{a}_2)^c$), set $L_{(\mu,\lambda_1,\cdots,\lambda_\sigma)}=\ZZ^+\bigcup \{0\}$ (resp. $L_{(\mu,\lambda_1,\cdots,\lambda_\sigma)}=\ZZ^-\bigcup \{0\}$). Set
\begin{align*}
 \mathbb{S}_\mu=\{(m,n)\in\ZZ^2: \mu(\textbf{a}_1)^{n}\mu(\textbf{a}_2)^{m}
 \leq e^{\emph{\textbf{l}}(|n|+|m|)}\}.
\end{align*}
As the first step proving the proposition, we need the following result, whose proof is similar to that of Lemma \ref{le:2}:
\begin{lemma}\label{le:13}  For any $\textbf{u}=u_1\cdots u_\sigma\in \mathcal{U}_{\mu,\sigma}$ (see \eqref{for:1} of Section \ref{sec:3}), where $u_i\in \mathfrak{g}_{\lambda_i,A}\cap\mathfrak{G}^1$, $1\leq i\leq \sigma$, any $n\in\ZZ$ and $m\in L_{(\mu,\lambda_1,\cdots,\lambda_\sigma)}$, we have
  \begin{align}\label{for:26}
 \Big\|(\text{Ad}_{z_2}|&_{J_\mu})^{m-1}(\text{Ad}_{z_1}|_{J_\mu})^{n-1}\big\langle \xi\circ (\mathfrak{z}_1^{-n}\mathfrak{z}_2^{-m}), \,\textbf{u}\psi\big\rangle\Big\|\notag\\
 &\leq Ce^{-\textbf{l}_2(|n|+|m|)}\norm{\xi}_{\sigma+1}\norm{\psi}_{\sigma+1}
\end{align}
where $\textbf{l}_2=\min\{2 \textbf{l},\, 2\sigma \textbf{l}^2\}$, for any $\xi,\,\psi\in J_\mu(\mathcal{O}^{\sigma+1})$.
\end{lemma}
\begin{proof} Set
\begin{align*}
 \mathcal{E}_{n,m}=(\text{Ad}&_{z_2}|_{J_\mu})^{m-1}(\text{Ad}_{z_1}|_{J_\mu})^{n-1}\langle \xi\circ (\mathfrak{z}_1^{-n}\mathfrak{z}_2^{-m}), \textbf{u}\psi\rangle.
\end{align*}
For any $n,\,m\in\ZZ$ we have
\begin{align}\label{for:444}
  \norm{\mathcal{E}_{n,m}} &\leq \norm{(\text{Ad}_{z_2}|_{J_\mu})^{m-1}}\norm{(\text{Ad}_{z_1}|_{J_\mu})^{n-1}}
  \Big\|\big\langle \xi\circ (\mathfrak{z}_1^{-n}\mathfrak{z}_2^{-m}), \,\textbf{u}\psi\big\rangle\Big\|\notag\\
  &\overset{\text{(a)}}{\leq} C\mu(\textbf{a}_1)^{n-1}\mu(\textbf{a}_2)^{m-1}e^{ \emph{\textbf{l}}_1(|n-1|+|m-1|)}\|\langle \xi\circ (\mathfrak{z}_1^{-n}\mathfrak{z}_2^{-m}), \textbf{u}\psi\rangle\|
  \end{align}
Here in $(a)$ we use \eqref{for:447} of Proposition \ref{po:9} and \eqref{for:435} of Corollary \ref{for:129}.

\smallskip

\noindent\emph{Step 1: We show that}: if $(m,n)\in \mathbb{S}_\mu$, then
\begin{align*}
 \norm{\mathcal{E}_{n,m}}\leq Ce^{-2 \emph{\textbf{l}}(|n|+|m|)}\norm{\xi}_{1}\norm{\psi}_{\sigma+1}.
\end{align*}
For any $(m,n)\in \mathbb{S}_\mu$, by using \eqref{for:444} and by using \eqref{for:403} of Proposition \ref{po:9} and \eqref{for:435} of Corollary \ref{for:129} we have
\begin{align}\label{for:55}
 \norm{\mathcal{E}_{n,m}}&\leq C\mu(\textbf{a}_1)^{n-1}\mu(\textbf{a}_2)^{m-1}e^{ \emph{\textbf{l}}_1(|n-1|+|m-1|)}\cdot e^{-4 \emph{\textbf{l}}(|n|+|m|)}\norm{\xi}_{1}\norm{\textbf{u}\psi}_{1}\notag\\
  &\overset{\text{(3)}}{\leq} C_1e^{2\emph{\textbf{l}}(|n-1|+|m-1|)}e^{-4 \emph{\textbf{l}}(|n|+|m|)}\norm{\xi}_{1}\norm{\textbf{u}\psi}_{1}\notag\\
  &\leq
  C_2e^{-2 \emph{\textbf{l}}(|n|+|m|)}\norm{\xi}_{1}\norm{\psi}_{\sigma+1}.
\end{align}
Here in $(1)$ we use the assumption for $\mathbb{S}_\mu$ and the fact that $\emph{\textbf{l}}_1=\emph{\textbf{l}}^2< \emph{\textbf{l}}$.

\smallskip
\noindent\emph{Step 2: We show that}: if $\mu=1$, then for any $n\in\ZZ$ and $m\in L_{(\mu,\lambda_1,\cdots,\lambda_\sigma)}$ we have
\begin{align*}
 \norm{\mathcal{E}_{n,m}}\leq Ce^{-2 \emph{\textbf{l}}(|n|+|m|)}\norm{\xi}_{1}\norm{\psi}_{\sigma+1}.
\end{align*}
We note that if $\mu=1$, then for any $n\in\ZZ$ and $m\in L_{(\mu,\lambda_1,\cdots,\lambda_\sigma)}$, $(m,n)\in \mathbb{S}_\mu$. Then the result follows from \emph{Step 1}.

\smallskip

\emph{Note}. From  now on, we only consider the case of $(m,n)\notin \mathbb{S}_\mu$. Otherwise, the result follows directly from \emph{Step 1}.

\smallskip
\noindent\emph{Step 3: We show that}: if $\mu\neq1$, then for any $n\in\ZZ$ and  $m\in L_{(\mu,\lambda_1,\cdots,\lambda_\sigma)}$ we have
\begin{align*}
 &\norm{\mathcal{E}_{n,m}}\leq Ce^{-4\sigma\emph{\textbf{l}}_1(|n|+|m|)}\norm{\xi}_{\sigma}\norm{\psi}.
\end{align*}
From the definition of $L_{(\mu,\lambda_1,\cdots,\lambda_\sigma)}$ we note that in this case
\begin{align*}
 \Big(\big(\Pi_{i=1}^\sigma\lambda_i(\textbf{a}_2)\big)\mu(\textbf{a}_2)^{-c}\Big)^{m}>1,
\end{align*}
which further gives
\begin{align}\label{for:16}
\Pi_{i=1}^\sigma\big(\lambda_i&(\textbf{a}_1)^{n}\lambda_i(\textbf{a}_2)^{m}\big)=\mu(\textbf{a}_1)^{cn}\mu(\textbf{a}_2)^{cm}
\Big(\big(\Pi_{i=1}^\sigma\lambda_i(\textbf{a}_2)\big)\mu(\textbf{a}_2)^{-c}\Big)^{m}\notag\\
&\geq \mu(\textbf{a}_1)^{cn}\mu(\textbf{a}_2)^{cm}\overset{\text{(1)}}{\geq} \mu(\textbf{a}_1)^{5\sigma \emph{\textbf{l}} n}\mu(\textbf{a}_2)^{5\sigma \emph{\textbf{l}}m}.
\end{align}
Here in $(1)$ we use \eqref{for:445} and
note that $\mu(\textbf{a}_1)^{n}\mu(\textbf{a}_2)^{m}>1$ since $(m,n)\notin \mathbb{S}_\mu$.

Hence, we have
\begin{align}\label{for:21}
  \norm{\mathcal{E}_{n,m}}&\overset{\text{(1)}}{\leq} C\mu(\textbf{a}_1)^{n-1}\mu(\textbf{a}_2)^{m-1}e^{\emph{\textbf{l}}_1(|n-1|+|m-1|)}\notag\\
  &\cdot\big\|u_\sigma\cdots u_1(\xi\circ (\mathfrak{z}_1^{-n}\mathfrak{z}_2^{-m}))\big\|\norm{\psi}\notag\\
  &\overset{\text{(2)}}{=} C\mu(\textbf{a}_1)^{n-1}\mu(\textbf{a}_2)^{m-1}e^{\emph{\textbf{l}}_1(|n-1|+|m-1|)}\notag\\
  &\cdot(\Pi_{i=1}^\sigma\norm{\text{Ad}_{\mathfrak{z}_1^{-n}\mathfrak{z}_2^{-m}}u_i})\big\|(\textbf{u}_{m,n}\xi)\circ (\mathfrak{z}_1^{-n}\mathfrak{z}_2^{-m})\big \|\norm{\psi}\notag\\
  &\overset{\text{(3)}}{\leq} C_1\mu(\textbf{a}_1)^{n}\mu(\textbf{a}_2)^{m}e^{\emph{\textbf{l}}_1(|n|+|m|)}\notag\\
  &\cdot \Pi_{i=1}^\sigma\big(\lambda_i(\textbf{a}_1)^{-n}\lambda_i(\textbf{a}_2)^{-m}e^{\emph{\textbf{l}}_1(|n|+|m|)}\big)\norm{\xi}_{\sigma}\norm{\psi}\notag\\
  &\overset{\text{(4)}}{\leq} C_1\mu(\textbf{a}_1)^{n}\mu(\textbf{a}_2)^{m}e^{\emph{\textbf{l}}_1(|n|+|m|)}\notag\\
  &\cdot \mu(\textbf{a}_1)^{-5\sigma \emph{\textbf{l}}n}\mu(\textbf{a}_2)^{-5\sigma \emph{\textbf{l}}m}e^{\sigma \emph{\textbf{l}}_1(|n|+|m|)}\norm{\xi}_{\sigma}\norm{\psi}\notag\\
  &=C_1\big(\mu(\textbf{a}_1)^{-n}\mu(\textbf{a}_2)^{-m}\big)^{5\sigma \emph{\textbf{l}}-1}e^{(\sigma+1) \emph{\textbf{l}}_1(|n|+|m|)}\norm{\xi}_{\sigma}\norm{\psi}\notag\\
  &\overset{\text{(5)}}{\leq}
  C_{1}e^{-(5\sigma \emph{\textbf{l}}-1)\emph{\textbf{l}}(|n|+|m|)}e^{(\sigma+1) \emph{\textbf{l}}_1(|n|+|m|)}\norm{\xi}_{\sigma}\norm{\psi}\notag\\
  &\overset{\text{(6)}}{\leq}
  C_{1}e^{-2\sigma \emph{\textbf{l}}_1(|n|+|m|)}\norm{\xi}_{\sigma}\norm{\psi}
\end{align}
Here in $(1)$ we use \eqref{for:444} and cauchy schwarz inequality; in $(2)$
\begin{align*}
 \textbf{u}_{m,n}=(\Pi_{i=1}^\sigma\norm{\text{Ad}_{\mathfrak{z}_1^{-n}\mathfrak{z}_2^{-m}}u_i})^{-1}\text{Ad}_{\mathfrak{z}_1^{-n}\mathfrak{z}_2^{-m}}
 (u_\sigma)\cdots\text{Ad}_{\mathfrak{z}_1^{-n}\mathfrak{z}_2^{-m}}(u_1)
\end{align*}
in $(3)$ we use \eqref{for:447} of Proposition \ref{po:9} and \eqref{for:435} of Corollary \ref{for:129}, as well as  $u_i\in \mathfrak{g}_{\lambda_i,A}\cap\mathfrak{G}^1$, $1\leq i\leq \sigma$; in $(4)$ we use \eqref{for:16}; in $(5)$ we use $(m,n)\notin \mathbb{S}_1$, which implies that
\begin{align*}
\mu(\textbf{a}_1)^{-n}\mu(\textbf{a}_2)^{-m}
 \leq e^{-\emph{\textbf{l}}(|n|+|m|)};
\end{align*}
in $(6)$ we recall \eqref{for:27} of Section \ref{for:432}:
\begin{align*}
 e^{-(5\sigma \emph{\textbf{l}}-1) \emph{\textbf{l}}}e^{(\sigma+1) \emph{\textbf{l}}_1}\leq  e^{-4\sigma \emph{\textbf{l}}^2}e^{2\sigma\emph{\textbf{l}}_1}
 =e^{-2\sigma \emph{\textbf{l}}_1}.
\end{align*}

\eqref{for:55} and \eqref{for:21} imply the result;  Hence we finish the proof.
\end{proof}

\subsection{Proof of Proposition \ref{po:1}}

We fix $l>0$ ($l$ is determined in \emph{Step 6}, see \eqref{for:363}).  To simplify notations, we denote
\begin{align*}
B_i=(\text{Ad}_{z_i}|_{J_\mu}),\,\,i=1,2\quad\text{and}\quad \Phi_{m}=\sum_{n=-l}^l B_2^{m}B_1^{n-1}\theta_\mu\circ (\mathfrak{z}_1^{-n}\mathfrak{z}_2^{-m}).
\end{align*}
Suppose $\textbf{u}\in \mathcal{U}_{\mu,\sigma}$ and $\psi\in J_{\mu}(\mathcal{O}^{\sigma+1})$.

Typically distinctions between the proofs of Proposition \ref{po:1} and Proposition \ref{po:2} (see Section \ref{sec:44})
 are minimal, and they usually appear at
the level of notations.

\smallskip
\noindent\emph{Step 1: We show that}
\begin{align}\label{for:233}
 \big\|\langle\mathfrak{B}_{\mu,z_1,\textbf{u},\ZZ}(\theta_\mu),\,\psi\rangle\big\|\leq \big\|\langle \Phi_{0},\, \textbf{u}\psi\rangle\big\|+Ce^{-3l\emph{\textbf{l}}}\norm{\theta_\mu}_\sigma\norm{\psi}_{\sigma+1}.
\end{align}
We have
\begin{align*}
&\big\|\langle\mathfrak{B}_{\mu,z_1,\textbf{u},\ZZ}(\theta_\mu),\,\psi\rangle\big\|\notag\\
    &=\Big\|\sum_n\big\langle (\text{Ad}_{z_1}|_{J_\mu})^{n-1}(\theta_\mu\circ \mathfrak{z}_1^{-n}),\, \textbf{u}\psi\big\rangle\Big\|\notag\\
    &\leq \big\|\langle \Phi_{0},\, \textbf{u}\psi\rangle\big\|+\big\|\sum_{|n|>l}\langle B_1^{n-1}\theta_\mu\circ \mathfrak{z}_1^{-n},\, \textbf{u}\psi\rangle\big\|\notag\\
    &\overset{\text{(1)}}{\leq} \big\|\langle \Phi_{0},\, \textbf{u}\psi\rangle\big\|+\sum_{|n|>l}Ce^{-3\emph{\textbf{l}}|n|}\norm{\theta_\mu}_\sigma\norm{\psi}_{\sigma+1}\notag\\
    &\leq \big\|\langle \Phi_{0},\, \textbf{u}\psi\rangle\big\|+C_1e^{-3l\emph{\textbf{l}}}\norm{\theta_\mu}_\sigma\norm{\psi}_{\sigma+1}.
  \end{align*}
Here in $(1)$ we use \eqref{for:448} of Lemma \ref{le:12}. Then we finish the proof.

\smallskip
\noindent\emph{Step 2: We show that}
\begin{align}\label{for:60}
\big\|\langle \Phi_{0},\, \textbf{u}\psi\rangle\big\|\leq\sum_{m\in L_{(\mu,\lambda_1,\cdots,\lambda_\sigma)}}\big\|\langle \Phi_{m}-\Phi_{m+1}, \textbf{u}\psi\rangle\big\|.
\end{align}
For any $m\in L_{(\mu,\lambda_1,\cdots,\lambda_\sigma)}$,  it follows from Lemma \ref{le:13} that
\begin{align}
 \big\|\langle \Phi_{m},\, \textbf{u}\psi\rangle\big\|&\leq\sum_{n} Ce^{-\emph{\textbf{l}}_2(|n|+|m|)}\norm{\xi}_{\sigma}\norm{\psi}_{\sigma+1}\notag\\
  &\leq C_1e^{-\emph{\textbf{l}}_2|m|}\norm{\xi}_{\sigma}\norm{\psi}_{\sigma+1}.
\end{align}
Letting $|m|\to\infty$, we get
\begin{align*}
 \big\|\langle \Phi_{m},\, \textbf{u}\psi\rangle\big\|\to 0,\qquad\text{as }|m|\to\infty.
\end{align*}
This implies \eqref{for:60}.

\smallskip
\noindent\emph{Step 3: We show that}
\begin{align}\label{for:451}
 \sum_{m\in L_{(\mu,\lambda_1,\cdots,\lambda_\sigma)}}\Big\|\big\langle \Phi_{m}-\Phi_{m+1}, \textbf{u}\psi\big\rangle\Big\|\leq \mathcal{E}_1+\mathcal{E}_2,
\end{align}
where
\begin{gather*}
 \mathcal{E}_1=\sum_{m\in L_{\mu}}\Big\|\langle B_2^{m}\Psi_{l}\circ \mathfrak{z}_2^{-(m+1)}, \,\textbf{u}\psi\rangle\Big\| \qquad \text{and}\\
 \mathcal{E}_2=\sum_{m\in L_{\mu}}\Big\|\big\langle B_2^{m}\Psi_{l}\circ \mathfrak{z}_2^{-(m+1)}-(\Phi_{m}-\Phi_{m+1}), \,\textbf{u}\psi\big\rangle\Big\|.
\end{gather*}
Here $\Psi_{l}$ is defined as follows:
\begin{align}\label{for:3}
 \Psi_{l}=B_1^{-(l+1)}\varphi_\mu\circ \mathfrak{z}_1^{l+1} -B_1^{l}\varphi_\mu\circ \mathfrak{z}_1^{-l}-\sum_{n=-l}^l B_1^{n-1}\omega_\mu \circ \mathfrak{z}_1^{-n}.
\end{align}
From equation \ref{for:15}, we have
\begin{align}\label{for:4}
 &\sum_{n=-\ell}^\ell B_1^{n-1}\theta_\mu\circ(z_2z_1^{-n})-\sum_{n=-\ell}^\ell B_1^{n-1}B_2\theta_\mu\circ z_1^{-n}=\Psi_{\ell}
\end{align}
The proof is the same as that of \emph{Step 3} in Section \ref{sec:44} if we replace $\pi_{u}(f_1\circ \tau^{-1})\psi$ by $\textbf{u}\psi$ and replace $z_i$ by $\mathfrak{z}_i$, $i=1,\,2$.

\smallskip
\noindent\emph{Step 4: We show that}
\begin{align}\label{cor:5}
 \mathcal{E}_1\leq C(\norm{\omega_\mu}_{\sigma+1}+e^{-l\emph{\textbf{l}}_2}\norm{\varphi_\mu}_{\sigma+1})\norm{\psi}_{\sigma+1}.
\end{align}
The proof is the same as that of \emph{Step 4} in Section \ref{sec:44} if we replace $\pi_{u}(f_1\circ \tau^{-1})\psi$ by $\textbf{u}\psi$, replace $z_i$ by $\mathfrak{z}_i$, $i=1,\,2$, replace Lemma \ref{le:2} by Lemma \ref{le:13}, replace $\gamma_3$ by $\emph{\textbf{l}}_2$ and replace $\kappa$ by $\sigma+1$.

\smallskip
\noindent\emph{Step 5: We show that}
\begin{align}
\mathcal{E}_2&\leq \sum_{m\in L_{(\mu,\lambda_1,\cdots,\lambda_\sigma)}}\big\|\langle Y_{1,m}, \,\textbf{u}\psi\rangle\big\|+\sum_{m\in L_{(\mu,\lambda_1,\cdots,\lambda_\sigma)}}\big\|\langle Y_{2,m},\, \textbf{u}\psi\rangle\big\|
\end{align}
where
\begin{gather*}
 Y_{1,m}=\sum_{n=-l}^l B_2^{m}B_1^{n-1}(\theta_\mu\circ s_n-\theta_\mu)\circ(\mathfrak{z}_1^{-n}\mathfrak{z}_2^{-m}),\\
 Y_{2,m}=\sum_{n=-l}^l B_2^{m+1}B_1^{n-1}\big((\text{Ad}_{k_n}|_{J_\mu})-I\big)\theta_\mu\circ (\mathfrak{z}_1^{-n}\mathfrak{z}_2^{-(m+1)}), \qquad\text{and}\\
 s_n=\mathfrak{z}_2\mathfrak{z}_1^{-n}\mathfrak{z}_2^{-1}\mathfrak{z}_1^n,\quad k_n=z_1^{-(n-1)}z_2^{-1}z_1^{n-1}z_2.
\end{gather*}
By using the form of $\Psi_{\ell}$ in \eqref{for:4}, we have
\begin{align*}
  B_2^{m}\Psi_{\ell}\circ \mathfrak{z}_2^{-(m+1)}&-(\Phi_{m}-\Phi_{m+1})=Y_{1,m}-Y_{2,m}.
\end{align*}
This implies the result.

\smallskip
\noindent\emph{Step 6:}  Let
\begin{align*}
 c_{\textbf{a}_1}=\max_{\mu\in \Delta_A}|\log\mu(\textbf{a}_1)|\quad\text{and}\quad \norm{\log[z_1,z_2]}=x^2.
\end{align*}
We can assume $x$ is sufficiently small such that $x\leq \delta$ (see \eqref{for:35} of Section \ref{sec:3}). We can choose $l>0$ such that
\begin{align*}
  e^{l }=x^{-\frac{5}{8}c_{\textbf{a}_1}^{-1}}.
\end{align*}
The next result studies the commutator relations between
$z_2^{\pm}$ and $z_1^{n}$. We show that for any $|n|\leq l$
\begin{align}\label{le:6}
\max\{\textbf{\emph{d}}(s_n,e),\,\textbf{\emph{d}}(k_n,e)\}\leq x.
\end{align}
\smallskip
We can write $z_2z_1=\exp(\mathfrak{n})z_1z_2$. Then $\norm{\mathfrak{n}}=x^2$. We note that
\begin{align*}
  z_2z_1z_2^{-1}&=\exp(\mathfrak{n})z_1\Rightarrow z_2z_1^{|n|}z_2^{-1}=(\exp(\mathfrak{n})z_1)^{|n|}\\
  &\Rightarrow z_2z_1^{|n|}z_2^{-1}z_1^{-|n|}=\exp(\mathfrak{n})\exp(\text{Ad}_{z_1}\mathfrak{n})\cdots \exp(\text{Ad}_{z_1^{|n|-1}}\mathfrak{n}).
\end{align*}
For any $0\leq j\leq |n|-1\leq l-1$ we have
\begin{align}\label{for:449}
 \norm{\text{Ad}_{z_1^{j}}\mathfrak{n}}\overset{\text{(1)}}{\leq} Ce^{lc_{\textbf{a}_1}}e^{l\emph{\textbf{l}}}\norm{\mathfrak{n}}\overset{\text{(2)}}{\leq} Ce^{l c_{\textbf{a}_1}\frac{6}{5}}x^2\overset{\text{(3)}}{\leq} Cx^{\frac{5}{4}}<x.
\end{align}
Here in $(1)$ we use \eqref{for:108} of Section \ref{sec:3}; in $(2)$ from \eqref{for:63} and \eqref{for:110} we see that $\emph{\textbf{l}}< \frac{1}{5}c_{\textbf{a}_1}$;
in $(3)$ we note that $x$ is sufficiently small.

Hence we have
\begin{align}\label{for:450}
  \textbf{\emph{d}}(z_2&z_1^{|n|}z_2^{-1}z_1^{-|n|},\,e)\overset{\text{(1)}}{\leq} \sum_{j=0}^{|n|-1} \textbf{\emph{d}}(\exp\big(\text{Ad}_{z_1^j}\mathfrak{n}\big),e)
  \overset{\text{(2)}}{\leq} C\sum_{j=0}^{|n|-1} \norm{\text{Ad}_{z_1^{j}}\mathfrak{n}}\notag\\
  &\overset{\text{(3)}}{\leq}  C_1l\cdot x^{\frac{5}{4}} \leq C_2x^{\frac{8}{7}}.
\end{align}
Here in $(1)$ we use right invariance of $\textbf{\emph{d}}$ (see );  in $(2)$ we use \eqref{for:35} of Section \ref{sec:3}
and note that $\norm{\text{Ad}_{z_1^{j}}\mathfrak{n}}< x\leq \delta$ (see \eqref{for:449});
in $(3)$ we use \eqref{for:449}.

Since $z_2z_1^{|n|}z_2^{-1}z_1^{-|n|}=(z_1^{|n|}z_2z_1^{-|n|}z_2^{-1})^{-1}$ and $x$ is sufficiently small, by using \eqref{for:35} of Section \ref{sec:3}, we see that \eqref{for:450} also implies that
\begin{align*}
 \textbf{\emph{d}}(z_1^{|n|}z_2z_1^{-|n|}z_2^{-1},\,e)\leq Cx^{\frac{8}{7}}.
\end{align*}
We also note that
\begin{align*}
  z_2z_1z_2^{-1}&=\exp(\mathfrak{n})z_1\Rightarrow z_2z_1^{-|n|}z_2^{-1}=(\exp(\mathfrak{n})z_1)^{-|n|}\\
  &\Rightarrow z_2z_1^{-|n|}z_2^{-1}\mathfrak{n}z_1^{|n|}=\exp(-dz_1^{-1}\mathfrak{n})\exp(-dz_1^{-2}\mathfrak{n})\cdots \exp(dz_1^{-|n|}\mathfrak{n}).
\end{align*}
By using the same reasoning, we can also show that
\begin{align*}
 \textbf{\emph{d}}(z_2z_1^{-|n|}z_2^{-1}z_1^{|n|},\,e)\leq Cx^{\frac{8}{7}}\quad\text{and}\quad\textbf{\emph{d}}(z_1^{-|n|}z_2z_1^{|n|}z_2^{-1},\,e)\leq Cx^{\frac{8}{7}}.
\end{align*}
Since $z_i$ and $\mathfrak{z}_i$ are conjugate, $i=1,\,2$, the above discussion implies the result.

\smallskip
\noindent\emph{Step 7: We show that }
\begin{align}\label{le:11}
  \mathcal{E}_2\leq Cx^{\frac{8}{7}}\|\theta_\mu\|_{\sigma+2}\norm{\psi}_{\sigma+1}.
\end{align}
The proof is the same as that of \emph{Step 7} in Section \ref{sec:44} if we replace $\pi_{u}(f_1\circ \tau^{-1})\psi$ by $\textbf{u}\psi$, replace $L_{\mu}$ by $L_{(\mu,\lambda_1,\cdots,\lambda_\sigma)}$, replace Lemma \ref{le:2} by Lemma \ref{le:13}, replace $\gamma_3$ by $\emph{\textbf{l}}_2$ and replace $\kappa$ by $\sigma+1$.

\smallskip

\emph{Conclusion}: It follows from \eqref{for:233}, \eqref{for:60}, \eqref{for:451}, \eqref{cor:5} and \eqref{le:11} that
\begin{align*}
    &\big\|\langle\mathfrak{B}_{\mu,z_1,\textbf{u},\ZZ}(\theta_\mu),\,\psi\rangle\big\|\notag\\
    &\leq Ce^{-3l\emph{\textbf{l}}}\norm{\theta_\mu}_\sigma\norm{\psi}_{\sigma+1}+C(\norm{\omega_\mu}_{\sigma+1}+e^{-l\emph{\textbf{l}}_2}\norm{\varphi_\mu}_{\sigma+1})\norm{\psi}_{\sigma+1}\\
    &+Cx^{\frac{8}{7}}\|\theta_\mu\|_{\sigma+2}\norm{\psi}_{\sigma+1}\\
     &\overset{\text{(1)}}{\leq} C(x^2)^{c_0}\norm{\theta_\mu}_\sigma\norm{\psi}_{\sigma+1}+C(\norm{\omega_\mu}_{\sigma+1}+(x^2)^{c_0}\norm{\varphi_\mu}_{\sigma+1})\norm{\psi}_{\sigma+1}\\
    &+C(x^2)^{c_0}\|\theta_\mu\|_{\sigma+2}\norm{\psi}_{\sigma+1}.
  \end{align*}
Here in $(1)$ we note that $e^{l }=x^{-\frac{5}{8}c_{\textbf{a}_1}^{-1}}$ and
let $c_0=\min\{\frac{4}{7},\,\frac{15}{16}\emph{\textbf{l}}c_{\textbf{a}_1}^{-1},\,\frac{5}{16}\emph{\textbf{l}}_2c_{\textbf{a}_1}^{-1}\}$. Hence we finish the proof.



\end{document}